%% file: Main_singleStor_TPS_OMG.tex
\documentclass[journal,twocolumn]{IEEEtran}

\input{defs.tex}

\usepackage{scalefnt}

\begin{document}

\title{Online Modified Greedy Algorithm\\ for Storage Control under Uncertainty}
\author{Junjie~Qin,~\IEEEmembership{Student Member, IEEE,}
        Yinlam~Chow,~\IEEEmembership{Student Member, IEEE,}
        Jiyan~Yang,~\IEEEmembership{Student Member, IEEE,}
        and~Ram~Rajagopal,~\IEEEmembership{Member, IEEE}%
\thanks{This research was supported in part by the Satre Family fellowship, and in part by the Tomkat Center for Sustainable Energy.}% <-this % stops a space
%\thanks{Y. Zhao is with the Dept. of Electrical Engineering, Stanford University, Stanford, CA, 94305 USA, and with the Dept. of
%Electrical Engineering Princeton University, Princeton, NJ, 08544 USA, (e-mail: yuez@princeton.edu).}%
\thanks{J. Qin, Y. Chow and J. Yang are with the Institute for Computational and Mathematical Engineering, Stanford University, Stanford, CA, 94305 USA, (e-mail: \{jqin,ychow,jiyan\}@stanford.edu).}%
\thanks{R. Rajagopal is with the Department of Civil and Environmental Engineering, Stanford University, Stanford, CA 94305 USA (e-mail: ramr@stanford.edu).}%
%\thanks{A. Goldsmith is with the Dept. of Electrical Engineering, Stanford University, Stanford, CA 94305 USA (e-mail: andrea@stanford.edu).}%
%\thanks{H. V. Poor is with the Dept. of Electrical Engineering, Princeton University, Princeton, NJ 08544 USA (e-mail: poor@princeton.edu).}
}

\maketitle
 \addtolength{\textfloatsep}{-5mm}
 \addtolength{\dblfloatsep}{-10mm}
 \addtolength{\belowcaptionskip}{-5mm}

\begin{abstract}
This paper studies the general problem of operating energy storage under uncertainty. 
%The kind of storage under consideration is generalized, that is, any form of power system devices and entities that can be suitably modelled and controlled using a discrete-time storage model defined rigorously in this paper. Noteworthy examples include conventional energy storage such as batteries and fly-wheels, and virtual storage such as pre-emptive deferrable loads and aggregations of thermostatistically controlled loads. 
Two fundamental sources of uncertainty are considered, namely the uncertainty in the unexpected fluctuation of the net demand process and the uncertainty in the locational marginal prices. We propose a very simple algorithm termed Online Modified Greedy (OMG) algorithm for this problem. A stylized analysis for the algorithm is performed, which shows that comparing to the optimal cost of the corresponding stochastic control problem, the sub-optimality of OMG is controlled by an easily computable bound. This suggests that, albeit simple, OMG is guaranteed to have good performance in cases when the bound is small. Meanwhile, OMG together with the sub-optimality bound can be used  to provide a lower bound for the optimal cost. Such a lower bound can be valuable in evaluating other heuristic algorithms. For the latter cases, a semidefinite program is derived to minimize the sub-optimality bound of OMG. Numerical experiments are conducted to verify our theoretical analysis and to demonstrate the use of the algorithm.
\end{abstract}


\input{introduction-conf.tex}
\input{problemFormulation-conf.tex}

\input{algorithm-conf.tex}
\input{numericalExample-conf.tex}

%\input{conclusion-conf.tex}
%\input{singleStorageProofDraft.tex}

%\section*{Acknowledgments}
%This work was supported in part by TomKat Center for Sustainable Energy,  in part by a Satre Family Fellowship, and in part by a Powell Foundation Fellowship. The authors wish to thank Walter Murray for useful comments.

\bibliography{jqin}
\bibliographystyle{IEEEtran}
%\clearpage
\appendices

\input{appendix_1bus.tex}
\input{appendix_markov_process.tex}

\end{document}

%% file: defs.tex
\usepackage{amsfonts}
\usepackage{standalone}
\usepackage{amsmath}
\usepackage{fix-cm}
\usepackage[T1]{fontenc}

\usepackage{bm}
\usepackage{psfrag}
\usepackage{pgf}
\usepackage{enumerate}
\usepackage{amssymb}
\usepackage{bbm}
\usepackage{color}
\usepackage{algorithm}
\usepackage{algpseudocode}
\usepackage{subfigure}
\usepackage{url}
\def\ie{{\it i.e.}}
\def\eg{{\it e.g.}}
\algnewcommand\algorithmicinput{\textbf{Input:}}
\algnewcommand\INPUT{\item[\algorithmicinput]}
\algnewcommand\Offline{\item[\textbf{Offline-Phase:}]}
\algnewcommand\Online{\item[\textbf{Online-Phase:}]}
\def\nn{\nonumber}
\def\expec{\mathbbm{E}}

\def\defeq{\triangleq}

\newcommand{\st}{\textup{subject to}}
\DeclareMathOperator{\minimize}{\text{\textnormal{minimize}}}

\newtheorem{assumption}{Assumption}
\newtheorem{theorem}{Theorem}
\newtheorem{lemma}{Lemma}

\newtheorem{remark}{Remark}
\newtheorem{definition}{Definition}

\newtheorem{example}{Example}

\def\b{\vv{b}}
\def\u{\vv{u}}
\def\v{\vv{v}}

\def\f{\vv{f}}
\def\q{\vv{q}}
\def\p{\vv{p}}
\def\d{\vv{d}}
\def\h{\vv{h}}

\def\b{s}
\def\u{u}

\def\f{f}

\def\bmax{S^{\max}}
\def\bmin{S^{\min}}
\def\bdot{S^{(\cdot)}}
\def\umax{U^{\max}}
\def\umin{U^{\min}}
\def\udot{U^{(\cdot)}}

\def\upos{\u^{+}}
\def\uneg{\u^{-}}
\def\p{p}
\def\q{q}
\def\d{\delta}
\def\dr{\d^\mathrm{R}}
\def\la{\lambda}
\def\muC{\mu^\mathrm{C}}
\def\muD{\mu^\mathrm{D}}

\newcommand{\pos}[1]{\left(#1\right)^+}
\renewcommand{\neg}[1]{\left(#1\right)^-} %x = \pos{x} - \neg{x}
\def\minimize{\mbox{minimize  }}
\def\st{\mbox{subject to  }}

\def\h{h}
\def\hC{\h^\mathrm{C}}
\def\hD{\h^\mathrm{D}}

\def\g{g}
\def\MWh{\text{MWh}}
\def\MW{\text{MW}}
\def\hr{\text{h}}
\newcommand{\opttag}[1]{\mbox{\bf#1  }}
\def\F{\mathcal{F}}

\def\dt{\Delta t}
\def\pmin{\p^{\min}}
\def\pmax{\p^{\max}}
\def\J{J}
\newcommand{\Jk}[1]{\J_\mathrm{P#1}}
\newcommand{\Jkt}[2]{\J_{\mathrm{P#1},#2}}

\newcommand{\Pk}[1]{\text{\bf P#1}}

\def\ustat{\u^{\mathrm{stat}}}
\def\fstat{\f^{\mathrm{stat}}}

\def\upistat{\u^{\pi,\mathrm{stat}}}
\def\fpistat{\f^{\pi,\mathrm{stat}}}

\def\vpi{\v^\pi}
\def\upi{\u^\pi}
\def\fpi{\f^\pi}

\def\vpistar{\v^{\pi,\star}}
\def\upistar{\u^{\pi,\star}}
\def\fpistar{\f^{\pi,\star}}

\def\bs{\widehat \b}
\def\ks{\Gamma}
\def\ksmin{\ks^{\min}}
\def\ksmax{\ks^{\max}}
\def\W{W}
\def\Wmax{\W^{\max}}
\def\v{v}

\def\vpistat{\v^{\pi,\mathrm{stat}}}
\def\vstat{\v^{\mathrm{stat}}}

\def\uhat{\u^{\mathrm{ol}}}
\def\vhat{\v^{\mathrm{ol}}}

\def\upihat{\u^{\pi,\mathrm{ol}}}

\def\fpihat{\f^{\pi,\mathrm{ol}}}

\def\Dl{\underline{D}}
\def\Du{\overline{D}}
\def\M{M}
\def\Mone{\M^{\u}}
\def\Mtwo{\M^{\b}}

\def\Mfour{\mathcal M}
\def\N{\eta}
\def\None{\N^{\u}}
\def\Ntwo{\N^{\b}}
\def\L{L}
\def\Ls{\Delta}

\def\DT{\Delta T}
\def\R{N_{\omega_0}}

\def\dmin{\d^{\min}}
\def\dmax{\d^{\max}}
\def\X{X}
\def\Xumin{\X^{\min,\u}}
\def\Xumax{\X^{\max,\u}}
\def\Xbmin{\X^{\min,\b}}
\def\Xbmax{\X^{\max,\b}}
\def\Xudot{\X^{(\cdot),\u}}

\def\Xbdot{\X^{(\cdot),\b}}

\newcommand{\ylmod}[1]{{#1}}

\def\Tday{\mathcal{T}^\mathrm{Day}}

\def\prl{\tilde \p}
\def\drl{\tilde \d}
\def\grl{\tilde \g}
\def\Sm{\mathbb{S}}

%% file: introduction-conf.tex
\section{Introduction}
%To ensure a sustainable energy future, deep penetration of renewable energy generation is essential. Renewable energy resources, such as wind and solar, are intrinsically variable. Uncertainties associated with these intermittent and volatile resources pose a significant challenge to their integration into the existing grid infrastructure \cite{NRELWest2010}. More flexibility, especially in shifting energy supply and/or demand across time and network, are desired to cope with the increased uncertainties.

Energy storage provides the functionality of shifting energy across time. A vast array of technologies, such as  batteries, flywheels, pumped-hydro, and compressed air energy storages, are available for such a purpose \cite{Denholm2010}.
 %and have been improving in terms of both their technical performance and reduced costs 
%Albeit concerns still exist about the capital and installation costs of large scale energy storage systems, certain amount of storage is viewed as necessary to facilitate renewable integration. Due to these considerations, new market rules have been developed for a class of energy storage resources by NYISO (the New York Independent System Operator) \cite{StorNYISO2010}, and aggressive targets for storage installation have been set recently by the state of California \cite{caStorTarget}.
Furthermore, flexible or controllable demand provides another ubiquitous source of storage. Deferrable loads -- including many thermal loads, loads of internet data-centers and loads corresponding to charging electric vehicles (EVs) over certain time intervals \cite{thermalStor1993, ThatteXieStorValue2012s} -- can be interpreted as \emph{storage of demand} \cite{ObRACC2013}. % in contrast to conventional \emph{storage of energy}, based on the equivalence that storing into the storage of energy defers the use of the corresponding energy, while storing into the storage of demand defers the fulfillment of the corresponding demand \cite{ObRACC2013}.
Other controllable loads which can possibly be shifted to an earlier or later time, such as thermostatically controlled loads (TCLs), may be modeled and controlled as a storage with negative lower bound and positive upper bound on the storage level \cite{Callaway2009, HaoSanandajiPoollaVincent2013}.
%\ylmod{are categorized as \emph{generalized storage}. Mathematically, the dynamics of generalized storage has the same form as the dynamics of storage of energy, but with different storage parameters.}%can be viewed as a kind of \emph{generalized storage},
%with the same mathematical form for the dynamics as storage of energy, yet different specifications for some of the storage parameters \cite{Callaway2009, HaoSanandajiPoollaVincent2013}.
These forms of storage enable inter-temporal shifting of excess energy supply and/or demand, and significantly reduce the reserve requirement and thus system costs. %\footnote{We focus on reserve energy use assuming the reserve capacity has been purchased in advance.} that is prepared for the uncertainties in renewable generation. This further lowers the economic cost and greenhouse gas emission in electric grid operations. }%due to the uncertainty introduced by renewable generation,  thus lowering both the economic costs and greenhouse gas emission of operating the electric grid.

The problem of optimal storage operation under various sources of uncertainty remains challenging.
Two categories of approaches have been proposed in the literature. The first category is based on exploiting structures of specific problem instances, usually using dynamic programming (DP). These structural results are valuable in providing insights about the system, and often lead to analytical solution of these problem instances. For example, analytical solutions to optimal storage arbitrage with stochastic prices have been derived in  \cite{QinPESGM12} without storage ramping constraints, and in \cite{MITrampStor} with ramping constraints. Problems of using energy storage to minimize energy imbalance are studied in various contexts; see \cite{SuEGTPS, RLDSACC} for reducing reserve energy requirements in power system dispatch, \cite{BitarRACC_colocated, Powell} for operating storage co-located with a wind farm, \cite{IBMload, DataCenter} for operating storage co-located with end-user demands, and \cite{StorDRLongbo} for storage with demand response.
However, such approaches rely heavily on specific assumptions of the type of storage, the form of the cost function, and the distribution of uncertain parameters. Generalizing analytical results to other specifications and more complex settings is usually difficult. 

In many cases, DP can also lead to efficient computational methods, notably algorithms based on value iteration, policy iteration or linear programming. For storage operation problems, as the state space, action space and disturbance space are all continuous,  approximations based on discretization \cite{chow1991optimal} or simulation \cite{Jain:2010:SOM:1836310.1836495} are needed. Although error bounds are available for these approximations, the computational cost of these methods usually grow exponentially with the dimensionality of the problem. This phenomenon, known as \emph{curses of dimensionality}, makes DP based computational methods not well suited for some instances of the storage control problems. More importantly, implementing DP based approaches requires full information of the probability distribution of the stochastic parameters, which may not be readily available. 

%and consequently this approach is mostly used to analyze single storage systems.
%Most of these studies (i.e. \cite{SuEGTPS, RLDSACC, BitarRACC_colocated, Powell, IBMload}) used stochastic dynamic programming to study problem structures or to derive analytical solutions. 

The other category is using heuristic algorithms, such as Model Predictive Control (MPC) \cite{XieEtAlWindStorMPC} and look-ahead policies \cite{NRELStorValue2013},  to identify sub-optimal storage control rules. Usually based on deterministic (convex) optimization, these approaches can be easily applied to general networks. The major drawback is that these approaches usually do not have any performance guarantee. Consequently,  it lacks theoretical justification for implementing them in real systems.
Examples of this category can be found in \cite{XieEtAlWindStorMPC} and references therein.
%A Model Predictive Control (MPC) \ylmod{approach} % based approach
% for controlling battery energy storage systems was developed in \cite{XieEtAlWindStorMPC}, where the performance was validated only via numerical simulation.

%Energy storage operation in the electric network under uncertainties is less well studied. With quadratic cost functions and simplified system constraints, optimal affine policies were derived in \cite{ZTMStorAllerton} without inequality constraints and in \cite{QRPESNetStor} with inequality constraints.

%The notion of using controllable loads as a form of virtual storage is not new, yet analytical characterization of the equivalent storage model and the use of it as a proxy to derive algorithms for controlling such loads are recent. A storage interpretation of deferrable loads and a method to calculate the corresponding storage capacity were given in \cite{ObRACC2013}. Modeling the power demand from a population of TCLs as storage was discussed in \cite{Callaway2009}, whose model was then encapsulated by the notion of generalized battery models in \cite{HaoSanandajiPoollaVincent2013}.

This work aims at designing online deterministic optimizations that solve the stochastic control problem with provable guarantees. It contributes to the existing literature in the following ways. 
First, we formalize the notion of {\it generalized storage} as a dynamic model that captures a variety of power system components which provide the functionality of storage. 
Second, we formulate the problem of  optimal storage operation under uncertainty as a stochastic control problem with general cost functions, and provide examples of applications that can be encapsulated by such a formulation. 
Third, we develop an online modified greedy (OMG) algorithm for this problem,
 %\rev{\footnote{\rev{Although closely related to the classical Lyapunov theory} \rev{for stability, the theory and techniques of Lyapunov optimization are relatively recent. See \cite{NeelyBook} for more details.}}}
and derive performance guarantees in the form of sub-optimality bounds for the algorithm. The OMG algorithm is very simple as it only requires solving a deterministic optimization in each step, and it  needs a very little amount of information regarding the probability distribution of the stochastic parameters. 
The sub-optimality bounds are not only of theoretical interests, but also suggests the use of OMG in many cases where accurate methods such as those based on DP are not applicable. 
Furthermore, these bounds are useful in evaluating the performance of other sub-optimal algorithms when the optimal costs are difficult to compute. They can also be used to estimate the maximum cost reduction that can be achieved by {\it any} storage control policies, thus provides understandings for the limit of a certain storage system.
To the best of our knowledge, this is the first algorithm with provable guarantees for the \emph{general} storage operation problem with both stochastic price and demand.  

Our methodology is built upon on the theory of Lyapunov optimization \cite{NeelyBook}, which was developed for queueing networks and has been applied to the context of energy storage control in recent work including \cite{DataCenter}, \cite{StorDRLongbo}, \cite{lyap1} and \cite{lyap2}. Different from these work, which analyze specific setups for storage operation, we aim to provide a general framework where the storage, co-located with a controllable resource and any stochastic uncontrollable resource, can be operated to minimize an arbitrary convex cost function. To achieve this goal, we have introduced a much more general storage model which i) captures energy dissipation over time, ii)  requires minimal assumptions in terms of the storage parameters to model \eg, storage of demand and TCLs, and iii) allows charging and discharging energy losses. 
In contrast, most of the existing work analyzes ideal energy storages without any of the above features, with the exception that \cite{lyap2} models charging energy losses. Modeling these features leads to a different online program, requires a new analysis for the algorithm, and results in different sub-optimality bounds. In particular, the new bounds developed in this paper scale very differently with the storage capacity compared to the bounds appeared in the prior work since we have captured the effect that large storage can lose more energy due to energy dissipation.
Preliminary results related to this paper appeared in \cite{QCYR:acm}.
This paper significantly generalizes \cite{QCYR:acm} by modeling additional controllable devices connected to the bus, dealing with general convex cost functions instead of piecewise linear costs, developing examples and analytical solutions for the online program to facilitate implementation, and conducting new case studies.
%Among many differences between  \cite{QCYR:acm} and this paper, the journal version models additional controllable devices connected to the bus, deals with general convex cost functions instead of piecewise linear costs, develop examples and analytical solutions for the online program to facilitate implementation, and conducts realistic numerical examples.

%Our methodology is closely related to that of \cite{DataCenter}, where the focus is on solving the problem of operating an idealized energy storage (with no energy dissipation over time, and no charging/discharging conversion loss) at data-centers. Our objective is to provide an algorithm to operate generalized storage network in a wide range of different settings. This requires an extended or a new analysis in the following aspects. From the modeling perspective, in order to capture applications such as deferrable loads and TCLs, we do not assume storage level is non-negative, instead, we only assume each storage is feasible (see Assumption~\ref{assume:feas} for more details).  Furthermore, modeling the dissipation of energy over time leads to a new sub-optimality bound; the bound in \cite{DataCenter} becomes a special case of our bound when the dissipation factor (or storage efficiency) is one.  A semidefinite program  is constructed to decide parameters of the algorithm in order to minimize the sub-optimality bound. Finally, the aspect of power network appears to be completely new.  

The rest of the paper is organized as follows. Section 2 formulates the problem of operating a generalized storage under uncertainty. Section 3 gives the online algorithm and states the performance guarantee. %Section 4 analyzes the single bus case in detail with a generalized storage, and 
%Section 5 provides a summary of results for general storage networks. 
Numerical examples are then given in Section 4. Section 5 concludes the paper.

%% file: problemFormulation-conf.tex
\section{Problem Formulation}\label{sec:problem}
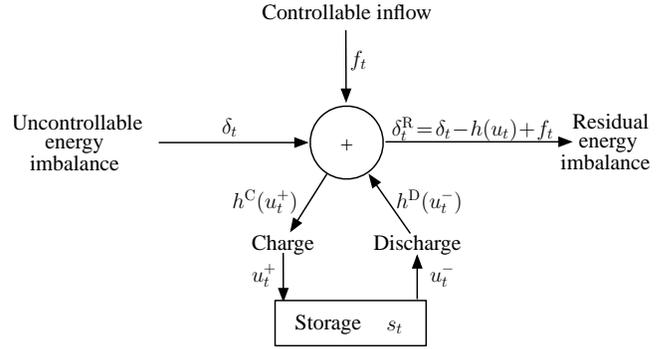
\begin{figure}[htbp]
\centering
\scalebox{0.55}{\scalefont{1.5} \input{./fig/diagram-singleStroEPS3.pgf}}
\caption{System diagram.}
\label{fig:diagSingle}
\end{figure}

Working with slotted time, we use $t$ as the index for an arbitrary time period and denote the constant length of each time period by $\dt$.
Using $\dt$, we can convert from power units (\eg, MW) to energy units (\eg, MWh) and vice versa with ease.\footnote{We work with real power in this paper. Incorporating reactive power and more detailed power flow model with storage is an important future direction.} 
%As $\dt$ may not be convenient numbers such as one hour or one minute,
For convenience and assuming a proper conversion, we work with energy units in this paper, albeit many power system quantities are conventionally specified in power units.
The system diagram is depicted in Figure~\ref{fig:diagSingle}.
%\begin{figure}[htbp]
%\centerline{
%\includegraphics[width = 0.40\textwidth]{../fig/diagSingleStor}}
%\caption{Diagram of a single-bus storage system. (TODO: remove or update network inflow (note that the network inflow can be positive or negative); another way will be to keep the network inflow as a decision variable--purchase and sell to the bulk grid; perhaps a more vivid representation of the cause of energy imbalance.)}
%\label{fig:diagSingle}
%\end{figure}

\subsection{Generalized Storage}\label{sec:prob:singleBusModel}
We start by describing a \emph{generalized storage} model, which is specified by the following elements:
\begin{itemize}%[leftmargin=0.15in]
\item The \emph{storage level} or State of Charge (SoC) $\b_t$ summarizes the status  of the storage at time period $t$. If $\b_t \ge 0$, it represents the amount of energy in storage; if $\b_t\le 0$, $-\b_t$ can represent the amount of currently deferred (and not fulfilled) demand.  It satisfies $\b _t \in [\bmin, \bmax]$, where $\bmax$ is the storage capacity, and $\bmin$ is the minimum allowed storage level. %It will be made clear later that we can assume without loss of generality $\bmin \ge 0$.
\item The \emph{storage operation} $\u_t$ summarizes the charging (when $\u_t\ge 0$) and discharging (when $u_t\le 0$) operations of the storage. It satisfies charging and discharging ramping constraints, \ie, $\u_t \in [\umin, \umax]$, where $\umin (\le 0)$ is the negation of the maximum discharge within each time period, and $\umax (\ge 0)$ is the maximum charge within each time period. We also use $\upos_t = \max(u_t,0)$ and $\uneg_t = \max(-u_t,0)$ to denote the charging and discharging operations, respectively.
\item The \emph{storage conversion function} $\h$ maps the storage operation $\u_t$ into its effect on the bus. In particular, it is composed of two functions, namely the \emph{charging conversion function} $\hC$, and the \emph{discharging conversion function} $\hD$, such that  $\hC(\upos_t)$ is the amount of energy  drawn from the bus due to $\upos_t$ amount of charge, and $\hD(\uneg_t)$ is the amount of energy that is injected into the bus due to $\uneg_t$ amount of discharge, whence
\[\h(\u_t) \defeq \hC(\upos_t)- \hD(\uneg_t)\] is the energy drawn from the bus by the storage.
\item The \emph{storage dynamics} is then
\begin{equation}\label{eq:storDyn}
\b_{t+1} = \la \b_t + \u _t,
\end{equation}
where $\la\in (0,1]$ is the \emph{storage efficiency} which models the loss over time even if there is no storage operation.
\end{itemize}
We provide the definition of a generalized storage as follows.
\begin{definition}\label{assume:feas}
For $t= 1, 2, \dots$, the controlled dynamic system with state $\b_t \in [\bmin, \bmax]$, control $\u_t \in [\umin, \umax]$, and dynamics $\b_{t+1} = \la \b_t + \u_t$ is deemed a \emph{generalized storage model} if
the set of parameters $\Sm = \{\la, \bmin, \bmax, \umin, \umax\}$ satisfies the following conditions:
\begin{itemize}
\item (feasibility) $\la \bmin + \umax \ge \bmin$ and $\la \bmax + \umin \le \bmax$;
\item (controllability) $\la \bmax + \umax \ge \bmax$ and $\la \bmin + \umin \le \bmin$.
\end{itemize}
In addition, the effect of the storage operation on the bus is captured by the conversion function $h$.
\end{definition}

The feasibility and controllability conditions can be interpreted as follows. Feasibility means that
starting from any feasible storage level, there exists a feasible storage operation such that the storage level in the next time period is feasible. Every storage system must satisfy the feasibility condition.
 Controllability requires that starting from any feasible storage level, there exists a sequence of feasible storage operations to reach any feasible storage level in a finite number of time periods. The linear nature of the dynamics \eqref{eq:storDyn} reduces the controllability requirements to the inequalities shown in Definition~\ref{assume:feas}, which hold for all practical storage systems except for pathological cases. Apparently, controllability implies feasibility. It will become clear that the feasibility condition is crucial in proving various results in this paper; it is often used in place of the positive storage level condition which does not hold for generalized storage models. The controllability condition is mostly introduced to simplify the presentation; see \cite{QCYR:report} for more details regarding how to relax it.

A few examples of generalized storage models are provided below.
\begin{example}[Storage of Energy]\label{eg:soe}
Storage of energy can be modeled as a generalized storage with  $\bmax \ge \bmin \ge 0$. Here $\umin$ and $\umax$ correspond to the power rating of the storage, up to a multiple of the length of each time period $\dt$.  By setting $\hC(\upos_t) = (1/\muC)\upos_t$, and $\hD(\uneg_t) = \muD \uneg_t$, one models the energy loss during charging and discharging operations. Here $\muC\in (0,1]$ is the charging efficiency; $\muD\in (0,1]$ is the discharging efficiency; and the round-trip efficiency of the energy storage is $\muC \muD$. For instance, based on the information from \cite{QinPESGM12}, a sodium sulfur (NaS) battery
and a compressed air energy storage (CAES) can be modeled with parameters shown in Table~\ref{tab:soe}.
\begin{table}[htbp]\label{tab:soe}
 \caption{Parameters for energy storage in Example 1. Here $\dt = 1$\textnormal{h}, $\umin = -\umax$, and $\muD = \muC$. }
 \center
\begin{tabular}{lccccccccccc}\hline
 & $\bmin$  & $\bmax$ & $\umax$ &$\muC$  & $\la$ \\
  \hline
NaS & $0 \MWh$  & $100 \MWh$ & $10 \MW \cdot 1 \hr$  &$0.85$  & $0.97$ \\
CAES & $0 \MWh$  & $3000 \MWh$ & $300 \MW \cdot 1 \hr$  &$0.85$  & $1.00$ \\\hline
\end{tabular}
\end{table}
\end{example}
\begin{example}[Storage of Demand]\label{eg:sod}
Pre-emptive deferrable loads may be modeled as storage of demand, with $-\b_t$ corresponding to the accumulated deferred (but not yet fulfilled) load up to time $t$ , and with $\u_t$ corresponding to the amount of load to defer/fulfill in time period $t$. We have $\bmin \le \bmax \le 0$ in this case. Storage of demand differs from storage of energy in the sense that it has to be discharged before charging is allowed.
%Setting $\hC(\upos_t) = \upos_t$ and $\hD(\uneg_t) = \uneg_t$ converts the charging/discharging of load to injection/withdraw of energy.
The conversion function can usually be set to $\h(\u_t) = \u_t$,
and generally $\la = 1$ in deferrable load related applications.
\end{example}
\begin{example}[Battery Model for Aggregation of TCLs]\label{eg:gbm}
It is shown recently that an aggregation of TCLs may be modeled as a generalized battery \cite{HaoSanandajiPoollaVincent2013}. With a linear approximation, a discrete time version of such a model can be cast into our framework by setting $\bmax \ge 0$ representing the maximum amount of virtual energy storage that can be obtained by pre-cooling without affecting the comfort level of the users. By a symmetric argument, $\bmin = - \bmax$ . Other storage parameters can be set properly according to Definition 1 of \cite{HaoSanandajiPoollaVincent2013}, and we have $\la \le 1$ to model energy dissipation.
\end{example}

\subsection{System Model and Cost Functions}
The generalized storage is connected to a bus together with several other system components.
For time period $t$, the local \emph{uncontrollable energy imbalance}, denoted by $\d_t$, is defined to be the difference between the uncontrollable local generation, such as energy generated by solar panel or priorly dispatched generators, and the demand. The sign convention is such that $\d_t \le 0$ ($\d_t>0$) represents a net demand (supply) at the bus.
Due to the limited predictability, both the local generation and demand can be stochastic, and therefore $\d_t$ is stochastic in general. The bus could be connected to another controllable component/device such as a standby generator or motor, from (to) which the energy inflow (outflow) is denoted by $f_t \ge 0$ ($f_t <0$) and we have $f_t \in \F$ for all $t$ where $\F$ is a convex and compact set.

%[OPTIONAL]
%
%The following more stringent assumptions are optional for our model and not essential for establishing the performance of our algorithm. Yet under these assumption, our analysis is tight in certain sense.

%Notice that storage variables $\b_t$ and $\u_t$ are described in terms of the stored object that can possibly be energy or demand. However, the suitably defined storage conversion functions map the corresponding charge or discharge into quantities in terms of energy,
The \emph{residual energy imbalance}, after accounting for the controllable inflow and storage operation, is then given by:
\begin{equation}\label{eq:reEI}
\dr_t \defeq\d_t - \h(\u_t) + f_t = \d_t - \hC(\upos_t) + \hD(\uneg_t) + f_t,
\end{equation}
which represents the overall output of the sub-system under consideration. Such energy imbalance may be matched by energy inflow/outflow from the main grid, at certain cost. Let
\begin{equation}
\g_t\defeq g_t(\u_t, f_t, \d_t, \p_t)
\end{equation}
be a convex cost function\footnote{Report \cite{QCYR:report} discusses how and to what extent the convexity requirement can be reduced.} for time period $t$, where $\p_t$ is a stochastic price parameter modeling for example the locational marginal price (LMP) at the bus. Different functional forms of $\g_t$ encode different uses of the storage. We provide the functional forms of $\g_t$ for the two fundamental use cases of the storage, namely, 
to exploit the inter-temporal differences in prices and to balance the unexpected fluctuations in net demand across time periods. We also provide another example where these two effects are somewhat combined.

\def\pl{\p^\mathrm{LMP}}
\begin{example}[Arbitrage]\label{eg:arbi}
%Another important application of storage is to arbitrage in the electricity markets.
Third-party owned storage devices may be used to arbitrage price variations in the electricity spot market.
Consider the case that the bus is only connected to a storage, \ie, $\d_t = 0$ and $f_t = 0$.
For arbitrage purpose and given a stochastic sequence of locational marginal prices $\{\p_t: t\ge 1\}$, the following cost function may be used%The problem of maximizing the expected arbitrage profit using storage operations can be cast as an instance of \eqref{prob:singleBusGeneral}, with the cost function (\ie, negative profit) given by:
\begin{equation}\label{eq:arbicost}
\g_t = -\p_t \dr_t = \p_t (\hC(\upos_t) - \hD(\uneg_t)),
\end{equation}
to characterize the negation of the stage-wise profit earned by storage operations.
\end{example}
\begin{example}[Balancing/Regulation]\label{eg:ba}
Storage may be used by the system operator or ancillary service providers 
to minimize residual energy imbalance  given by some stochastic net demand $\{\d_t: t\ge 1\}$ process. Typical cost functions penalize the positive and negative residual energy imbalance differently, and may have different penalties at different time periods, \eg, to model the different consequences of load shedding at different times of each day. The problem of optimal storage control for such a purpose can be modeled by problem~\eqref{prob:singleBusGeneral} with the cost function
\begin{equation} \label{eq:bacost}
\g_t = q^+_t \pos{\dr_t} + q^-_t \neg{\dr_t},
\end{equation}
where $q^+_t$ and $q^-_t$ are the penalties for each unit of positive and negative residual energy imbalance at time period $t$, respectively. %This cost is clearly a special case of \eqref{eq:single BusGeneralCost}.
%with $L=2$, $\p(t,1) = q^+_t$, $\p(t,2) = q^-_t$, $\alI(1) = 1$, $\alI(2) = -1$, $$......
\end{example}

\begin{example}[Storage Co-Located with Stochastic Generation] \label{eg:co} 
For storage co-located with a wind farm or an end-user, 
it can be the case that both the net energy imbalances and the prices are stochastic. %A set of notable examples include operating storage co-located with stochastic renewable generation such as wind while facing a stochastic locational marginal price sequence.
Applications of this type can be cast into our framework using $\{\d_t: t\ge 1\}$ to model the stochastic generation or demand process, and $\{\p_t: t\ge 1\}$  to model the stochastic prices. A possible cost function is
\begin{equation}\label{eq:cocost}
\g_t = \p_t \neg{\dr_t},
\end{equation}
where the excessive supplied energy is curtailed with no cost/benefit, and the excessive demand is supplied via buying energy from the market at stochastic price $\p_t$.
%This class of problems typically involves the optimization of both the storage operation and \ylmod{certain contract (or commitment) level adjustment} for the stochastic generation or demand. The sub-problem of storage operation (for fixed contract level) can usually be posed as \eqref{prob:singleBusGeneral}, \ylmod{where stochastic process $\d_t$ models the deviation of the stochastic generation (demand) from the contract level, and stochastic processes $\p(t,\ell)$, model the locational marginal prices. }
\end{example}

%We consider the following \emph{stochastic piecewise linear cost function} for each fixed bus
%\begin{align}\label{eq:singleBusGeneralCost}
%\g_t &= \sum_{\ell = 1}^L \p(t, \ell) \Big(\alI(\ell) \d_t - \alC(\ell) \hC(\upos_t) \\
%&\quad + \alD(\ell) \hD(\uneg_t) + \alF(\ell) f_t + \alConst(t,\ell)\Big)^+, \nn
%\end{align}
%where the parameter $\p(t, \ell)$ is in general stochastic, and follows a  prescribed probability law,  and $\alI(\ell)$, $\alC(\ell)$, $\alD(\ell)$, $\alF(\ell)$ and $\alConst(t, \ell)$ are constants, for each $\ell = 1, \dots, L$ and $t$. This cost function serves as a generalization of positive (and/or negative) part cost function of the residual energy imbalance, and it encapsulates many applications of storage as shown in Section~\ref{sec:prob:applications}. Our analysis applies to a more general class of cost functions; see Appendix~\ref{appendix:convexCost} for more details.
%

\subsection{Optimal Storage Operation Problem}\label{sec:prob:applications}
In case that all the stochastic parameters are known ahead of time, the optimization of the storage operation (possibly together with the controllable inflow) can be written as
\begin{subequations}\label{prob:singleBusDeterministic}
\begin{align}
\minimize \quad &  (1/T)  \sum_{t=1}^T \g_t \\
\st \quad & \b_{t+1} = \la \b_t + \u_t, \label{eq:ca}\\
        & \bmin \le \b_t \le \bmax, \label{eq:cb}\\
        & \umin \le \u_t \le \umax,\label{eq:cc}\\
        & f_t \in \F, \label{eq:cd}
\end{align}
\end{subequations}
where the optimization variables are $\u_t$ and $f_t$ for $t= 1, \dots T$, and the initial state $\b_1 \in [\bmin, \bmax]$ has an arbitrary given value. In the formulation above,  $T$ is 
the number of time periods that is considered for the storage operation problem. Although engineering practices often use a $T$ that corresponds to a relatively short time period (\eg, solving the problem for each week or month with the storage being operated every 5 minute to 1 hour), it leads to a loss of optimality, \ie, increased system cost, by using a $T$ that is less than the  \emph{decision horizon} \cite{DecisionHorizon1988} of the problem. Here the decision horizon, roughly speaking, is a $T$ such that the information in stage $T+1$ would not affect the optimal solution of the problem in the first $T$ stages. 
Since calculating the exact decision horizon under stochastic settings is not always possible, using a larger $T$ is usually preferable. 
%Using a decision horizon shorter than the actual total number of periods derived from the life span of the storage may simplify the problem; however under stochastic settings, identifying such a $T$ is challenging. Thus we work with $T$ derived from the life span of the storage under consideration. For example, for a battery with 3 year life span and operated each hour, we expected $T$ to be on the order of $10^4$.

Due to the fact that $\g_t$ depends on stochastic parameters $\d_t$ and $\p_t$ whose realizations are not known ahead of time, problem~\eqref{prob:singleBusDeterministic} is not well defined. In a \emph{risk neutral} setting, one may instead solve
\begin{subequations}\label{prob:singleBusGeneral}
\begin{align}
\minimize \quad &  (1/T) \expec \Big[ \sum_{t=1}^T \g_t \Big]\\
\st \quad & \eqref{eq:ca}, \eqref{eq:cb}, \eqref{eq:cc}, \eqref{eq:cd},
\end{align}
\end{subequations}
where the expectation is taken over the possible realizations of $\d_t$ and $\p_t$ for $t = 1, \dots, T$, and the goal is to identify optimal policies which are functions that map information available at stage $t$ to the optimal actions $\u_t$ and $\f_t$\footnote{\ylmod{Notation: In this paper, we denote control policies and actions $(\u_t,f_t)$ with the same set of variables. To differentiate, we use $(\upi_t,\fpi_t)$ to denote the corresponding control policy that induces action $(\u_t,\f_t)$ at time $t$.}}. The following challenges must be resolved in order to derive a practical algorithm for problem formulation~\eqref{prob:singleBusGeneral}. (i) Probability distributions of $\d_t$ and $\p_t$ are required for evaluating the objective function. This requires probabilistic forecasts for a long horizon, which often is practically infeasible. (ii) The exact offline optimal solution of problem~\eqref{prob:singleBusGeneral} is characterized by the Bellman's recursion \cite{bertsekas2007dynamic}, which is computationally intractable for problems with continuous variables such as \eqref{prob:singleBusGeneral}. No general solution exists for the aforementioned challenges; thus certain approximations are necessary. Usually, one has to seek a good tradeoff between the simplicity and the performance of the algorithm. In the remaining of this paper, we provide a very simple algorithm that has provable performance guarantees.

%% file: fig/diagram-singleStroEPS3.pgf
% Created by Eps2pgf 0.7.0 (build on 2008-08-24) on Sun Apr 27 02:47:25 PDT 2014
\begin{pgfpicture}
\pgfpathmoveto{\pgfqpoint{0cm}{0cm}}
\pgfpathlineto{\pgfqpoint{15.416cm}{0cm}}
\pgfpathlineto{\pgfqpoint{15.416cm}{8.431cm}}
\pgfpathlineto{\pgfqpoint{0cm}{8.431cm}}
\pgfpathclose
\pgfusepath{clip}
\begin{pgfscope}
\end{pgfscope}
\begin{pgfscope}
\begin{pgfscope}
\pgfpathmoveto{\pgfqpoint{0cm}{0cm}}
\pgfpathlineto{\pgfqpoint{15.416cm}{0cm}}
\pgfpathlineto{\pgfqpoint{15.416cm}{8.431cm}}
\pgfpathlineto{\pgfqpoint{0cm}{8.431cm}}
\pgfpathclose
\pgfusepath{clip}
\definecolor{eps2pgf_color}{gray}{1}\pgfsetstrokecolor{eps2pgf_color}\pgfsetfillcolor{eps2pgf_color}
\pgfpathmoveto{\pgfqpoint{-0.053cm}{15.117cm}}
\pgfpathlineto{\pgfqpoint{20.267cm}{15.117cm}}
\pgfpathlineto{\pgfqpoint{20.267cm}{-10.742cm}}
\pgfpathlineto{\pgfqpoint{-0.053cm}{-10.742cm}}
\pgfpathclose
\pgfusepath{fill}
\pgfpathmoveto{\pgfqpoint{8.667cm}{5.633cm}}
\pgfpathcurveto{\pgfqpoint{9.011cm}{5.289cm}}{\pgfqpoint{9.011cm}{4.73cm}}{\pgfqpoint{8.667cm}{4.386cm}}
\pgfpathcurveto{\pgfqpoint{8.323cm}{4.041cm}}{\pgfqpoint{7.764cm}{4.041cm}}{\pgfqpoint{7.42cm}{4.386cm}}
\pgfpathcurveto{\pgfqpoint{7.075cm}{4.73cm}}{\pgfqpoint{7.075cm}{5.289cm}}{\pgfqpoint{7.42cm}{5.633cm}}
\pgfpathcurveto{\pgfqpoint{7.764cm}{5.977cm}}{\pgfqpoint{8.323cm}{5.977cm}}{\pgfqpoint{8.667cm}{5.633cm}}
\pgfusepath{fill}
\pgfsetdash{}{0cm}
\pgfsetlinewidth{0.353mm}
\pgfsetroundcap
\pgfsetroundjoin
\definecolor{eps2pgf_color}{rgb}{0,0,0}\pgfsetstrokecolor{eps2pgf_color}\pgfsetfillcolor{eps2pgf_color}
\pgfpathmoveto{\pgfqpoint{8.667cm}{5.633cm}}
\pgfpathcurveto{\pgfqpoint{9.011cm}{5.289cm}}{\pgfqpoint{9.011cm}{4.73cm}}{\pgfqpoint{8.667cm}{4.386cm}}
\pgfpathcurveto{\pgfqpoint{8.323cm}{4.041cm}}{\pgfqpoint{7.764cm}{4.041cm}}{\pgfqpoint{7.42cm}{4.386cm}}
\pgfpathcurveto{\pgfqpoint{7.075cm}{4.73cm}}{\pgfqpoint{7.075cm}{5.289cm}}{\pgfqpoint{7.42cm}{5.633cm}}
\pgfpathcurveto{\pgfqpoint{7.764cm}{5.977cm}}{\pgfqpoint{8.323cm}{5.977cm}}{\pgfqpoint{8.667cm}{5.633cm}}
\pgfusepath{stroke}
\definecolor{eps2pgf_color}{rgb}{0,0,0}\pgfsetstrokecolor{eps2pgf_color}\pgfsetfillcolor{eps2pgf_color}
\pgftext[x=8.049cm,y=4.943cm,rotate=0]{\selectfont{+}}
\pgfsetdash{}{0cm}
\definecolor{eps2pgf_color}{rgb}{0,0,0}\pgfsetstrokecolor{eps2pgf_color}\pgfsetfillcolor{eps2pgf_color}
\pgfpathmoveto{\pgfqpoint{3.511cm}{5.009cm}}
\pgfpathlineto{\pgfqpoint{6.794cm}{5.009cm}}
\pgfusepath{stroke}
\pgfpathmoveto{\pgfqpoint{7.077cm}{5.009cm}}
\pgfpathlineto{\pgfqpoint{6.794cm}{5.115cm}}
\pgfpathlineto{\pgfqpoint{6.794cm}{4.904cm}}
\pgfpathclose
\pgfusepath{fill}
\pgfsetdash{}{0cm}
\pgfsetbuttcap
\pgfsetmiterjoin
\pgfpathmoveto{\pgfqpoint{7.077cm}{5.009cm}}
\pgfpathlineto{\pgfqpoint{6.794cm}{5.115cm}}
\pgfpathlineto{\pgfqpoint{6.794cm}{4.904cm}}
\pgfpathclose
\pgfusepath{stroke}
\pgfsetdash{}{0cm}
\pgfsetroundcap
\pgfsetroundjoin
\pgfpathmoveto{\pgfqpoint{8.992cm}{5.027cm}}
\pgfpathlineto{\pgfqpoint{13.176cm}{5.027cm}}
\pgfusepath{stroke}
\pgfpathmoveto{\pgfqpoint{13.458cm}{5.027cm}}
\pgfpathlineto{\pgfqpoint{13.176cm}{5.133cm}}
\pgfpathlineto{\pgfqpoint{13.176cm}{4.921cm}}
\pgfpathclose
\pgfusepath{fill}
\pgfsetdash{}{0cm}
\pgfsetbuttcap
\pgfsetmiterjoin
\pgfpathmoveto{\pgfqpoint{13.458cm}{5.027cm}}
\pgfpathlineto{\pgfqpoint{13.176cm}{5.133cm}}
\pgfpathlineto{\pgfqpoint{13.176cm}{4.921cm}}
\pgfpathclose
\pgfusepath{stroke}
\pgfsetdash{}{0cm}
\pgfsetroundcap
\pgfsetroundjoin
\pgfpathmoveto{\pgfqpoint{8.022cm}{7.714cm}}
\pgfpathlineto{\pgfqpoint{8.028cm}{6.299cm}}
\pgfusepath{stroke}
\pgfpathmoveto{\pgfqpoint{8.029cm}{6.017cm}}
\pgfpathlineto{\pgfqpoint{8.134cm}{6.3cm}}
\pgfpathlineto{\pgfqpoint{7.922cm}{6.299cm}}
\pgfpathclose
\pgfusepath{fill}
\pgfsetdash{}{0cm}
\pgfsetbuttcap
\pgfsetmiterjoin
\pgfpathmoveto{\pgfqpoint{8.029cm}{6.017cm}}
\pgfpathlineto{\pgfqpoint{8.134cm}{6.3cm}}
\pgfpathlineto{\pgfqpoint{7.922cm}{6.299cm}}
\pgfpathclose
\pgfusepath{stroke}
\pgfsetdash{}{0cm}
\pgfsetroundcap
\pgfsetroundjoin
\pgfpathmoveto{\pgfqpoint{7.564cm}{4.248cm}}
\pgfpathlineto{\pgfqpoint{6.851cm}{3.117cm}}
\pgfusepath{stroke}
\pgfpathmoveto{\pgfqpoint{6.7cm}{2.878cm}}
\pgfpathlineto{\pgfqpoint{6.94cm}{3.061cm}}
\pgfpathlineto{\pgfqpoint{6.761cm}{3.173cm}}
\pgfpathclose
\pgfusepath{fill}
\pgfsetdash{}{0cm}
\pgfsetbuttcap
\pgfsetmiterjoin
\pgfpathmoveto{\pgfqpoint{6.7cm}{2.878cm}}
\pgfpathlineto{\pgfqpoint{6.94cm}{3.061cm}}
\pgfpathlineto{\pgfqpoint{6.761cm}{3.173cm}}
\pgfpathclose
\pgfusepath{stroke}
\definecolor{eps2pgf_color}{rgb}{0,0,0}\pgfsetstrokecolor{eps2pgf_color}\pgfsetfillcolor{eps2pgf_color}
\pgftext[x=6.509cm,y=2.539cm,rotate=0]{\selectfont{Charge}}
\pgfsetdash{}{0cm}
\pgfsetroundcap
\pgfsetroundjoin
\definecolor{eps2pgf_color}{rgb}{0,0,0}\pgfsetstrokecolor{eps2pgf_color}\pgfsetfillcolor{eps2pgf_color}
\pgfpathmoveto{\pgfqpoint{9.574cm}{2.821cm}}
\pgfpathlineto{\pgfqpoint{8.759cm}{3.986cm}}
\pgfusepath{stroke}
\pgfpathmoveto{\pgfqpoint{8.597cm}{4.217cm}}
\pgfpathlineto{\pgfqpoint{8.673cm}{3.925cm}}
\pgfpathlineto{\pgfqpoint{8.846cm}{4.047cm}}
\pgfpathclose
\pgfusepath{fill}
\pgfsetdash{}{0cm}
\pgfsetbuttcap
\pgfsetmiterjoin
\pgfpathmoveto{\pgfqpoint{8.597cm}{4.217cm}}
\pgfpathlineto{\pgfqpoint{8.673cm}{3.925cm}}
\pgfpathlineto{\pgfqpoint{8.846cm}{4.047cm}}
\pgfpathclose
\pgfusepath{stroke}
\definecolor{eps2pgf_color}{rgb}{0,0,0}\pgfsetstrokecolor{eps2pgf_color}\pgfsetfillcolor{eps2pgf_color}
\pgftext[x=9.754cm,y=2.539cm,rotate=0]{\selectfont{Discharge}}
\pgfsetdash{}{0cm}
\pgfsetroundcap
\pgfsetroundjoin
\definecolor{eps2pgf_color}{rgb}{0,0,0}\pgfsetstrokecolor{eps2pgf_color}\pgfsetfillcolor{eps2pgf_color}
\pgfpathmoveto{\pgfqpoint{6.509cm}{2.328cm}}
\pgfpathlineto{\pgfqpoint{6.511cm}{1.476cm}}
\pgfusepath{stroke}
\pgfpathmoveto{\pgfqpoint{6.512cm}{1.193cm}}
\pgfpathlineto{\pgfqpoint{6.617cm}{1.476cm}}
\pgfpathlineto{\pgfqpoint{6.405cm}{1.475cm}}
\pgfpathclose
\pgfusepath{fill}
\pgfsetdash{}{0cm}
\pgfsetbuttcap
\pgfsetmiterjoin
\pgfpathmoveto{\pgfqpoint{6.512cm}{1.193cm}}
\pgfpathlineto{\pgfqpoint{6.617cm}{1.476cm}}
\pgfpathlineto{\pgfqpoint{6.405cm}{1.475cm}}
\pgfpathclose
\pgfusepath{stroke}
\pgfsetdash{}{0cm}
\pgfsetroundcap
\pgfsetroundjoin
\pgfpathmoveto{\pgfqpoint{9.747cm}{1.234cm}}
\pgfpathlineto{\pgfqpoint{9.747cm}{1.978cm}}
\pgfusepath{stroke}
\pgfpathmoveto{\pgfqpoint{9.747cm}{2.261cm}}
\pgfpathlineto{\pgfqpoint{9.641cm}{1.978cm}}
\pgfpathlineto{\pgfqpoint{9.853cm}{1.978cm}}
\pgfpathclose
\pgfusepath{fill}
\pgfsetdash{}{0cm}
\pgfsetbuttcap
\pgfsetmiterjoin
\pgfpathmoveto{\pgfqpoint{9.747cm}{2.261cm}}
\pgfpathlineto{\pgfqpoint{9.641cm}{1.978cm}}
\pgfpathlineto{\pgfqpoint{9.853cm}{1.978cm}}
\pgfpathclose
\pgfusepath{stroke}
\definecolor{eps2pgf_color}{rgb}{1,1,1}\pgfsetstrokecolor{eps2pgf_color}\pgfsetfillcolor{eps2pgf_color}
\pgfpathmoveto{\pgfqpoint{6.312cm}{1.188cm}}
\pgfpathlineto{\pgfqpoint{9.963cm}{1.188cm}}
\pgfpathlineto{\pgfqpoint{9.963cm}{0.095cm}}
\pgfpathlineto{\pgfqpoint{6.312cm}{0.095cm}}
\pgfpathclose
\pgfusepath{fill}
\pgfsetdash{}{0cm}
\pgfsetroundcap
\pgfsetroundjoin
\definecolor{eps2pgf_color}{rgb}{0,0,0}\pgfsetstrokecolor{eps2pgf_color}\pgfsetfillcolor{eps2pgf_color}
\pgfpathmoveto{\pgfqpoint{6.312cm}{1.188cm}}
\pgfpathlineto{\pgfqpoint{9.963cm}{1.188cm}}
\pgfpathlineto{\pgfqpoint{9.963cm}{0.095cm}}
\pgfpathlineto{\pgfqpoint{6.312cm}{0.095cm}}
\pgfpathclose
\pgfusepath{stroke}
\definecolor{eps2pgf_color}{rgb}{0,0,0}\pgfsetstrokecolor{eps2pgf_color}\pgfsetfillcolor{eps2pgf_color}
\pgftext[x=9.132cm,y=0.538cm,rotate=0]{\selectfont{         $\b_t$}}
\pgftext[x=7.599cm,y=0.609cm,rotate=0]{\selectfont{Storage}}
\pgftext[x=8.049cm,y=8.177cm,rotate=0]{\selectfont{Controllable inflow}}
\pgftext[x=1.600cm,y=5.497cm,rotate=0]{\selectfont{Uncontrollable }}
\pgftext[x=1.466cm,y=4.921cm,rotate=0]{\selectfont{energy}}
\pgftext[x=1.466cm,y=4.509cm,rotate=0]{\selectfont{imbalance}}
\pgftext[x=14.429cm,y=5.526cm,rotate=0]{\selectfont{Residual}}
\pgftext[x=14.432cm,y=4.95cm,rotate=0]{\selectfont{energy }}
\pgftext[x=14.300cm,y=4.538cm,rotate=0]{\selectfont{imbalance}}
\pgftext[x=8.328cm,y=6.991cm,rotate=0]{\selectfont{$f_t$}}
\pgftext[x=5.217cm,y=5.378cm,rotate=0]{\selectfont{$\d_t$}}
\pgftext[x=6.050cm,y=1.791cm,rotate=0]{\selectfont{$\u^+_t$}}
\pgftext[x=10.355cm,y=1.791cm,rotate=0]{\selectfont{$\u^-_t$}}
\pgftext[x=6.071cm,y=3.586cm,rotate=0]{\selectfont{$\hC(\u^+_t)$}}
\pgftext[x=10.039cm,y=3.586cm,rotate=0]{\selectfont{$\hD(\u^-_t)$}}
\pgftext[x=11.087cm,y=5.34cm,rotate=0]{\selectfont{$\dr_t\! =\! \d_t\!-\!h(u_t)\!+\!f_t$}}
\end{pgfscope}
\end{pgfscope}
\end{pgfpicture}

%% file: algorithm-conf.tex
\section{The Online Modified Greedy Algorithm}\label{sec:onlineALG}
\subsection{Algorithm}
Among algorithms that have been proposed to solve problem \eqref{prob:singleBusGeneral}, the greedy (or myopic) algorithm is one of the simplest. In an online setting where at the beginning of each time period $t$ the realizations of the stochastic parameters, $\drl_t$ and $\prl_t$, are revealed to the operator, the \emph{greedy algorithm} solves
\begin{subequations}\label{prob:greedy}
\begin{align}
\minimize \quad &   \grl_t= \g_t(\u_t, \f_t, \drl_t, \prl_t) \\
\st \quad
        & \bmin \le \la \b_t + \u_t \le \bmax, \label{eq:cbgre}\\
        & \umin \le \u_t \le \umax,\\
        & \f_t \in \F,
\end{align}
\end{subequations}
where  the optimization variables are $\u_t$ and $\f_t$. Other than rare cases, the greedy algorithm is sub-optimal for problem~\eqref{prob:singleBusGeneral}, and the level of sub-optimality is usually difficult to characterize. In the reminder of this section, we show that a slight modification of \eqref{prob:greedy} renders an algorithm that comes with provable bounds to optimality.

The algorithm, termed the online modified greedy (OMG) algorithm, is composed of an offline and online phase. Next we describe the input data to the algorithm and each phase. \vspace{.2cm}

\noindent{\bf Input Data.}
Other than data specifying the storage model ($\Sm$ and $h$),
 %and the cost functions\footnote{Due to the greedy nature of the algorithm, this requirement can be reduced such that the cost function $\g_t(\cdot)$ needs only to become known to the operator at the beginning of time period $t$. } (the functional form of $\g_t$), 
OMG  requires two more parameters regarding the cost functions, denoted by $\Dl g$ and $\Du g$ which are defined as follows.
\begin{definition}
Let $y\defeq (f, \d, \p)$. For function $\phi_t(u, y) \defeq \g_t(\u, f, \d, \p)$ that is convex (but not necessarily differentiable) in $\u$, a real number $\alpha$ is called a (partial) subgradient of $\phi_t$ with respect to argument $\u$ at given $(\u, y)$ if $\phi_t(u',y) \ge \phi_t(u, y) + \alpha (u'-u)$ for all $u' \in [\umin, \umax]$. The set of all subgradients at $(\u, y)$, denoted by $\partial_u \phi_t(u,y)$, is called the (partial) subdifferential of $\phi_t(u, y)$ with respect to $\u$ at $(\u,y)$. Denote $\mathcal{U} \defeq [\umin,\umax]$, $\mathcal{Y} \defeq \F \times [\dmin, \dmax] \times [\pmin, \pmax]$, $\mathbbm{Z}_+ \defeq \{1,2,\dots\}$, where $[\dmin, \dmax]$ and $[\pmin, \pmax]$ are the compact supports for $\d_t$ and $\p_t$, respectively. Define the set
\[
D \g \defeq \bigcup_{(t,u,y) \in \mathbbm{Z}_{+}\times\mathcal{U} \times \mathcal{Y}} \partial_u \phi_t(\u, y),
\]
and let real numbers $\Dl g$ and $\Du g$ be defined such that
\begin{equation}
\Dl\g \le \inf D\g \le \sup D\g \le \Du\g.
\end{equation}
That is, $\Dl \g$ and $\Du \g$ are a lower bound and an upper bound of the subgradient of $\phi_t$ over  its (compact) domain and over all time periods, respectively.
\end{definition}

The quantities $\Dl \g$ and $\Du \g$ partially characterize how sensitive the cost is in perturbation of storage operation.
It will be shown later that a smaller $\Du \g - \Dl \g$ leads to a tighter sub-optimality bound of our algorithm, so that if possible one should select $\Dl g = \inf Dg$ and $\Du g = \sup Dg$. We demonstrate the procedure of calculating $\Du \g$ and $\Dl \g$ for cost functions discussed in Examples~\ref{eg:arbi},~\ref{eg:ba} and~\ref{eg:co} under the simplification that the conversion  function $h$ is the identity mapping, \ie, $h(\u) = \u$. 
\begin{example}[Calculate $\Dl \g$ and $\Du \g$]
(i) For the arbitrage cost function~\eqref{eq:arbicost}, we have
\[
\partial_u \g_t(\u, \p_t) = \{\p_t\} \mbox{ and } D\g = [\pmin, \pmax].
\]
 Thus one can set $\Dl \g = \pmin$ and $\Du \g =  \pmax$.

(ii) For the balancing cost~\eqref{eq:bacost}, if for example the penalty rate is homogeneous across time (\ie, $q^+_t \equiv q^+ \ge 0$, $\q^-_t \equiv q^- \ge 0$)\footnote{We also assume the feasible set is such that both $\dr_t >0$ and $\dr_t <0$ are possible for certain (but not necessarily the same) $t$ and $(\u_t, \f_t, \d_t)$.} 
%and loss of load is penalized with a higher rate than curtail of generation (\ie, $q^+ \ge q^- \ge 0$), 
, then it is easy to check that
$ D\g = [-q^+, q^-]$,
and so $\Dl \g = -q^+$ and $\Du \g =  q^-$.

(iii) For the cost function~\eqref{eq:cocost} and positive prices ($\pmax \ge \pmin \ge 0$), one can use
$\Dl \g= 0$ and $\Du\g = \pmax$.
\end{example}

For more general cost functions, one may obtain $\Du \g$ and $\Dl \g$ by solving certain optimization problems.

\begin{remark}[Distribution-Free Method] The OMG algorithm is a distribution-free method in the sense that almost no information regarding the distribution of the stochastic parameters $\d_t$ and $\p_t$ is required. The only exception is when calculating $\Du\g$ and $\Dl \g$, the supports of $\d_t$ and $\p_t$ may be needed. But compared to the entire distribution functions, it is much easier to estimate the supports of the stochastic parameters  from historical data.
\end{remark}

\begin{remark}[Determine the Supports for $\d_t$ and $\p_t$]
The supports for $\d_t$ and $\p_t$ may be determined based on the physical parameters of the system. For instance, if $\d_t$ models the wind power generation process, then $\dmin$ and $\dmax$ may be determined using the minimal possible wind generation (which is $0$ in many cases) and the nameplate capacity for the wind farm, respectively; if $\p_t$ models the locational marginal prices at the bus, then it can be bounded using an estimate of the maximal marginal cost of generation.  Another possible approach is to estimate the supports using the forecasts of $\d_t$ and $\p_t$, which in turn are based on historical observation of the processes. Techniques that are used to determine the uncertainty sets for robust optimization can be used here; interested readers are referred to \cite{bertsimas2009constructing} for more details. As in general a smaller $\Du g - \Dl g$ leads to better performance guarantees, it is beneficial to  obtain a tight estimate for the supports of the stochastic parameters. 
\end{remark}

\vspace{.2cm}

\noindent{\bf Offline Phase.}
The algorithm depends on two algorithmic parameters, namely a shift parameter $\ks$ and a weight parameter $\W$, that should be selected offline. Any pair $(\ks, \W)$ satisfies the following conditions can be used\footnote{Discussions of the intuitions behinds the algorithmic parameters are deferred to the part describing the online phase of the algorithm. The conditions on $\ks$ and $\W$ follow from the feasibility requirement of the algorithm; see Appendix A for more details.}:
\begin{align}
\ksmin\le & \ks \le \ksmax, \label{eq:ksbounds}\\
0 < & \W  \le \Wmax, \label{eq:Wbounds}
\end{align}
where
\begin{equation}\label{eq:ineq_1}
\ksmin \defeq \frac{1}{\la} \left(-\W \Dl \g + \umax - \bmax\right),
\end{equation}
\begin{equation}\label{eq:ineq_2}
\ksmax \defeq \frac{1}{\la} \left(-\W \Du \g - \bmin+ \umin\right),\end{equation}
and
\begin{equation}\label{eq:W_max}
\Wmax \defeq \frac{(\bmax - \bmin) -(\umax - \umin)}{\Du \g - \Dl \g}  .
\end{equation}
Note that the interval for $\W$ in \eqref{eq:Wbounds} is well-defined under a mild condition (see the next subsection for more details), and the interval for $\ks$ in \eqref{eq:ksbounds} is always well-defined. It will be clear later that the sub-optimality bound depends on the choice of $(\ks, \W)$. Here we provide two possible ways for selecting these parameters.
% detailed discussion of both approaches are delayed to the end of this section.
\begin{itemize}
\item The \emph{maximum weight} approach (\texttt{maxW}): Setting $\W = \Wmax$, one reduces the interval in \eqref{eq:ksbounds} to a singleton ($\ksmin = \ksmax$) and
\begin{equation}\label{eq:ksifWmax}
\ks  = \frac{\Dl \g(\bmin - \umin) -\Du \g(\bmax - \umax)}{\la(\Du \g - \Dl \g)}.
\end{equation}
Using this parameter configuration in a sense sets OMG to be the ``greediest'' in the range of admissible parameter specifications. 
\item The \emph{minimum sub-optimality bound} approach (\texttt{minS}): It turns out that the sub-optimality bound of OMG as a function of $(\ks, \W)$ can be minimized using a semidefinite program reformulation (see Lemma~\ref{SDP_P3_PO} in the next section). Empirical results show that using the bound minimizing $(\ks, \W)$, one often obtains better lower bounds for the optimal costs. Thus this is the recommended approach if one runs the OMG algorithm for the purpose of evaluating other algorithms. It is not necessarily the case that the actual algorithm performance with this choice of algorithmic parameters is optimized -- minimizing the sub-optimality bound is not equivalent to minimizing the actual sub-optimality.
\end{itemize}

\begin{remark}
For ideal storage ($\la =1$), the maximum weight and minimum sub-optimality bound approaches coincide.
\end{remark}\vspace{.2cm}

\noindent{\bf Online Phase.} At the beginning of each time period $t$, the OMG algorithm solves the following modified version of program \eqref{prob:greedy},
\begin{subequations}\label{prob:P3_alg}
\begin{align}
 \minimize  &  \la (\b_t+\ks) \u_t + \W \grl_t\\
\st  %& \b_{t+1} = \la_v \b_t + \u_t, \label{P1:dynamics}\\
        %& \bmin \le \b_t \le \bmax, \\
        %& \bmin - \la \b_t \le \u_t \le \bmax - \la \b_t, \label{P1:bbounds-u}\\
        & \umin \le \u_t \le \umax, \\
        &  f_t \in \F,
\end{align}
\end{subequations}
for the storage operation $\u_t$ and controllable inflow $f_t$. Comparing the above optimization~\eqref{prob:P3_alg} to optimization~\eqref{prob:greedy}, one notices two modifications. The first modification is in the objective function. Instead of directly optimizing the cost at the current time period, the OMG algorithm optimizes a weighted combination of the stage-wise cost and a linear term of $\u_t$ depending on the shifted storage level $\b_t +\ks$. Here the weight parameter $W$ decides the importance of the original cost in this weighted combination, while the shift parameter $\ks$ defines the shifted state given the original state $\b_t$.  
Roughly speaking, the shifted state $\b_t+\ks$ belongs to an interval $[\bmin+\ks, \bmax+\ks]$ which usually contains $0$.
If the storage level is relatively high, the shifted state is greater than $0$, such that the state-dependent term (\ie, $\la (\b_t + \ks) \u_t$) encourages a negative $\u_t$ (discharge) to minimize the weighted sum. As a result, the storage level in the next time period will be brought down. On the other hand, if the storage level is relatively low, the shifted state is smaller than $0$, such that the state-dependent term encourages a positive $\u_t$ (charge) and consequently the next stage storage level is increased. These two effects together help to hedge against uncertainty by maintaining a storage level somewhere in the middle of the feasible interval.
 The second modification is the deletion of the constraint~\eqref{eq:cbgre}. We will show later that by selecting $(\ks, \W)$ satisfying conditions~\eqref{eq:ksbounds} and~\eqref{eq:Wbounds}, the constraint~\eqref{eq:cbgre} holds automatically. However, for the purpose of robustness (considering the possibility of feeding incorrect parameters to the algorithm), one can optionally add the constraint~\eqref{eq:cbgre} to~\eqref{prob:P3_alg}.

In case that $f_t = 0$, the online optimization usually can be solved analytically. This leads to further simplification of the implementation. Assuming $h$ is the identity mapping, we work out the analytical solutions of~\eqref{prob:P3_alg} with the cost functions given in Examples~\ref{eg:arbi} and~\ref{eg:ba}.
\begin{example}[Analytical Solutions of the Online Program]
(i) For the arbitrage cost function~\eqref{eq:arbicost}, the optimal storage operation $\u^\star_t$ is as follows:%\footnote{In the last case, any $\u_t \in [\umin, \umax]$ is optimal. We choose $\u^\star_t = 0$ for convenience. Similar treatment is used in the sequel.}:
\[
\u^\star_t = \begin{cases}
\umin &\mbox{if } \b_t > (\W\p_t/\la) - \ks,\\
\umax &\mbox{if } \b_t \le (\W\p_t/\la) - \ks.
%0 &\mbox{if } \b_t = (\W\p_t/\la) - \ks.
\end{cases}
\]
\end{example}

(ii) For the balancing cost function~\eqref{eq:bacost}, the optimal storage operation is
\[
\u^\star_t = \begin{cases}
\umin &\mbox{if } \b_t > (\W q^-_t/\la) - \ks,\\
\umax &\mbox{if } \b_t < (-\W q^+_t/\la) - \ks,\\
\Pi_{\mathcal{U}}(-\d_t) &\mbox{if } (-\W q^+_t/\la) - \ks \le \b_t \le (\W q^-_t/\la) - \ks,
\end{cases}
\]
where $\Pi_{\mathcal{U}}(\cdot)$ is the (Euclidean) projection operator for the feasible set of storage operation $\mathcal{U} = [\umin,\umax]$, \ie, $\Pi_{\mathcal{U}}(-\d_t) = \min \left(\max(-\d_t, \umin), \umax\right)$.

%For the balancing cost function~\eqref{eq:bacost}, depending on the realization of the energy imbalance $\d_t$, there are three cases.
%\begin{itemize}
%\item If $-\d_t \le \umin$, the optimal storage operation is
%\[
%\u^\star_t = \begin{cases}
%\umin &\mbox{if } \b_t > (-\W q^+_t/\la) - \ks,\\
%\umax &\mbox{if } \b_t < (-\W q^+_t/\la) - \ks,\\
%0 &\mbox{if } \b_t = (-\W q^+_t/\la) - \ks.
%\end{cases}
%\]
%\item If $-\d_t \ge \umax$, the optimal storage operation is
%\[
%\u^\star_t = \begin{cases}
%\umin &\mbox{if } \b_t > (\W q^-_t/\la) - \ks,\\
%\umax &\mbox{if } \b_t < (\W q^-_t/\la) - \ks,\\
%0 &\mbox{if } \b_t = (\W q^-_t/\la) - \ks.
%\end{cases}
%\]
%\item If $\umin <-\d_t < \umax$, the optimal storage operation is
%\[
%\u^\star_t = \begin{cases}
%\umin &\mbox{if } \b_t > (\W q^-_t/\la) - \ks,\\
%\umax &\mbox{if } \b_t < (-\W q^+_t/\la) - \ks,\\
%-\d_t &\mbox{if } (-\W q^+_t/\la) - \ks \le \b_t \le (\W q^-_t/\la) - \ks.
%\end{cases}
%\]
%\end{itemize}

We close this subsection by summarizing the algorithm in a compact form (Algorithm~\ref{alg}).
\begin{algorithm}[htbp]
\caption{Online Modified Greedy Algorithm}
\begin{algorithmic}
\INPUT $\Dl \g$, $\Du \g$, $\Sm$, $h$, and the functional form of $\g_t$.
\Offline Determine $(\ks,\W)$ using either the maximum weight or minimum sub-optimality bound approaches.
\Online
\For {each time period $t$}
    \State Observe realizations of $\d_t$ and $\p_t$ and solve~\eqref{prob:P3_alg}.
%    \State Solve the following deterministic optimization for storage operation $\u_t$:
%\begin{subequations}\label{prob:P3_alg}
%\begin{align}
% \minimize  &  \la (\b_t+\ks) \u_t + \W \g_t\\
%\st  %& \b_{t+1} = \la_v \b_t + \u_t, \label{P1:dynamics}\\
%        %& \bmin \le \b_t \le \bmax, \\
%        %& \bmin - \la \b_t \le \u_t \le \bmax - \la \b_t, \label{P1:bbounds-u}\\
%        & \umin \le \u_t \le \umax,
%\end{align}
%\end{subequations}
\EndFor
\end{algorithmic}\label{alg}
\end{algorithm}
%1. State the modified optimization and discuss some intuition. 2. Stress the fact that it's extremely simple to implement. 3. Given analytical solution to guild the implementation.
\subsection{Analysis of the Algorithm Performance}
We proceed by providing a stylized analysis for the algorithm performance.
\begin{assumption}\label{ass:1}
The following assumptions are in force for the analysis in this section.
\begin{enumerate}[{\bf A}1]
\item Infinite horizon: The horizon length $T$ approaches to infinity.
\item IID disturbance: The imbalance process $\{\d_t: t\ge 1\}$ is independent and identically distributed (i.i.d.) across $t$ and is supported on a compact interval $[\dmin, \dmax]$. Similarly, the process $\{\p_t: t \ge 1\}$ is i.i.d. across $t$ and is supported on a compact interval $[\pmin, \pmax]$. Here $\d_t$ and $\p_t$ may be correlated. 
\item Frequent acting: The storage parameters satisfy $\umax-\umin < \bmax -\bmin$.
\end{enumerate}
\end{assumption}
Here {\bf A}1 and {\bf A}2 are technical assumptions introduced to simplify the exposition. 
An extra term of $O(1/T)$ appears in the sub-optimality bound when {\bf A}1 is relaxed.\footnote{See Remark~\ref{rk:ft} for some additional discussions.}
For $T$ on the order of  $10^3$ (which is, \eg, corresponding to operating the storage every 30 minutes for a month or every 5 minutes for a week) or larger, this term is negligible.  The bounds in this section may not be accurate for applications with truly small $T$. Appendix B discusses how to reduce {\bf A}2. Under these two assumptions, the storage operation problem can be cast as an infinite horizon average cost stochastic optimal control problem in the following form
\begin{subequations}\label{prob:singleBusGeneralinf}
\begin{align}
\minimize \quad &  \lim_{T \to \infty}  (1/T) \expec \Big[ \sum_{t=1}^T \g_t \Big]\\
\st \quad &\eqref{eq:ca}, \eqref{eq:cb}, \eqref{eq:cc}, \eqref{eq:cd},
\end{align}
\end{subequations}
where we aim to find a control policy that maps the information available up to each of the stages to \ylmod{control actions} that minimizes the expected average cost  and satisfies all the constraints for each time period $t$.

Assumption {\bf A}3 appears to be a restriction on the physical parameters of the storage model. It states that the range of feasible storage control $\umax -\umin$ is smaller than the range of storage levels $\bmax -\bmin$, \ie, the ramping limits of the storage is relatively small compared to the storage capacity. This is, nevertheless, not completely true as the designer of the storage controller usually also has the freedom to select the frequency of the controller in a range of possible values. More specifically, for a fixed storage system, it has a certain storage capacity (\eg, energy rating in unit of MWh, and \ie, $\bmax-\bmin$ in our notation) and certain charging/discharging ramping capacity (\eg, power rating in unit of MW, and denoted by $r^+$ and $r^-$ for charging and discharging rate, respectively). We have $\umax = r^+ \dt$, $\umin = -r^- \dt$, and therefore $\umax -\umin = (r^++r^-)\dt$ can be made smaller than $\bmax -\bmin$ as long as the frequency of the controller is high enough (or the length of each time period $\dt$ is small enough).

Define \ylmod{$J(\upi,\fpi)$} as the value (or total cost) function of~\eqref{prob:singleBusGeneral} induced by the \ylmod{sequence of control policies $\{(\upi_t,\fpi_t),\,t\geq 1\}$ and $J^\star=J(\upistar, \fpistar)$} as the optimal value of the average cost stochastic control problem with \ylmod{$\{(\upistar_t,\fpistar_t),\,t\geq 1\}$} being the corresponding optimal \ylmod{sequence of control policies}. Sometimes we also use the notation $J(\upi)$ when the $\fpi$ sequence is clear from the context.
We are ready to state the main theorem regarding the performance of the OMG algorithm.
\begin{theorem} [Performance]\label{thm:perf_lyap}
The \ylmod{control policy} sequence \ylmod{$(\upihat, \fpihat) \defeq\{(\upihat_t, \fpihat_t), t\ge 1\}$} generated by the OMG algorithm is feasible with respect to all constraints of~\eqref{prob:singleBusGeneral} and its sub-optimality  is bounded by $\M(\ks)/\W$, that is
\begin{equation}\label{eq:1busiidperf_bdd}
J^\star \le  \ylmod{J(\upihat,\fpihat)} \le J^\star + \M(\ks)/\W,
\end{equation}
%Let $\uhat_t$ be the control action chosen by solving \textbf{P3} and
%\[
%\begin{split}
%\Jk{3}^\star_t=&\sum_{\ell = 1}^L \p(t, \ell) \Bigg(\alI(\ell) \d_t - \alC(\ell)\hC\left( \pos{\uhat_t}\right) \\
%&\quad + \alD(\ell) \hD\left(\neg{\uhat_t}\right) + \alConst(\ell)\Bigg)^+.
%\end{split}
%\]
%Then
% \begin{equation} \label{eq:perf_bdd}
%  \lim_{T \to \infty} \frac{1}{T}\expec \left[ \sum_{t=1}^{T} \Jk{3}^\star_t\right]  \leq \Jk{1}^\star + \frac{\M(\ks) }{\W}.
%\end{equation}
where
\begin{align*}
&\M(\ks)=\Mone(\ks)+\la(1-\la)\Mtwo(\ks),\\
&\Mone(\ks) =\! \frac{1}{2}\! \max\left(\! \left(\umin\!+(1-\la)\ks\right)^2\!\!,\left(\umax\!+(1-\la)\ks\right)^2\! \right)\!,\\
&\Mtwo(\ks)= \max\left( \left(\bmin+\ks\right)^2,\left(\bmax+\ks\right)^2 \right).
\end{align*}
\end{theorem}
%\begin{IEEEproof}
%A quadratic Lyapunov function is constructed. The relation between the Lyapunov drift and the objective function of the online optimization problem in \eqref{prob:P3_alg} is exploited, which in turn relates to the objective function of the average cost stochastic control problem in \eqref{prob:singleBusGeneral}. Appendix~\ref{sec:app:1bus} contains the whole proof.
%\end{IEEEproof}

The theorem above guarantees that the cost of the OMG algorithm is bounded above by $J^\star + \M(\ks)/\W$. 
The proof of the theorem is relegated to Appendix A. 
The sub-optimality bound $\M(\ks)/\W$ reduces to a much simpler form if $\la = 1$.
\begin{remark}[Sub-Optimality Bound, $\la =1 $]\label{remark:subo:la1}
For a storage with $\la = 1$, we have \[\M \defeq \M(\ks) = (1/2) \max((\umin)^2, (\umax)^2),\] and the online algorithm is no worse than $\M/\W$ sub-optimal.
 In this case, one would optimize the performance by setting
\[
\W=\Wmax = \frac{(\bmax - \bmin)-(\umax-\umin)}{\Du \g\! - \Dl \g},
\]
and the corresponding
interval $[\ksmin,\ksmax]$ is a singleton with $\ksmin=\ksmax$ being the expression displayed in \eqref{eq:ksifWmax}.
%\end{remark}
%\begin{remark}[Optimality for Efficient Unbounded Storage]
Let $\bmax - \bmin = \rho (\umax - \umin)$. Suppose $|\umax| = |\umin|$. For ideal storage ($\la = 1$),  the sub-optimality bound is
\[
\frac{\M}{\W}= \frac{(1/2)(\Du \g - \Dl g)(\umax)^2}{(\bmax- \bmin) - (\umax - \umin)} = \frac{\Du \g - \Dl g}{4(\rho - 1)} \umax.
\]
For fixed $\umax$, as storage capacity increases, \ie,  $\rho \to \infty$, the sub-optimality $(\M/\W) \to 0$. That is, OMG is near-optimal for ideal storage with small ramping limits and a large capacity.  On the other hand, if $\umax$ and $\bmax$ increases with their ratio $\rho$ fixed, the bound increases linearly with $\umax$. %Note that these two types of scaling correspond to different physical scenarios.
\end{remark}

%Based on the above theorem, we have shown that the online control policy, derived from Lyapunov optimization, is $O( \M(\ks)/{\W})-$optimal. First for $\la=1$, it is easy to see that this online control policy is $O(1/\W)-$optimal. In this case, the performance is optimized by choosing $\W=\Wmax$. Next, in the following theorem, we will find the value of $(\ks,\W)$ that optimizes the sub-optimal bound, for the general case with $\la\in(0,1)$, using convex optimization.
For the remaining  case $\la \in (0,1)$, the sub-optimality bound is no longer monotone in $W$ as choosing a smaller $W$ can lead to a larger interval $[\ksmin, \ksmax]$ potentially containing a $\ks$ which in turn leads to smaller $\Mone(\ks)$ and $\Mtwo(\ks)$ values.
Thus it requires solving an optimization program to identify the bound-minimizing parameter pair $(\ks, \W)$. In the next result, we state a semidefinite program to find $(\ks^\star, \W^\star)$ that solves the following  parameter optimization program
\begin{subequations}\label{P3:PO}
\begin{align*}
\opttag{PO:}\minimize & \quad \M(\ks)/\W\\
\st & \quad  \ksmin\le \ks \le \ksmax,\,\, 0< \W \le \Wmax,
\end{align*}
\end{subequations}
where the optimization variables are $\ks$ and $W$. 
%In the current form, this program appears to be non-convex. The next result reformulates {\bf PO} into a semidefinite program. 

\begin{lemma}[Semidefinite Reformulation of {\bf PO}]\label{SDP_P3_PO}
Let symmetric positive definite matrices $\Xumin$, $\Xumax$, $\Xbmin$ and $\Xbmax$ be defined as follows
%\begin{equation}
%\Xumin\!\! = \begin{bmatrix}
%\None&\!\!\umin+(1-\la)\ks\\
%*& 2\W
%\end{bmatrix},
%\Xumax\!\! = \begin{bmatrix}
%\None&\!\!\umax+(1-\la)\ks\\
%*& 2\W
%\end{bmatrix},
%\end{equation}
%\begin{equation}
%\Xbmin\!\! = \begin{bmatrix}
%\Ntwo&\!\!\bmin+\ks\\
%*& \W
%\end{bmatrix},
%\Xbmax\!\! = \begin{bmatrix}
%\Ntwo&\!\!\bmax+\ks\\
%*& \W
%\end{bmatrix},
%\end{equation}
\begin{equation*}
\!\!\Xudot\!\! = \!\!\begin{bmatrix}
\None&\!\!\udot+(1-\la)\ks\\
*& 2\W
\end{bmatrix}, \,\,
\Xbdot\!\! = \!\!\begin{bmatrix}
\Ntwo&\!\!\bdot+\ks\\
*& \W
\end{bmatrix},\!\!
\end{equation*}
where $(\cdot)$ can be either $\max$ or $\min$, and $\None$ and $\Ntwo$ are auxilliary variables. Then {\bf PO} can be solved via the following semidefinite program
\begin{subequations}\label{prob:sdp}
\begin{align}
\!\!\emph{\minimize} \quad & \None+\la(1-\la)\Ntwo \!\!\\
\!\!\emph{\st} \quad & \ksmin\le  \ks \le \ksmax,\,\, 0 <  \W  \le \Wmax, \!\!\\
& \Xumin,\Xumax, \Xbmin, \Xbmax \succeq 0, \!\!
\end{align}
\end{subequations}
where the optimization variables are $\W$, $\ks$, $\None$, $\Ntwo$, $\Xumin$, $\Xumax$, $\Xbmin$ and $\Xbmax$, and
 $\ksmin$ and $\ksmax$ are linear functions of $\W$ as defined in \eqref{eq:ineq_1} and \eqref{eq:ineq_2}.
\end{lemma}
This lemma provides us an efficient way to evaluate the minimum sub-optimality bound over all the algorithmic parameter choices. In the next example, we compare the minimum sub-optimality bounds for the case with $\la = 1$  and that for the case with $\la <1$. 
\begin{example}[Scaling of Sub-Optimality Bounds]
While the performance bounds in Theorem 1 holds for any instance of generalized storage models, it is useful to understand how the bound varies with the parameters of the storage system. For simplicity, we consider the balancing cost function with $\Du g = 1$ and $\Dl g =-1$, and storage systems with $\bmin = 0$ and $\umin = -\umax$.\footnote{Section~\ref{sec:baiid} will consider a similar setup. With the balancing cost~\eqref{costba}, the bounds calculated in this example can be physically interpreted as the average imbalance per unit. }
Motivated by discussions in Remark 4, we consider the following two sets of scenarios.
\begin{itemize}
\item Increasing the storage capacity $\bmax$ with a fixed $\umax/\bmax$ ratio: This set of scenarios can model e.g. a storage system consisting of $n$ identical battery modules with a common $\umax_i/\bmax_i$ ratio for battery $i$, $i=1,\dots,n$, whose system-wise charging and discharging circuit capacity is not constraining. As a demonstration, in this example, we fix the $\umax/\bmax=0.1$.
\item Increasing the storage capacity $\bmax$ with a fixed $\umax$: This set of scenarios can model e.g. a storage system consisting of $n$ identical battery modules and whose charging and discharging limits are determined by the shared system-wise charging/discharging circuit ratings instead of the intrinsic charging/discharging rates of each of the battery modules. In this example, we fix $\umax=0.01$. 
\end{itemize}
Figure~\ref{fig:x1} shows that the sub-optimality bound grows linearly with the storage capacity in the first set of scenarios for both $\la = 1$ and $\la <1$, and that larger $\la$ leads to smaller bounds. Choosing the algorithmic parameters using the SDP proposed in Lemma 1 (\texttt{minS}) leads to smaller bounds compared to the max weight heuristic (\texttt{maxW}) and the improvement is more significant when $\la$ is smaller.
Figure~\ref{fig:x2} depicts the bounds in the second set of scenarios, where it is shown that for $\la = 1$, increasing the storage capacity with fixed $\umax$ drives the sub-optimality bound to zero as predicted by Remark 4. However, the behavior of the bounds for $\la <1$ is very different in this set of scenarios due to the fact that larger storage capacity implies potentially more energy dissipation over time. As such, the sub-optimality bounds for both \texttt{minS} and \texttt{maxW} in fact grow with the storage capacity in a nonlinear fashion. Figure~\ref{fig:x3} plots the bounds amortized by the corresponding storage capacity. %\footnote{The amortized sub-optimality bounds in the setting of $n$ batteries mentioned above scale similarly.}.
For storage with energy dissipation, instead of approaching zero, the amortized sub-optimality decreases with the storage capacity and approaches a positive constant which increases with $(1-\la)$. 
\begin{figure}[htbp]
\centering\hspace{-1.em}%
\input{./fig/case1_1_1.tex}\hspace{-.2cm}%
\input{./fig/case1_2_1.tex}
\caption{The sub-optimality bound increases with $\bmax$ linearly when $\umax/\bmax$ ratio is fixed.}
\label{fig:x1}
\end{figure}
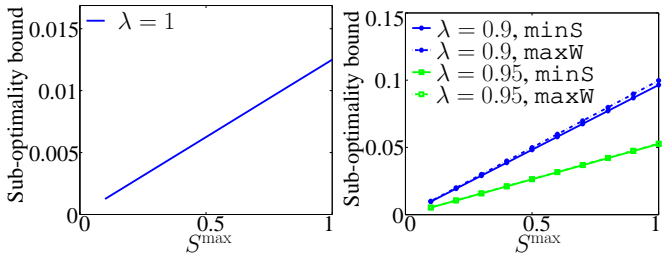
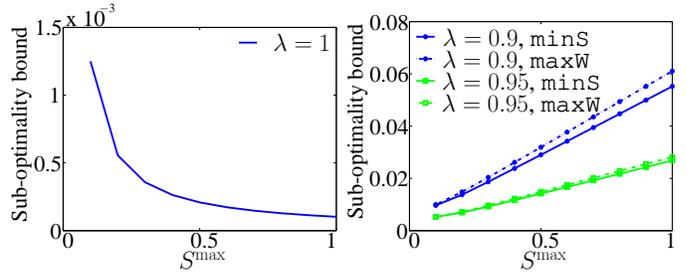
\begin{figure}[htbp]
\centering\hspace{-1.em}%
\input{./fig/case2_1_1.tex}\hspace{-.2cm}%
\input{./fig/case2_2_1.tex}
\caption{The sub-optimality bound decreases with $\bmax$ when $\la =1$ but increases with $\bmax$ when $\la <1$, given that $\umax$ is fixed.}
\label{fig:x2}
\end{figure}
\begin{figure}[htbp]
\centering\hspace{-1.em}%
\input{./fig/case2_1_2.tex}\hspace{-.2cm}%
\input{./fig/case2_2_2.tex}
\caption{The sub-optimality bound amortized by $\bmax$ decreases with $\bmax$ when $\umax$ is fixed.}
\label{fig:x3}
\end{figure}
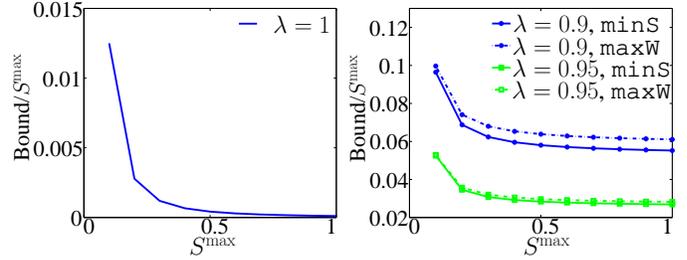

\end{example}

The surprising difference in the left and right panels of Figure~\ref{fig:x2} suggests the importance of modeling the energy dissipation in real-world applications.
\begin{remark}[Practical Guarantees of Lyapunov Methods]
Prior studies  \cite{DataCenter, StorDRLongbo, lyap1, lyap2}, which do not consider energy dissipation over time, have advocated the use of Lyapunov type methods for large storage based on the scaling shown in the left panel of Figure~\ref{fig:x2}. However, when energy dissipation is considered, the sub-optimality bound in fact grows with the storage capacity. Thus it is unclear that whether Lyapunov type methods are more suitable for large storage systems than smaller ones when there is energy dissipation. Furthermore, even for systems with a tiny amount of energy dissipation, it is very important to gauge the performance of Lyapunov methods using bounds for $\la <1$ as the bounds for $\la=1$ may substantially underestimate the sub-optimality especially for storage with a large capacity. 
\end{remark}

%\begin{IEEEproof}
%The result follows from Schur complement. See Appendix~\ref{sec:app:1bus} for details.
%\end{IEEEproof}

We close this section by discussing an implication of the performance theorem.
%\begin{remark}[Optimality at the Fast-Acting Limit]
%%Let the length of each time period be $\Delta t$.
%At the limit $\Delta t \to 0$, the online algorithm is optimal. Indeed, as discussed in Section~\ref{sec:problem}, both $|\umin|$ and $|\umax|$ are linear in $\Delta t$, such that $|\umax| \to 0$ and $|\umin|\to 0 $ as $\Delta t\to 0$. Meanwhile, $\la \to 1$ as $\Delta t \to 0$. So by Remark~\ref{remark:subo:la1}, it is easy to verify that the sub-optimality $\M/\W$ converges to zero as $\Delta t \to 0$.
%\end{remark}

\begin{remark}[Value of Storage and Percentage Cost Savings]\label{remark:vos}
In all applications including those discussed in Example~\ref{eg:arbi},~\ref{eg:ba}, and~\ref{eg:co},
the Operational Value of Storage (VoS) is broadly defined as the savings in the long term system cost due to storage operation. Such an index is usually calculated by assuming storage is operated optimally. In stochastic environments, the optimal system cost with storage operation is hard to obtain in general settings. Consider the case that $f_t = 0$. In our notations, let $\u^{\pi,\mathrm{ns}}$ denote the \ylmod{control policy sequence} $\{\upi_t: \upi_t = 0,  t\ge 1 \}$ which corresponds to no storage operation. Then
\[
\ylmod{\mathrm{VoS} = J(\u^{\pi,\mathrm{ns}}) - J^\star,}
\]
and it can be estimated by the interval
\[
\ylmod{\left[J(\u^{\pi,\mathrm{ns}})\! -\! J(\upihat),\,\,  J(\u^{\pi,\mathrm{ns}})\! -\! J(\upihat) \!+ \!\frac{\M}{\W} \right].}
\]
Additionally, for a storage operation \ylmod{control policy sequence $\upi$}, the percentage cost savings due to storage can then be defined by $(J(\u^{\pi,\mathrm{ns}}) - J(\upi))/J(\u^{\pi,\mathrm{ns}})$. An upper bound of this for any storage control policy can be obtained via $(J(\u^{\pi,\mathrm{ns}}) - J(\upihat) + \M/\W)/J(\u^{\pi,\mathrm{ns}})$, which to an extent summarizes the limit of a storage system in providing cost reduction.
\end{remark}

%% file: fig/case1_1_1.tex
\scalebox{0.235}{\scalefont{2} \input{./fig/case1_1_1.pgf}}

%% file: fig/case1_1_1.pgf
% Created by Eps2pgf 0.7.0 (build on 2008-08-24) on Mon Apr 06 01:29:46 PDT 2015
\begin{pgfpicture}
\pgfpathmoveto{\pgfqpoint{0.635cm}{6.315cm}}
\pgfpathlineto{\pgfqpoint{19.473cm}{6.315cm}}
\pgfpathlineto{\pgfqpoint{19.473cm}{20.849cm}}
\pgfpathlineto{\pgfqpoint{0.635cm}{20.849cm}}
\pgfpathclose
\pgfusepath{clip}
\begin{pgfscope}
\begin{pgfscope}
\pgfpathmoveto{\pgfqpoint{0.635cm}{20.864cm}}
\pgfpathlineto{\pgfqpoint{0.635cm}{6.344cm}}
\pgfpathlineto{\pgfqpoint{19.491cm}{6.344cm}}
\pgfpathlineto{\pgfqpoint{19.491cm}{20.864cm}}
\pgfpathclose
\pgfusepath{clip}
\definecolor{eps2pgf_color}{gray}{1}\pgfsetstrokecolor{eps2pgf_color}\pgfsetfillcolor{eps2pgf_color}
\pgfpathmoveto{\pgfqpoint{0.635cm}{21.59cm}}
\pgfpathlineto{\pgfqpoint{0.635cm}{6.341cm}}
\pgfpathlineto{\pgfqpoint{20.99cm}{6.341cm}}
\pgfpathlineto{\pgfqpoint{20.99cm}{21.59cm}}
\pgfpathclose
\pgfusepath{fill}
\pgfpathmoveto{\pgfqpoint{4.777cm}{8.614cm}}
\pgfpathlineto{\pgfqpoint{4.777cm}{20.446cm}}
\pgfpathlineto{\pgfqpoint{19.053cm}{20.446cm}}
\pgfpathlineto{\pgfqpoint{19.053cm}{8.614cm}}
\pgfpathclose
\pgfseteorule\pgfusepath{fill}\pgfsetnonzerorule
\pgfsetdash{}{0cm}
\pgfsetlinewidth{0.176mm}
\pgfsetroundjoin
\pgfpathmoveto{\pgfqpoint{4.777cm}{8.614cm}}
\pgfpathlineto{\pgfqpoint{4.777cm}{20.446cm}}
\pgfpathlineto{\pgfqpoint{19.053cm}{20.446cm}}
\pgfpathlineto{\pgfqpoint{19.053cm}{8.614cm}}
\pgfpathlineto{\pgfqpoint{4.777cm}{8.614cm}}
\pgfusepath{stroke}
\pgfsetdash{}{0cm}
\definecolor{eps2pgf_color}{gray}{0}\pgfsetstrokecolor{eps2pgf_color}\pgfsetfillcolor{eps2pgf_color}
\pgfpathmoveto{\pgfqpoint{4.777cm}{8.614cm}}
\pgfpathlineto{\pgfqpoint{19.053cm}{8.614cm}}
\pgfusepath{stroke}
\pgfsetdash{}{0cm}
\pgfpathmoveto{\pgfqpoint{4.777cm}{20.446cm}}
\pgfpathlineto{\pgfqpoint{19.053cm}{20.446cm}}
\pgfusepath{stroke}
\pgfsetdash{}{0cm}
\pgfpathmoveto{\pgfqpoint{4.777cm}{8.614cm}}
\pgfpathlineto{\pgfqpoint{4.777cm}{20.446cm}}
\pgfusepath{stroke}
\pgfsetdash{}{0cm}
\pgfpathmoveto{\pgfqpoint{19.053cm}{8.614cm}}
\pgfpathlineto{\pgfqpoint{19.053cm}{20.446cm}}
\pgfusepath{stroke}
\pgfsetdash{}{0cm}
\pgfpathmoveto{\pgfqpoint{4.777cm}{8.614cm}}
\pgfpathlineto{\pgfqpoint{19.053cm}{8.614cm}}
\pgfusepath{stroke}
\pgfsetdash{}{0cm}
\pgfpathmoveto{\pgfqpoint{4.777cm}{8.614cm}}
\pgfpathlineto{\pgfqpoint{4.777cm}{20.446cm}}
\pgfusepath{stroke}
\pgfsetdash{}{0cm}
\pgfpathmoveto{\pgfqpoint{4.777cm}{8.614cm}}
\pgfpathlineto{\pgfqpoint{4.777cm}{8.758cm}}
\pgfusepath{stroke}
\pgfsetdash{}{0cm}
\pgfpathmoveto{\pgfqpoint{4.777cm}{20.446cm}}
\pgfpathlineto{\pgfqpoint{4.777cm}{20.305cm}}
\pgfusepath{stroke}
\pgftext[x=4.777cm,y=7.891cm+.2cm,rotate=0]{  \fontsize{36}{36.14}\selectfont{ {0}}}
\pgfsetdash{}{0cm}
\pgfpathmoveto{\pgfqpoint{11.915cm}{8.614cm}}
\pgfpathlineto{\pgfqpoint{11.915cm}{8.758cm}}
\pgfusepath{stroke}
\pgfsetdash{}{0cm}
\pgfpathmoveto{\pgfqpoint{11.915cm}{20.446cm}}
\pgfpathlineto{\pgfqpoint{11.915cm}{20.305cm}}
\pgfusepath{stroke}
\pgftext[x=11.913cm,y=7.891cm+.2cm,rotate=0]{  \fontsize{36}{36.14}\selectfont{ {0.5}}}
\pgfsetdash{}{0cm}
\pgfpathmoveto{\pgfqpoint{19.053cm}{8.614cm}}
\pgfpathlineto{\pgfqpoint{19.053cm}{8.758cm}}
\pgfusepath{stroke}
\pgfsetdash{}{0cm}
\pgfpathmoveto{\pgfqpoint{19.053cm}{20.446cm}}
\pgfpathlineto{\pgfqpoint{19.053cm}{20.305cm}}
\pgfusepath{stroke}
\pgftext[x=19.002cm-.4cm,y=7.901cm+.2cm,rotate=0]{  \fontsize{36}{36.14}\selectfont{ {1}}}
\pgfsetdash{}{0cm}
\pgfpathmoveto{\pgfqpoint{4.777cm}{8.614cm}}
\pgfpathlineto{\pgfqpoint{4.918cm}{8.614cm}}
\pgfusepath{stroke}
\pgfsetdash{}{0cm}
\pgfpathmoveto{\pgfqpoint{19.053cm}{8.614cm}}
\pgfpathlineto{\pgfqpoint{18.909cm}{8.614cm}}
\pgfusepath{stroke}
\pgftext[x=4.381cm-.3cm,y=8.582cm,rotate=0]{  \fontsize{36}{36.14}\selectfont{ {0}}}
\pgfsetdash{}{0cm}
\pgfpathmoveto{\pgfqpoint{4.777cm}{12.121cm}}
\pgfpathlineto{\pgfqpoint{4.918cm}{12.121cm}}
\pgfusepath{stroke}
\pgfsetdash{}{0cm}
\pgfpathmoveto{\pgfqpoint{19.053cm}{12.121cm}}
\pgfpathlineto{\pgfqpoint{18.909cm}{12.121cm}}
\pgfusepath{stroke}
\pgftext[x=3.35cm-.3cm,y=12.089cm,rotate=0]{  \fontsize{36}{36.14}\selectfont{ {0.005}}}
\pgfsetdash{}{0cm}
\pgfpathmoveto{\pgfqpoint{4.777cm}{15.628cm}}
\pgfpathlineto{\pgfqpoint{4.918cm}{15.628cm}}
\pgfusepath{stroke}
\pgfsetdash{}{0cm}
\pgfpathmoveto{\pgfqpoint{19.053cm}{15.628cm}}
\pgfpathlineto{\pgfqpoint{18.909cm}{15.628cm}}
\pgfusepath{stroke}
\pgftext[x=3.562cm-.3cm,y=15.596cm,rotate=0]{  \fontsize{36}{36.14}\selectfont{ {0.01}}}
\pgfsetdash{}{0cm}
\pgfpathmoveto{\pgfqpoint{4.777cm}{19.132cm}}
\pgfpathlineto{\pgfqpoint{4.918cm}{19.132cm}}
\pgfusepath{stroke}
\pgfsetdash{}{0cm}
\pgfpathmoveto{\pgfqpoint{19.053cm}{19.132cm}}
\pgfpathlineto{\pgfqpoint{18.909cm}{19.132cm}}
\pgfusepath{stroke}
\pgftext[x=3.35cm-.3cm,y=19.1cm,rotate=0]{  \fontsize{36}{36.14}\selectfont{ {0.015}}}
\pgfsetdash{}{0cm}
\pgfpathmoveto{\pgfqpoint{4.777cm}{8.614cm}}
\pgfpathlineto{\pgfqpoint{19.053cm}{8.614cm}}
\pgfusepath{stroke}
\pgfsetdash{}{0cm}
\pgfpathmoveto{\pgfqpoint{4.777cm}{20.446cm}}
\pgfpathlineto{\pgfqpoint{19.053cm}{20.446cm}}
\pgfusepath{stroke}
\pgfsetdash{}{0cm}
\pgfpathmoveto{\pgfqpoint{4.777cm}{8.614cm}}
\pgfpathlineto{\pgfqpoint{4.777cm}{20.446cm}}
\pgfusepath{stroke}
\pgfsetdash{}{0cm}
\pgfpathmoveto{\pgfqpoint{19.053cm}{8.614cm}}
\pgfpathlineto{\pgfqpoint{19.053cm}{20.446cm}}
\pgfusepath{stroke}
\begin{pgfscope}
\pgfpathmoveto{\pgfqpoint{4.777cm}{20.446cm}}
\pgfpathlineto{\pgfqpoint{4.777cm}{8.611cm}}
\pgfpathlineto{\pgfqpoint{19.056cm}{8.611cm}}
\pgfpathlineto{\pgfqpoint{19.056cm}{20.446cm}}
\pgfpathclose
\pgfusepath{clip}
\pgfsetdash{}{0cm}
\pgfsetlinewidth{1.058mm}
\definecolor{eps2pgf_color}{rgb}{0,0,1}\pgfsetstrokecolor{eps2pgf_color}\pgfsetfillcolor{eps2pgf_color}
\pgfpathmoveto{\pgfqpoint{6.203cm}{9.493cm}}
\pgfpathlineto{\pgfqpoint{7.632cm}{10.369cm}}
\pgfpathlineto{\pgfqpoint{9.058cm}{11.245cm}}
\pgfpathlineto{\pgfqpoint{10.486cm}{12.121cm}}
\pgfpathlineto{\pgfqpoint{11.915cm}{12.997cm}}
\pgfpathlineto{\pgfqpoint{13.341cm}{13.873cm}}
\pgfpathlineto{\pgfqpoint{14.77cm}{14.752cm}}
\pgfpathlineto{\pgfqpoint{16.195cm}{15.628cm}}
\pgfpathlineto{\pgfqpoint{17.624cm}{16.504cm}}
\pgfpathlineto{\pgfqpoint{19.053cm}{17.38cm}}
\pgfusepath{stroke}
\end{pgfscope}
\pgftext[x=11.933cm,y=6.968cm,rotate=0]{  \fontsize{40}{50}{\selectfont $\bmax$}}
\pgftext[x=1.257cm,y=14.53cm,rotate=90]{  \fontsize{36}{36.14}\selectfont{ {Sub-optimality bound}}}
\pgftext[x=4.727cm,y=8.487cm,rotate=0]{\fontsize{10.04}{12.04}\selectfont{ { }}}
\pgftext[x=19.006cm,y=20.323cm,rotate=0]{\fontsize{10.04}{12.04}\selectfont{ { }}}
\pgftext[x=8.249cm,y=19.581cm,rotate=0]{  \fontsize{40}{36.14}\selectfont{ {$\la =1$}}}
\begin{pgfscope}
\pgfpathmoveto{\pgfqpoint{4.951cm}{20.27cm}}
\pgfpathlineto{\pgfqpoint{4.951cm}{18.897cm}}
\pgfpathlineto{\pgfqpoint{10.498cm}{18.897cm}}
\pgfpathlineto{\pgfqpoint{10.498cm}{20.27cm}}
\pgfpathclose
\pgfusepath{clip}
\pgfsetdash{}{0cm}
\pgfsetlinewidth{1.058mm}
\definecolor{eps2pgf_color}{rgb}{0,0,1}\pgfsetstrokecolor{eps2pgf_color}\pgfsetfillcolor{eps2pgf_color}
\pgfpathmoveto{\pgfqpoint{5.153cm}{19.585cm}}
\pgfpathlineto{\pgfqpoint{6.182cm}{19.585cm}}
\pgfusepath{stroke}
\end{pgfscope}
\pgfsetdash{}{0cm}
\definecolor{eps2pgf_color}{rgb}{0,0,1}\pgfsetstrokecolor{eps2pgf_color}\pgfsetfillcolor{eps2pgf_color}
\pgfusepath{stroke}
\end{pgfscope}
\end{pgfscope}
\end{pgfpicture}

%% file: fig/case1_2_1.tex
\scalebox{0.235}{\scalefont{2} \input{./fig/case1_2_1.pgf}}

%% file: fig/case1_2_1.pgf
% Created by Eps2pgf 0.7.0 (build on 2008-08-24) on Mon Apr 06 15:10:32 PDT 2015
\begin{pgfpicture}
\pgfpathmoveto{\pgfqpoint{0.917cm}{6.562cm}}
\pgfpathlineto{\pgfqpoint{19.638cm}{6.562cm}}
\pgfpathlineto{\pgfqpoint{19.638cm}{21.378cm}}
\pgfpathlineto{\pgfqpoint{0.917cm}{21.378cm}}
\pgfpathclose
\pgfusepath{clip}
\begin{pgfscope}
\begin{pgfscope}
\pgfpathmoveto{\pgfqpoint{0.917cm}{21.378cm}}
\pgfpathlineto{\pgfqpoint{19.67cm}{21.378cm}}
\pgfpathlineto{\pgfqpoint{19.67cm}{6.582cm}}
\pgfpathlineto{\pgfqpoint{0.917cm}{6.582cm}}
\pgfpathclose
\pgfusepath{clip}
\begin{pgfscope}
\definecolor{eps2pgf_color}{gray}{1}\pgfsetstrokecolor{eps2pgf_color}\pgfsetfillcolor{eps2pgf_color}
\pgfpathmoveto{\pgfqpoint{0.917cm}{21.378cm}}
\pgfpathlineto{\pgfqpoint{19.673cm}{21.378cm}}
\pgfpathlineto{\pgfqpoint{19.673cm}{6.579cm}}
\pgfpathlineto{\pgfqpoint{0.917cm}{6.579cm}}
\pgfpathclose
\pgfusepath{fill}
\end{pgfscope}
\definecolor{eps2pgf_color}{gray}{1}\pgfsetstrokecolor{eps2pgf_color}\pgfsetfillcolor{eps2pgf_color}
\pgfpathmoveto{\pgfqpoint{4.471cm}{8.852cm}}
\pgfpathlineto{\pgfqpoint{4.471cm}{20.27cm}}
\pgfpathlineto{\pgfqpoint{18.794cm}{20.27cm}}
\pgfpathlineto{\pgfqpoint{18.794cm}{8.852cm}}
\pgfpathclose
\pgfseteorule\pgfusepath{fill}\pgfsetnonzerorule
\pgfsetdash{}{0cm}
\pgfsetlinewidth{0.176mm}
\pgfsetroundjoin
\pgfpathmoveto{\pgfqpoint{4.471cm}{8.852cm}}
\pgfpathlineto{\pgfqpoint{4.471cm}{20.27cm}}
\pgfpathlineto{\pgfqpoint{18.794cm}{20.27cm}}
\pgfpathlineto{\pgfqpoint{18.794cm}{8.852cm}}
\pgfpathlineto{\pgfqpoint{4.471cm}{8.852cm}}
\pgfusepath{stroke}
\pgfsetdash{}{0cm}
\definecolor{eps2pgf_color}{gray}{0}\pgfsetstrokecolor{eps2pgf_color}\pgfsetfillcolor{eps2pgf_color}
\pgfpathmoveto{\pgfqpoint{4.471cm}{8.852cm}}
\pgfpathlineto{\pgfqpoint{18.794cm}{8.852cm}}
\pgfusepath{stroke}
\pgfsetdash{}{0cm}
\pgfpathmoveto{\pgfqpoint{4.471cm}{20.27cm}}
\pgfpathlineto{\pgfqpoint{18.794cm}{20.27cm}}
\pgfusepath{stroke}
\pgfsetdash{}{0cm}
\pgfpathmoveto{\pgfqpoint{4.471cm}{8.852cm}}
\pgfpathlineto{\pgfqpoint{4.471cm}{20.27cm}}
\pgfusepath{stroke}
\pgfsetdash{}{0cm}
\pgfpathmoveto{\pgfqpoint{18.794cm}{8.852cm}}
\pgfpathlineto{\pgfqpoint{18.794cm}{20.27cm}}
\pgfusepath{stroke}
\pgfsetdash{}{0cm}
\pgfpathmoveto{\pgfqpoint{4.471cm}{8.852cm}}
\pgfpathlineto{\pgfqpoint{18.794cm}{8.852cm}}
\pgfusepath{stroke}
\pgfsetdash{}{0cm}
\pgfpathmoveto{\pgfqpoint{4.471cm}{8.852cm}}
\pgfpathlineto{\pgfqpoint{4.471cm}{20.27cm}}
\pgfusepath{stroke}
\pgfsetdash{}{0cm}
\pgfpathmoveto{\pgfqpoint{4.471cm}{8.852cm}}
\pgfpathlineto{\pgfqpoint{4.471cm}{8.996cm}}
\pgfusepath{stroke}
\pgfsetdash{}{0cm}
\pgfpathmoveto{\pgfqpoint{4.471cm}{20.273cm}}
\pgfpathlineto{\pgfqpoint{4.471cm}{20.129cm}}
\pgfusepath{stroke}
\pgftext[x=4.472cm,y=8.129cm+.2cm,rotate=0]{ \fontsize{36}{36.14}\selectfont{ {0}}}
\pgfsetdash{}{0cm}
\pgfpathmoveto{\pgfqpoint{11.633cm}{8.852cm}}
\pgfpathlineto{\pgfqpoint{11.633cm}{8.996cm}}
\pgfusepath{stroke}
\pgfsetdash{}{0cm}
\pgfpathmoveto{\pgfqpoint{11.633cm}{20.273cm}}
\pgfpathlineto{\pgfqpoint{11.633cm}{20.129cm}}
\pgfusepath{stroke}
\pgftext[x=11.631cm,y=8.129cm+.2cm,rotate=0]{ \fontsize{36}{36.14}\selectfont{ {0.5}}}
\pgfsetdash{}{0cm}
\pgfpathmoveto{\pgfqpoint{18.794cm}{8.852cm}}
\pgfpathlineto{\pgfqpoint{18.794cm}{8.996cm}}
\pgfusepath{stroke}
\pgfsetdash{}{0cm}
\pgfpathmoveto{\pgfqpoint{18.794cm}{20.273cm}}
\pgfpathlineto{\pgfqpoint{18.794cm}{20.129cm}}
\pgfusepath{stroke}
\pgftext[x=18.744cm-.4cm,y=8.139cm+.2cm,rotate=0]{ \fontsize{36}{36.14}\selectfont{ {1}}}
\pgfsetdash{}{0cm}
\pgfpathmoveto{\pgfqpoint{4.471cm}{8.852cm}}
\pgfpathlineto{\pgfqpoint{4.613cm}{8.852cm}}
\pgfusepath{stroke}
\pgfsetdash{}{0cm}
\pgfpathmoveto{\pgfqpoint{18.794cm}{8.852cm}}
\pgfpathlineto{\pgfqpoint{18.65cm}{8.852cm}}
\pgfusepath{stroke}
\pgftext[x=4.075cm-.3cm,y=8.82cm,rotate=0]{ \fontsize{36}{36.14}\selectfont{ {0}}}
\pgfsetdash{}{0cm}
\pgfpathmoveto{\pgfqpoint{4.471cm}{12.659cm}}
\pgfpathlineto{\pgfqpoint{4.613cm}{12.659cm}}
\pgfusepath{stroke}
\pgfsetdash{}{0cm}
\pgfpathmoveto{\pgfqpoint{18.794cm}{12.659cm}}
\pgfpathlineto{\pgfqpoint{18.65cm}{12.659cm}}
\pgfusepath{stroke}
\pgftext[x=3.338cm-.3cm,y=12.627cm,rotate=0]{ \fontsize{36}{36.14}\selectfont{ {0.05}}}
\pgfsetdash{}{0cm}
\pgfpathmoveto{\pgfqpoint{4.471cm}{16.466cm}}
\pgfpathlineto{\pgfqpoint{4.613cm}{16.466cm}}
\pgfusepath{stroke}
\pgfsetdash{}{0cm}
\pgfpathmoveto{\pgfqpoint{18.794cm}{16.466cm}}
\pgfpathlineto{\pgfqpoint{18.65cm}{16.466cm}}
\pgfusepath{stroke}
\pgftext[x=3.55cm-.3cm,y=16.434cm,rotate=0]{ \fontsize{36}{36.14}\selectfont{ {0.1}}}
\pgfsetdash{}{0cm}
\pgfpathmoveto{\pgfqpoint{4.471cm}{20.273cm}}
\pgfpathlineto{\pgfqpoint{4.613cm}{20.273cm}}
\pgfusepath{stroke}
\pgfsetdash{}{0cm}
\pgfpathmoveto{\pgfqpoint{18.794cm}{20.273cm}}
\pgfpathlineto{\pgfqpoint{18.65cm}{20.273cm}}
\pgfusepath{stroke}
\pgftext[x=3.338cm-.3cm,y=20.241cm-.3cm,rotate=0]{ \fontsize{36}{36.14}\selectfont{ {0.15}}}
\pgfsetdash{}{0cm}
\pgfpathmoveto{\pgfqpoint{4.471cm}{8.852cm}}
\pgfpathlineto{\pgfqpoint{18.794cm}{8.852cm}}
\pgfusepath{stroke}
\pgfsetdash{}{0cm}
\pgfpathmoveto{\pgfqpoint{4.471cm}{20.27cm}}
\pgfpathlineto{\pgfqpoint{18.794cm}{20.27cm}}
\pgfusepath{stroke}
\pgfsetdash{}{0cm}
\pgfpathmoveto{\pgfqpoint{4.471cm}{8.852cm}}
\pgfpathlineto{\pgfqpoint{4.471cm}{20.27cm}}
\pgfusepath{stroke}
\pgfsetdash{}{0cm}
\pgfpathmoveto{\pgfqpoint{18.794cm}{8.852cm}}
\pgfpathlineto{\pgfqpoint{18.794cm}{20.27cm}}
\pgfusepath{stroke}
\begin{pgfscope}
\pgfpathmoveto{\pgfqpoint{4.471cm}{20.27cm}}
\pgfpathlineto{\pgfqpoint{18.797cm}{20.27cm}}
\pgfpathlineto{\pgfqpoint{18.797cm}{8.849cm}}
\pgfpathlineto{\pgfqpoint{4.471cm}{8.849cm}}
\pgfpathclose
\pgfusepath{clip}
\pgfsetdash{}{0cm}
\pgfsetlinewidth{1.058mm}
\definecolor{eps2pgf_color}{rgb}{0,0,1}\pgfsetstrokecolor{eps2pgf_color}\pgfsetfillcolor{eps2pgf_color}
\pgfpathmoveto{\pgfqpoint{5.903cm}{9.587cm}}
\pgfpathlineto{\pgfqpoint{7.335cm}{10.322cm}}
\pgfpathlineto{\pgfqpoint{8.767cm}{11.057cm}}
\pgfpathlineto{\pgfqpoint{10.198cm}{11.789cm}}
\pgfpathlineto{\pgfqpoint{11.633cm}{12.524cm}}
\pgfpathlineto{\pgfqpoint{13.065cm}{13.259cm}}
\pgfpathlineto{\pgfqpoint{14.496cm}{13.991cm}}
\pgfpathlineto{\pgfqpoint{15.928cm}{14.726cm}}
\pgfpathlineto{\pgfqpoint{17.36cm}{15.46cm}}
\pgfpathlineto{\pgfqpoint{18.794cm}{16.192cm}}
\pgfusepath{stroke}
\end{pgfscope}
\pgfsetdash{}{0cm}
\pgfsetlinewidth{1.058mm}
\definecolor{eps2pgf_color}{rgb}{0,0,1}\pgfsetstrokecolor{eps2pgf_color}\pgfsetfillcolor{eps2pgf_color}
\pgfpathmoveto{\pgfqpoint{5.762cm}{9.587cm}}
\pgfpathlineto{\pgfqpoint{6.044cm}{9.587cm}}
\pgfusepath{stroke}
\pgfsetdash{}{0cm}
\pgfpathmoveto{\pgfqpoint{5.903cm}{9.728cm}}
\pgfpathlineto{\pgfqpoint{5.903cm}{9.446cm}}
\pgfusepath{stroke}
\pgfsetdash{}{0cm}
\pgfpathmoveto{\pgfqpoint{7.194cm}{10.322cm}}
\pgfpathlineto{\pgfqpoint{7.476cm}{10.322cm}}
\pgfusepath{stroke}
\pgfsetdash{}{0cm}
\pgfpathmoveto{\pgfqpoint{7.335cm}{10.463cm}}
\pgfpathlineto{\pgfqpoint{7.335cm}{10.181cm}}
\pgfusepath{stroke}
\pgfsetdash{}{0cm}
\pgfpathmoveto{\pgfqpoint{8.625cm}{11.057cm}}
\pgfpathlineto{\pgfqpoint{8.908cm}{11.057cm}}
\pgfusepath{stroke}
\pgfsetdash{}{0cm}
\pgfpathmoveto{\pgfqpoint{8.767cm}{11.198cm}}
\pgfpathlineto{\pgfqpoint{8.767cm}{10.916cm}}
\pgfusepath{stroke}
\pgfsetdash{}{0cm}
\pgfpathmoveto{\pgfqpoint{10.057cm}{11.789cm}}
\pgfpathlineto{\pgfqpoint{10.339cm}{11.789cm}}
\pgfusepath{stroke}
\pgfsetdash{}{0cm}
\pgfpathmoveto{\pgfqpoint{10.198cm}{11.93cm}}
\pgfpathlineto{\pgfqpoint{10.198cm}{11.648cm}}
\pgfusepath{stroke}
\pgfsetdash{}{0cm}
\pgfpathmoveto{\pgfqpoint{11.492cm}{12.524cm}}
\pgfpathlineto{\pgfqpoint{11.774cm}{12.524cm}}
\pgfusepath{stroke}
\pgfsetdash{}{0cm}
\pgfpathmoveto{\pgfqpoint{11.633cm}{12.665cm}}
\pgfpathlineto{\pgfqpoint{11.633cm}{12.382cm}}
\pgfusepath{stroke}
\pgfsetdash{}{0cm}
\pgfpathmoveto{\pgfqpoint{12.923cm}{13.259cm}}
\pgfpathlineto{\pgfqpoint{13.206cm}{13.259cm}}
\pgfusepath{stroke}
\pgfsetdash{}{0cm}
\pgfpathmoveto{\pgfqpoint{13.065cm}{13.4cm}}
\pgfpathlineto{\pgfqpoint{13.065cm}{13.117cm}}
\pgfusepath{stroke}
\pgfsetdash{}{0cm}
\pgfpathmoveto{\pgfqpoint{14.355cm}{13.991cm}}
\pgfpathlineto{\pgfqpoint{14.637cm}{13.991cm}}
\pgfusepath{stroke}
\pgfsetdash{}{0cm}
\pgfpathmoveto{\pgfqpoint{14.496cm}{14.132cm}}
\pgfpathlineto{\pgfqpoint{14.496cm}{13.849cm}}
\pgfusepath{stroke}
\pgfsetdash{}{0cm}
\pgfpathmoveto{\pgfqpoint{15.787cm}{14.726cm}}
\pgfpathlineto{\pgfqpoint{16.069cm}{14.726cm}}
\pgfusepath{stroke}
\pgfsetdash{}{0cm}
\pgfpathmoveto{\pgfqpoint{15.928cm}{14.867cm}}
\pgfpathlineto{\pgfqpoint{15.928cm}{14.584cm}}
\pgfusepath{stroke}
\pgfsetdash{}{0cm}
\pgfpathmoveto{\pgfqpoint{17.218cm}{15.46cm}}
\pgfpathlineto{\pgfqpoint{17.501cm}{15.46cm}}
\pgfusepath{stroke}
\pgfsetdash{}{0cm}
\pgfpathmoveto{\pgfqpoint{17.36cm}{15.602cm}}
\pgfpathlineto{\pgfqpoint{17.36cm}{15.319cm}}
\pgfusepath{stroke}
\pgfsetdash{}{0cm}
\pgfpathmoveto{\pgfqpoint{18.653cm}{16.192cm}}
\pgfpathlineto{\pgfqpoint{18.935cm}{16.192cm}}
\pgfusepath{stroke}
\pgfsetdash{}{0cm}
\pgfpathmoveto{\pgfqpoint{18.794cm}{16.334cm}}
\pgfpathlineto{\pgfqpoint{18.794cm}{16.051cm}}
\pgfusepath{stroke}
\pgfsetdash{}{0cm}
\pgfpathmoveto{\pgfqpoint{5.806cm}{9.684cm}}
\pgfpathlineto{\pgfqpoint{6cm}{9.49cm}}
\pgfusepath{stroke}
\pgfsetdash{}{0cm}
\pgfpathmoveto{\pgfqpoint{6cm}{9.684cm}}
\pgfpathlineto{\pgfqpoint{5.806cm}{9.49cm}}
\pgfusepath{stroke}
\pgfsetdash{}{0cm}
\pgfpathmoveto{\pgfqpoint{7.238cm}{10.419cm}}
\pgfpathlineto{\pgfqpoint{7.432cm}{10.225cm}}
\pgfusepath{stroke}
\pgfsetdash{}{0cm}
\pgfpathmoveto{\pgfqpoint{7.432cm}{10.419cm}}
\pgfpathlineto{\pgfqpoint{7.238cm}{10.225cm}}
\pgfusepath{stroke}
\pgfsetdash{}{0cm}
\pgfpathmoveto{\pgfqpoint{8.67cm}{11.154cm}}
\pgfpathlineto{\pgfqpoint{8.864cm}{10.96cm}}
\pgfusepath{stroke}
\pgfsetdash{}{0cm}
\pgfpathmoveto{\pgfqpoint{8.864cm}{11.154cm}}
\pgfpathlineto{\pgfqpoint{8.67cm}{10.96cm}}
\pgfusepath{stroke}
\pgfsetdash{}{0cm}
\pgfpathmoveto{\pgfqpoint{10.101cm}{11.886cm}}
\pgfpathlineto{\pgfqpoint{10.295cm}{11.692cm}}
\pgfusepath{stroke}
\pgfsetdash{}{0cm}
\pgfpathmoveto{\pgfqpoint{10.295cm}{11.886cm}}
\pgfpathlineto{\pgfqpoint{10.101cm}{11.692cm}}
\pgfusepath{stroke}
\pgfsetdash{}{0cm}
\pgfpathmoveto{\pgfqpoint{11.536cm}{12.621cm}}
\pgfpathlineto{\pgfqpoint{11.73cm}{12.427cm}}
\pgfusepath{stroke}
\pgfsetdash{}{0cm}
\pgfpathmoveto{\pgfqpoint{11.73cm}{12.621cm}}
\pgfpathlineto{\pgfqpoint{11.536cm}{12.427cm}}
\pgfusepath{stroke}
\pgfsetdash{}{0cm}
\pgfpathmoveto{\pgfqpoint{12.968cm}{13.356cm}}
\pgfpathlineto{\pgfqpoint{13.162cm}{13.162cm}}
\pgfusepath{stroke}
\pgfsetdash{}{0cm}
\pgfpathmoveto{\pgfqpoint{13.162cm}{13.356cm}}
\pgfpathlineto{\pgfqpoint{12.968cm}{13.162cm}}
\pgfusepath{stroke}
\pgfsetdash{}{0cm}
\pgfpathmoveto{\pgfqpoint{14.399cm}{14.088cm}}
\pgfpathlineto{\pgfqpoint{14.593cm}{13.894cm}}
\pgfusepath{stroke}
\pgfsetdash{}{0cm}
\pgfpathmoveto{\pgfqpoint{14.593cm}{14.088cm}}
\pgfpathlineto{\pgfqpoint{14.399cm}{13.894cm}}
\pgfusepath{stroke}
\pgfsetdash{}{0cm}
\pgfpathmoveto{\pgfqpoint{15.831cm}{14.823cm}}
\pgfpathlineto{\pgfqpoint{16.025cm}{14.629cm}}
\pgfusepath{stroke}
\pgfsetdash{}{0cm}
\pgfpathmoveto{\pgfqpoint{16.025cm}{14.823cm}}
\pgfpathlineto{\pgfqpoint{15.831cm}{14.629cm}}
\pgfusepath{stroke}
\pgfsetdash{}{0cm}
\pgfpathmoveto{\pgfqpoint{17.263cm}{15.558cm}}
\pgfpathlineto{\pgfqpoint{17.457cm}{15.363cm}}
\pgfusepath{stroke}
\pgfsetdash{}{0cm}
\pgfpathmoveto{\pgfqpoint{17.457cm}{15.558cm}}
\pgfpathlineto{\pgfqpoint{17.263cm}{15.363cm}}
\pgfusepath{stroke}
\pgfsetdash{}{0cm}
\pgfpathmoveto{\pgfqpoint{18.697cm}{16.29cm}}
\pgfpathlineto{\pgfqpoint{18.891cm}{16.095cm}}
\pgfusepath{stroke}
\pgfsetdash{}{0cm}
\pgfpathmoveto{\pgfqpoint{18.891cm}{16.29cm}}
\pgfpathlineto{\pgfqpoint{18.697cm}{16.095cm}}
\pgfusepath{stroke}
\begin{pgfscope}
\pgfpathmoveto{\pgfqpoint{4.471cm}{20.27cm}}
\pgfpathlineto{\pgfqpoint{18.797cm}{20.27cm}}
\pgfpathlineto{\pgfqpoint{18.797cm}{8.849cm}}
\pgfpathlineto{\pgfqpoint{4.471cm}{8.849cm}}
\pgfpathclose
\pgfusepath{clip}
\pgfsetdash{{0.018cm}{0.141cm}{0.212cm}{0.141cm}}{0cm}
\pgfpathmoveto{\pgfqpoint{5.903cm}{9.613cm}}
\pgfpathlineto{\pgfqpoint{7.335cm}{10.372cm}}
\pgfpathlineto{\pgfqpoint{8.767cm}{11.13cm}}
\pgfpathlineto{\pgfqpoint{10.198cm}{11.889cm}}
\pgfpathlineto{\pgfqpoint{11.633cm}{12.647cm}}
\pgfpathlineto{\pgfqpoint{13.065cm}{13.406cm}}
\pgfpathlineto{\pgfqpoint{14.496cm}{14.164cm}}
\pgfpathlineto{\pgfqpoint{15.928cm}{14.925cm}}
\pgfpathlineto{\pgfqpoint{17.36cm}{15.684cm}}
\pgfpathlineto{\pgfqpoint{18.794cm}{16.442cm}}
\pgfusepath{stroke}
\end{pgfscope}
\pgfsetdash{}{0cm}
\pgfpathmoveto{\pgfqpoint{5.762cm}{9.613cm}}
\pgfpathlineto{\pgfqpoint{6.044cm}{9.613cm}}
\pgfusepath{stroke}
\pgfsetdash{}{0cm}
\pgfpathmoveto{\pgfqpoint{5.903cm}{9.754cm}}
\pgfpathlineto{\pgfqpoint{5.903cm}{9.472cm}}
\pgfusepath{stroke}
\pgfsetdash{}{0cm}
\pgfpathmoveto{\pgfqpoint{7.194cm}{10.372cm}}
\pgfpathlineto{\pgfqpoint{7.476cm}{10.372cm}}
\pgfusepath{stroke}
\pgfsetdash{}{0cm}
\pgfpathmoveto{\pgfqpoint{7.335cm}{10.513cm}}
\pgfpathlineto{\pgfqpoint{7.335cm}{10.231cm}}
\pgfusepath{stroke}
\pgfsetdash{}{0cm}
\pgfpathmoveto{\pgfqpoint{8.625cm}{11.13cm}}
\pgfpathlineto{\pgfqpoint{8.908cm}{11.13cm}}
\pgfusepath{stroke}
\pgfsetdash{}{0cm}
\pgfpathmoveto{\pgfqpoint{8.767cm}{11.271cm}}
\pgfpathlineto{\pgfqpoint{8.767cm}{10.989cm}}
\pgfusepath{stroke}
\pgfsetdash{}{0cm}
\pgfpathmoveto{\pgfqpoint{10.057cm}{11.889cm}}
\pgfpathlineto{\pgfqpoint{10.339cm}{11.889cm}}
\pgfusepath{stroke}
\pgfsetdash{}{0cm}
\pgfpathmoveto{\pgfqpoint{10.198cm}{12.03cm}}
\pgfpathlineto{\pgfqpoint{10.198cm}{11.747cm}}
\pgfusepath{stroke}
\pgfsetdash{}{0cm}
\pgfpathmoveto{\pgfqpoint{11.492cm}{12.647cm}}
\pgfpathlineto{\pgfqpoint{11.774cm}{12.647cm}}
\pgfusepath{stroke}
\pgfsetdash{}{0cm}
\pgfpathmoveto{\pgfqpoint{11.633cm}{12.788cm}}
\pgfpathlineto{\pgfqpoint{11.633cm}{12.506cm}}
\pgfusepath{stroke}
\pgfsetdash{}{0cm}
\pgfpathmoveto{\pgfqpoint{12.923cm}{13.406cm}}
\pgfpathlineto{\pgfqpoint{13.206cm}{13.406cm}}
\pgfusepath{stroke}
\pgfsetdash{}{0cm}
\pgfpathmoveto{\pgfqpoint{13.065cm}{13.547cm}}
\pgfpathlineto{\pgfqpoint{13.065cm}{13.264cm}}
\pgfusepath{stroke}
\pgfsetdash{}{0cm}
\pgfpathmoveto{\pgfqpoint{14.355cm}{14.164cm}}
\pgfpathlineto{\pgfqpoint{14.637cm}{14.164cm}}
\pgfusepath{stroke}
\pgfsetdash{}{0cm}
\pgfpathmoveto{\pgfqpoint{14.496cm}{14.305cm}}
\pgfpathlineto{\pgfqpoint{14.496cm}{14.023cm}}
\pgfusepath{stroke}
\pgfsetdash{}{0cm}
\pgfpathmoveto{\pgfqpoint{15.787cm}{14.925cm}}
\pgfpathlineto{\pgfqpoint{16.069cm}{14.925cm}}
\pgfusepath{stroke}
\pgfsetdash{}{0cm}
\pgfpathmoveto{\pgfqpoint{15.928cm}{15.067cm}}
\pgfpathlineto{\pgfqpoint{15.928cm}{14.784cm}}
\pgfusepath{stroke}
\pgfsetdash{}{0cm}
\pgfpathmoveto{\pgfqpoint{17.218cm}{15.684cm}}
\pgfpathlineto{\pgfqpoint{17.501cm}{15.684cm}}
\pgfusepath{stroke}
\pgfsetdash{}{0cm}
\pgfpathmoveto{\pgfqpoint{17.36cm}{15.825cm}}
\pgfpathlineto{\pgfqpoint{17.36cm}{15.543cm}}
\pgfusepath{stroke}
\pgfsetdash{}{0cm}
\pgfpathmoveto{\pgfqpoint{18.653cm}{16.442cm}}
\pgfpathlineto{\pgfqpoint{18.935cm}{16.442cm}}
\pgfusepath{stroke}
\pgfsetdash{}{0cm}
\pgfpathmoveto{\pgfqpoint{18.794cm}{16.583cm}}
\pgfpathlineto{\pgfqpoint{18.794cm}{16.301cm}}
\pgfusepath{stroke}
\pgfsetdash{}{0cm}
\pgfpathmoveto{\pgfqpoint{5.806cm}{9.71cm}}
\pgfpathlineto{\pgfqpoint{6cm}{9.516cm}}
\pgfusepath{stroke}
\pgfsetdash{}{0cm}
\pgfpathmoveto{\pgfqpoint{6cm}{9.71cm}}
\pgfpathlineto{\pgfqpoint{5.806cm}{9.516cm}}
\pgfusepath{stroke}
\pgfsetdash{}{0cm}
\pgfpathmoveto{\pgfqpoint{7.238cm}{10.469cm}}
\pgfpathlineto{\pgfqpoint{7.432cm}{10.275cm}}
\pgfusepath{stroke}
\pgfsetdash{}{0cm}
\pgfpathmoveto{\pgfqpoint{7.432cm}{10.469cm}}
\pgfpathlineto{\pgfqpoint{7.238cm}{10.275cm}}
\pgfusepath{stroke}
\pgfsetdash{}{0cm}
\pgfpathmoveto{\pgfqpoint{8.67cm}{11.227cm}}
\pgfpathlineto{\pgfqpoint{8.864cm}{11.033cm}}
\pgfusepath{stroke}
\pgfsetdash{}{0cm}
\pgfpathmoveto{\pgfqpoint{8.864cm}{11.227cm}}
\pgfpathlineto{\pgfqpoint{8.67cm}{11.033cm}}
\pgfusepath{stroke}
\pgfsetdash{}{0cm}
\pgfpathmoveto{\pgfqpoint{10.101cm}{11.986cm}}
\pgfpathlineto{\pgfqpoint{10.295cm}{11.792cm}}
\pgfusepath{stroke}
\pgfsetdash{}{0cm}
\pgfpathmoveto{\pgfqpoint{10.295cm}{11.986cm}}
\pgfpathlineto{\pgfqpoint{10.101cm}{11.792cm}}
\pgfusepath{stroke}
\pgfsetdash{}{0cm}
\pgfpathmoveto{\pgfqpoint{11.536cm}{12.744cm}}
\pgfpathlineto{\pgfqpoint{11.73cm}{12.55cm}}
\pgfusepath{stroke}
\pgfsetdash{}{0cm}
\pgfpathmoveto{\pgfqpoint{11.73cm}{12.744cm}}
\pgfpathlineto{\pgfqpoint{11.536cm}{12.55cm}}
\pgfusepath{stroke}
\pgfsetdash{}{0cm}
\pgfpathmoveto{\pgfqpoint{12.968cm}{13.503cm}}
\pgfpathlineto{\pgfqpoint{13.162cm}{13.309cm}}
\pgfusepath{stroke}
\pgfsetdash{}{0cm}
\pgfpathmoveto{\pgfqpoint{13.162cm}{13.503cm}}
\pgfpathlineto{\pgfqpoint{12.968cm}{13.309cm}}
\pgfusepath{stroke}
\pgfsetdash{}{0cm}
\pgfpathmoveto{\pgfqpoint{14.399cm}{14.261cm}}
\pgfpathlineto{\pgfqpoint{14.593cm}{14.067cm}}
\pgfusepath{stroke}
\pgfsetdash{}{0cm}
\pgfpathmoveto{\pgfqpoint{14.593cm}{14.261cm}}
\pgfpathlineto{\pgfqpoint{14.399cm}{14.067cm}}
\pgfusepath{stroke}
\pgfsetdash{}{0cm}
\pgfpathmoveto{\pgfqpoint{15.831cm}{15.022cm}}
\pgfpathlineto{\pgfqpoint{16.025cm}{14.828cm}}
\pgfusepath{stroke}
\pgfsetdash{}{0cm}
\pgfpathmoveto{\pgfqpoint{16.025cm}{15.022cm}}
\pgfpathlineto{\pgfqpoint{15.831cm}{14.828cm}}
\pgfusepath{stroke}
\pgfsetdash{}{0cm}
\pgfpathmoveto{\pgfqpoint{17.263cm}{15.781cm}}
\pgfpathlineto{\pgfqpoint{17.457cm}{15.587cm}}
\pgfusepath{stroke}
\pgfsetdash{}{0cm}
\pgfpathmoveto{\pgfqpoint{17.457cm}{15.781cm}}
\pgfpathlineto{\pgfqpoint{17.263cm}{15.587cm}}
\pgfusepath{stroke}
\pgfsetdash{}{0cm}
\pgfpathmoveto{\pgfqpoint{18.697cm}{16.539cm}}
\pgfpathlineto{\pgfqpoint{18.891cm}{16.345cm}}
\pgfusepath{stroke}
\pgfsetdash{}{0cm}
\pgfpathmoveto{\pgfqpoint{18.891cm}{16.539cm}}
\pgfpathlineto{\pgfqpoint{18.697cm}{16.345cm}}
\pgfusepath{stroke}
\begin{pgfscope}
\pgfpathmoveto{\pgfqpoint{4.471cm}{20.27cm}}
\pgfpathlineto{\pgfqpoint{18.797cm}{20.27cm}}
\pgfpathlineto{\pgfqpoint{18.797cm}{8.849cm}}
\pgfpathlineto{\pgfqpoint{4.471cm}{8.849cm}}
\pgfpathclose
\pgfusepath{clip}
\pgfsetdash{}{0cm}
\definecolor{eps2pgf_color}{rgb}{0,1,0}\pgfsetstrokecolor{eps2pgf_color}\pgfsetfillcolor{eps2pgf_color}
\pgfpathmoveto{\pgfqpoint{5.903cm}{9.252cm}}
\pgfpathlineto{\pgfqpoint{7.335cm}{9.651cm}}
\pgfpathlineto{\pgfqpoint{8.767cm}{10.051cm}}
\pgfpathlineto{\pgfqpoint{10.198cm}{10.451cm}}
\pgfpathlineto{\pgfqpoint{11.633cm}{10.851cm}}
\pgfpathlineto{\pgfqpoint{13.065cm}{11.251cm}}
\pgfpathlineto{\pgfqpoint{14.496cm}{11.65cm}}
\pgfpathlineto{\pgfqpoint{15.928cm}{12.05cm}}
\pgfpathlineto{\pgfqpoint{17.36cm}{12.45cm}}
\pgfpathlineto{\pgfqpoint{18.794cm}{12.85cm}}
\pgfusepath{stroke}
\end{pgfscope}
\pgfsetdash{}{0cm}
\pgfsetmiterjoin
\definecolor{eps2pgf_color}{rgb}{0,1,0}\pgfsetstrokecolor{eps2pgf_color}\pgfsetfillcolor{eps2pgf_color}
\pgfpathmoveto{\pgfqpoint{5.791cm}{9.363cm}}
\pgfpathlineto{\pgfqpoint{6.015cm}{9.363cm}}
\pgfpathlineto{\pgfqpoint{6.015cm}{9.14cm}}
\pgfpathlineto{\pgfqpoint{5.791cm}{9.14cm}}
\pgfpathlineto{\pgfqpoint{5.791cm}{9.363cm}}
\pgfpathclose
\pgfusepath{stroke}
\pgfsetdash{}{0cm}
\pgfpathmoveto{\pgfqpoint{7.223cm}{9.763cm}}
\pgfpathlineto{\pgfqpoint{7.447cm}{9.763cm}}
\pgfpathlineto{\pgfqpoint{7.447cm}{9.54cm}}
\pgfpathlineto{\pgfqpoint{7.223cm}{9.54cm}}
\pgfpathlineto{\pgfqpoint{7.223cm}{9.763cm}}
\pgfpathclose
\pgfusepath{stroke}
\pgfsetdash{}{0cm}
\pgfpathmoveto{\pgfqpoint{8.655cm}{10.163cm}}
\pgfpathlineto{\pgfqpoint{8.878cm}{10.163cm}}
\pgfpathlineto{\pgfqpoint{8.878cm}{9.94cm}}
\pgfpathlineto{\pgfqpoint{8.655cm}{9.94cm}}
\pgfpathlineto{\pgfqpoint{8.655cm}{10.163cm}}
\pgfpathclose
\pgfusepath{stroke}
\pgfsetdash{}{0cm}
\pgfpathmoveto{\pgfqpoint{10.087cm}{10.563cm}}
\pgfpathlineto{\pgfqpoint{10.31cm}{10.563cm}}
\pgfpathlineto{\pgfqpoint{10.31cm}{10.339cm}}
\pgfpathlineto{\pgfqpoint{10.087cm}{10.339cm}}
\pgfpathlineto{\pgfqpoint{10.087cm}{10.563cm}}
\pgfpathclose
\pgfusepath{stroke}
\pgfsetdash{}{0cm}
\pgfpathmoveto{\pgfqpoint{11.521cm}{10.963cm}}
\pgfpathlineto{\pgfqpoint{11.745cm}{10.963cm}}
\pgfpathlineto{\pgfqpoint{11.745cm}{10.739cm}}
\pgfpathlineto{\pgfqpoint{11.521cm}{10.739cm}}
\pgfpathlineto{\pgfqpoint{11.521cm}{10.963cm}}
\pgfpathclose
\pgfusepath{stroke}
\pgfsetdash{}{0cm}
\pgfpathmoveto{\pgfqpoint{12.953cm}{11.362cm}}
\pgfpathlineto{\pgfqpoint{13.176cm}{11.362cm}}
\pgfpathlineto{\pgfqpoint{13.176cm}{11.139cm}}
\pgfpathlineto{\pgfqpoint{12.953cm}{11.139cm}}
\pgfpathlineto{\pgfqpoint{12.953cm}{11.362cm}}
\pgfpathclose
\pgfusepath{stroke}
\pgfsetdash{}{0cm}
\pgfpathmoveto{\pgfqpoint{14.385cm}{11.762cm}}
\pgfpathlineto{\pgfqpoint{14.608cm}{11.762cm}}
\pgfpathlineto{\pgfqpoint{14.608cm}{11.539cm}}
\pgfpathlineto{\pgfqpoint{14.385cm}{11.539cm}}
\pgfpathlineto{\pgfqpoint{14.385cm}{11.762cm}}
\pgfpathclose
\pgfusepath{stroke}
\pgfsetdash{}{0cm}
\pgfpathmoveto{\pgfqpoint{15.816cm}{12.162cm}}
\pgfpathlineto{\pgfqpoint{16.04cm}{12.162cm}}
\pgfpathlineto{\pgfqpoint{16.04cm}{11.939cm}}
\pgfpathlineto{\pgfqpoint{15.816cm}{11.939cm}}
\pgfpathlineto{\pgfqpoint{15.816cm}{12.162cm}}
\pgfpathclose
\pgfusepath{stroke}
\pgfsetdash{}{0cm}
\pgfpathmoveto{\pgfqpoint{17.248cm}{12.562cm}}
\pgfpathlineto{\pgfqpoint{17.471cm}{12.562cm}}
\pgfpathlineto{\pgfqpoint{17.471cm}{12.338cm}}
\pgfpathlineto{\pgfqpoint{17.248cm}{12.338cm}}
\pgfpathlineto{\pgfqpoint{17.248cm}{12.562cm}}
\pgfpathclose
\pgfusepath{stroke}
\pgfsetdash{}{0cm}
\pgfpathmoveto{\pgfqpoint{18.683cm}{12.962cm}}
\pgfpathlineto{\pgfqpoint{18.906cm}{12.962cm}}
\pgfpathlineto{\pgfqpoint{18.906cm}{12.738cm}}
\pgfpathlineto{\pgfqpoint{18.683cm}{12.738cm}}
\pgfpathlineto{\pgfqpoint{18.683cm}{12.962cm}}
\pgfpathclose
\pgfusepath{stroke}
\begin{pgfscope}
\pgfpathmoveto{\pgfqpoint{4.471cm}{20.27cm}}
\pgfpathlineto{\pgfqpoint{18.797cm}{20.27cm}}
\pgfpathlineto{\pgfqpoint{18.797cm}{8.849cm}}
\pgfpathlineto{\pgfqpoint{4.471cm}{8.849cm}}
\pgfpathclose
\pgfusepath{clip}
\pgfsetdash{{0.018cm}{0.141cm}{0.212cm}{0.141cm}}{0cm}
\pgfpathmoveto{\pgfqpoint{5.903cm}{9.255cm}}
\pgfpathlineto{\pgfqpoint{7.335cm}{9.657cm}}
\pgfpathlineto{\pgfqpoint{8.767cm}{10.06cm}}
\pgfpathlineto{\pgfqpoint{10.198cm}{10.463cm}}
\pgfpathlineto{\pgfqpoint{11.633cm}{10.866cm}}
\pgfpathlineto{\pgfqpoint{13.065cm}{11.265cm}}
\pgfpathlineto{\pgfqpoint{14.496cm}{11.668cm}}
\pgfpathlineto{\pgfqpoint{15.928cm}{12.071cm}}
\pgfpathlineto{\pgfqpoint{17.36cm}{12.474cm}}
\pgfpathlineto{\pgfqpoint{18.794cm}{12.876cm}}
\pgfusepath{stroke}
\end{pgfscope}
\pgfsetdash{}{0cm}
\pgfpathmoveto{\pgfqpoint{5.791cm}{9.366cm}}
\pgfpathlineto{\pgfqpoint{6.015cm}{9.366cm}}
\pgfpathlineto{\pgfqpoint{6.015cm}{9.143cm}}
\pgfpathlineto{\pgfqpoint{5.791cm}{9.143cm}}
\pgfpathlineto{\pgfqpoint{5.791cm}{9.366cm}}
\pgfpathclose
\pgfusepath{stroke}
\pgfsetdash{}{0cm}
\pgfpathmoveto{\pgfqpoint{7.223cm}{9.769cm}}
\pgfpathlineto{\pgfqpoint{7.447cm}{9.769cm}}
\pgfpathlineto{\pgfqpoint{7.447cm}{9.546cm}}
\pgfpathlineto{\pgfqpoint{7.223cm}{9.546cm}}
\pgfpathlineto{\pgfqpoint{7.223cm}{9.769cm}}
\pgfpathclose
\pgfusepath{stroke}
\pgfsetdash{}{0cm}
\pgfpathmoveto{\pgfqpoint{8.655cm}{10.172cm}}
\pgfpathlineto{\pgfqpoint{8.878cm}{10.172cm}}
\pgfpathlineto{\pgfqpoint{8.878cm}{9.948cm}}
\pgfpathlineto{\pgfqpoint{8.655cm}{9.948cm}}
\pgfpathlineto{\pgfqpoint{8.655cm}{10.172cm}}
\pgfpathclose
\pgfusepath{stroke}
\pgfsetdash{}{0cm}
\pgfpathmoveto{\pgfqpoint{10.087cm}{10.575cm}}
\pgfpathlineto{\pgfqpoint{10.31cm}{10.575cm}}
\pgfpathlineto{\pgfqpoint{10.31cm}{10.351cm}}
\pgfpathlineto{\pgfqpoint{10.087cm}{10.351cm}}
\pgfpathlineto{\pgfqpoint{10.087cm}{10.575cm}}
\pgfpathclose
\pgfusepath{stroke}
\pgfsetdash{}{0cm}
\pgfpathmoveto{\pgfqpoint{11.521cm}{10.977cm}}
\pgfpathlineto{\pgfqpoint{11.745cm}{10.977cm}}
\pgfpathlineto{\pgfqpoint{11.745cm}{10.754cm}}
\pgfpathlineto{\pgfqpoint{11.521cm}{10.754cm}}
\pgfpathlineto{\pgfqpoint{11.521cm}{10.977cm}}
\pgfpathclose
\pgfusepath{stroke}
\pgfsetdash{}{0cm}
\pgfpathmoveto{\pgfqpoint{12.953cm}{11.377cm}}
\pgfpathlineto{\pgfqpoint{13.176cm}{11.377cm}}
\pgfpathlineto{\pgfqpoint{13.176cm}{11.154cm}}
\pgfpathlineto{\pgfqpoint{12.953cm}{11.154cm}}
\pgfpathlineto{\pgfqpoint{12.953cm}{11.377cm}}
\pgfpathclose
\pgfusepath{stroke}
\pgfsetdash{}{0cm}
\pgfpathmoveto{\pgfqpoint{14.385cm}{11.78cm}}
\pgfpathlineto{\pgfqpoint{14.608cm}{11.78cm}}
\pgfpathlineto{\pgfqpoint{14.608cm}{11.556cm}}
\pgfpathlineto{\pgfqpoint{14.385cm}{11.556cm}}
\pgfpathlineto{\pgfqpoint{14.385cm}{11.78cm}}
\pgfpathclose
\pgfusepath{stroke}
\pgfsetdash{}{0cm}
\pgfpathmoveto{\pgfqpoint{15.816cm}{12.183cm}}
\pgfpathlineto{\pgfqpoint{16.04cm}{12.183cm}}
\pgfpathlineto{\pgfqpoint{16.04cm}{11.959cm}}
\pgfpathlineto{\pgfqpoint{15.816cm}{11.959cm}}
\pgfpathlineto{\pgfqpoint{15.816cm}{12.183cm}}
\pgfpathclose
\pgfusepath{stroke}
\pgfsetdash{}{0cm}
\pgfpathmoveto{\pgfqpoint{17.248cm}{12.585cm}}
\pgfpathlineto{\pgfqpoint{17.471cm}{12.585cm}}
\pgfpathlineto{\pgfqpoint{17.471cm}{12.362cm}}
\pgfpathlineto{\pgfqpoint{17.248cm}{12.362cm}}
\pgfpathlineto{\pgfqpoint{17.248cm}{12.585cm}}
\pgfpathclose
\pgfusepath{stroke}
\pgfsetdash{}{0cm}
\pgfpathmoveto{\pgfqpoint{18.683cm}{12.988cm}}
\pgfpathlineto{\pgfqpoint{18.906cm}{12.988cm}}
\pgfpathlineto{\pgfqpoint{18.906cm}{12.765cm}}
\pgfpathlineto{\pgfqpoint{18.683cm}{12.765cm}}
\pgfpathlineto{\pgfqpoint{18.683cm}{12.988cm}}
\pgfpathclose
\pgfusepath{stroke}
\begin{pgfscope}
\pgfpathmoveto{\pgfqpoint{4.471cm}{20.27cm}}
\pgfpathlineto{\pgfqpoint{18.797cm}{20.27cm}}
\pgfpathlineto{\pgfqpoint{18.797cm}{8.849cm}}
\pgfpathlineto{\pgfqpoint{4.471cm}{8.849cm}}
\pgfpathclose
\pgfusepath{clip}
\end{pgfscope}
\definecolor{eps2pgf_color}{gray}{0}\pgfsetstrokecolor{eps2pgf_color}\pgfsetfillcolor{eps2pgf_color}
\pgftext[x=11.651cm,y=7.206cm,rotate=0]{ \fontsize{40}{36.14}\selectfont{ {$\bmax$}}}
\pgftext[x=1.539cm,y=14.56cm,rotate=90]{ \fontsize{36}{36.14}\selectfont{ {Sub-optimality bound}}}
\pgftext[x=4.421cm,y=8.725cm,rotate=0]{\fontsize{10.04}{12.04}\selectfont{ { }}}
\pgftext[x=9.91cm+.2cm,y=19.309cm,rotate=0]{ \fontsize{40}{36.14}\selectfont{ {$\la = 0.9$, \texttt{minS}}}}
\begin{pgfscope}
\pgfpathmoveto{\pgfqpoint{4.648cm}{20.094cm}}
\pgfpathlineto{\pgfqpoint{14.875cm}{20.094cm}}
\pgfpathlineto{\pgfqpoint{14.875cm}{14.981cm}}
\pgfpathlineto{\pgfqpoint{4.648cm}{14.981cm}}
\pgfpathclose
\pgfusepath{clip}
\pgfsetdash{}{0cm}
\definecolor{eps2pgf_color}{rgb}{0,0,1}\pgfsetstrokecolor{eps2pgf_color}\pgfsetfillcolor{eps2pgf_color}
\pgfpathmoveto{\pgfqpoint{4.854cm}{19.415cm}}
\pgfpathlineto{\pgfqpoint{5.897cm}{19.415cm}}
\pgfusepath{stroke}
\begin{pgfscope}
\pgfpathmoveto{\pgfqpoint{4.98cm}{19.808cm}}
\pgfpathlineto{\pgfqpoint{5.771cm}{19.808cm}}
\pgfpathlineto{\pgfqpoint{5.771cm}{19.018cm}}
\pgfpathlineto{\pgfqpoint{4.98cm}{19.018cm}}
\pgfpathclose
\pgfusepath{clip}
\pgfsetdash{}{0cm}
\pgfpathmoveto{\pgfqpoint{5.233cm}{19.415cm}}
\pgfpathlineto{\pgfqpoint{5.515cm}{19.415cm}}
\pgfusepath{stroke}
\pgfsetdash{}{0cm}
\pgfpathmoveto{\pgfqpoint{5.374cm}{19.556cm}}
\pgfpathlineto{\pgfqpoint{5.374cm}{19.273cm}}
\pgfusepath{stroke}
\pgfsetdash{}{0cm}
\pgfpathmoveto{\pgfqpoint{5.277cm}{19.512cm}}
\pgfpathlineto{\pgfqpoint{5.471cm}{19.318cm}}
\pgfusepath{stroke}
\pgfsetdash{}{0cm}
\pgfpathmoveto{\pgfqpoint{5.471cm}{19.512cm}}
\pgfpathlineto{\pgfqpoint{5.277cm}{19.318cm}}
\pgfusepath{stroke}
\end{pgfscope}
\end{pgfscope}
\pgftext[x=10.087cm,y=18.059cm,rotate=0]{ \fontsize{40}{36.14}\selectfont{ {$\la = 0.9$, \texttt{maxW}}}}
\begin{pgfscope}
\pgfpathmoveto{\pgfqpoint{4.648cm}{20.094cm}}
\pgfpathlineto{\pgfqpoint{14.875cm}{20.094cm}}
\pgfpathlineto{\pgfqpoint{14.875cm}{14.981cm}}
\pgfpathlineto{\pgfqpoint{4.648cm}{14.981cm}}
\pgfpathclose
\pgfusepath{clip}
\pgfsetdash{{0.018cm}{0.141cm}{0.212cm}{0.141cm}}{0cm}
\definecolor{eps2pgf_color}{rgb}{0,0,1}\pgfsetstrokecolor{eps2pgf_color}\pgfsetfillcolor{eps2pgf_color}
\pgfpathmoveto{\pgfqpoint{4.854cm}{18.165cm}}
\pgfpathlineto{\pgfqpoint{5.897cm}{18.165cm}}
\pgfusepath{stroke}
\begin{pgfscope}
\pgfpathmoveto{\pgfqpoint{4.98cm}{18.559cm}}
\pgfpathlineto{\pgfqpoint{5.771cm}{18.559cm}}
\pgfpathlineto{\pgfqpoint{5.771cm}{17.768cm}}
\pgfpathlineto{\pgfqpoint{4.98cm}{17.768cm}}
\pgfpathclose
\pgfusepath{clip}
\pgfsetdash{}{0cm}
\pgfpathmoveto{\pgfqpoint{5.233cm}{18.165cm}}
\pgfpathlineto{\pgfqpoint{5.515cm}{18.165cm}}
\pgfusepath{stroke}
\pgfsetdash{}{0cm}
\pgfpathmoveto{\pgfqpoint{5.374cm}{18.306cm}}
\pgfpathlineto{\pgfqpoint{5.374cm}{18.024cm}}
\pgfusepath{stroke}
\pgfsetdash{}{0cm}
\pgfpathmoveto{\pgfqpoint{5.277cm}{18.262cm}}
\pgfpathlineto{\pgfqpoint{5.471cm}{18.068cm}}
\pgfusepath{stroke}
\pgfsetdash{}{0cm}
\pgfpathmoveto{\pgfqpoint{5.471cm}{18.262cm}}
\pgfpathlineto{\pgfqpoint{5.277cm}{18.068cm}}
\pgfusepath{stroke}
\end{pgfscope}
\end{pgfscope}
\pgftext[x=10.205cm+.2cm,y=16.81cm,rotate=0]{ \fontsize{40}{36.14}\selectfont{ {$\la = 0.95$, \texttt{minS}}}}
\begin{pgfscope}
\pgfpathmoveto{\pgfqpoint{4.648cm}{20.094cm}}
\pgfpathlineto{\pgfqpoint{14.875cm}{20.094cm}}
\pgfpathlineto{\pgfqpoint{14.875cm}{14.981cm}}
\pgfpathlineto{\pgfqpoint{4.648cm}{14.981cm}}
\pgfpathclose
\pgfusepath{clip}
\pgfsetdash{}{0cm}
\definecolor{eps2pgf_color}{rgb}{0,1,0}\pgfsetstrokecolor{eps2pgf_color}\pgfsetfillcolor{eps2pgf_color}
\pgfpathmoveto{\pgfqpoint{4.854cm}{16.916cm}}
\pgfpathlineto{\pgfqpoint{5.897cm}{16.916cm}}
\pgfusepath{stroke}
\begin{pgfscope}
\pgfpathmoveto{\pgfqpoint{4.98cm}{17.31cm}}
\pgfpathlineto{\pgfqpoint{5.771cm}{17.31cm}}
\pgfpathlineto{\pgfqpoint{5.771cm}{16.519cm}}
\pgfpathlineto{\pgfqpoint{4.98cm}{16.519cm}}
\pgfpathclose
\pgfusepath{clip}
\pgfsetdash{}{0cm}
\pgfpathmoveto{\pgfqpoint{5.262cm}{17.027cm}}
\pgfpathlineto{\pgfqpoint{5.486cm}{17.027cm}}
\pgfpathlineto{\pgfqpoint{5.486cm}{16.804cm}}
\pgfpathlineto{\pgfqpoint{5.262cm}{16.804cm}}
\pgfpathlineto{\pgfqpoint{5.262cm}{17.027cm}}
\pgfpathclose
\pgfusepath{stroke}
\end{pgfscope}
\end{pgfscope}
\pgftext[x=10.381cm,y=15.56cm,rotate=0]{ \fontsize{40}{36.14}\selectfont{ {$\la = 0.95$, \texttt{maxW}}}}
\begin{pgfscope}
\pgfpathmoveto{\pgfqpoint{4.648cm}{20.094cm}}
\pgfpathlineto{\pgfqpoint{14.875cm}{20.094cm}}
\pgfpathlineto{\pgfqpoint{14.875cm}{14.981cm}}
\pgfpathlineto{\pgfqpoint{4.648cm}{14.981cm}}
\pgfpathclose
\pgfusepath{clip}
\pgfsetdash{{0.018cm}{0.141cm}{0.212cm}{0.141cm}}{0cm}
\definecolor{eps2pgf_color}{rgb}{0,1,0}\pgfsetstrokecolor{eps2pgf_color}\pgfsetfillcolor{eps2pgf_color}
\pgfpathmoveto{\pgfqpoint{4.854cm}{15.666cm}}
\pgfpathlineto{\pgfqpoint{5.897cm}{15.666cm}}
\pgfusepath{stroke}
\begin{pgfscope}
\pgfpathmoveto{\pgfqpoint{4.98cm}{16.06cm}}
\pgfpathlineto{\pgfqpoint{5.771cm}{16.06cm}}
\pgfpathlineto{\pgfqpoint{5.771cm}{15.269cm}}
\pgfpathlineto{\pgfqpoint{4.98cm}{15.269cm}}
\pgfpathclose
\pgfusepath{clip}
\pgfsetdash{}{0cm}
\pgfpathmoveto{\pgfqpoint{5.262cm}{15.778cm}}
\pgfpathlineto{\pgfqpoint{5.486cm}{15.778cm}}
\pgfpathlineto{\pgfqpoint{5.486cm}{15.555cm}}
\pgfpathlineto{\pgfqpoint{5.262cm}{15.555cm}}
\pgfpathlineto{\pgfqpoint{5.262cm}{15.778cm}}
\pgfpathclose
\pgfusepath{stroke}
\end{pgfscope}
\end{pgfscope}
\pgfsetdash{}{0cm}
\pgfsetlinewidth{0.176mm}
\definecolor{eps2pgf_color}{rgb}{0,1,0}\pgfsetstrokecolor{eps2pgf_color}\pgfsetfillcolor{eps2pgf_color}
\pgfusepath{stroke}
\end{pgfscope}
\end{pgfscope}
\end{pgfpicture}

%% file: fig/case2_1_1.tex
\scalebox{0.235}{\scalefont{2} \input{./fig/case2_1_1.pgf}}

%% file: fig/case2_1_1.pgf
% Created by Eps2pgf 0.7.0 (build on 2008-08-24) on Mon Apr 06 01:29:13 PDT 2015
\begin{pgfpicture}
\pgfpathmoveto{\pgfqpoint{0.635cm}{6.315cm}}
\pgfpathlineto{\pgfqpoint{19.473cm}{6.315cm}}
\pgfpathlineto{\pgfqpoint{19.473cm}{21.731cm}}
\pgfpathlineto{\pgfqpoint{0.635cm}{21.731cm}}
\pgfpathclose
\pgfusepath{clip}
\begin{pgfscope}
\begin{pgfscope}
\pgfpathmoveto{\pgfqpoint{0.635cm}{21.731cm}}
\pgfpathlineto{\pgfqpoint{0.635cm}{6.344cm}}
\pgfpathlineto{\pgfqpoint{19.491cm}{6.344cm}}
\pgfpathlineto{\pgfqpoint{19.491cm}{21.731cm}}
\pgfpathclose
\pgfusepath{clip}
\definecolor{eps2pgf_color}{gray}{1}\pgfsetstrokecolor{eps2pgf_color}\pgfsetfillcolor{eps2pgf_color}
\pgfpathmoveto{\pgfqpoint{0.635cm}{21.59cm}}
\pgfpathlineto{\pgfqpoint{0.635cm}{6.341cm}}
\pgfpathlineto{\pgfqpoint{20.99cm}{6.341cm}}
\pgfpathlineto{\pgfqpoint{20.99cm}{21.59cm}}
\pgfpathclose
\pgfusepath{fill}
\pgfpathmoveto{\pgfqpoint{3.601cm}{8.614cm}}
\pgfpathlineto{\pgfqpoint{3.601cm}{20.12cm}}
\pgfpathlineto{\pgfqpoint{19.053cm}{20.12cm}}
\pgfpathlineto{\pgfqpoint{19.053cm}{8.614cm}}
\pgfpathclose
\pgfseteorule\pgfusepath{fill}\pgfsetnonzerorule
\pgfsetdash{}{0cm}
\pgfsetlinewidth{0.176mm}
\pgfsetroundjoin
\pgfpathmoveto{\pgfqpoint{3.601cm}{8.614cm}}
\pgfpathlineto{\pgfqpoint{3.601cm}{20.12cm}}
\pgfpathlineto{\pgfqpoint{19.053cm}{20.12cm}}
\pgfpathlineto{\pgfqpoint{19.053cm}{8.614cm}}
\pgfpathlineto{\pgfqpoint{3.601cm}{8.614cm}}
\pgfusepath{stroke}
\pgfsetdash{}{0cm}
\definecolor{eps2pgf_color}{gray}{0}\pgfsetstrokecolor{eps2pgf_color}\pgfsetfillcolor{eps2pgf_color}
\pgfpathmoveto{\pgfqpoint{3.601cm}{8.614cm}}
\pgfpathlineto{\pgfqpoint{19.053cm}{8.614cm}}
\pgfusepath{stroke}
\pgfsetdash{}{0cm}
\pgfpathmoveto{\pgfqpoint{3.601cm}{20.12cm}}
\pgfpathlineto{\pgfqpoint{19.053cm}{20.12cm}}
\pgfusepath{stroke}
\pgfsetdash{}{0cm}
\pgfpathmoveto{\pgfqpoint{3.601cm}{8.614cm}}
\pgfpathlineto{\pgfqpoint{3.601cm}{20.12cm}}
\pgfusepath{stroke}
\pgfsetdash{}{0cm}
\pgfpathmoveto{\pgfqpoint{19.053cm}{8.614cm}}
\pgfpathlineto{\pgfqpoint{19.053cm}{20.12cm}}
\pgfusepath{stroke}
\pgfsetdash{}{0cm}
\pgfpathmoveto{\pgfqpoint{3.601cm}{8.614cm}}
\pgfpathlineto{\pgfqpoint{19.053cm}{8.614cm}}
\pgfusepath{stroke}
\pgfsetdash{}{0cm}
\pgfpathmoveto{\pgfqpoint{3.601cm}{8.614cm}}
\pgfpathlineto{\pgfqpoint{3.601cm}{20.12cm}}
\pgfusepath{stroke}
\pgfsetdash{}{0cm}
\pgfpathmoveto{\pgfqpoint{3.601cm}{8.614cm}}
\pgfpathlineto{\pgfqpoint{3.601cm}{8.769cm}}
\pgfusepath{stroke}
\pgfsetdash{}{0cm}
\pgfpathmoveto{\pgfqpoint{3.601cm}{20.12cm}}
\pgfpathlineto{\pgfqpoint{3.601cm}{19.967cm}}
\pgfusepath{stroke}
\pgftext[x=3.602cm,y=7.891cm+.2cm,rotate=0]{  \fontsize{36}{36.14}\selectfont{ {0}}}
\pgfsetdash{}{0cm}
\pgfpathmoveto{\pgfqpoint{11.327cm}{8.614cm}}
\pgfpathlineto{\pgfqpoint{11.327cm}{8.769cm}}
\pgfusepath{stroke}
\pgfsetdash{}{0cm}
\pgfpathmoveto{\pgfqpoint{11.327cm}{20.12cm}}
\pgfpathlineto{\pgfqpoint{11.327cm}{19.967cm}}
\pgfusepath{stroke}
\pgftext[x=11.325cm,y=7.891cm+.2cm,rotate=0]{  \fontsize{36}{36.14}\selectfont{ {0.5}}}
\pgfsetdash{}{0cm}
\pgfpathmoveto{\pgfqpoint{19.053cm}{8.614cm}}
\pgfpathlineto{\pgfqpoint{19.053cm}{8.769cm}}
\pgfusepath{stroke}
\pgfsetdash{}{0cm}
\pgfpathmoveto{\pgfqpoint{19.053cm}{20.12cm}}
\pgfpathlineto{\pgfqpoint{19.053cm}{19.967cm}}
\pgfusepath{stroke}
\pgftext[x=19.002cm-.5cm,y=7.901cm+.2cm,rotate=0]{  \fontsize{36}{36.14}\selectfont{ {1}}}
\pgfsetdash{}{0cm}
\pgfpathmoveto{\pgfqpoint{3.601cm}{8.614cm}}
\pgfpathlineto{\pgfqpoint{3.754cm}{8.614cm}}
\pgfusepath{stroke}
\pgfsetdash{}{0cm}
\pgfpathmoveto{\pgfqpoint{19.053cm}{8.614cm}}
\pgfpathlineto{\pgfqpoint{18.897cm}{8.614cm}}
\pgfusepath{stroke}
\pgftext[x=3.205cm-.3cm,y=8.582cm,rotate=0]{  \fontsize{36}{36.14}\selectfont{ {0}}}
\pgfsetdash{}{0cm}
\pgfpathmoveto{\pgfqpoint{3.601cm}{12.45cm}}
\pgfpathlineto{\pgfqpoint{3.754cm}{12.45cm}}
\pgfusepath{stroke}
\pgfsetdash{}{0cm}
\pgfpathmoveto{\pgfqpoint{19.053cm}{12.45cm}}
\pgfpathlineto{\pgfqpoint{18.897cm}{12.45cm}}
\pgfusepath{stroke}
\pgftext[x=2.761cm-.3cm,y=12.418cm,rotate=0]{  \fontsize{36}{36.14}\selectfont{ {0.5}}}
\pgfsetdash{}{0cm}
\pgfpathmoveto{\pgfqpoint{3.601cm}{16.287cm}}
\pgfpathlineto{\pgfqpoint{3.754cm}{16.287cm}}
\pgfusepath{stroke}
\pgfsetdash{}{0cm}
\pgfpathmoveto{\pgfqpoint{19.053cm}{16.287cm}}
\pgfpathlineto{\pgfqpoint{18.897cm}{16.287cm}}
\pgfusepath{stroke}
\pgftext[x=3.154cm-.3cm,y=16.265cm,rotate=0]{  \fontsize{36}{36.14}\selectfont{ {1}}}
\pgfsetdash{}{0cm}
\pgfpathmoveto{\pgfqpoint{3.601cm}{20.12cm}}
\pgfpathlineto{\pgfqpoint{3.754cm}{20.12cm}}
\pgfusepath{stroke}
\pgfsetdash{}{0cm}
\pgfpathmoveto{\pgfqpoint{19.053cm}{20.12cm}}
\pgfpathlineto{\pgfqpoint{18.897cm}{20.12cm}}
\pgfusepath{stroke}
\pgftext[x=2.795cm-.3cm,y=20.088cm,rotate=0]{  \fontsize{36}{36.14}\selectfont{ {1.5}}}
\pgftext[x=4.588cm,y=20.641cm,rotate=0]{  \fontsize{36}{36.14}\selectfont{ {x 10}}}
\pgftext[x=6.004cm,y=21.176cm,rotate=0]{\fontsize{20.08}{24.09}\selectfont{ {$-$3}}}
\pgfsetdash{}{0cm}
\pgfpathmoveto{\pgfqpoint{3.601cm}{8.614cm}}
\pgfpathlineto{\pgfqpoint{19.053cm}{8.614cm}}
\pgfusepath{stroke}
\pgfsetdash{}{0cm}
\pgfpathmoveto{\pgfqpoint{3.601cm}{20.12cm}}
\pgfpathlineto{\pgfqpoint{19.053cm}{20.12cm}}
\pgfusepath{stroke}
\pgfsetdash{}{0cm}
\pgfpathmoveto{\pgfqpoint{3.601cm}{8.614cm}}
\pgfpathlineto{\pgfqpoint{3.601cm}{20.12cm}}
\pgfusepath{stroke}
\pgfsetdash{}{0cm}
\pgfpathmoveto{\pgfqpoint{19.053cm}{8.614cm}}
\pgfpathlineto{\pgfqpoint{19.053cm}{20.12cm}}
\pgfusepath{stroke}
\begin{pgfscope}
\pgfpathmoveto{\pgfqpoint{3.601cm}{20.12cm}}
\pgfpathlineto{\pgfqpoint{3.601cm}{8.611cm}}
\pgfpathlineto{\pgfqpoint{19.056cm}{8.611cm}}
\pgfpathlineto{\pgfqpoint{19.056cm}{20.12cm}}
\pgfpathclose
\pgfusepath{clip}
\pgfsetdash{}{0cm}
\pgfsetlinewidth{1.058mm}
\definecolor{eps2pgf_color}{rgb}{0,0,1}\pgfsetstrokecolor{eps2pgf_color}\pgfsetfillcolor{eps2pgf_color}
\pgfpathmoveto{\pgfqpoint{5.145cm}{18.203cm}}
\pgfpathlineto{\pgfqpoint{6.691cm}{12.876cm}}
\pgfpathlineto{\pgfqpoint{8.234cm}{11.354cm}}
\pgfpathlineto{\pgfqpoint{9.781cm}{10.633cm}}
\pgfpathlineto{\pgfqpoint{11.327cm}{10.213cm}}
\pgfpathlineto{\pgfqpoint{12.871cm}{9.937cm}}
\pgfpathlineto{\pgfqpoint{14.417cm}{9.743cm}}
\pgfpathlineto{\pgfqpoint{15.96cm}{9.598cm}}
\pgfpathlineto{\pgfqpoint{17.507cm}{9.487cm}}
\pgfpathlineto{\pgfqpoint{19.053cm}{9.399cm}}
\pgfusepath{stroke}
\end{pgfscope}
\pgftext[x=11.345cm,y=6.968cm,rotate=0]{  \fontsize{40}{50}{\selectfont $\bmax$}}
\pgftext[x=1.257cm,y=14.366cm,rotate=90]{  \fontsize{36}{36.14}\selectfont{ {Sub-optimality bound}}}
\pgftext[x=3.551cm,y=8.487cm,rotate=0]{\fontsize{10.04}{12.04}\selectfont{ { }}}
\pgftext[x=19.006cm,y=19.997cm,rotate=0]{\fontsize{10.04}{12.04}\selectfont{ { }}}
\pgftext[x=16.633cm,y=19.254cm,rotate=0]{  \fontsize{40}{36.14}\selectfont{ {  $\la = 1$}}}
\begin{pgfscope}
\pgfpathmoveto{\pgfqpoint{13.335cm}{19.944cm}}
\pgfpathlineto{\pgfqpoint{13.335cm}{18.571cm}}
\pgfpathlineto{\pgfqpoint{18.882cm}{18.571cm}}
\pgfpathlineto{\pgfqpoint{18.882cm}{19.944cm}}
\pgfpathclose
\pgfusepath{clip}
\pgfsetdash{}{0cm}
\pgfsetlinewidth{1.058mm}
\definecolor{eps2pgf_color}{rgb}{0,0,1}\pgfsetstrokecolor{eps2pgf_color}\pgfsetfillcolor{eps2pgf_color}
\pgfpathmoveto{\pgfqpoint{13.538cm}{19.259cm}}
\pgfpathlineto{\pgfqpoint{14.567cm}{19.259cm}}
\pgfusepath{stroke}
\end{pgfscope}
\pgfsetdash{}{0cm}
\definecolor{eps2pgf_color}{rgb}{0,0,1}\pgfsetstrokecolor{eps2pgf_color}\pgfsetfillcolor{eps2pgf_color}
\pgfusepath{stroke}
\end{pgfscope}
\end{pgfscope}
\end{pgfpicture}

%% file: fig/case2_2_1.tex
\scalebox{0.235}{\scalefont{2} \input{./fig/case2_2_1.pgf}}

%% file: fig/case2_2_1.pgf
% Created by Eps2pgf 0.7.0 (build on 2008-08-24) on Mon Apr 06 01:29:30 PDT 2015
\begin{pgfpicture}
\pgfpathmoveto{\pgfqpoint{0.635cm}{6.315cm}}
\pgfpathlineto{\pgfqpoint{19.579cm}{6.315cm}}
\pgfpathlineto{\pgfqpoint{19.579cm}{21.167cm}}
\pgfpathlineto{\pgfqpoint{0.635cm}{21.167cm}}
\pgfpathclose
\pgfusepath{clip}
\begin{pgfscope}
\begin{pgfscope}
\pgfpathmoveto{\pgfqpoint{0.635cm}{21.178cm}}
\pgfpathlineto{\pgfqpoint{0.635cm}{6.344cm}}
\pgfpathlineto{\pgfqpoint{19.591cm}{6.344cm}}
\pgfpathlineto{\pgfqpoint{19.591cm}{21.178cm}}
\pgfpathclose
\pgfusepath{clip}
\definecolor{eps2pgf_color}{gray}{1}\pgfsetstrokecolor{eps2pgf_color}\pgfsetfillcolor{eps2pgf_color}
\pgfpathmoveto{\pgfqpoint{0.635cm}{21.59cm}}
\pgfpathlineto{\pgfqpoint{0.635cm}{6.341cm}}
\pgfpathlineto{\pgfqpoint{20.99cm}{6.341cm}}
\pgfpathlineto{\pgfqpoint{20.99cm}{21.59cm}}
\pgfpathclose
\pgfusepath{fill}
\pgfpathmoveto{\pgfqpoint{4.189cm}{8.614cm}}
\pgfpathlineto{\pgfqpoint{4.189cm}{20.446cm}}
\pgfpathlineto{\pgfqpoint{19.053cm}{20.446cm}}
\pgfpathlineto{\pgfqpoint{19.053cm}{8.614cm}}
\pgfpathclose
\pgfseteorule\pgfusepath{fill}\pgfsetnonzerorule
\pgfsetdash{}{0cm}
\pgfsetlinewidth{0.176mm}
\pgfsetroundjoin
\pgfpathmoveto{\pgfqpoint{4.189cm}{8.614cm}}
\pgfpathlineto{\pgfqpoint{4.189cm}{20.446cm}}
\pgfpathlineto{\pgfqpoint{19.053cm}{20.446cm}}
\pgfpathlineto{\pgfqpoint{19.053cm}{8.614cm}}
\pgfpathlineto{\pgfqpoint{4.189cm}{8.614cm}}
\pgfusepath{stroke}
\pgfsetdash{}{0cm}
\definecolor{eps2pgf_color}{gray}{0}\pgfsetstrokecolor{eps2pgf_color}\pgfsetfillcolor{eps2pgf_color}
\pgfpathmoveto{\pgfqpoint{4.189cm}{8.614cm}}
\pgfpathlineto{\pgfqpoint{19.053cm}{8.614cm}}
\pgfusepath{stroke}
\pgfsetdash{}{0cm}
\pgfpathmoveto{\pgfqpoint{4.189cm}{20.446cm}}
\pgfpathlineto{\pgfqpoint{19.053cm}{20.446cm}}
\pgfusepath{stroke}
\pgfsetdash{}{0cm}
\pgfpathmoveto{\pgfqpoint{4.189cm}{8.614cm}}
\pgfpathlineto{\pgfqpoint{4.189cm}{20.446cm}}
\pgfusepath{stroke}
\pgfsetdash{}{0cm}
\pgfpathmoveto{\pgfqpoint{19.053cm}{8.614cm}}
\pgfpathlineto{\pgfqpoint{19.053cm}{20.446cm}}
\pgfusepath{stroke}
\pgfsetdash{}{0cm}
\pgfpathmoveto{\pgfqpoint{4.189cm}{8.614cm}}
\pgfpathlineto{\pgfqpoint{19.053cm}{8.614cm}}
\pgfusepath{stroke}
\pgfsetdash{}{0cm}
\pgfpathmoveto{\pgfqpoint{4.189cm}{8.614cm}}
\pgfpathlineto{\pgfqpoint{4.189cm}{20.446cm}}
\pgfusepath{stroke}
\pgfsetdash{}{0cm}
\pgfpathmoveto{\pgfqpoint{4.189cm}{8.614cm}}
\pgfpathlineto{\pgfqpoint{4.189cm}{8.764cm}}
\pgfusepath{stroke}
\pgfsetdash{}{0cm}
\pgfpathmoveto{\pgfqpoint{4.189cm}{20.446cm}}
\pgfpathlineto{\pgfqpoint{4.189cm}{20.299cm}}
\pgfusepath{stroke}
\pgftext[x=4.189cm,y=7.891cm+.2cm,rotate=0]{  \fontsize{36}{36.14}\selectfont{ {0}}}
\pgfsetdash{}{0cm}
\pgfpathmoveto{\pgfqpoint{11.621cm}{8.614cm}}
\pgfpathlineto{\pgfqpoint{11.621cm}{8.764cm}}
\pgfusepath{stroke}
\pgfsetdash{}{0cm}
\pgfpathmoveto{\pgfqpoint{11.621cm}{20.446cm}}
\pgfpathlineto{\pgfqpoint{11.621cm}{20.299cm}}
\pgfusepath{stroke}
\pgftext[x=11.619cm,y=7.891cm+.2cm,rotate=0]{  \fontsize{36}{36.14}\selectfont{ {0.5}}}
\pgfsetdash{}{0cm}
\pgfpathmoveto{\pgfqpoint{19.053cm}{8.614cm}}
\pgfpathlineto{\pgfqpoint{19.053cm}{8.764cm}}
\pgfusepath{stroke}
\pgfsetdash{}{0cm}
\pgfpathmoveto{\pgfqpoint{19.053cm}{20.446cm}}
\pgfpathlineto{\pgfqpoint{19.053cm}{20.299cm}}
\pgfusepath{stroke}
\pgftext[x=19.002cm-.3cm,y=7.901cm+.2cm,rotate=0]{  \fontsize{36}{36.14}\selectfont{ {1}}}
\pgfsetdash{}{0cm}
\pgfpathmoveto{\pgfqpoint{4.189cm}{8.614cm}}
\pgfpathlineto{\pgfqpoint{4.336cm}{8.614cm}}
\pgfusepath{stroke}
\pgfsetdash{}{0cm}
\pgfpathmoveto{\pgfqpoint{19.053cm}{8.614cm}}
\pgfpathlineto{\pgfqpoint{18.903cm}{8.614cm}}
\pgfusepath{stroke}
\pgftext[x=3.793cm-.3cm,y=8.582cm,rotate=0]{  \fontsize{36}{36.14}\selectfont{ {0}}}
\pgfsetdash{}{0cm}
\pgfpathmoveto{\pgfqpoint{4.189cm}{11.574cm}}
\pgfpathlineto{\pgfqpoint{4.336cm}{11.574cm}}
\pgfusepath{stroke}
\pgfsetdash{}{0cm}
\pgfpathmoveto{\pgfqpoint{19.053cm}{11.574cm}}
\pgfpathlineto{\pgfqpoint{18.903cm}{11.574cm}}
\pgfusepath{stroke}
\pgftext[x=3.052cm-.3cm,y=11.542cm,rotate=0]{  \fontsize{36}{36.14}\selectfont{ {0.02}}}
\pgfsetdash{}{0cm}
\pgfpathmoveto{\pgfqpoint{4.189cm}{14.532cm}}
\pgfpathlineto{\pgfqpoint{4.336cm}{14.532cm}}
\pgfusepath{stroke}
\pgfsetdash{}{0cm}
\pgfpathmoveto{\pgfqpoint{19.053cm}{14.532cm}}
\pgfpathlineto{\pgfqpoint{18.903cm}{14.532cm}}
\pgfusepath{stroke}
\pgftext[x=3.06cm-.3cm,y=14.5cm,rotate=0]{  \fontsize{36}{36.14}\selectfont{ {0.04}}}
\pgfsetdash{}{0cm}
\pgfpathmoveto{\pgfqpoint{4.189cm}{17.489cm}}
\pgfpathlineto{\pgfqpoint{4.336cm}{17.489cm}}
\pgfusepath{stroke}
\pgfsetdash{}{0cm}
\pgfpathmoveto{\pgfqpoint{19.053cm}{17.489cm}}
\pgfpathlineto{\pgfqpoint{18.903cm}{17.489cm}}
\pgfusepath{stroke}
\pgftext[x=3.058cm-.3cm,y=17.457cm,rotate=0]{  \fontsize{36}{36.14}\selectfont{ {0.06}}}
\pgfsetdash{}{0cm}
\pgfpathmoveto{\pgfqpoint{4.189cm}{20.446cm}}
\pgfpathlineto{\pgfqpoint{4.336cm}{20.446cm}}
\pgfusepath{stroke}
\pgfsetdash{}{0cm}
\pgfpathmoveto{\pgfqpoint{19.053cm}{20.446cm}}
\pgfpathlineto{\pgfqpoint{18.903cm}{20.446cm}}
\pgfusepath{stroke}
\pgftext[x=3.057cm-.3cm,y=20.414cm,rotate=0]{  \fontsize{36}{36.14}\selectfont{ {0.08}}}
\pgfsetdash{}{0cm}
\pgfpathmoveto{\pgfqpoint{4.189cm}{8.614cm}}
\pgfpathlineto{\pgfqpoint{19.053cm}{8.614cm}}
\pgfusepath{stroke}
\pgfsetdash{}{0cm}
\pgfpathmoveto{\pgfqpoint{4.189cm}{20.446cm}}
\pgfpathlineto{\pgfqpoint{19.053cm}{20.446cm}}
\pgfusepath{stroke}
\pgfsetdash{}{0cm}
\pgfpathmoveto{\pgfqpoint{4.189cm}{8.614cm}}
\pgfpathlineto{\pgfqpoint{4.189cm}{20.446cm}}
\pgfusepath{stroke}
\pgfsetdash{}{0cm}
\pgfpathmoveto{\pgfqpoint{19.053cm}{8.614cm}}
\pgfpathlineto{\pgfqpoint{19.053cm}{20.446cm}}
\pgfusepath{stroke}
\begin{pgfscope}
\pgfpathmoveto{\pgfqpoint{4.189cm}{20.446cm}}
\pgfpathlineto{\pgfqpoint{4.189cm}{8.611cm}}
\pgfpathlineto{\pgfqpoint{19.056cm}{8.611cm}}
\pgfpathlineto{\pgfqpoint{19.056cm}{20.446cm}}
\pgfpathclose
\pgfusepath{clip}
\pgfsetdash{}{0cm}
\pgfsetlinewidth{1.058mm}
\definecolor{eps2pgf_color}{rgb}{0,0,1}\pgfsetstrokecolor{eps2pgf_color}\pgfsetfillcolor{eps2pgf_color}
\pgfpathmoveto{\pgfqpoint{5.674cm}{10.042cm}}
\pgfpathlineto{\pgfqpoint{7.161cm}{10.648cm}}
\pgfpathlineto{\pgfqpoint{8.646cm}{11.38cm}}
\pgfpathlineto{\pgfqpoint{10.134cm}{12.138cm}}
\pgfpathlineto{\pgfqpoint{11.621cm}{12.906cm}}
\pgfpathlineto{\pgfqpoint{13.106cm}{13.679cm}}
\pgfpathlineto{\pgfqpoint{14.593cm}{14.452cm}}
\pgfpathlineto{\pgfqpoint{16.078cm}{15.228cm}}
\pgfpathlineto{\pgfqpoint{17.565cm}{16.007cm}}
\pgfpathlineto{\pgfqpoint{19.053cm}{16.783cm}}
\pgfusepath{stroke}
\end{pgfscope}
\begin{pgfscope}
\pgfpathmoveto{\pgfqpoint{5.28cm}{17.177cm}}
\pgfpathlineto{\pgfqpoint{5.28cm}{9.646cm}}
\pgfpathlineto{\pgfqpoint{19.45cm}{9.646cm}}
\pgfpathlineto{\pgfqpoint{19.45cm}{17.177cm}}
\pgfpathclose
\pgfusepath{clip}
\pgfsetdash{}{0cm}
\pgfsetlinewidth{1.058mm}
\definecolor{eps2pgf_color}{rgb}{0,0,1}\pgfsetstrokecolor{eps2pgf_color}\pgfsetfillcolor{eps2pgf_color}
\pgfpathmoveto{\pgfqpoint{5.533cm}{10.042cm}}
\pgfpathlineto{\pgfqpoint{5.815cm}{10.042cm}}
\pgfusepath{stroke}
\pgfsetdash{}{0cm}
\pgfpathmoveto{\pgfqpoint{5.674cm}{10.184cm}}
\pgfpathlineto{\pgfqpoint{5.674cm}{9.901cm}}
\pgfusepath{stroke}
\pgfsetdash{}{0cm}
\pgfpathmoveto{\pgfqpoint{7.02cm}{10.648cm}}
\pgfpathlineto{\pgfqpoint{7.302cm}{10.648cm}}
\pgfusepath{stroke}
\pgfsetdash{}{0cm}
\pgfpathmoveto{\pgfqpoint{7.161cm}{10.789cm}}
\pgfpathlineto{\pgfqpoint{7.161cm}{10.507cm}}
\pgfusepath{stroke}
\pgfsetdash{}{0cm}
\pgfpathmoveto{\pgfqpoint{8.505cm}{11.38cm}}
\pgfpathlineto{\pgfqpoint{8.787cm}{11.38cm}}
\pgfusepath{stroke}
\pgfsetdash{}{0cm}
\pgfpathmoveto{\pgfqpoint{8.646cm}{11.521cm}}
\pgfpathlineto{\pgfqpoint{8.646cm}{11.239cm}}
\pgfusepath{stroke}
\pgfsetdash{}{0cm}
\pgfpathmoveto{\pgfqpoint{9.992cm}{12.138cm}}
\pgfpathlineto{\pgfqpoint{10.275cm}{12.138cm}}
\pgfusepath{stroke}
\pgfsetdash{}{0cm}
\pgfpathmoveto{\pgfqpoint{10.134cm}{12.28cm}}
\pgfpathlineto{\pgfqpoint{10.134cm}{11.997cm}}
\pgfusepath{stroke}
\pgfsetdash{}{0cm}
\pgfpathmoveto{\pgfqpoint{11.48cm}{12.906cm}}
\pgfpathlineto{\pgfqpoint{11.762cm}{12.906cm}}
\pgfusepath{stroke}
\pgfsetdash{}{0cm}
\pgfpathmoveto{\pgfqpoint{11.621cm}{13.047cm}}
\pgfpathlineto{\pgfqpoint{11.621cm}{12.765cm}}
\pgfusepath{stroke}
\pgfsetdash{}{0cm}
\pgfpathmoveto{\pgfqpoint{12.965cm}{13.679cm}}
\pgfpathlineto{\pgfqpoint{13.247cm}{13.679cm}}
\pgfusepath{stroke}
\pgfsetdash{}{0cm}
\pgfpathmoveto{\pgfqpoint{13.106cm}{13.82cm}}
\pgfpathlineto{\pgfqpoint{13.106cm}{13.538cm}}
\pgfusepath{stroke}
\pgfsetdash{}{0cm}
\pgfpathmoveto{\pgfqpoint{14.452cm}{14.452cm}}
\pgfpathlineto{\pgfqpoint{14.734cm}{14.452cm}}
\pgfusepath{stroke}
\pgfsetdash{}{0cm}
\pgfpathmoveto{\pgfqpoint{14.593cm}{14.593cm}}
\pgfpathlineto{\pgfqpoint{14.593cm}{14.311cm}}
\pgfusepath{stroke}
\pgfsetdash{}{0cm}
\pgfpathmoveto{\pgfqpoint{15.937cm}{15.228cm}}
\pgfpathlineto{\pgfqpoint{16.219cm}{15.228cm}}
\pgfusepath{stroke}
\pgfsetdash{}{0cm}
\pgfpathmoveto{\pgfqpoint{16.078cm}{15.369cm}}
\pgfpathlineto{\pgfqpoint{16.078cm}{15.087cm}}
\pgfusepath{stroke}
\pgfsetdash{}{0cm}
\pgfpathmoveto{\pgfqpoint{17.424cm}{16.007cm}}
\pgfpathlineto{\pgfqpoint{17.707cm}{16.007cm}}
\pgfusepath{stroke}
\pgfsetdash{}{0cm}
\pgfpathmoveto{\pgfqpoint{17.565cm}{16.148cm}}
\pgfpathlineto{\pgfqpoint{17.565cm}{15.866cm}}
\pgfusepath{stroke}
\pgfsetdash{}{0cm}
\pgfpathmoveto{\pgfqpoint{18.912cm}{16.783cm}}
\pgfpathlineto{\pgfqpoint{19.194cm}{16.783cm}}
\pgfusepath{stroke}
\pgfsetdash{}{0cm}
\pgfpathmoveto{\pgfqpoint{19.053cm}{16.925cm}}
\pgfpathlineto{\pgfqpoint{19.053cm}{16.642cm}}
\pgfusepath{stroke}
\pgfsetdash{}{0cm}
\pgfpathmoveto{\pgfqpoint{5.577cm}{10.139cm}}
\pgfpathlineto{\pgfqpoint{5.771cm}{9.945cm}}
\pgfusepath{stroke}
\pgfsetdash{}{0cm}
\pgfpathmoveto{\pgfqpoint{5.771cm}{10.139cm}}
\pgfpathlineto{\pgfqpoint{5.577cm}{9.945cm}}
\pgfusepath{stroke}
\pgfsetdash{}{0cm}
\pgfpathmoveto{\pgfqpoint{7.064cm}{10.745cm}}
\pgfpathlineto{\pgfqpoint{7.258cm}{10.551cm}}
\pgfusepath{stroke}
\pgfsetdash{}{0cm}
\pgfpathmoveto{\pgfqpoint{7.258cm}{10.745cm}}
\pgfpathlineto{\pgfqpoint{7.064cm}{10.551cm}}
\pgfusepath{stroke}
\pgfsetdash{}{0cm}
\pgfpathmoveto{\pgfqpoint{8.549cm}{11.477cm}}
\pgfpathlineto{\pgfqpoint{8.743cm}{11.283cm}}
\pgfusepath{stroke}
\pgfsetdash{}{0cm}
\pgfpathmoveto{\pgfqpoint{8.743cm}{11.477cm}}
\pgfpathlineto{\pgfqpoint{8.549cm}{11.283cm}}
\pgfusepath{stroke}
\pgfsetdash{}{0cm}
\pgfpathmoveto{\pgfqpoint{10.037cm}{12.236cm}}
\pgfpathlineto{\pgfqpoint{10.231cm}{12.041cm}}
\pgfusepath{stroke}
\pgfsetdash{}{0cm}
\pgfpathmoveto{\pgfqpoint{10.231cm}{12.236cm}}
\pgfpathlineto{\pgfqpoint{10.037cm}{12.041cm}}
\pgfusepath{stroke}
\pgfsetdash{}{0cm}
\pgfpathmoveto{\pgfqpoint{11.524cm}{13.003cm}}
\pgfpathlineto{\pgfqpoint{11.718cm}{12.809cm}}
\pgfusepath{stroke}
\pgfsetdash{}{0cm}
\pgfpathmoveto{\pgfqpoint{11.718cm}{13.003cm}}
\pgfpathlineto{\pgfqpoint{11.524cm}{12.809cm}}
\pgfusepath{stroke}
\pgfsetdash{}{0cm}
\pgfpathmoveto{\pgfqpoint{13.009cm}{13.776cm}}
\pgfpathlineto{\pgfqpoint{13.203cm}{13.582cm}}
\pgfusepath{stroke}
\pgfsetdash{}{0cm}
\pgfpathmoveto{\pgfqpoint{13.203cm}{13.776cm}}
\pgfpathlineto{\pgfqpoint{13.009cm}{13.582cm}}
\pgfusepath{stroke}
\pgfsetdash{}{0cm}
\pgfpathmoveto{\pgfqpoint{14.496cm}{14.549cm}}
\pgfpathlineto{\pgfqpoint{14.69cm}{14.355cm}}
\pgfusepath{stroke}
\pgfsetdash{}{0cm}
\pgfpathmoveto{\pgfqpoint{14.69cm}{14.549cm}}
\pgfpathlineto{\pgfqpoint{14.496cm}{14.355cm}}
\pgfusepath{stroke}
\pgfsetdash{}{0cm}
\pgfpathmoveto{\pgfqpoint{15.981cm}{15.325cm}}
\pgfpathlineto{\pgfqpoint{16.175cm}{15.131cm}}
\pgfusepath{stroke}
\pgfsetdash{}{0cm}
\pgfpathmoveto{\pgfqpoint{16.175cm}{15.325cm}}
\pgfpathlineto{\pgfqpoint{15.981cm}{15.131cm}}
\pgfusepath{stroke}
\pgfsetdash{}{0cm}
\pgfpathmoveto{\pgfqpoint{17.468cm}{16.104cm}}
\pgfpathlineto{\pgfqpoint{17.662cm}{15.91cm}}
\pgfusepath{stroke}
\pgfsetdash{}{0cm}
\pgfpathmoveto{\pgfqpoint{17.662cm}{16.104cm}}
\pgfpathlineto{\pgfqpoint{17.468cm}{15.91cm}}
\pgfusepath{stroke}
\pgfsetdash{}{0cm}
\pgfpathmoveto{\pgfqpoint{18.956cm}{16.88cm}}
\pgfpathlineto{\pgfqpoint{19.15cm}{16.686cm}}
\pgfusepath{stroke}
\pgfsetdash{}{0cm}
\pgfpathmoveto{\pgfqpoint{19.15cm}{16.88cm}}
\pgfpathlineto{\pgfqpoint{18.956cm}{16.686cm}}
\pgfusepath{stroke}
\end{pgfscope}
\begin{pgfscope}
\pgfpathmoveto{\pgfqpoint{4.189cm}{20.446cm}}
\pgfpathlineto{\pgfqpoint{4.189cm}{8.611cm}}
\pgfpathlineto{\pgfqpoint{19.056cm}{8.611cm}}
\pgfpathlineto{\pgfqpoint{19.056cm}{20.446cm}}
\pgfpathclose
\pgfusepath{clip}
\pgfsetdash{{0.018cm}{0.141cm}{0.212cm}{0.141cm}}{0cm}
\pgfsetlinewidth{1.058mm}
\definecolor{eps2pgf_color}{rgb}{0,0,1}\pgfsetstrokecolor{eps2pgf_color}\pgfsetfillcolor{eps2pgf_color}
\pgfpathmoveto{\pgfqpoint{5.674cm}{10.089cm}}
\pgfpathlineto{\pgfqpoint{7.161cm}{10.807cm}}
\pgfpathlineto{\pgfqpoint{8.646cm}{11.633cm}}
\pgfpathlineto{\pgfqpoint{10.134cm}{12.48cm}}
\pgfpathlineto{\pgfqpoint{11.621cm}{13.335cm}}
\pgfpathlineto{\pgfqpoint{13.106cm}{14.193cm}}
\pgfpathlineto{\pgfqpoint{14.593cm}{15.055cm}}
\pgfpathlineto{\pgfqpoint{16.078cm}{15.919cm}}
\pgfpathlineto{\pgfqpoint{17.565cm}{16.783cm}}
\pgfpathlineto{\pgfqpoint{19.053cm}{17.648cm}}
\pgfusepath{stroke}
\end{pgfscope}
\begin{pgfscope}
\pgfpathmoveto{\pgfqpoint{5.28cm}{18.042cm}}
\pgfpathlineto{\pgfqpoint{5.28cm}{9.693cm}}
\pgfpathlineto{\pgfqpoint{19.45cm}{9.693cm}}
\pgfpathlineto{\pgfqpoint{19.45cm}{18.042cm}}
\pgfpathclose
\pgfusepath{clip}
\pgfsetdash{}{0cm}
\pgfsetlinewidth{1.058mm}
\definecolor{eps2pgf_color}{rgb}{0,0,1}\pgfsetstrokecolor{eps2pgf_color}\pgfsetfillcolor{eps2pgf_color}
\pgfpathmoveto{\pgfqpoint{5.533cm}{10.089cm}}
\pgfpathlineto{\pgfqpoint{5.815cm}{10.089cm}}
\pgfusepath{stroke}
\pgfsetdash{}{0cm}
\pgfpathmoveto{\pgfqpoint{5.674cm}{10.231cm}}
\pgfpathlineto{\pgfqpoint{5.674cm}{9.948cm}}
\pgfusepath{stroke}
\pgfsetdash{}{0cm}
\pgfpathmoveto{\pgfqpoint{7.02cm}{10.807cm}}
\pgfpathlineto{\pgfqpoint{7.302cm}{10.807cm}}
\pgfusepath{stroke}
\pgfsetdash{}{0cm}
\pgfpathmoveto{\pgfqpoint{7.161cm}{10.948cm}}
\pgfpathlineto{\pgfqpoint{7.161cm}{10.666cm}}
\pgfusepath{stroke}
\pgfsetdash{}{0cm}
\pgfpathmoveto{\pgfqpoint{8.505cm}{11.633cm}}
\pgfpathlineto{\pgfqpoint{8.787cm}{11.633cm}}
\pgfusepath{stroke}
\pgfsetdash{}{0cm}
\pgfpathmoveto{\pgfqpoint{8.646cm}{11.774cm}}
\pgfpathlineto{\pgfqpoint{8.646cm}{11.492cm}}
\pgfusepath{stroke}
\pgfsetdash{}{0cm}
\pgfpathmoveto{\pgfqpoint{9.992cm}{12.48cm}}
\pgfpathlineto{\pgfqpoint{10.275cm}{12.48cm}}
\pgfusepath{stroke}
\pgfsetdash{}{0cm}
\pgfpathmoveto{\pgfqpoint{10.134cm}{12.621cm}}
\pgfpathlineto{\pgfqpoint{10.134cm}{12.338cm}}
\pgfusepath{stroke}
\pgfsetdash{}{0cm}
\pgfpathmoveto{\pgfqpoint{11.48cm}{13.335cm}}
\pgfpathlineto{\pgfqpoint{11.762cm}{13.335cm}}
\pgfusepath{stroke}
\pgfsetdash{}{0cm}
\pgfpathmoveto{\pgfqpoint{11.621cm}{13.476cm}}
\pgfpathlineto{\pgfqpoint{11.621cm}{13.194cm}}
\pgfusepath{stroke}
\pgfsetdash{}{0cm}
\pgfpathmoveto{\pgfqpoint{12.965cm}{14.193cm}}
\pgfpathlineto{\pgfqpoint{13.247cm}{14.193cm}}
\pgfusepath{stroke}
\pgfsetdash{}{0cm}
\pgfpathmoveto{\pgfqpoint{13.106cm}{14.335cm}}
\pgfpathlineto{\pgfqpoint{13.106cm}{14.052cm}}
\pgfusepath{stroke}
\pgfsetdash{}{0cm}
\pgfpathmoveto{\pgfqpoint{14.452cm}{15.055cm}}
\pgfpathlineto{\pgfqpoint{14.734cm}{15.055cm}}
\pgfusepath{stroke}
\pgfsetdash{}{0cm}
\pgfpathmoveto{\pgfqpoint{14.593cm}{15.196cm}}
\pgfpathlineto{\pgfqpoint{14.593cm}{14.914cm}}
\pgfusepath{stroke}
\pgfsetdash{}{0cm}
\pgfpathmoveto{\pgfqpoint{15.937cm}{15.919cm}}
\pgfpathlineto{\pgfqpoint{16.219cm}{15.919cm}}
\pgfusepath{stroke}
\pgfsetdash{}{0cm}
\pgfpathmoveto{\pgfqpoint{16.078cm}{16.06cm}}
\pgfpathlineto{\pgfqpoint{16.078cm}{15.778cm}}
\pgfusepath{stroke}
\pgfsetdash{}{0cm}
\pgfpathmoveto{\pgfqpoint{17.424cm}{16.783cm}}
\pgfpathlineto{\pgfqpoint{17.707cm}{16.783cm}}
\pgfusepath{stroke}
\pgfsetdash{}{0cm}
\pgfpathmoveto{\pgfqpoint{17.565cm}{16.925cm}}
\pgfpathlineto{\pgfqpoint{17.565cm}{16.642cm}}
\pgfusepath{stroke}
\pgfsetdash{}{0cm}
\pgfpathmoveto{\pgfqpoint{18.912cm}{17.648cm}}
\pgfpathlineto{\pgfqpoint{19.194cm}{17.648cm}}
\pgfusepath{stroke}
\pgfsetdash{}{0cm}
\pgfpathmoveto{\pgfqpoint{19.053cm}{17.789cm}}
\pgfpathlineto{\pgfqpoint{19.053cm}{17.507cm}}
\pgfusepath{stroke}
\pgfsetdash{}{0cm}
\pgfpathmoveto{\pgfqpoint{5.577cm}{10.186cm}}
\pgfpathlineto{\pgfqpoint{5.771cm}{9.992cm}}
\pgfusepath{stroke}
\pgfsetdash{}{0cm}
\pgfpathmoveto{\pgfqpoint{5.771cm}{10.186cm}}
\pgfpathlineto{\pgfqpoint{5.577cm}{9.992cm}}
\pgfusepath{stroke}
\pgfsetdash{}{0cm}
\pgfpathmoveto{\pgfqpoint{7.064cm}{10.904cm}}
\pgfpathlineto{\pgfqpoint{7.258cm}{10.71cm}}
\pgfusepath{stroke}
\pgfsetdash{}{0cm}
\pgfpathmoveto{\pgfqpoint{7.258cm}{10.904cm}}
\pgfpathlineto{\pgfqpoint{7.064cm}{10.71cm}}
\pgfusepath{stroke}
\pgfsetdash{}{0cm}
\pgfpathmoveto{\pgfqpoint{8.549cm}{11.73cm}}
\pgfpathlineto{\pgfqpoint{8.743cm}{11.536cm}}
\pgfusepath{stroke}
\pgfsetdash{}{0cm}
\pgfpathmoveto{\pgfqpoint{8.743cm}{11.73cm}}
\pgfpathlineto{\pgfqpoint{8.549cm}{11.536cm}}
\pgfusepath{stroke}
\pgfsetdash{}{0cm}
\pgfpathmoveto{\pgfqpoint{10.037cm}{12.577cm}}
\pgfpathlineto{\pgfqpoint{10.231cm}{12.382cm}}
\pgfusepath{stroke}
\pgfsetdash{}{0cm}
\pgfpathmoveto{\pgfqpoint{10.231cm}{12.577cm}}
\pgfpathlineto{\pgfqpoint{10.037cm}{12.382cm}}
\pgfusepath{stroke}
\pgfsetdash{}{0cm}
\pgfpathmoveto{\pgfqpoint{11.524cm}{13.432cm}}
\pgfpathlineto{\pgfqpoint{11.718cm}{13.238cm}}
\pgfusepath{stroke}
\pgfsetdash{}{0cm}
\pgfpathmoveto{\pgfqpoint{11.718cm}{13.432cm}}
\pgfpathlineto{\pgfqpoint{11.524cm}{13.238cm}}
\pgfusepath{stroke}
\pgfsetdash{}{0cm}
\pgfpathmoveto{\pgfqpoint{13.009cm}{14.29cm}}
\pgfpathlineto{\pgfqpoint{13.203cm}{14.096cm}}
\pgfusepath{stroke}
\pgfsetdash{}{0cm}
\pgfpathmoveto{\pgfqpoint{13.203cm}{14.29cm}}
\pgfpathlineto{\pgfqpoint{13.009cm}{14.096cm}}
\pgfusepath{stroke}
\pgfsetdash{}{0cm}
\pgfpathmoveto{\pgfqpoint{14.496cm}{15.152cm}}
\pgfpathlineto{\pgfqpoint{14.69cm}{14.958cm}}
\pgfusepath{stroke}
\pgfsetdash{}{0cm}
\pgfpathmoveto{\pgfqpoint{14.69cm}{15.152cm}}
\pgfpathlineto{\pgfqpoint{14.496cm}{14.958cm}}
\pgfusepath{stroke}
\pgfsetdash{}{0cm}
\pgfpathmoveto{\pgfqpoint{15.981cm}{16.016cm}}
\pgfpathlineto{\pgfqpoint{16.175cm}{15.822cm}}
\pgfusepath{stroke}
\pgfsetdash{}{0cm}
\pgfpathmoveto{\pgfqpoint{16.175cm}{16.016cm}}
\pgfpathlineto{\pgfqpoint{15.981cm}{15.822cm}}
\pgfusepath{stroke}
\pgfsetdash{}{0cm}
\pgfpathmoveto{\pgfqpoint{17.468cm}{16.88cm}}
\pgfpathlineto{\pgfqpoint{17.662cm}{16.686cm}}
\pgfusepath{stroke}
\pgfsetdash{}{0cm}
\pgfpathmoveto{\pgfqpoint{17.662cm}{16.88cm}}
\pgfpathlineto{\pgfqpoint{17.468cm}{16.686cm}}
\pgfusepath{stroke}
\pgfsetdash{}{0cm}
\pgfpathmoveto{\pgfqpoint{18.956cm}{17.745cm}}
\pgfpathlineto{\pgfqpoint{19.15cm}{17.551cm}}
\pgfusepath{stroke}
\pgfsetdash{}{0cm}
\pgfpathmoveto{\pgfqpoint{19.15cm}{17.745cm}}
\pgfpathlineto{\pgfqpoint{18.956cm}{17.551cm}}
\pgfusepath{stroke}
\end{pgfscope}
\begin{pgfscope}
\pgfpathmoveto{\pgfqpoint{4.189cm}{20.446cm}}
\pgfpathlineto{\pgfqpoint{4.189cm}{8.611cm}}
\pgfpathlineto{\pgfqpoint{19.056cm}{8.611cm}}
\pgfpathlineto{\pgfqpoint{19.056cm}{20.446cm}}
\pgfpathclose
\pgfusepath{clip}
\pgfsetdash{}{0cm}
\pgfsetlinewidth{1.058mm}
\definecolor{eps2pgf_color}{rgb}{0,1,0}\pgfsetstrokecolor{eps2pgf_color}\pgfsetfillcolor{eps2pgf_color}
\pgfpathmoveto{\pgfqpoint{5.674cm}{9.393cm}}
\pgfpathlineto{\pgfqpoint{7.161cm}{9.637cm}}
\pgfpathlineto{\pgfqpoint{8.646cm}{9.978cm}}
\pgfpathlineto{\pgfqpoint{10.134cm}{10.339cm}}
\pgfpathlineto{\pgfqpoint{11.621cm}{10.71cm}}
\pgfpathlineto{\pgfqpoint{13.106cm}{11.083cm}}
\pgfpathlineto{\pgfqpoint{14.593cm}{11.456cm}}
\pgfpathlineto{\pgfqpoint{16.078cm}{11.833cm}}
\pgfpathlineto{\pgfqpoint{17.565cm}{12.209cm}}
\pgfpathlineto{\pgfqpoint{19.053cm}{12.588cm}}
\pgfusepath{stroke}
\end{pgfscope}
\begin{pgfscope}
\pgfpathmoveto{\pgfqpoint{5.28cm}{12.982cm}}
\pgfpathlineto{\pgfqpoint{5.28cm}{8.996cm}}
\pgfpathlineto{\pgfqpoint{19.45cm}{8.996cm}}
\pgfpathlineto{\pgfqpoint{19.45cm}{12.982cm}}
\pgfpathclose
\pgfusepath{clip}
\pgfsetdash{}{0cm}
\pgfsetlinewidth{1.058mm}
\pgfsetmiterjoin
\definecolor{eps2pgf_color}{rgb}{0,1,0}\pgfsetstrokecolor{eps2pgf_color}\pgfsetfillcolor{eps2pgf_color}
\pgfpathmoveto{\pgfqpoint{5.562cm}{9.504cm}}
\pgfpathlineto{\pgfqpoint{5.786cm}{9.504cm}}
\pgfpathlineto{\pgfqpoint{5.786cm}{9.281cm}}
\pgfpathlineto{\pgfqpoint{5.562cm}{9.281cm}}
\pgfpathlineto{\pgfqpoint{5.562cm}{9.504cm}}
\pgfpathclose
\pgfusepath{stroke}
\pgfsetdash{}{0cm}
\pgfpathmoveto{\pgfqpoint{7.05cm}{9.748cm}}
\pgfpathlineto{\pgfqpoint{7.273cm}{9.748cm}}
\pgfpathlineto{\pgfqpoint{7.273cm}{9.525cm}}
\pgfpathlineto{\pgfqpoint{7.05cm}{9.525cm}}
\pgfpathlineto{\pgfqpoint{7.05cm}{9.748cm}}
\pgfpathclose
\pgfusepath{stroke}
\pgfsetdash{}{0cm}
\pgfpathmoveto{\pgfqpoint{8.534cm}{10.089cm}}
\pgfpathlineto{\pgfqpoint{8.758cm}{10.089cm}}
\pgfpathlineto{\pgfqpoint{8.758cm}{9.866cm}}
\pgfpathlineto{\pgfqpoint{8.534cm}{9.866cm}}
\pgfpathlineto{\pgfqpoint{8.534cm}{10.089cm}}
\pgfpathclose
\pgfusepath{stroke}
\pgfsetdash{}{0cm}
\pgfpathmoveto{\pgfqpoint{10.022cm}{10.451cm}}
\pgfpathlineto{\pgfqpoint{10.245cm}{10.451cm}}
\pgfpathlineto{\pgfqpoint{10.245cm}{10.228cm}}
\pgfpathlineto{\pgfqpoint{10.022cm}{10.228cm}}
\pgfpathlineto{\pgfqpoint{10.022cm}{10.451cm}}
\pgfpathclose
\pgfusepath{stroke}
\pgfsetdash{}{0cm}
\pgfpathmoveto{\pgfqpoint{11.509cm}{10.821cm}}
\pgfpathlineto{\pgfqpoint{11.733cm}{10.821cm}}
\pgfpathlineto{\pgfqpoint{11.733cm}{10.598cm}}
\pgfpathlineto{\pgfqpoint{11.509cm}{10.598cm}}
\pgfpathlineto{\pgfqpoint{11.509cm}{10.821cm}}
\pgfpathclose
\pgfusepath{stroke}
\pgfsetdash{}{0cm}
\pgfpathmoveto{\pgfqpoint{12.994cm}{11.195cm}}
\pgfpathlineto{\pgfqpoint{13.217cm}{11.195cm}}
\pgfpathlineto{\pgfqpoint{13.217cm}{10.971cm}}
\pgfpathlineto{\pgfqpoint{12.994cm}{10.971cm}}
\pgfpathlineto{\pgfqpoint{12.994cm}{11.195cm}}
\pgfpathclose
\pgfusepath{stroke}
\pgfsetdash{}{0cm}
\pgfpathmoveto{\pgfqpoint{14.482cm}{11.568cm}}
\pgfpathlineto{\pgfqpoint{14.705cm}{11.568cm}}
\pgfpathlineto{\pgfqpoint{14.705cm}{11.345cm}}
\pgfpathlineto{\pgfqpoint{14.482cm}{11.345cm}}
\pgfpathlineto{\pgfqpoint{14.482cm}{11.568cm}}
\pgfpathclose
\pgfusepath{stroke}
\pgfsetdash{}{0cm}
\pgfpathmoveto{\pgfqpoint{15.966cm}{11.944cm}}
\pgfpathlineto{\pgfqpoint{16.19cm}{11.944cm}}
\pgfpathlineto{\pgfqpoint{16.19cm}{11.721cm}}
\pgfpathlineto{\pgfqpoint{15.966cm}{11.721cm}}
\pgfpathlineto{\pgfqpoint{15.966cm}{11.944cm}}
\pgfpathclose
\pgfusepath{stroke}
\pgfsetdash{}{0cm}
\pgfpathmoveto{\pgfqpoint{17.454cm}{12.321cm}}
\pgfpathlineto{\pgfqpoint{17.677cm}{12.321cm}}
\pgfpathlineto{\pgfqpoint{17.677cm}{12.097cm}}
\pgfpathlineto{\pgfqpoint{17.454cm}{12.097cm}}
\pgfpathlineto{\pgfqpoint{17.454cm}{12.321cm}}
\pgfpathclose
\pgfusepath{stroke}
\pgfsetdash{}{0cm}
\pgfpathmoveto{\pgfqpoint{18.941cm}{12.7cm}}
\pgfpathlineto{\pgfqpoint{19.165cm}{12.7cm}}
\pgfpathlineto{\pgfqpoint{19.165cm}{12.477cm}}
\pgfpathlineto{\pgfqpoint{18.941cm}{12.477cm}}
\pgfpathlineto{\pgfqpoint{18.941cm}{12.7cm}}
\pgfpathclose
\pgfusepath{stroke}
\end{pgfscope}
\begin{pgfscope}
\pgfpathmoveto{\pgfqpoint{4.189cm}{20.446cm}}
\pgfpathlineto{\pgfqpoint{4.189cm}{8.611cm}}
\pgfpathlineto{\pgfqpoint{19.056cm}{8.611cm}}
\pgfpathlineto{\pgfqpoint{19.056cm}{20.446cm}}
\pgfpathclose
\pgfusepath{clip}
\pgfsetdash{{0.018cm}{0.141cm}{0.212cm}{0.141cm}}{0cm}
\pgfsetlinewidth{1.058mm}
\definecolor{eps2pgf_color}{rgb}{0,1,0}\pgfsetstrokecolor{eps2pgf_color}\pgfsetfillcolor{eps2pgf_color}
\pgfpathmoveto{\pgfqpoint{5.674cm}{9.396cm}}
\pgfpathlineto{\pgfqpoint{7.161cm}{9.672cm}}
\pgfpathlineto{\pgfqpoint{8.646cm}{10.037cm}}
\pgfpathlineto{\pgfqpoint{10.134cm}{10.419cm}}
\pgfpathlineto{\pgfqpoint{11.621cm}{10.807cm}}
\pgfpathlineto{\pgfqpoint{13.106cm}{11.201cm}}
\pgfpathlineto{\pgfqpoint{14.593cm}{11.595cm}}
\pgfpathlineto{\pgfqpoint{16.078cm}{11.992cm}}
\pgfpathlineto{\pgfqpoint{17.565cm}{12.388cm}}
\pgfpathlineto{\pgfqpoint{19.053cm}{12.785cm}}
\pgfusepath{stroke}
\end{pgfscope}
\begin{pgfscope}
\pgfpathmoveto{\pgfqpoint{5.28cm}{13.179cm}}
\pgfpathlineto{\pgfqpoint{5.28cm}{8.999cm}}
\pgfpathlineto{\pgfqpoint{19.45cm}{8.999cm}}
\pgfpathlineto{\pgfqpoint{19.45cm}{13.179cm}}
\pgfpathclose
\pgfusepath{clip}
\pgfsetdash{}{0cm}
\pgfsetlinewidth{1.058mm}
\pgfsetmiterjoin
\definecolor{eps2pgf_color}{rgb}{0,1,0}\pgfsetstrokecolor{eps2pgf_color}\pgfsetfillcolor{eps2pgf_color}
\pgfpathmoveto{\pgfqpoint{5.562cm}{9.507cm}}
\pgfpathlineto{\pgfqpoint{5.786cm}{9.507cm}}
\pgfpathlineto{\pgfqpoint{5.786cm}{9.284cm}}
\pgfpathlineto{\pgfqpoint{5.562cm}{9.284cm}}
\pgfpathlineto{\pgfqpoint{5.562cm}{9.507cm}}
\pgfpathclose
\pgfusepath{stroke}
\pgfsetdash{}{0cm}
\pgfpathmoveto{\pgfqpoint{7.05cm}{9.784cm}}
\pgfpathlineto{\pgfqpoint{7.273cm}{9.784cm}}
\pgfpathlineto{\pgfqpoint{7.273cm}{9.56cm}}
\pgfpathlineto{\pgfqpoint{7.05cm}{9.56cm}}
\pgfpathlineto{\pgfqpoint{7.05cm}{9.784cm}}
\pgfpathclose
\pgfusepath{stroke}
\pgfsetdash{}{0cm}
\pgfpathmoveto{\pgfqpoint{8.534cm}{10.148cm}}
\pgfpathlineto{\pgfqpoint{8.758cm}{10.148cm}}
\pgfpathlineto{\pgfqpoint{8.758cm}{9.925cm}}
\pgfpathlineto{\pgfqpoint{8.534cm}{9.925cm}}
\pgfpathlineto{\pgfqpoint{8.534cm}{10.148cm}}
\pgfpathclose
\pgfusepath{stroke}
\pgfsetdash{}{0cm}
\pgfpathmoveto{\pgfqpoint{10.022cm}{10.53cm}}
\pgfpathlineto{\pgfqpoint{10.245cm}{10.53cm}}
\pgfpathlineto{\pgfqpoint{10.245cm}{10.307cm}}
\pgfpathlineto{\pgfqpoint{10.022cm}{10.307cm}}
\pgfpathlineto{\pgfqpoint{10.022cm}{10.53cm}}
\pgfpathclose
\pgfusepath{stroke}
\pgfsetdash{}{0cm}
\pgfpathmoveto{\pgfqpoint{11.509cm}{10.918cm}}
\pgfpathlineto{\pgfqpoint{11.733cm}{10.918cm}}
\pgfpathlineto{\pgfqpoint{11.733cm}{10.695cm}}
\pgfpathlineto{\pgfqpoint{11.509cm}{10.695cm}}
\pgfpathlineto{\pgfqpoint{11.509cm}{10.918cm}}
\pgfpathclose
\pgfusepath{stroke}
\pgfsetdash{}{0cm}
\pgfpathmoveto{\pgfqpoint{12.994cm}{11.312cm}}
\pgfpathlineto{\pgfqpoint{13.217cm}{11.312cm}}
\pgfpathlineto{\pgfqpoint{13.217cm}{11.089cm}}
\pgfpathlineto{\pgfqpoint{12.994cm}{11.089cm}}
\pgfpathlineto{\pgfqpoint{12.994cm}{11.312cm}}
\pgfpathclose
\pgfusepath{stroke}
\pgfsetdash{}{0cm}
\pgfpathmoveto{\pgfqpoint{14.482cm}{11.706cm}}
\pgfpathlineto{\pgfqpoint{14.705cm}{11.706cm}}
\pgfpathlineto{\pgfqpoint{14.705cm}{11.483cm}}
\pgfpathlineto{\pgfqpoint{14.482cm}{11.483cm}}
\pgfpathlineto{\pgfqpoint{14.482cm}{11.706cm}}
\pgfpathclose
\pgfusepath{stroke}
\pgfsetdash{}{0cm}
\pgfpathmoveto{\pgfqpoint{15.966cm}{12.103cm}}
\pgfpathlineto{\pgfqpoint{16.19cm}{12.103cm}}
\pgfpathlineto{\pgfqpoint{16.19cm}{11.88cm}}
\pgfpathlineto{\pgfqpoint{15.966cm}{11.88cm}}
\pgfpathlineto{\pgfqpoint{15.966cm}{12.103cm}}
\pgfpathclose
\pgfusepath{stroke}
\pgfsetdash{}{0cm}
\pgfpathmoveto{\pgfqpoint{17.454cm}{12.5cm}}
\pgfpathlineto{\pgfqpoint{17.677cm}{12.5cm}}
\pgfpathlineto{\pgfqpoint{17.677cm}{12.277cm}}
\pgfpathlineto{\pgfqpoint{17.454cm}{12.277cm}}
\pgfpathlineto{\pgfqpoint{17.454cm}{12.5cm}}
\pgfpathclose
\pgfusepath{stroke}
\pgfsetdash{}{0cm}
\pgfpathmoveto{\pgfqpoint{18.941cm}{12.897cm}}
\pgfpathlineto{\pgfqpoint{19.165cm}{12.897cm}}
\pgfpathlineto{\pgfqpoint{19.165cm}{12.674cm}}
\pgfpathlineto{\pgfqpoint{18.941cm}{12.674cm}}
\pgfpathlineto{\pgfqpoint{18.941cm}{12.897cm}}
\pgfpathclose
\pgfusepath{stroke}
\end{pgfscope}
\begin{pgfscope}
\pgfpathmoveto{\pgfqpoint{4.189cm}{20.446cm}}
\pgfpathlineto{\pgfqpoint{4.189cm}{8.611cm}}
\pgfpathlineto{\pgfqpoint{19.056cm}{8.611cm}}
\pgfpathlineto{\pgfqpoint{19.056cm}{20.446cm}}
\pgfpathclose
\pgfusepath{clip}
\end{pgfscope}
\pgftext[x=11.639cm,y=6.968cm,rotate=0]{  \fontsize{40}{50}{\selectfont $\bmax$}}
\pgftext[x=1.257cm,y=14.53cm,rotate=90]{  \fontsize{36}{36.14}\selectfont{ {Sub-optimality bound}}}
\pgftext[x=4.139cm,y=8.487cm,rotate=0]{\fontsize{10.04}{12.04}\selectfont{ { }}}
\pgftext[x=19.006cm,y=20.323cm,rotate=0]{\fontsize{10.04}{12.04}\selectfont{ { }}}
\pgftext[x=8.916cm+1.2cm,y=19.485cm,rotate=0]{  \fontsize{40}{36.14}\selectfont{ {$\la = 0.9$, \texttt{minS} }}}
\begin{pgfscope}
\pgfpathmoveto{\pgfqpoint{4.366cm}{20.27cm}}
\pgfpathlineto{\pgfqpoint{4.366cm}{15.158cm}}
\pgfpathlineto{\pgfqpoint{13.179cm}{15.158cm}}
\pgfpathlineto{\pgfqpoint{13.179cm}{20.27cm}}
\pgfpathclose
\pgfusepath{clip}
\pgfsetdash{}{0cm}
\pgfsetlinewidth{1.058mm}
\definecolor{eps2pgf_color}{rgb}{0,0,1}\pgfsetstrokecolor{eps2pgf_color}\pgfsetfillcolor{eps2pgf_color}
\pgfpathmoveto{\pgfqpoint{4.571cm}{19.591cm}}
\pgfpathlineto{\pgfqpoint{5.609cm}{19.591cm}}
\pgfusepath{stroke}
\begin{pgfscope}
\pgfpathmoveto{\pgfqpoint{4.698cm}{19.985cm}}
\pgfpathlineto{\pgfqpoint{4.698cm}{19.194cm}}
\pgfpathlineto{\pgfqpoint{5.489cm}{19.194cm}}
\pgfpathlineto{\pgfqpoint{5.489cm}{19.985cm}}
\pgfpathclose
\pgfusepath{clip}
\pgfsetdash{}{0cm}
\pgfpathmoveto{\pgfqpoint{4.951cm}{19.591cm}}
\pgfpathlineto{\pgfqpoint{5.233cm}{19.591cm}}
\pgfusepath{stroke}
\pgfsetdash{}{0cm}
\pgfpathmoveto{\pgfqpoint{5.092cm}{19.732cm}}
\pgfpathlineto{\pgfqpoint{5.092cm}{19.45cm}}
\pgfusepath{stroke}
\pgfsetdash{}{0cm}
\pgfpathmoveto{\pgfqpoint{4.995cm}{19.688cm}}
\pgfpathlineto{\pgfqpoint{5.189cm}{19.494cm}}
\pgfusepath{stroke}
\pgfsetdash{}{0cm}
\pgfpathmoveto{\pgfqpoint{5.189cm}{19.688cm}}
\pgfpathlineto{\pgfqpoint{4.995cm}{19.494cm}}
\pgfusepath{stroke}
\end{pgfscope}
\end{pgfscope}
\pgftext[x=9.093cm+1cm,y=18.236cm,rotate=0]{  \fontsize{40}{36.14}\selectfont{ {$\la = 0.9$, \texttt{maxW}  }}}
\begin{pgfscope}
\pgfpathmoveto{\pgfqpoint{4.366cm}{20.27cm}}
\pgfpathlineto{\pgfqpoint{4.366cm}{15.158cm}}
\pgfpathlineto{\pgfqpoint{13.179cm}{15.158cm}}
\pgfpathlineto{\pgfqpoint{13.179cm}{20.27cm}}
\pgfpathclose
\pgfusepath{clip}
\pgfsetdash{{0.018cm}{0.141cm}{0.212cm}{0.141cm}}{0cm}
\pgfsetlinewidth{1.058mm}
\definecolor{eps2pgf_color}{rgb}{0,0,1}\pgfsetstrokecolor{eps2pgf_color}\pgfsetfillcolor{eps2pgf_color}
\pgfpathmoveto{\pgfqpoint{4.571cm}{18.342cm}}
\pgfpathlineto{\pgfqpoint{5.609cm}{18.342cm}}
\pgfusepath{stroke}
\begin{pgfscope}
\pgfpathmoveto{\pgfqpoint{4.698cm}{18.735cm}}
\pgfpathlineto{\pgfqpoint{4.698cm}{17.945cm}}
\pgfpathlineto{\pgfqpoint{5.489cm}{17.945cm}}
\pgfpathlineto{\pgfqpoint{5.489cm}{18.735cm}}
\pgfpathclose
\pgfusepath{clip}
\pgfsetdash{}{0cm}
\pgfpathmoveto{\pgfqpoint{4.951cm}{18.342cm}}
\pgfpathlineto{\pgfqpoint{5.233cm}{18.342cm}}
\pgfusepath{stroke}
\pgfsetdash{}{0cm}
\pgfpathmoveto{\pgfqpoint{5.092cm}{18.483cm}}
\pgfpathlineto{\pgfqpoint{5.092cm}{18.2cm}}
\pgfusepath{stroke}
\pgfsetdash{}{0cm}
\pgfpathmoveto{\pgfqpoint{4.995cm}{18.439cm}}
\pgfpathlineto{\pgfqpoint{5.189cm}{18.244cm}}
\pgfusepath{stroke}
\pgfsetdash{}{0cm}
\pgfpathmoveto{\pgfqpoint{5.189cm}{18.439cm}}
\pgfpathlineto{\pgfqpoint{4.995cm}{18.244cm}}
\pgfusepath{stroke}
\end{pgfscope}
\end{pgfscope}
\pgftext[x=9.211cm+1.2cm,y=16.986cm,rotate=0]{  \fontsize{40}{36.14}\selectfont{ {$\la = 0.95$, \texttt{minS} }}}
\begin{pgfscope}
\pgfpathmoveto{\pgfqpoint{4.366cm}{20.27cm}}
\pgfpathlineto{\pgfqpoint{4.366cm}{15.158cm}}
\pgfpathlineto{\pgfqpoint{13.179cm}{15.158cm}}
\pgfpathlineto{\pgfqpoint{13.179cm}{20.27cm}}
\pgfpathclose
\pgfusepath{clip}
\pgfsetdash{}{0cm}
\pgfsetlinewidth{1.058mm}
\definecolor{eps2pgf_color}{rgb}{0,1,0}\pgfsetstrokecolor{eps2pgf_color}\pgfsetfillcolor{eps2pgf_color}
\pgfpathmoveto{\pgfqpoint{4.571cm}{17.092cm}}
\pgfpathlineto{\pgfqpoint{5.609cm}{17.092cm}}
\pgfusepath{stroke}
\begin{pgfscope}
\pgfpathmoveto{\pgfqpoint{4.698cm}{17.486cm}}
\pgfpathlineto{\pgfqpoint{4.698cm}{16.695cm}}
\pgfpathlineto{\pgfqpoint{5.489cm}{16.695cm}}
\pgfpathlineto{\pgfqpoint{5.489cm}{17.486cm}}
\pgfpathclose
\pgfusepath{clip}
\pgfsetdash{}{0cm}
\pgfsetmiterjoin
\pgfpathmoveto{\pgfqpoint{4.98cm}{17.204cm}}
\pgfpathlineto{\pgfqpoint{5.203cm}{17.204cm}}
\pgfpathlineto{\pgfqpoint{5.203cm}{16.98cm}}
\pgfpathlineto{\pgfqpoint{4.98cm}{16.98cm}}
\pgfpathlineto{\pgfqpoint{4.98cm}{17.204cm}}
\pgfpathclose
\pgfusepath{stroke}
\end{pgfscope}
\end{pgfscope}
\pgftext[x=9.387cm+1cm,y=15.737cm,rotate=0]{  \fontsize{40}{36.14}\selectfont{ {$\la = 0.95$, \texttt{maxW}  }}}
\begin{pgfscope}
\pgfpathmoveto{\pgfqpoint{4.366cm}{20.27cm}}
\pgfpathlineto{\pgfqpoint{4.366cm}{15.158cm}}
\pgfpathlineto{\pgfqpoint{13.179cm}{15.158cm}}
\pgfpathlineto{\pgfqpoint{13.179cm}{20.27cm}}
\pgfpathclose
\pgfusepath{clip}
\pgfsetdash{{0.018cm}{0.141cm}{0.212cm}{0.141cm}}{0cm}
\pgfsetlinewidth{1.058mm}
\definecolor{eps2pgf_color}{rgb}{0,1,0}\pgfsetstrokecolor{eps2pgf_color}\pgfsetfillcolor{eps2pgf_color}
\pgfpathmoveto{\pgfqpoint{4.571cm}{15.843cm}}
\pgfpathlineto{\pgfqpoint{5.609cm}{15.843cm}}
\pgfusepath{stroke}
\begin{pgfscope}
\pgfpathmoveto{\pgfqpoint{4.698cm}{16.237cm}}
\pgfpathlineto{\pgfqpoint{4.698cm}{15.446cm}}
\pgfpathlineto{\pgfqpoint{5.489cm}{15.446cm}}
\pgfpathlineto{\pgfqpoint{5.489cm}{16.237cm}}
\pgfpathclose
\pgfusepath{clip}
\pgfsetdash{}{0cm}
\pgfsetmiterjoin
\pgfpathmoveto{\pgfqpoint{4.98cm}{15.954cm}}
\pgfpathlineto{\pgfqpoint{5.203cm}{15.954cm}}
\pgfpathlineto{\pgfqpoint{5.203cm}{15.731cm}}
\pgfpathlineto{\pgfqpoint{4.98cm}{15.731cm}}
\pgfpathlineto{\pgfqpoint{4.98cm}{15.954cm}}
\pgfpathclose
\pgfusepath{stroke}
\end{pgfscope}
\end{pgfscope}
\pgfsetdash{}{0cm}
\definecolor{eps2pgf_color}{rgb}{0,1,0}\pgfsetstrokecolor{eps2pgf_color}\pgfsetfillcolor{eps2pgf_color}
\pgfusepath{stroke}
\end{pgfscope}
\end{pgfscope}
\end{pgfpicture}

%% file: fig/case2_1_2.tex
\scalebox{0.235}{\scalefont{2} \input{./fig/case2_1_2.pgf}}

%% file: fig/case2_1_2.pgf
% Created by Eps2pgf 0.7.0 (build on 2008-08-24) on Mon Apr 06 01:29:22 PDT 2015
\begin{pgfpicture}
\pgfpathmoveto{\pgfqpoint{0.635cm}{6.315cm}}
\pgfpathlineto{\pgfqpoint{19.473cm}{6.315cm}}
\pgfpathlineto{\pgfqpoint{19.473cm}{21.167cm}}
\pgfpathlineto{\pgfqpoint{0.635cm}{21.167cm}}
\pgfpathclose
\pgfusepath{clip}
\begin{pgfscope}
\begin{pgfscope}
\pgfpathmoveto{\pgfqpoint{0.635cm}{21.178cm}}
\pgfpathlineto{\pgfqpoint{0.635cm}{6.344cm}}
\pgfpathlineto{\pgfqpoint{19.491cm}{6.344cm}}
\pgfpathlineto{\pgfqpoint{19.491cm}{21.178cm}}
\pgfpathclose
\pgfusepath{clip}
\definecolor{eps2pgf_color}{gray}{1}\pgfsetstrokecolor{eps2pgf_color}\pgfsetfillcolor{eps2pgf_color}
\pgfpathmoveto{\pgfqpoint{0.635cm}{21.59cm}}
\pgfpathlineto{\pgfqpoint{0.635cm}{6.341cm}}
\pgfpathlineto{\pgfqpoint{20.99cm}{6.341cm}}
\pgfpathlineto{\pgfqpoint{20.99cm}{21.59cm}}
\pgfpathclose
\pgfusepath{fill}
\pgfpathmoveto{\pgfqpoint{4.777cm}{8.614cm}}
\pgfpathlineto{\pgfqpoint{4.777cm}{20.449cm}}
\pgfpathlineto{\pgfqpoint{19.053cm}{20.449cm}}
\pgfpathlineto{\pgfqpoint{19.053cm}{8.614cm}}
\pgfpathclose
\pgfseteorule\pgfusepath{fill}\pgfsetnonzerorule
\pgfsetdash{}{0cm}
\pgfsetlinewidth{0.176mm}
\pgfsetroundjoin
\pgfpathmoveto{\pgfqpoint{4.777cm}{8.614cm}}
\pgfpathlineto{\pgfqpoint{4.777cm}{20.449cm}}
\pgfpathlineto{\pgfqpoint{19.053cm}{20.449cm}}
\pgfpathlineto{\pgfqpoint{19.053cm}{8.614cm}}
\pgfpathlineto{\pgfqpoint{4.777cm}{8.614cm}}
\pgfusepath{stroke}
\pgfsetdash{}{0cm}
\definecolor{eps2pgf_color}{gray}{0}\pgfsetstrokecolor{eps2pgf_color}\pgfsetfillcolor{eps2pgf_color}
\pgfpathmoveto{\pgfqpoint{4.777cm}{8.614cm}}
\pgfpathlineto{\pgfqpoint{19.053cm}{8.614cm}}
\pgfusepath{stroke}
\pgfsetdash{}{0cm}
\pgfpathmoveto{\pgfqpoint{4.777cm}{20.449cm}}
\pgfpathlineto{\pgfqpoint{19.053cm}{20.449cm}}
\pgfusepath{stroke}
\pgfsetdash{}{0cm}
\pgfpathmoveto{\pgfqpoint{4.777cm}{8.614cm}}
\pgfpathlineto{\pgfqpoint{4.777cm}{20.449cm}}
\pgfusepath{stroke}
\pgfsetdash{}{0cm}
\pgfpathmoveto{\pgfqpoint{19.053cm}{8.614cm}}
\pgfpathlineto{\pgfqpoint{19.053cm}{20.449cm}}
\pgfusepath{stroke}
\pgfsetdash{}{0cm}
\pgfpathmoveto{\pgfqpoint{4.777cm}{8.614cm}}
\pgfpathlineto{\pgfqpoint{19.053cm}{8.614cm}}
\pgfusepath{stroke}
\pgfsetdash{}{0cm}
\pgfpathmoveto{\pgfqpoint{4.777cm}{8.614cm}}
\pgfpathlineto{\pgfqpoint{4.777cm}{20.449cm}}
\pgfusepath{stroke}
\pgfsetdash{}{0cm}
\pgfpathmoveto{\pgfqpoint{4.777cm}{8.614cm}}
\pgfpathlineto{\pgfqpoint{4.777cm}{8.758cm}}
\pgfusepath{stroke}
\pgfsetdash{}{0cm}
\pgfpathmoveto{\pgfqpoint{4.777cm}{20.446cm}}
\pgfpathlineto{\pgfqpoint{4.777cm}{20.305cm}}
\pgfusepath{stroke}
\pgftext[x=4.777cm,y=7.891cm+.2cm,rotate=0]{  \fontsize{36}{36.14}\selectfont{ {0}}}
\pgfsetdash{}{0cm}
\pgfpathmoveto{\pgfqpoint{11.915cm}{8.614cm}}
\pgfpathlineto{\pgfqpoint{11.915cm}{8.758cm}}
\pgfusepath{stroke}
\pgfsetdash{}{0cm}
\pgfpathmoveto{\pgfqpoint{11.915cm}{20.446cm}}
\pgfpathlineto{\pgfqpoint{11.915cm}{20.305cm}}
\pgfusepath{stroke}
\pgftext[x=11.913cm,y=7.891cm+.2cm,rotate=0]{  \fontsize{36}{36.14}\selectfont{ {0.5}}}
\pgfsetdash{}{0cm}
\pgfpathmoveto{\pgfqpoint{19.053cm}{8.614cm}}
\pgfpathlineto{\pgfqpoint{19.053cm}{8.758cm}}
\pgfusepath{stroke}
\pgfsetdash{}{0cm}
\pgfpathmoveto{\pgfqpoint{19.053cm}{20.446cm}}
\pgfpathlineto{\pgfqpoint{19.053cm}{20.305cm}}
\pgfusepath{stroke}
\pgftext[x=19.002cm-.5cm,y=7.901cm+.2cm,rotate=0]{  \fontsize{36}{36.14}\selectfont{ {1}}}
\pgfsetdash{}{0cm}
\pgfpathmoveto{\pgfqpoint{4.777cm}{8.614cm}}
\pgfpathlineto{\pgfqpoint{4.918cm}{8.614cm}}
\pgfusepath{stroke}
\pgfsetdash{}{0cm}
\pgfpathmoveto{\pgfqpoint{19.053cm}{8.614cm}}
\pgfpathlineto{\pgfqpoint{18.909cm}{8.614cm}}
\pgfusepath{stroke}
\pgftext[x=4.381cm-.3cm,y=8.582cm,rotate=0]{  \fontsize{36}{36.14}\selectfont{ {0}}}
\pgfsetdash{}{0cm}
\pgfpathmoveto{\pgfqpoint{4.777cm}{12.559cm}}
\pgfpathlineto{\pgfqpoint{4.918cm}{12.559cm}}
\pgfusepath{stroke}
\pgfsetdash{}{0cm}
\pgfpathmoveto{\pgfqpoint{19.053cm}{12.559cm}}
\pgfpathlineto{\pgfqpoint{18.909cm}{12.559cm}}
\pgfusepath{stroke}
\pgftext[x=3.35cm-.3cm,y=12.527cm,rotate=0]{  \fontsize{36}{36.14}\selectfont{ {0.005}}}
\pgfsetdash{}{0cm}
\pgfpathmoveto{\pgfqpoint{4.777cm}{16.504cm}}
\pgfpathlineto{\pgfqpoint{4.918cm}{16.504cm}}
\pgfusepath{stroke}
\pgfsetdash{}{0cm}
\pgfpathmoveto{\pgfqpoint{19.053cm}{16.504cm}}
\pgfpathlineto{\pgfqpoint{18.909cm}{16.504cm}}
\pgfusepath{stroke}
\pgftext[x=3.562cm-.3cm,y=16.472cm,rotate=0]{  \fontsize{36}{36.14}\selectfont{ {0.01}}}
\pgfsetdash{}{0cm}
\pgfpathmoveto{\pgfqpoint{4.777cm}{20.446cm}}
\pgfpathlineto{\pgfqpoint{4.918cm}{20.446cm}}
\pgfusepath{stroke}
\pgfsetdash{}{0cm}
\pgfpathmoveto{\pgfqpoint{19.053cm}{20.446cm}}
\pgfpathlineto{\pgfqpoint{18.909cm}{20.446cm}}
\pgfusepath{stroke}
\pgftext[x=3.35cm-.3cm,y=20.414cm,rotate=0]{  \fontsize{36}{36.14}\selectfont{ {0.015}}}
\pgfsetdash{}{0cm}
\pgfpathmoveto{\pgfqpoint{4.777cm}{8.614cm}}
\pgfpathlineto{\pgfqpoint{19.053cm}{8.614cm}}
\pgfusepath{stroke}
\pgfsetdash{}{0cm}
\pgfpathmoveto{\pgfqpoint{4.777cm}{20.449cm}}
\pgfpathlineto{\pgfqpoint{19.053cm}{20.449cm}}
\pgfusepath{stroke}
\pgfsetdash{}{0cm}
\pgfpathmoveto{\pgfqpoint{4.777cm}{8.614cm}}
\pgfpathlineto{\pgfqpoint{4.777cm}{20.449cm}}
\pgfusepath{stroke}
\pgfsetdash{}{0cm}
\pgfpathmoveto{\pgfqpoint{19.053cm}{8.614cm}}
\pgfpathlineto{\pgfqpoint{19.053cm}{20.449cm}}
\pgfusepath{stroke}
\begin{pgfscope}
\pgfpathmoveto{\pgfqpoint{4.777cm}{20.446cm}}
\pgfpathlineto{\pgfqpoint{4.777cm}{8.611cm}}
\pgfpathlineto{\pgfqpoint{19.056cm}{8.611cm}}
\pgfpathlineto{\pgfqpoint{19.056cm}{20.446cm}}
\pgfpathclose
\pgfusepath{clip}
\pgfsetdash{}{0cm}
\pgfsetlinewidth{1.058mm}
\definecolor{eps2pgf_color}{rgb}{0,0,1}\pgfsetstrokecolor{eps2pgf_color}\pgfsetfillcolor{eps2pgf_color}
\pgfpathmoveto{\pgfqpoint{6.203cm}{18.477cm}}
\pgfpathlineto{\pgfqpoint{7.632cm}{10.807cm}}
\pgfpathlineto{\pgfqpoint{9.058cm}{9.554cm}}
\pgfpathlineto{\pgfqpoint{10.486cm}{9.134cm}}
\pgfpathlineto{\pgfqpoint{11.915cm}{8.943cm}}
\pgfpathlineto{\pgfqpoint{13.341cm}{8.843cm}}
\pgfpathlineto{\pgfqpoint{14.77cm}{8.781cm}}
\pgfpathlineto{\pgfqpoint{16.195cm}{8.743cm}}
\pgfpathlineto{\pgfqpoint{17.624cm}{8.714cm}}
\pgfpathlineto{\pgfqpoint{19.053cm}{8.696cm}}
\pgfusepath{stroke}
\end{pgfscope}
\pgftext[x=11.933cm,y=6.968cm,rotate=0]{  \fontsize{40}{50}{\selectfont $\bmax$}}
\pgftext[x=1.249cm,y=14.553cm,rotate=90]{  \fontsize{36}{36.14}\selectfont{ {Bound/$\bmax$}}}
\pgftext[x=4.727cm,y=8.487cm,rotate=0]{\fontsize{10.04}{12.04}\selectfont{ { }}}
\pgftext[x=19.006cm,y=20.323cm,rotate=0]{\fontsize{10.04}{12.04}\selectfont{ { }}}
\pgftext[x=16.633cm,y=19.581cm,rotate=0]{  \fontsize{40}{36.14}\selectfont{ {  $\la = 1$}}}
\begin{pgfscope}
\pgfpathmoveto{\pgfqpoint{13.335cm}{20.27cm}}
\pgfpathlineto{\pgfqpoint{13.335cm}{18.897cm}}
\pgfpathlineto{\pgfqpoint{18.882cm}{18.897cm}}
\pgfpathlineto{\pgfqpoint{18.882cm}{20.27cm}}
\pgfpathclose
\pgfusepath{clip}
\pgfsetdash{}{0cm}
\pgfsetlinewidth{1.058mm}
\definecolor{eps2pgf_color}{rgb}{0,0,1}\pgfsetstrokecolor{eps2pgf_color}\pgfsetfillcolor{eps2pgf_color}
\pgfpathmoveto{\pgfqpoint{13.538cm}{19.585cm}}
\pgfpathlineto{\pgfqpoint{14.567cm}{19.585cm}}
\pgfusepath{stroke}
\end{pgfscope}
\pgfsetdash{}{0cm}
\definecolor{eps2pgf_color}{rgb}{0,0,1}\pgfsetstrokecolor{eps2pgf_color}\pgfsetfillcolor{eps2pgf_color}
\pgfusepath{stroke}
\end{pgfscope}
\end{pgfscope}
\end{pgfpicture}

%% file: fig/case2_2_2.tex
\scalebox{0.235}{\scalefont{2} \input{./fig/case2_2_2.pgf}}

%% file: fig/case2_2_2.pgf
% Created by Eps2pgf 0.7.0 (build on 2008-08-24) on Mon Apr 06 15:10:12 PDT 2015
\begin{pgfpicture}
\pgfpathmoveto{\pgfqpoint{0.635cm}{6.315cm}}
\pgfpathlineto{\pgfqpoint{19.579cm}{6.315cm}}
\pgfpathlineto{\pgfqpoint{19.579cm}{20.849cm}}
\pgfpathlineto{\pgfqpoint{0.635cm}{20.849cm}}
\pgfpathclose
\pgfusepath{clip}
\begin{pgfscope}
\begin{pgfscope}
\pgfpathmoveto{\pgfqpoint{0.635cm}{20.867cm}}
\pgfpathlineto{\pgfqpoint{0.635cm}{6.344cm}}
\pgfpathlineto{\pgfqpoint{19.591cm}{6.344cm}}
\pgfpathlineto{\pgfqpoint{19.591cm}{20.867cm}}
\pgfpathclose
\pgfusepath{clip}
\definecolor{eps2pgf_color}{gray}{1}\pgfsetstrokecolor{eps2pgf_color}\pgfsetfillcolor{eps2pgf_color}
\pgfpathmoveto{\pgfqpoint{0.635cm}{21.59cm}}
\pgfpathlineto{\pgfqpoint{0.635cm}{6.341cm}}
\pgfpathlineto{\pgfqpoint{20.99cm}{6.341cm}}
\pgfpathlineto{\pgfqpoint{20.99cm}{21.59cm}}
\pgfpathclose
\pgfusepath{fill}
\pgfpathmoveto{\pgfqpoint{4.189cm}{8.614cm}}
\pgfpathlineto{\pgfqpoint{4.189cm}{20.449cm}}
\pgfpathlineto{\pgfqpoint{19.053cm}{20.449cm}}
\pgfpathlineto{\pgfqpoint{19.053cm}{8.614cm}}
\pgfpathclose
\pgfseteorule\pgfusepath{fill}\pgfsetnonzerorule
\pgfsetdash{}{0cm}
\pgfsetlinewidth{0.176mm}
\pgfsetroundjoin
\pgfpathmoveto{\pgfqpoint{4.189cm}{8.614cm}}
\pgfpathlineto{\pgfqpoint{4.189cm}{20.449cm}}
\pgfpathlineto{\pgfqpoint{19.053cm}{20.449cm}}
\pgfpathlineto{\pgfqpoint{19.053cm}{8.614cm}}
\pgfpathlineto{\pgfqpoint{4.189cm}{8.614cm}}
\pgfusepath{stroke}
\pgfsetdash{}{0cm}
\definecolor{eps2pgf_color}{gray}{0}\pgfsetstrokecolor{eps2pgf_color}\pgfsetfillcolor{eps2pgf_color}
\pgfpathmoveto{\pgfqpoint{4.189cm}{8.614cm}}
\pgfpathlineto{\pgfqpoint{19.053cm}{8.614cm}}
\pgfusepath{stroke}
\pgfsetdash{}{0cm}
\pgfpathmoveto{\pgfqpoint{4.189cm}{20.449cm}}
\pgfpathlineto{\pgfqpoint{19.053cm}{20.449cm}}
\pgfusepath{stroke}
\pgfsetdash{}{0cm}
\pgfpathmoveto{\pgfqpoint{4.189cm}{8.614cm}}
\pgfpathlineto{\pgfqpoint{4.189cm}{20.449cm}}
\pgfusepath{stroke}
\pgfsetdash{}{0cm}
\pgfpathmoveto{\pgfqpoint{19.053cm}{8.614cm}}
\pgfpathlineto{\pgfqpoint{19.053cm}{20.449cm}}
\pgfusepath{stroke}
\pgfsetdash{}{0cm}
\pgfpathmoveto{\pgfqpoint{4.189cm}{8.614cm}}
\pgfpathlineto{\pgfqpoint{19.053cm}{8.614cm}}
\pgfusepath{stroke}
\pgfsetdash{}{0cm}
\pgfpathmoveto{\pgfqpoint{4.189cm}{8.614cm}}
\pgfpathlineto{\pgfqpoint{4.189cm}{20.449cm}}
\pgfusepath{stroke}
\pgfsetdash{}{0cm}
\pgfpathmoveto{\pgfqpoint{4.189cm}{8.614cm}}
\pgfpathlineto{\pgfqpoint{4.189cm}{8.764cm}}
\pgfusepath{stroke}
\pgfsetdash{}{0cm}
\pgfpathmoveto{\pgfqpoint{4.189cm}{20.446cm}}
\pgfpathlineto{\pgfqpoint{4.189cm}{20.299cm}}
\pgfusepath{stroke}
\pgftext[x=4.189cm,y=7.891cm+.2cm,rotate=0]{ \fontsize{36}{36.14}\selectfont{ {0}}}
\pgfsetdash{}{0cm}
\pgfpathmoveto{\pgfqpoint{11.621cm}{8.614cm}}
\pgfpathlineto{\pgfqpoint{11.621cm}{8.764cm}}
\pgfusepath{stroke}
\pgfsetdash{}{0cm}
\pgfpathmoveto{\pgfqpoint{11.621cm}{20.446cm}}
\pgfpathlineto{\pgfqpoint{11.621cm}{20.299cm}}
\pgfusepath{stroke}
\pgftext[x=11.619cm,y=7.891cm+.2cm,rotate=0]{ \fontsize{36}{36.14}\selectfont{ {0.5}}}
\pgfsetdash{}{0cm}
\pgfpathmoveto{\pgfqpoint{19.053cm}{8.614cm}}
\pgfpathlineto{\pgfqpoint{19.053cm}{8.764cm}}
\pgfusepath{stroke}
\pgfsetdash{}{0cm}
\pgfpathmoveto{\pgfqpoint{19.053cm}{20.446cm}}
\pgfpathlineto{\pgfqpoint{19.053cm}{20.299cm}}
\pgfusepath{stroke}
\pgftext[x=19.002cm-.5cm,y=7.901cm+.2cm,rotate=0]{ \fontsize{36}{36.14}\selectfont{ {1}}}
\pgfsetdash{}{0cm}
\pgfpathmoveto{\pgfqpoint{4.189cm}{8.614cm}}
\pgfpathlineto{\pgfqpoint{4.336cm}{8.614cm}}
\pgfusepath{stroke}
\pgfsetdash{}{0cm}
\pgfpathmoveto{\pgfqpoint{19.053cm}{8.614cm}}
\pgfpathlineto{\pgfqpoint{18.903cm}{8.614cm}}
\pgfusepath{stroke}
\pgftext[x=3.052cm-.3cm,y=8.582cm,rotate=0]{ \fontsize{36}{36.14}\selectfont{ {0.02}}}
\pgfsetdash{}{0cm}
\pgfpathmoveto{\pgfqpoint{4.189cm}{10.766cm}}
\pgfpathlineto{\pgfqpoint{4.336cm}{10.766cm}}
\pgfusepath{stroke}
\pgfsetdash{}{0cm}
\pgfpathmoveto{\pgfqpoint{19.053cm}{10.766cm}}
\pgfpathlineto{\pgfqpoint{18.903cm}{10.766cm}}
\pgfusepath{stroke}
\pgftext[x=3.06cm-.3cm,y=10.734cm,rotate=0]{ \fontsize{36}{36.14}\selectfont{ {0.04}}}
\pgfsetdash{}{0cm}
\pgfpathmoveto{\pgfqpoint{4.189cm}{12.918cm}}
\pgfpathlineto{\pgfqpoint{4.336cm}{12.918cm}}
\pgfusepath{stroke}
\pgfsetdash{}{0cm}
\pgfpathmoveto{\pgfqpoint{19.053cm}{12.918cm}}
\pgfpathlineto{\pgfqpoint{18.903cm}{12.918cm}}
\pgfusepath{stroke}
\pgftext[x=3.058cm-.3cm,y=12.886cm,rotate=0]{ \fontsize{36}{36.14}\selectfont{ {0.06}}}
\pgfsetdash{}{0cm}
\pgfpathmoveto{\pgfqpoint{4.189cm}{15.069cm}}
\pgfpathlineto{\pgfqpoint{4.336cm}{15.069cm}}
\pgfusepath{stroke}
\pgfsetdash{}{0cm}
\pgfpathmoveto{\pgfqpoint{19.053cm}{15.069cm}}
\pgfpathlineto{\pgfqpoint{18.903cm}{15.069cm}}
\pgfusepath{stroke}
\pgftext[x=3.057cm-.3cm,y=15.038cm,rotate=0]{ \fontsize{36}{36.14}\selectfont{ {0.08}}}
\pgfsetdash{}{0cm}
\pgfpathmoveto{\pgfqpoint{4.189cm}{17.221cm}}
\pgfpathlineto{\pgfqpoint{4.336cm}{17.221cm}}
\pgfusepath{stroke}
\pgfsetdash{}{0cm}
\pgfpathmoveto{\pgfqpoint{19.053cm}{17.221cm}}
\pgfpathlineto{\pgfqpoint{18.903cm}{17.221cm}}
\pgfusepath{stroke}
\pgftext[x=3.267cm-.3cm,y=17.189cm,rotate=0]{ \fontsize{36}{36.14}\selectfont{ {0.1}}}
\pgfsetdash{}{0cm}
\pgfpathmoveto{\pgfqpoint{4.189cm}{19.373cm}}
\pgfpathlineto{\pgfqpoint{4.336cm}{19.373cm}}
\pgfusepath{stroke}
\pgfsetdash{}{0cm}
\pgfpathmoveto{\pgfqpoint{19.053cm}{19.373cm}}
\pgfpathlineto{\pgfqpoint{18.903cm}{19.373cm}}
\pgfusepath{stroke}
\pgftext[x=3.052cm-.3cm,y=19.341cm,rotate=0]{ \fontsize{36}{36.14}\selectfont{ {0.12}}}
\pgfsetdash{}{0cm}
\pgfpathmoveto{\pgfqpoint{4.189cm}{8.614cm}}
\pgfpathlineto{\pgfqpoint{19.053cm}{8.614cm}}
\pgfusepath{stroke}
\pgfsetdash{}{0cm}
\pgfpathmoveto{\pgfqpoint{4.189cm}{20.449cm}}
\pgfpathlineto{\pgfqpoint{19.053cm}{20.449cm}}
\pgfusepath{stroke}
\pgfsetdash{}{0cm}
\pgfpathmoveto{\pgfqpoint{4.189cm}{8.614cm}}
\pgfpathlineto{\pgfqpoint{4.189cm}{20.449cm}}
\pgfusepath{stroke}
\pgfsetdash{}{0cm}
\pgfpathmoveto{\pgfqpoint{19.053cm}{8.614cm}}
\pgfpathlineto{\pgfqpoint{19.053cm}{20.449cm}}
\pgfusepath{stroke}
\begin{pgfscope}
\pgfpathmoveto{\pgfqpoint{4.189cm}{20.449cm}}
\pgfpathlineto{\pgfqpoint{4.189cm}{8.611cm}}
\pgfpathlineto{\pgfqpoint{19.056cm}{8.611cm}}
\pgfpathlineto{\pgfqpoint{19.056cm}{20.449cm}}
\pgfpathclose
\pgfusepath{clip}
\pgfsetdash{}{0cm}
\pgfsetlinewidth{1.058mm}
\definecolor{eps2pgf_color}{rgb}{0,0,1}\pgfsetstrokecolor{eps2pgf_color}\pgfsetfillcolor{eps2pgf_color}
\pgfpathmoveto{\pgfqpoint{5.674cm}{16.836cm}}
\pgfpathlineto{\pgfqpoint{7.161cm}{13.858cm}}
\pgfpathlineto{\pgfqpoint{8.646cm}{13.17cm}}
\pgfpathlineto{\pgfqpoint{10.134cm}{12.871cm}}
\pgfpathlineto{\pgfqpoint{11.621cm}{12.706cm}}
\pgfpathlineto{\pgfqpoint{13.106cm}{12.6cm}}
\pgfpathlineto{\pgfqpoint{14.593cm}{12.529cm}}
\pgfpathlineto{\pgfqpoint{16.078cm}{12.477cm}}
\pgfpathlineto{\pgfqpoint{17.565cm}{12.435cm}}
\pgfpathlineto{\pgfqpoint{19.053cm}{12.403cm}}
\pgfusepath{stroke}
\end{pgfscope}
\pgfsetdash{}{0cm}
\pgfsetlinewidth{1.058mm}
\definecolor{eps2pgf_color}{rgb}{0,0,1}\pgfsetstrokecolor{eps2pgf_color}\pgfsetfillcolor{eps2pgf_color}
\pgfpathmoveto{\pgfqpoint{5.533cm}{16.836cm}}
\pgfpathlineto{\pgfqpoint{5.815cm}{16.836cm}}
\pgfusepath{stroke}
\pgfsetdash{}{0cm}
\pgfpathmoveto{\pgfqpoint{5.674cm}{16.977cm}}
\pgfpathlineto{\pgfqpoint{5.674cm}{16.695cm}}
\pgfusepath{stroke}
\pgfsetdash{}{0cm}
\pgfpathmoveto{\pgfqpoint{7.02cm}{13.858cm}}
\pgfpathlineto{\pgfqpoint{7.302cm}{13.858cm}}
\pgfusepath{stroke}
\pgfsetdash{}{0cm}
\pgfpathmoveto{\pgfqpoint{7.161cm}{13.999cm}}
\pgfpathlineto{\pgfqpoint{7.161cm}{13.717cm}}
\pgfusepath{stroke}
\pgfsetdash{}{0cm}
\pgfpathmoveto{\pgfqpoint{8.505cm}{13.17cm}}
\pgfpathlineto{\pgfqpoint{8.787cm}{13.17cm}}
\pgfusepath{stroke}
\pgfsetdash{}{0cm}
\pgfpathmoveto{\pgfqpoint{8.646cm}{13.311cm}}
\pgfpathlineto{\pgfqpoint{8.646cm}{13.029cm}}
\pgfusepath{stroke}
\pgfsetdash{}{0cm}
\pgfpathmoveto{\pgfqpoint{9.992cm}{12.871cm}}
\pgfpathlineto{\pgfqpoint{10.275cm}{12.871cm}}
\pgfusepath{stroke}
\pgfsetdash{}{0cm}
\pgfpathmoveto{\pgfqpoint{10.134cm}{13.012cm}}
\pgfpathlineto{\pgfqpoint{10.134cm}{12.729cm}}
\pgfusepath{stroke}
\pgfsetdash{}{0cm}
\pgfpathmoveto{\pgfqpoint{11.48cm}{12.706cm}}
\pgfpathlineto{\pgfqpoint{11.762cm}{12.706cm}}
\pgfusepath{stroke}
\pgfsetdash{}{0cm}
\pgfpathmoveto{\pgfqpoint{11.621cm}{12.847cm}}
\pgfpathlineto{\pgfqpoint{11.621cm}{12.565cm}}
\pgfusepath{stroke}
\pgfsetdash{}{0cm}
\pgfpathmoveto{\pgfqpoint{12.965cm}{12.6cm}}
\pgfpathlineto{\pgfqpoint{13.247cm}{12.6cm}}
\pgfusepath{stroke}
\pgfsetdash{}{0cm}
\pgfpathmoveto{\pgfqpoint{13.106cm}{12.741cm}}
\pgfpathlineto{\pgfqpoint{13.106cm}{12.459cm}}
\pgfusepath{stroke}
\pgfsetdash{}{0cm}
\pgfpathmoveto{\pgfqpoint{14.452cm}{12.529cm}}
\pgfpathlineto{\pgfqpoint{14.734cm}{12.529cm}}
\pgfusepath{stroke}
\pgfsetdash{}{0cm}
\pgfpathmoveto{\pgfqpoint{14.593cm}{12.671cm}}
\pgfpathlineto{\pgfqpoint{14.593cm}{12.388cm}}
\pgfusepath{stroke}
\pgfsetdash{}{0cm}
\pgfpathmoveto{\pgfqpoint{15.937cm}{12.477cm}}
\pgfpathlineto{\pgfqpoint{16.219cm}{12.477cm}}
\pgfusepath{stroke}
\pgfsetdash{}{0cm}
\pgfpathmoveto{\pgfqpoint{16.078cm}{12.618cm}}
\pgfpathlineto{\pgfqpoint{16.078cm}{12.335cm}}
\pgfusepath{stroke}
\pgfsetdash{}{0cm}
\pgfpathmoveto{\pgfqpoint{17.424cm}{12.435cm}}
\pgfpathlineto{\pgfqpoint{17.707cm}{12.435cm}}
\pgfusepath{stroke}
\pgfsetdash{}{0cm}
\pgfpathmoveto{\pgfqpoint{17.565cm}{12.577cm}}
\pgfpathlineto{\pgfqpoint{17.565cm}{12.294cm}}
\pgfusepath{stroke}
\pgfsetdash{}{0cm}
\pgfpathmoveto{\pgfqpoint{18.912cm}{12.403cm}}
\pgfpathlineto{\pgfqpoint{19.194cm}{12.403cm}}
\pgfusepath{stroke}
\pgfsetdash{}{0cm}
\pgfpathmoveto{\pgfqpoint{19.053cm}{12.544cm}}
\pgfpathlineto{\pgfqpoint{19.053cm}{12.262cm}}
\pgfusepath{stroke}
\pgfsetdash{}{0cm}
\pgfpathmoveto{\pgfqpoint{5.577cm}{16.933cm}}
\pgfpathlineto{\pgfqpoint{5.771cm}{16.739cm}}
\pgfusepath{stroke}
\pgfsetdash{}{0cm}
\pgfpathmoveto{\pgfqpoint{5.771cm}{16.933cm}}
\pgfpathlineto{\pgfqpoint{5.577cm}{16.739cm}}
\pgfusepath{stroke}
\pgfsetdash{}{0cm}
\pgfpathmoveto{\pgfqpoint{7.064cm}{13.955cm}}
\pgfpathlineto{\pgfqpoint{7.258cm}{13.761cm}}
\pgfusepath{stroke}
\pgfsetdash{}{0cm}
\pgfpathmoveto{\pgfqpoint{7.258cm}{13.955cm}}
\pgfpathlineto{\pgfqpoint{7.064cm}{13.761cm}}
\pgfusepath{stroke}
\pgfsetdash{}{0cm}
\pgfpathmoveto{\pgfqpoint{8.549cm}{13.267cm}}
\pgfpathlineto{\pgfqpoint{8.743cm}{13.073cm}}
\pgfusepath{stroke}
\pgfsetdash{}{0cm}
\pgfpathmoveto{\pgfqpoint{8.743cm}{13.267cm}}
\pgfpathlineto{\pgfqpoint{8.549cm}{13.073cm}}
\pgfusepath{stroke}
\pgfsetdash{}{0cm}
\pgfpathmoveto{\pgfqpoint{10.037cm}{12.968cm}}
\pgfpathlineto{\pgfqpoint{10.231cm}{12.773cm}}
\pgfusepath{stroke}
\pgfsetdash{}{0cm}
\pgfpathmoveto{\pgfqpoint{10.231cm}{12.968cm}}
\pgfpathlineto{\pgfqpoint{10.037cm}{12.773cm}}
\pgfusepath{stroke}
\pgfsetdash{}{0cm}
\pgfpathmoveto{\pgfqpoint{11.524cm}{12.803cm}}
\pgfpathlineto{\pgfqpoint{11.718cm}{12.609cm}}
\pgfusepath{stroke}
\pgfsetdash{}{0cm}
\pgfpathmoveto{\pgfqpoint{11.718cm}{12.803cm}}
\pgfpathlineto{\pgfqpoint{11.524cm}{12.609cm}}
\pgfusepath{stroke}
\pgfsetdash{}{0cm}
\pgfpathmoveto{\pgfqpoint{13.009cm}{12.697cm}}
\pgfpathlineto{\pgfqpoint{13.203cm}{12.503cm}}
\pgfusepath{stroke}
\pgfsetdash{}{0cm}
\pgfpathmoveto{\pgfqpoint{13.203cm}{12.697cm}}
\pgfpathlineto{\pgfqpoint{13.009cm}{12.503cm}}
\pgfusepath{stroke}
\pgfsetdash{}{0cm}
\pgfpathmoveto{\pgfqpoint{14.496cm}{12.627cm}}
\pgfpathlineto{\pgfqpoint{14.69cm}{12.432cm}}
\pgfusepath{stroke}
\pgfsetdash{}{0cm}
\pgfpathmoveto{\pgfqpoint{14.69cm}{12.627cm}}
\pgfpathlineto{\pgfqpoint{14.496cm}{12.432cm}}
\pgfusepath{stroke}
\pgfsetdash{}{0cm}
\pgfpathmoveto{\pgfqpoint{15.981cm}{12.574cm}}
\pgfpathlineto{\pgfqpoint{16.175cm}{12.38cm}}
\pgfusepath{stroke}
\pgfsetdash{}{0cm}
\pgfpathmoveto{\pgfqpoint{16.175cm}{12.574cm}}
\pgfpathlineto{\pgfqpoint{15.981cm}{12.38cm}}
\pgfusepath{stroke}
\pgfsetdash{}{0cm}
\pgfpathmoveto{\pgfqpoint{17.468cm}{12.532cm}}
\pgfpathlineto{\pgfqpoint{17.662cm}{12.338cm}}
\pgfusepath{stroke}
\pgfsetdash{}{0cm}
\pgfpathmoveto{\pgfqpoint{17.662cm}{12.532cm}}
\pgfpathlineto{\pgfqpoint{17.468cm}{12.338cm}}
\pgfusepath{stroke}
\pgfsetdash{}{0cm}
\pgfpathmoveto{\pgfqpoint{18.956cm}{12.5cm}}
\pgfpathlineto{\pgfqpoint{19.15cm}{12.306cm}}
\pgfusepath{stroke}
\pgfsetdash{}{0cm}
\pgfpathmoveto{\pgfqpoint{19.15cm}{12.5cm}}
\pgfpathlineto{\pgfqpoint{18.956cm}{12.306cm}}
\pgfusepath{stroke}
\begin{pgfscope}
\pgfpathmoveto{\pgfqpoint{4.189cm}{20.449cm}}
\pgfpathlineto{\pgfqpoint{4.189cm}{8.611cm}}
\pgfpathlineto{\pgfqpoint{19.056cm}{8.611cm}}
\pgfpathlineto{\pgfqpoint{19.056cm}{20.449cm}}
\pgfpathclose
\pgfusepath{clip}
\pgfsetdash{{0.018cm}{0.141cm}{0.212cm}{0.141cm}}{0cm}
\pgfpathmoveto{\pgfqpoint{5.674cm}{17.186cm}}
\pgfpathlineto{\pgfqpoint{7.161cm}{14.434cm}}
\pgfpathlineto{\pgfqpoint{8.646cm}{13.776cm}}
\pgfpathlineto{\pgfqpoint{10.134cm}{13.488cm}}
\pgfpathlineto{\pgfqpoint{11.621cm}{13.329cm}}
\pgfpathlineto{\pgfqpoint{13.106cm}{13.226cm}}
\pgfpathlineto{\pgfqpoint{14.593cm}{13.156cm}}
\pgfpathlineto{\pgfqpoint{16.078cm}{13.103cm}}
\pgfpathlineto{\pgfqpoint{17.565cm}{13.065cm}}
\pgfpathlineto{\pgfqpoint{19.053cm}{13.032cm}}
\pgfusepath{stroke}
\end{pgfscope}
\pgfsetdash{}{0cm}
\pgfpathmoveto{\pgfqpoint{5.533cm}{17.186cm}}
\pgfpathlineto{\pgfqpoint{5.815cm}{17.186cm}}
\pgfusepath{stroke}
\pgfsetdash{}{0cm}
\pgfpathmoveto{\pgfqpoint{5.674cm}{17.327cm}}
\pgfpathlineto{\pgfqpoint{5.674cm}{17.045cm}}
\pgfusepath{stroke}
\pgfsetdash{}{0cm}
\pgfpathmoveto{\pgfqpoint{7.02cm}{14.434cm}}
\pgfpathlineto{\pgfqpoint{7.302cm}{14.434cm}}
\pgfusepath{stroke}
\pgfsetdash{}{0cm}
\pgfpathmoveto{\pgfqpoint{7.161cm}{14.576cm}}
\pgfpathlineto{\pgfqpoint{7.161cm}{14.293cm}}
\pgfusepath{stroke}
\pgfsetdash{}{0cm}
\pgfpathmoveto{\pgfqpoint{8.505cm}{13.776cm}}
\pgfpathlineto{\pgfqpoint{8.787cm}{13.776cm}}
\pgfusepath{stroke}
\pgfsetdash{}{0cm}
\pgfpathmoveto{\pgfqpoint{8.646cm}{13.917cm}}
\pgfpathlineto{\pgfqpoint{8.646cm}{13.635cm}}
\pgfusepath{stroke}
\pgfsetdash{}{0cm}
\pgfpathmoveto{\pgfqpoint{9.992cm}{13.488cm}}
\pgfpathlineto{\pgfqpoint{10.275cm}{13.488cm}}
\pgfusepath{stroke}
\pgfsetdash{}{0cm}
\pgfpathmoveto{\pgfqpoint{10.134cm}{13.629cm}}
\pgfpathlineto{\pgfqpoint{10.134cm}{13.347cm}}
\pgfusepath{stroke}
\pgfsetdash{}{0cm}
\pgfpathmoveto{\pgfqpoint{11.48cm}{13.329cm}}
\pgfpathlineto{\pgfqpoint{11.762cm}{13.329cm}}
\pgfusepath{stroke}
\pgfsetdash{}{0cm}
\pgfpathmoveto{\pgfqpoint{11.621cm}{13.47cm}}
\pgfpathlineto{\pgfqpoint{11.621cm}{13.188cm}}
\pgfusepath{stroke}
\pgfsetdash{}{0cm}
\pgfpathmoveto{\pgfqpoint{12.965cm}{13.226cm}}
\pgfpathlineto{\pgfqpoint{13.247cm}{13.226cm}}
\pgfusepath{stroke}
\pgfsetdash{}{0cm}
\pgfpathmoveto{\pgfqpoint{13.106cm}{13.367cm}}
\pgfpathlineto{\pgfqpoint{13.106cm}{13.085cm}}
\pgfusepath{stroke}
\pgfsetdash{}{0cm}
\pgfpathmoveto{\pgfqpoint{14.452cm}{13.156cm}}
\pgfpathlineto{\pgfqpoint{14.734cm}{13.156cm}}
\pgfusepath{stroke}
\pgfsetdash{}{0cm}
\pgfpathmoveto{\pgfqpoint{14.593cm}{13.297cm}}
\pgfpathlineto{\pgfqpoint{14.593cm}{13.015cm}}
\pgfusepath{stroke}
\pgfsetdash{}{0cm}
\pgfpathmoveto{\pgfqpoint{15.937cm}{13.103cm}}
\pgfpathlineto{\pgfqpoint{16.219cm}{13.103cm}}
\pgfusepath{stroke}
\pgfsetdash{}{0cm}
\pgfpathmoveto{\pgfqpoint{16.078cm}{13.244cm}}
\pgfpathlineto{\pgfqpoint{16.078cm}{12.962cm}}
\pgfusepath{stroke}
\pgfsetdash{}{0cm}
\pgfpathmoveto{\pgfqpoint{17.424cm}{13.065cm}}
\pgfpathlineto{\pgfqpoint{17.707cm}{13.065cm}}
\pgfusepath{stroke}
\pgfsetdash{}{0cm}
\pgfpathmoveto{\pgfqpoint{17.565cm}{13.206cm}}
\pgfpathlineto{\pgfqpoint{17.565cm}{12.923cm}}
\pgfusepath{stroke}
\pgfsetdash{}{0cm}
\pgfpathmoveto{\pgfqpoint{18.912cm}{13.032cm}}
\pgfpathlineto{\pgfqpoint{19.194cm}{13.032cm}}
\pgfusepath{stroke}
\pgfsetdash{}{0cm}
\pgfpathmoveto{\pgfqpoint{19.053cm}{13.173cm}}
\pgfpathlineto{\pgfqpoint{19.053cm}{12.891cm}}
\pgfusepath{stroke}
\pgfsetdash{}{0cm}
\pgfpathmoveto{\pgfqpoint{5.577cm}{17.283cm}}
\pgfpathlineto{\pgfqpoint{5.771cm}{17.089cm}}
\pgfusepath{stroke}
\pgfsetdash{}{0cm}
\pgfpathmoveto{\pgfqpoint{5.771cm}{17.283cm}}
\pgfpathlineto{\pgfqpoint{5.577cm}{17.089cm}}
\pgfusepath{stroke}
\pgfsetdash{}{0cm}
\pgfpathmoveto{\pgfqpoint{7.064cm}{14.532cm}}
\pgfpathlineto{\pgfqpoint{7.258cm}{14.337cm}}
\pgfusepath{stroke}
\pgfsetdash{}{0cm}
\pgfpathmoveto{\pgfqpoint{7.258cm}{14.532cm}}
\pgfpathlineto{\pgfqpoint{7.064cm}{14.337cm}}
\pgfusepath{stroke}
\pgfsetdash{}{0cm}
\pgfpathmoveto{\pgfqpoint{8.549cm}{13.873cm}}
\pgfpathlineto{\pgfqpoint{8.743cm}{13.679cm}}
\pgfusepath{stroke}
\pgfsetdash{}{0cm}
\pgfpathmoveto{\pgfqpoint{8.743cm}{13.873cm}}
\pgfpathlineto{\pgfqpoint{8.549cm}{13.679cm}}
\pgfusepath{stroke}
\pgfsetdash{}{0cm}
\pgfpathmoveto{\pgfqpoint{10.037cm}{13.585cm}}
\pgfpathlineto{\pgfqpoint{10.231cm}{13.391cm}}
\pgfusepath{stroke}
\pgfsetdash{}{0cm}
\pgfpathmoveto{\pgfqpoint{10.231cm}{13.585cm}}
\pgfpathlineto{\pgfqpoint{10.037cm}{13.391cm}}
\pgfusepath{stroke}
\pgfsetdash{}{0cm}
\pgfpathmoveto{\pgfqpoint{11.524cm}{13.426cm}}
\pgfpathlineto{\pgfqpoint{11.718cm}{13.232cm}}
\pgfusepath{stroke}
\pgfsetdash{}{0cm}
\pgfpathmoveto{\pgfqpoint{11.718cm}{13.426cm}}
\pgfpathlineto{\pgfqpoint{11.524cm}{13.232cm}}
\pgfusepath{stroke}
\pgfsetdash{}{0cm}
\pgfpathmoveto{\pgfqpoint{13.009cm}{13.323cm}}
\pgfpathlineto{\pgfqpoint{13.203cm}{13.129cm}}
\pgfusepath{stroke}
\pgfsetdash{}{0cm}
\pgfpathmoveto{\pgfqpoint{13.203cm}{13.323cm}}
\pgfpathlineto{\pgfqpoint{13.009cm}{13.129cm}}
\pgfusepath{stroke}
\pgfsetdash{}{0cm}
\pgfpathmoveto{\pgfqpoint{14.496cm}{13.253cm}}
\pgfpathlineto{\pgfqpoint{14.69cm}{13.059cm}}
\pgfusepath{stroke}
\pgfsetdash{}{0cm}
\pgfpathmoveto{\pgfqpoint{14.69cm}{13.253cm}}
\pgfpathlineto{\pgfqpoint{14.496cm}{13.059cm}}
\pgfusepath{stroke}
\pgfsetdash{}{0cm}
\pgfpathmoveto{\pgfqpoint{15.981cm}{13.2cm}}
\pgfpathlineto{\pgfqpoint{16.175cm}{13.006cm}}
\pgfusepath{stroke}
\pgfsetdash{}{0cm}
\pgfpathmoveto{\pgfqpoint{16.175cm}{13.2cm}}
\pgfpathlineto{\pgfqpoint{15.981cm}{13.006cm}}
\pgfusepath{stroke}
\pgfsetdash{}{0cm}
\pgfpathmoveto{\pgfqpoint{17.468cm}{13.162cm}}
\pgfpathlineto{\pgfqpoint{17.662cm}{12.968cm}}
\pgfusepath{stroke}
\pgfsetdash{}{0cm}
\pgfpathmoveto{\pgfqpoint{17.662cm}{13.162cm}}
\pgfpathlineto{\pgfqpoint{17.468cm}{12.968cm}}
\pgfusepath{stroke}
\pgfsetdash{}{0cm}
\pgfpathmoveto{\pgfqpoint{18.956cm}{13.129cm}}
\pgfpathlineto{\pgfqpoint{19.15cm}{12.935cm}}
\pgfusepath{stroke}
\pgfsetdash{}{0cm}
\pgfpathmoveto{\pgfqpoint{19.15cm}{13.129cm}}
\pgfpathlineto{\pgfqpoint{18.956cm}{12.935cm}}
\pgfusepath{stroke}
\begin{pgfscope}
\pgfpathmoveto{\pgfqpoint{4.189cm}{20.449cm}}
\pgfpathlineto{\pgfqpoint{4.189cm}{8.611cm}}
\pgfpathlineto{\pgfqpoint{19.056cm}{8.611cm}}
\pgfpathlineto{\pgfqpoint{19.056cm}{20.449cm}}
\pgfpathclose
\pgfusepath{clip}
\pgfsetdash{}{0cm}
\definecolor{eps2pgf_color}{rgb}{0,1,0}\pgfsetstrokecolor{eps2pgf_color}\pgfsetfillcolor{eps2pgf_color}
\pgfpathmoveto{\pgfqpoint{5.674cm}{12.112cm}}
\pgfpathlineto{\pgfqpoint{7.161cm}{10.181cm}}
\pgfpathlineto{\pgfqpoint{8.646cm}{9.769cm}}
\pgfpathlineto{\pgfqpoint{10.134cm}{9.601cm}}
\pgfpathlineto{\pgfqpoint{11.621cm}{9.51cm}}
\pgfpathlineto{\pgfqpoint{13.106cm}{9.454cm}}
\pgfpathlineto{\pgfqpoint{14.593cm}{9.416cm}}
\pgfpathlineto{\pgfqpoint{16.078cm}{9.39cm}}
\pgfpathlineto{\pgfqpoint{17.565cm}{9.369cm}}
\pgfpathlineto{\pgfqpoint{19.053cm}{9.352cm}}
\pgfusepath{stroke}
\end{pgfscope}
\pgfsetdash{}{0cm}
\pgfsetmiterjoin
\definecolor{eps2pgf_color}{rgb}{0,1,0}\pgfsetstrokecolor{eps2pgf_color}\pgfsetfillcolor{eps2pgf_color}
\pgfpathmoveto{\pgfqpoint{5.562cm}{12.224cm}}
\pgfpathlineto{\pgfqpoint{5.786cm}{12.224cm}}
\pgfpathlineto{\pgfqpoint{5.786cm}{12cm}}
\pgfpathlineto{\pgfqpoint{5.562cm}{12cm}}
\pgfpathlineto{\pgfqpoint{5.562cm}{12.224cm}}
\pgfpathclose
\pgfusepath{stroke}
\pgfsetdash{}{0cm}
\pgfpathmoveto{\pgfqpoint{7.05cm}{10.292cm}}
\pgfpathlineto{\pgfqpoint{7.273cm}{10.292cm}}
\pgfpathlineto{\pgfqpoint{7.273cm}{10.069cm}}
\pgfpathlineto{\pgfqpoint{7.05cm}{10.069cm}}
\pgfpathlineto{\pgfqpoint{7.05cm}{10.292cm}}
\pgfpathclose
\pgfusepath{stroke}
\pgfsetdash{}{0cm}
\pgfpathmoveto{\pgfqpoint{8.534cm}{9.881cm}}
\pgfpathlineto{\pgfqpoint{8.758cm}{9.881cm}}
\pgfpathlineto{\pgfqpoint{8.758cm}{9.657cm}}
\pgfpathlineto{\pgfqpoint{8.534cm}{9.657cm}}
\pgfpathlineto{\pgfqpoint{8.534cm}{9.881cm}}
\pgfpathclose
\pgfusepath{stroke}
\pgfsetdash{}{0cm}
\pgfpathmoveto{\pgfqpoint{10.022cm}{9.713cm}}
\pgfpathlineto{\pgfqpoint{10.245cm}{9.713cm}}
\pgfpathlineto{\pgfqpoint{10.245cm}{9.49cm}}
\pgfpathlineto{\pgfqpoint{10.022cm}{9.49cm}}
\pgfpathlineto{\pgfqpoint{10.022cm}{9.713cm}}
\pgfpathclose
\pgfusepath{stroke}
\pgfsetdash{}{0cm}
\pgfpathmoveto{\pgfqpoint{11.509cm}{9.622cm}}
\pgfpathlineto{\pgfqpoint{11.733cm}{9.622cm}}
\pgfpathlineto{\pgfqpoint{11.733cm}{9.399cm}}
\pgfpathlineto{\pgfqpoint{11.509cm}{9.399cm}}
\pgfpathlineto{\pgfqpoint{11.509cm}{9.622cm}}
\pgfpathclose
\pgfusepath{stroke}
\pgfsetdash{}{0cm}
\pgfpathmoveto{\pgfqpoint{12.994cm}{9.566cm}}
\pgfpathlineto{\pgfqpoint{13.217cm}{9.566cm}}
\pgfpathlineto{\pgfqpoint{13.217cm}{9.343cm}}
\pgfpathlineto{\pgfqpoint{12.994cm}{9.343cm}}
\pgfpathlineto{\pgfqpoint{12.994cm}{9.566cm}}
\pgfpathclose
\pgfusepath{stroke}
\pgfsetdash{}{0cm}
\pgfpathmoveto{\pgfqpoint{14.482cm}{9.528cm}}
\pgfpathlineto{\pgfqpoint{14.705cm}{9.528cm}}
\pgfpathlineto{\pgfqpoint{14.705cm}{9.305cm}}
\pgfpathlineto{\pgfqpoint{14.482cm}{9.305cm}}
\pgfpathlineto{\pgfqpoint{14.482cm}{9.528cm}}
\pgfpathclose
\pgfusepath{stroke}
\pgfsetdash{}{0cm}
\pgfpathmoveto{\pgfqpoint{15.966cm}{9.501cm}}
\pgfpathlineto{\pgfqpoint{16.19cm}{9.501cm}}
\pgfpathlineto{\pgfqpoint{16.19cm}{9.278cm}}
\pgfpathlineto{\pgfqpoint{15.966cm}{9.278cm}}
\pgfpathlineto{\pgfqpoint{15.966cm}{9.501cm}}
\pgfpathclose
\pgfusepath{stroke}
\pgfsetdash{}{0cm}
\pgfpathmoveto{\pgfqpoint{17.454cm}{9.481cm}}
\pgfpathlineto{\pgfqpoint{17.677cm}{9.481cm}}
\pgfpathlineto{\pgfqpoint{17.677cm}{9.257cm}}
\pgfpathlineto{\pgfqpoint{17.454cm}{9.257cm}}
\pgfpathlineto{\pgfqpoint{17.454cm}{9.481cm}}
\pgfpathclose
\pgfusepath{stroke}
\pgfsetdash{}{0cm}
\pgfpathmoveto{\pgfqpoint{18.941cm}{9.463cm}}
\pgfpathlineto{\pgfqpoint{19.165cm}{9.463cm}}
\pgfpathlineto{\pgfqpoint{19.165cm}{9.24cm}}
\pgfpathlineto{\pgfqpoint{18.941cm}{9.24cm}}
\pgfpathlineto{\pgfqpoint{18.941cm}{9.463cm}}
\pgfpathclose
\pgfusepath{stroke}
\begin{pgfscope}
\pgfpathmoveto{\pgfqpoint{4.189cm}{20.449cm}}
\pgfpathlineto{\pgfqpoint{4.189cm}{8.611cm}}
\pgfpathlineto{\pgfqpoint{19.056cm}{8.611cm}}
\pgfpathlineto{\pgfqpoint{19.056cm}{20.449cm}}
\pgfpathclose
\pgfusepath{clip}
\pgfsetdash{{0.018cm}{0.141cm}{0.212cm}{0.141cm}}{0cm}
\pgfpathmoveto{\pgfqpoint{5.674cm}{12.147cm}}
\pgfpathlineto{\pgfqpoint{7.161cm}{10.304cm}}
\pgfpathlineto{\pgfqpoint{8.646cm}{9.907cm}}
\pgfpathlineto{\pgfqpoint{10.134cm}{9.743cm}}
\pgfpathlineto{\pgfqpoint{11.621cm}{9.651cm}}
\pgfpathlineto{\pgfqpoint{13.106cm}{9.598cm}}
\pgfpathlineto{\pgfqpoint{14.593cm}{9.56cm}}
\pgfpathlineto{\pgfqpoint{16.078cm}{9.534cm}}
\pgfpathlineto{\pgfqpoint{17.565cm}{9.513cm}}
\pgfpathlineto{\pgfqpoint{19.053cm}{9.496cm}}
\pgfusepath{stroke}
\end{pgfscope}
\pgfsetdash{}{0cm}
\pgfpathmoveto{\pgfqpoint{5.562cm}{12.259cm}}
\pgfpathlineto{\pgfqpoint{5.786cm}{12.259cm}}
\pgfpathlineto{\pgfqpoint{5.786cm}{12.036cm}}
\pgfpathlineto{\pgfqpoint{5.562cm}{12.036cm}}
\pgfpathlineto{\pgfqpoint{5.562cm}{12.259cm}}
\pgfpathclose
\pgfusepath{stroke}
\pgfsetdash{}{0cm}
\pgfpathmoveto{\pgfqpoint{7.05cm}{10.416cm}}
\pgfpathlineto{\pgfqpoint{7.273cm}{10.416cm}}
\pgfpathlineto{\pgfqpoint{7.273cm}{10.192cm}}
\pgfpathlineto{\pgfqpoint{7.05cm}{10.192cm}}
\pgfpathlineto{\pgfqpoint{7.05cm}{10.416cm}}
\pgfpathclose
\pgfusepath{stroke}
\pgfsetdash{}{0cm}
\pgfpathmoveto{\pgfqpoint{8.534cm}{10.019cm}}
\pgfpathlineto{\pgfqpoint{8.758cm}{10.019cm}}
\pgfpathlineto{\pgfqpoint{8.758cm}{9.795cm}}
\pgfpathlineto{\pgfqpoint{8.534cm}{9.795cm}}
\pgfpathlineto{\pgfqpoint{8.534cm}{10.019cm}}
\pgfpathclose
\pgfusepath{stroke}
\pgfsetdash{}{0cm}
\pgfpathmoveto{\pgfqpoint{10.022cm}{9.854cm}}
\pgfpathlineto{\pgfqpoint{10.245cm}{9.854cm}}
\pgfpathlineto{\pgfqpoint{10.245cm}{9.631cm}}
\pgfpathlineto{\pgfqpoint{10.022cm}{9.631cm}}
\pgfpathlineto{\pgfqpoint{10.022cm}{9.854cm}}
\pgfpathclose
\pgfusepath{stroke}
\pgfsetdash{}{0cm}
\pgfpathmoveto{\pgfqpoint{11.509cm}{9.763cm}}
\pgfpathlineto{\pgfqpoint{11.733cm}{9.763cm}}
\pgfpathlineto{\pgfqpoint{11.733cm}{9.54cm}}
\pgfpathlineto{\pgfqpoint{11.509cm}{9.54cm}}
\pgfpathlineto{\pgfqpoint{11.509cm}{9.763cm}}
\pgfpathclose
\pgfusepath{stroke}
\pgfsetdash{}{0cm}
\pgfpathmoveto{\pgfqpoint{12.994cm}{9.71cm}}
\pgfpathlineto{\pgfqpoint{13.217cm}{9.71cm}}
\pgfpathlineto{\pgfqpoint{13.217cm}{9.487cm}}
\pgfpathlineto{\pgfqpoint{12.994cm}{9.487cm}}
\pgfpathlineto{\pgfqpoint{12.994cm}{9.71cm}}
\pgfpathclose
\pgfusepath{stroke}
\pgfsetdash{}{0cm}
\pgfpathmoveto{\pgfqpoint{14.482cm}{9.672cm}}
\pgfpathlineto{\pgfqpoint{14.705cm}{9.672cm}}
\pgfpathlineto{\pgfqpoint{14.705cm}{9.449cm}}
\pgfpathlineto{\pgfqpoint{14.482cm}{9.449cm}}
\pgfpathlineto{\pgfqpoint{14.482cm}{9.672cm}}
\pgfpathclose
\pgfusepath{stroke}
\pgfsetdash{}{0cm}
\pgfpathmoveto{\pgfqpoint{15.966cm}{9.646cm}}
\pgfpathlineto{\pgfqpoint{16.19cm}{9.646cm}}
\pgfpathlineto{\pgfqpoint{16.19cm}{9.422cm}}
\pgfpathlineto{\pgfqpoint{15.966cm}{9.422cm}}
\pgfpathlineto{\pgfqpoint{15.966cm}{9.646cm}}
\pgfpathclose
\pgfusepath{stroke}
\pgfsetdash{}{0cm}
\pgfpathmoveto{\pgfqpoint{17.454cm}{9.625cm}}
\pgfpathlineto{\pgfqpoint{17.677cm}{9.625cm}}
\pgfpathlineto{\pgfqpoint{17.677cm}{9.402cm}}
\pgfpathlineto{\pgfqpoint{17.454cm}{9.402cm}}
\pgfpathlineto{\pgfqpoint{17.454cm}{9.625cm}}
\pgfpathclose
\pgfusepath{stroke}
\pgfsetdash{}{0cm}
\pgfpathmoveto{\pgfqpoint{18.941cm}{9.607cm}}
\pgfpathlineto{\pgfqpoint{19.165cm}{9.607cm}}
\pgfpathlineto{\pgfqpoint{19.165cm}{9.384cm}}
\pgfpathlineto{\pgfqpoint{18.941cm}{9.384cm}}
\pgfpathlineto{\pgfqpoint{18.941cm}{9.607cm}}
\pgfpathclose
\pgfusepath{stroke}
\begin{pgfscope}
\pgfpathmoveto{\pgfqpoint{4.189cm}{20.449cm}}
\pgfpathlineto{\pgfqpoint{4.189cm}{8.611cm}}
\pgfpathlineto{\pgfqpoint{19.056cm}{8.611cm}}
\pgfpathlineto{\pgfqpoint{19.056cm}{20.449cm}}
\pgfpathclose
\pgfusepath{clip}
\end{pgfscope}
\definecolor{eps2pgf_color}{gray}{0}\pgfsetstrokecolor{eps2pgf_color}\pgfsetfillcolor{eps2pgf_color}
\pgftext[x=11.639cm,y=6.968cm,rotate=0]{ \fontsize{40}{36.14}\selectfont{ {$\bmax$}}}
\pgftext[x=1.257cm,y=14.53cm,rotate=90]{ \fontsize{36}{36.14}\selectfont{ {Bound/$\bmax$}}}
\pgftext[x=4.139cm,y=8.487cm,rotate=0]{\fontsize{10.04}{12.04}\selectfont{ { }}}
\pgftext[x=19.006cm,y=20.326cm,rotate=0]{\fontsize{10.04}{12.04}\selectfont{ { }}}
\pgftext[x=13.914cm+.2cm-.2cm,y=19.488cm,rotate=0]{ \fontsize{40}{36.14}\selectfont{ {$\la = 0.9$, \texttt{minS}}}}
\begin{pgfscope}
\pgfpathmoveto{\pgfqpoint{8.655cm}{20.273cm}}
\pgfpathlineto{\pgfqpoint{8.655cm}{15.161cm}}
\pgfpathlineto{\pgfqpoint{18.879cm}{15.161cm}}
\pgfpathlineto{\pgfqpoint{18.879cm}{20.273cm}}
\pgfpathclose
\pgfusepath{clip}
\pgfsetdash{}{0cm}
\definecolor{eps2pgf_color}{rgb}{0,0,1}\pgfsetstrokecolor{eps2pgf_color}\pgfsetfillcolor{eps2pgf_color}
\pgfpathmoveto{\pgfqpoint{8.861cm}{19.594cm}}
\pgfpathlineto{\pgfqpoint{9.901cm}{19.594cm}}
\pgfusepath{stroke}
\begin{pgfscope}
\pgfpathmoveto{\pgfqpoint{8.987cm}{19.988cm}}
\pgfpathlineto{\pgfqpoint{8.987cm}{19.197cm}}
\pgfpathlineto{\pgfqpoint{9.778cm}{19.197cm}}
\pgfpathlineto{\pgfqpoint{9.778cm}{19.988cm}}
\pgfpathclose
\pgfusepath{clip}
\pgfsetdash{}{0cm}
\pgfpathmoveto{\pgfqpoint{9.24cm}{19.594cm}}
\pgfpathlineto{\pgfqpoint{9.522cm}{19.594cm}}
\pgfusepath{stroke}
\pgfsetdash{}{0cm}
\pgfpathmoveto{\pgfqpoint{9.381cm}{19.735cm}}
\pgfpathlineto{\pgfqpoint{9.381cm}{19.453cm}}
\pgfusepath{stroke}
\pgfsetdash{}{0cm}
\pgfpathmoveto{\pgfqpoint{9.284cm}{19.691cm}}
\pgfpathlineto{\pgfqpoint{9.478cm}{19.497cm}}
\pgfusepath{stroke}
\pgfsetdash{}{0cm}
\pgfpathmoveto{\pgfqpoint{9.478cm}{19.691cm}}
\pgfpathlineto{\pgfqpoint{9.284cm}{19.497cm}}
\pgfusepath{stroke}
\end{pgfscope}
\end{pgfscope}
\pgftext[x=14.091cm-.2cm,y=18.239cm,rotate=0]{ \fontsize{40}{36.14}\selectfont{ {$\la = 0.9$, \texttt{maxW}}}}
\begin{pgfscope}
\pgfpathmoveto{\pgfqpoint{8.655cm}{20.273cm}}
\pgfpathlineto{\pgfqpoint{8.655cm}{15.161cm}}
\pgfpathlineto{\pgfqpoint{18.879cm}{15.161cm}}
\pgfpathlineto{\pgfqpoint{18.879cm}{20.273cm}}
\pgfpathclose
\pgfusepath{clip}
\pgfsetdash{{0.018cm}{0.141cm}{0.212cm}{0.141cm}}{0cm}
\definecolor{eps2pgf_color}{rgb}{0,0,1}\pgfsetstrokecolor{eps2pgf_color}\pgfsetfillcolor{eps2pgf_color}
\pgfpathmoveto{\pgfqpoint{8.861cm}{18.344cm}}
\pgfpathlineto{\pgfqpoint{9.901cm}{18.344cm}}
\pgfusepath{stroke}
\begin{pgfscope}
\pgfpathmoveto{\pgfqpoint{8.987cm}{18.738cm}}
\pgfpathlineto{\pgfqpoint{8.987cm}{17.948cm}}
\pgfpathlineto{\pgfqpoint{9.778cm}{17.948cm}}
\pgfpathlineto{\pgfqpoint{9.778cm}{18.738cm}}
\pgfpathclose
\pgfusepath{clip}
\pgfsetdash{}{0cm}
\pgfpathmoveto{\pgfqpoint{9.24cm}{18.344cm}}
\pgfpathlineto{\pgfqpoint{9.522cm}{18.344cm}}
\pgfusepath{stroke}
\pgfsetdash{}{0cm}
\pgfpathmoveto{\pgfqpoint{9.381cm}{18.486cm}}
\pgfpathlineto{\pgfqpoint{9.381cm}{18.203cm}}
\pgfusepath{stroke}
\pgfsetdash{}{0cm}
\pgfpathmoveto{\pgfqpoint{9.284cm}{18.441cm}}
\pgfpathlineto{\pgfqpoint{9.478cm}{18.247cm}}
\pgfusepath{stroke}
\pgfsetdash{}{0cm}
\pgfpathmoveto{\pgfqpoint{9.478cm}{18.441cm}}
\pgfpathlineto{\pgfqpoint{9.284cm}{18.247cm}}
\pgfusepath{stroke}
\end{pgfscope}
\end{pgfscope}
\pgftext[x=14.209cm+.2cm-.2cm,y=16.989cm,rotate=0]{ \fontsize{40}{36.14}\selectfont{ {$\la = 0.95$, \texttt{minS}}}}
\begin{pgfscope}
\pgfpathmoveto{\pgfqpoint{8.655cm}{20.273cm}}
\pgfpathlineto{\pgfqpoint{8.655cm}{15.161cm}}
\pgfpathlineto{\pgfqpoint{18.879cm}{15.161cm}}
\pgfpathlineto{\pgfqpoint{18.879cm}{20.273cm}}
\pgfpathclose
\pgfusepath{clip}
\pgfsetdash{}{0cm}
\definecolor{eps2pgf_color}{rgb}{0,1,0}\pgfsetstrokecolor{eps2pgf_color}\pgfsetfillcolor{eps2pgf_color}
\pgfpathmoveto{\pgfqpoint{8.861cm}{17.095cm}}
\pgfpathlineto{\pgfqpoint{9.901cm}{17.095cm}}
\pgfusepath{stroke}
\begin{pgfscope}
\pgfpathmoveto{\pgfqpoint{8.987cm}{17.489cm}}
\pgfpathlineto{\pgfqpoint{8.987cm}{16.698cm}}
\pgfpathlineto{\pgfqpoint{9.778cm}{16.698cm}}
\pgfpathlineto{\pgfqpoint{9.778cm}{17.489cm}}
\pgfpathclose
\pgfusepath{clip}
\pgfsetdash{}{0cm}
\pgfpathmoveto{\pgfqpoint{9.269cm}{17.207cm}}
\pgfpathlineto{\pgfqpoint{9.493cm}{17.207cm}}
\pgfpathlineto{\pgfqpoint{9.493cm}{16.983cm}}
\pgfpathlineto{\pgfqpoint{9.269cm}{16.983cm}}
\pgfpathlineto{\pgfqpoint{9.269cm}{17.207cm}}
\pgfpathclose
\pgfusepath{stroke}
\end{pgfscope}
\end{pgfscope}
\pgftext[x=14.385cm-.2cm,y=15.74cm,rotate=0]{ \fontsize{40}{36.14}\selectfont{ {$\la = 0.95$, \texttt{maxW}}}}
\begin{pgfscope}
\pgfpathmoveto{\pgfqpoint{8.655cm}{20.273cm}}
\pgfpathlineto{\pgfqpoint{8.655cm}{15.161cm}}
\pgfpathlineto{\pgfqpoint{18.879cm}{15.161cm}}
\pgfpathlineto{\pgfqpoint{18.879cm}{20.273cm}}
\pgfpathclose
\pgfusepath{clip}
\pgfsetdash{{0.018cm}{0.141cm}{0.212cm}{0.141cm}}{0cm}
\definecolor{eps2pgf_color}{rgb}{0,1,0}\pgfsetstrokecolor{eps2pgf_color}\pgfsetfillcolor{eps2pgf_color}
\pgfpathmoveto{\pgfqpoint{8.861cm}{15.846cm}}
\pgfpathlineto{\pgfqpoint{9.901cm}{15.846cm}}
\pgfusepath{stroke}
\begin{pgfscope}
\pgfpathmoveto{\pgfqpoint{8.987cm}{16.24cm}}
\pgfpathlineto{\pgfqpoint{8.987cm}{15.449cm}}
\pgfpathlineto{\pgfqpoint{9.778cm}{15.449cm}}
\pgfpathlineto{\pgfqpoint{9.778cm}{16.24cm}}
\pgfpathclose
\pgfusepath{clip}
\pgfsetdash{}{0cm}
\pgfpathmoveto{\pgfqpoint{9.269cm}{15.957cm}}
\pgfpathlineto{\pgfqpoint{9.493cm}{15.957cm}}
\pgfpathlineto{\pgfqpoint{9.493cm}{15.734cm}}
\pgfpathlineto{\pgfqpoint{9.269cm}{15.734cm}}
\pgfpathlineto{\pgfqpoint{9.269cm}{15.957cm}}
\pgfpathclose
\pgfusepath{stroke}
\end{pgfscope}
\end{pgfscope}
\pgfsetdash{}{0cm}
\pgfsetlinewidth{0.176mm}
\definecolor{eps2pgf_color}{rgb}{0,1,0}\pgfsetstrokecolor{eps2pgf_color}\pgfsetfillcolor{eps2pgf_color}
\pgfusepath{stroke}
\end{pgfscope}
\end{pgfscope}
\end{pgfpicture}

%% file: numericalExample-conf.tex
\section{Numerical Experiments}
\addtolength{\abovedisplayskip}{-2mm}
\addtolength{\belowdisplayskip}{-2mm}
\subsection{Balancing with IID Disturbance}\label{sec:baiid}
\begin{figure}[htbp]
\centering\hspace{-1.em}%
\input{./fig/I1.tex}\hspace{-.2cm}%
\input{./fig/I2.tex}
\caption{Algorithm performance with temporally homogeneous cost and ideal storage (left panel),  and temporally heterogeneous cost and non-ideal storage (right panel). The average costs represent average imbalance p.u. (\emph{cf.} cost \eqref{costba}) in the left panel and average penalty \eqref{eq:inhomo:cost} in the right panel, respectively.
%Note that, different from the setup in Remark~\ref{remark:subo:la1}, we scale $\umin$, $\umax$ together with $\bmax$ in this (and the following) numerical examples.}
}\label{fig:3}
\end{figure}
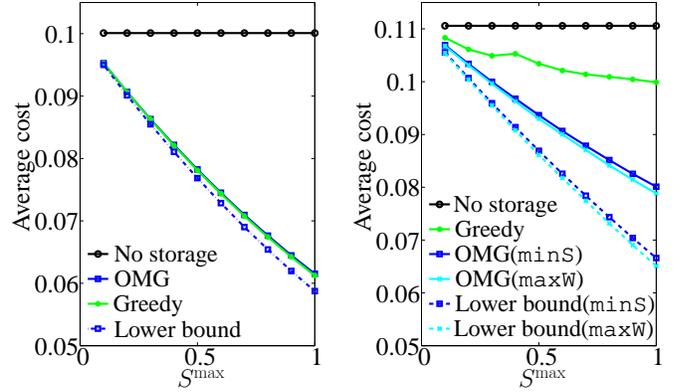

We first test our algorithm in a simple setting where the analytical solution for the optimal control policy is available, so that the algorithm performance can be compared against the true optimal costs. We consider the problem of using energy storage to minimize the energy imbalance as studied in \cite{SuEGTPS}, where it is shown that greedy storage operation is optimal if $\la = 1$ and if the following cost is considered
\begin{equation}\label{costba}
\g _t  = |\d_t - (1/\muC)\upos_t + \muD\uneg_t|.
\end{equation}
As in \cite{SuEGTPS}, we specify storage parameters in per unit, and  $\bmin = 0$. Let  $\muC = \muD = 1$ so that the parameterization of storage operation here is equivalent to that of \cite{SuEGTPS}. We assume each time period represents an hour, and $-\umin = \umax = (1/10) \bmax$.
In order to evaluate the performance, we simulate the $\d_t$ process by drawing i.i.d. samples from zero-mean Laplace distribution with standard deviation $\sigma_\d = 0.149$ per unit (p.u.) obtained from NREL data \cite{SuEGTPS}. The time horizon for the simulation is chosen to be $T=1000$. Figure~\ref{fig:3} (left panel) depicts the performance of OMG and the optimal cost $J^\star$ obtained from the greedy policy,  where it is shown that the costs of OMG are close to the optimal costs, and are better than what the (worst-case) sub-optimality bound predicts. \footnote{By an abuse of notation, in this section, we use $J^\star$ to denote the results from simulation, which are estimates of the true expectations.} 

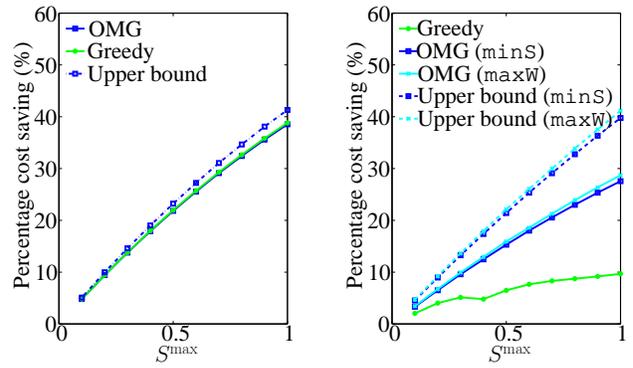
\begin{figure}[htbp]
\centering\hspace{-1.em}%
\input{./fig/I3.tex}
\input{./fig/I4.tex}
\caption{Percentage cost savings with temporally homogeneous cost and ideal storage (left panel),  and temporally heterogeneous cost and non-ideal storage (right panel).}
\label{fig:3pcs}
\end{figure}

A slight modification of the cost function would render a problem which does not have an analytical solution. Consider the setting where only unsatisfied demand is penalized with a higher penalty during the day ($7$ am to $7$ pm):
\begin{equation}\label{eq:inhomo:cost}
\!\!\g _t\!\! = \!\!
\begin{cases}
3\neg{\d_t\! -\! (\upos_t/\muC) + \muD\uneg_t}\!\!,\!\! &t\in \Tday,  \\
\neg{\d_t\! -\! (\upos_t/\muC) + \muD\uneg_t}\!\!, \!\!& \mbox{otherwise},
\end{cases}
\end{equation}
where $\Tday$ is the set of stages that corresponds to time points in the range of $7$ am to $7$ pm.
We run the same set of tests above, with the modification that now $\muC = \muD = 0.85$, and $\lambda = 0.9975$ (which corresponds to the NaS battery in Example~\ref{eg:soe} operated in 5 minute intervals). Note that the greedy policy is only a sub-optimal heuristic for this case. Figure~\ref{fig:3} (right panel) shows OMG performs significantly better than the greedy algorithm. The costs of our algorithm together with the lower bounds give narrow envelopes for the optimal average cost $J^\star$ in this setting, which can be used to evaluate the performance of other sub-optimal algorithms numerically. We have also shown the performance and lower bounds of the OMG algorithm with \texttt{minS} and \texttt{maxW} parameter settings. In this example, \texttt{minS} gives better lower bounds whereas \texttt{maxW} leads to lower costs. Figure~\ref{fig:3pcs} translates the cost numbers into the percentage cost savings of operating the storage (with various approaches) comparing to the no storage scenario. 

In both experiments, we also plot the costs of certainty equivalent/predictive storage control, whose solution can be shown to be $\u_t = 0 $ for all $t$. Consequently, the costs of such operation rule are the same as the system costs when there is no storage.

\subsection{Simulation with Real Price and Net Demand Data}
\begin{figure}[htbp]
\centering
\input{I4_p}
\input{I4_d}
\caption{Bar plots for hourly locational marginal price and forecast error data for a wind farm in PJM interconnection in January 2004. Power units have been converted to energy units.
}\label{fig:data}
\end{figure}
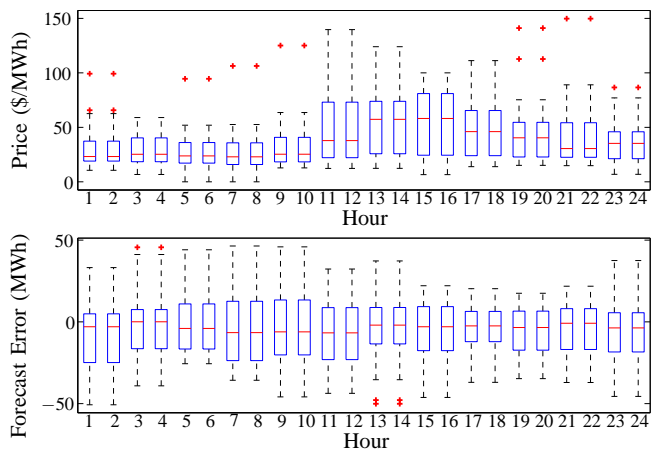

We consider a case where a storage is co-located with a wind farm. The wind farm operates the storage (i) to reduce wind power spillage caused by forecast errors, and (ii) to arbitrage price differences across different time periods.
The setting here is similar to Example~\ref{eg:co}, such that both the price and the net demand  are random. The stage-wise cost function is 
\[
 \g(t) = \p_t(\d_t - (1/\muC)\upos_t + \muD\uneg_t),
\]
where the $\{\p_t : t\ge 1\}$ and $\{\d_t: t\ge 1\}$ sequences are obtained from the LMP data from PJM interconnection and forecast error data from the NREL dataset \cite{NREL2010}  (Figure~\ref{fig:data}).  
We consider an ideal storage with capacity $\bmax = 5 \sigma_d$ and $\umax = - \umin = (1/20) \bmax$, where $\sigma_d = 20.1 $MWh is the empirical standard deviation of the wind power generation forecast error. The storage is operated every hour and the simulation is run for a month, \ie, $T = 360$. The average per stage cost without energy storage is $224.65$ \$, whereas the average per stage cost of greedy storage operation, OMG, and the offline clairvoyant optimal operation are $99.7\%$, $88.8\%$, and $75.7\%$ of the no storage cost, respectively. Here the offline clairvoyant optimal operation is calculated by solving a deterministic optimization assuming full knowledge of future $\d_t$ and $\p_t$ sequence, and is in general a loose lower bound of the optimal costs. The stochastic lower bound assuming i.i.d. disturbance suggests the minimal achievable per stage cost would be $83.2\%$ of the no storage cost. 

%I4_costarbi_smaxR5_mu1_la1_dt60_T360_sb1_caseID340_GR0.99739_LRS0.88859_LRW0.88859_OR0.7569.mat

\section{Conclusion}
In this paper, we formulate the problem of operating a generalized storage under uncertainty as a stochastic control problem. A very simple algorithm, termed online modified greedy algorithm, is proposed and analyzed. The sub-optimality of the algorithm is proved to be bounded by a function of the system parameters. The bound is efficiently computable and can be used to gauge the performance of the algorithm as well as to estimate the optimal cost. Numerical simulations are conducted to illustrate the use of the algorithm and to validate its effectiveness.

The following future directions are of interests for generalizing/improving the proposed method. 
i) The proposed algorithm does not require the knowledge of the full probability distributions of disturbances. While this may be advantageous when such information is not available, in case that it is available or partially available, extensions of the algorithm incorporating such information may generate a better storage control policy. 
ii) Our approach is easily generalizable to settings with multiple \emph{same-stage} variables, \ie, the controllable inflow can be a vector that lies in a given convex set. However, applications that also involves \emph{look-ahead} variables, such as those arising in the contexts that the storage is operated with a wind farm participating in the forward markets or that the storage itself participates in the forward markets, cannot directly be cast into our framework. Generalizing the algorithm for those contexts by e.g. incorporating ideas from \cite[Section 4.9.2]{NeelyBook} is an important future direction. iii) The current algorithm optimizes a single storage. Extending the algorithm to a setting with multiple storages that are connected via a power network will enable the algorithm to be applied to settings such as storage control in micro-grids. One possible way for such an extension is reported in \cite{7101867}.

%% file: fig/I1.tex
\scalebox{0.265}{\scalefont{2} \input{./fig/single_storage_I1.pgf}}

%% file: fig/single_storage_I1.pgf
% Created by Eps2pgf 0.7.0 (build on 2008-08-24) on Mon Apr 06 20:47:11 PDT 2015
\begin{pgfpicture}
\pgfpathmoveto{\pgfqpoint{2.293cm}{3.387cm}}
\pgfpathlineto{\pgfqpoint{19.226cm}{3.387cm}}
\pgfpathlineto{\pgfqpoint{19.226cm}{24.553cm}}
\pgfpathlineto{\pgfqpoint{2.293cm}{24.553cm}}
\pgfpathclose
\pgfusepath{clip}
\begin{pgfscope}
\begin{pgfscope}
\pgfpathmoveto{\pgfqpoint{2.293cm}{24.553cm}}
\pgfpathlineto{\pgfqpoint{19.229cm}{24.553cm}}
\pgfpathlineto{\pgfqpoint{19.229cm}{3.413cm}}
\pgfpathlineto{\pgfqpoint{2.293cm}{3.413cm}}
\pgfpathclose
\pgfusepath{clip}
\begin{pgfscope}
\definecolor{eps2pgf_color}{gray}{1}\pgfsetstrokecolor{eps2pgf_color}\pgfsetfillcolor{eps2pgf_color}
\pgfpathmoveto{\pgfqpoint{2.293cm}{24.553cm}}
\pgfpathlineto{\pgfqpoint{19.232cm}{24.553cm}}
\pgfpathlineto{\pgfqpoint{19.232cm}{3.41cm}}
\pgfpathlineto{\pgfqpoint{2.293cm}{3.41cm}}
\pgfpathclose
\pgfusepath{fill}
\end{pgfscope}
\definecolor{eps2pgf_color}{gray}{1}\pgfsetstrokecolor{eps2pgf_color}\pgfsetfillcolor{eps2pgf_color}
\pgfpathmoveto{\pgfqpoint{5.847cm}{5.739cm}}
\pgfpathlineto{\pgfqpoint{5.847cm}{22.969cm}}
\pgfpathlineto{\pgfqpoint{17.621cm}{22.969cm}}
\pgfpathlineto{\pgfqpoint{17.621cm}{5.739cm}}
\pgfpathclose
\pgfseteorule\pgfusepath{fill}\pgfsetnonzerorule
\pgfsetdash{}{0cm}
\pgfsetlinewidth{0.176mm}
\pgfsetroundjoin
\pgfpathmoveto{\pgfqpoint{5.847cm}{5.739cm}}
\pgfpathlineto{\pgfqpoint{5.847cm}{22.969cm}}
\pgfpathlineto{\pgfqpoint{17.621cm}{22.969cm}}
\pgfpathlineto{\pgfqpoint{17.621cm}{5.739cm}}
\pgfpathlineto{\pgfqpoint{5.847cm}{5.739cm}}
\pgfusepath{stroke}
\pgfsetdash{}{0cm}
\definecolor{eps2pgf_color}{gray}{0}\pgfsetstrokecolor{eps2pgf_color}\pgfsetfillcolor{eps2pgf_color}
\pgfpathmoveto{\pgfqpoint{5.847cm}{5.739cm}}
\pgfpathlineto{\pgfqpoint{17.621cm}{5.739cm}}
\pgfusepath{stroke}
\pgfsetdash{}{0cm}
\pgfpathmoveto{\pgfqpoint{5.847cm}{22.969cm}}
\pgfpathlineto{\pgfqpoint{17.621cm}{22.969cm}}
\pgfusepath{stroke}
\pgfsetdash{}{0cm}
\pgfpathmoveto{\pgfqpoint{5.847cm}{5.739cm}}
\pgfpathlineto{\pgfqpoint{5.847cm}{22.969cm}}
\pgfusepath{stroke}
\pgfsetdash{}{0cm}
\pgfpathmoveto{\pgfqpoint{17.621cm}{5.739cm}}
\pgfpathlineto{\pgfqpoint{17.621cm}{22.969cm}}
\pgfusepath{stroke}
\pgfsetdash{}{0cm}
\pgfpathmoveto{\pgfqpoint{5.847cm}{5.739cm}}
\pgfpathlineto{\pgfqpoint{17.621cm}{5.739cm}}
\pgfusepath{stroke}
\pgfsetdash{}{0cm}
\pgfpathmoveto{\pgfqpoint{5.847cm}{5.739cm}}
\pgfpathlineto{\pgfqpoint{5.847cm}{22.969cm}}
\pgfusepath{stroke}
\pgfsetdash{}{0cm}
\pgfpathmoveto{\pgfqpoint{5.847cm}{5.739cm}}
\pgfpathlineto{\pgfqpoint{5.847cm}{5.912cm}}
\pgfusepath{stroke}
\pgfsetdash{}{0cm}
\pgfpathmoveto{\pgfqpoint{5.847cm}{22.969cm}}
\pgfpathlineto{\pgfqpoint{5.847cm}{22.798cm}}
\pgfusepath{stroke}
\pgftext[x=5.848cm,y=5.016cm+.2cm,rotate=0]{\fontsize{36}{36.14}\selectfont{ {0}}}
\pgfsetdash{}{0cm}
\pgfpathmoveto{\pgfqpoint{11.733cm}{5.739cm}}
\pgfpathlineto{\pgfqpoint{11.733cm}{5.912cm}}
\pgfusepath{stroke}
\pgfsetdash{}{0cm}
\pgfpathmoveto{\pgfqpoint{11.733cm}{22.969cm}}
\pgfpathlineto{\pgfqpoint{11.733cm}{22.798cm}}
\pgfusepath{stroke}
\pgftext[x=11.731cm,y=5.016cm+.2cm,rotate=0]{\fontsize{36}{36.14}\selectfont{ {0.5}}}
\pgfsetdash{}{0cm}
\pgfpathmoveto{\pgfqpoint{17.621cm}{5.739cm}}
\pgfpathlineto{\pgfqpoint{17.621cm}{5.912cm}}
\pgfusepath{stroke}
\pgfsetdash{}{0cm}
\pgfpathmoveto{\pgfqpoint{17.621cm}{22.969cm}}
\pgfpathlineto{\pgfqpoint{17.621cm}{22.798cm}}
\pgfusepath{stroke}
\pgftext[x=17.571cm,y=5.026cm+.2cm,rotate=0]{\fontsize{36}{36.14}\selectfont{ {1}}}
\pgfsetdash{}{0cm}
\pgfpathmoveto{\pgfqpoint{5.847cm}{5.739cm}}
\pgfpathlineto{\pgfqpoint{6.018cm}{5.739cm}}
\pgfusepath{stroke}
\pgfsetdash{}{0cm}
\pgfpathmoveto{\pgfqpoint{17.621cm}{5.739cm}}
\pgfpathlineto{\pgfqpoint{17.448cm}{5.739cm}}
\pgfusepath{stroke}
\pgftext[x=4.714cm-.2cm,y=5.707cm,rotate=0]{\fontsize{36}{36.14}\selectfont{ {0.05}}}
\pgfsetdash{}{0cm}
\pgfpathmoveto{\pgfqpoint{5.847cm}{8.872cm}}
\pgfpathlineto{\pgfqpoint{6.018cm}{8.872cm}}
\pgfusepath{stroke}
\pgfsetdash{}{0cm}
\pgfpathmoveto{\pgfqpoint{17.621cm}{8.872cm}}
\pgfpathlineto{\pgfqpoint{17.448cm}{8.872cm}}
\pgfusepath{stroke}
\pgftext[x=4.716cm-.2cm,y=8.84cm,rotate=0]{\fontsize{36}{36.14}\selectfont{ {0.06}}}
\pgfsetdash{}{0cm}
\pgfpathmoveto{\pgfqpoint{5.847cm}{12.006cm}}
\pgfpathlineto{\pgfqpoint{6.018cm}{12.006cm}}
\pgfusepath{stroke}
\pgfsetdash{}{0cm}
\pgfpathmoveto{\pgfqpoint{17.621cm}{12.006cm}}
\pgfpathlineto{\pgfqpoint{17.448cm}{12.006cm}}
\pgfusepath{stroke}
\pgftext[x=4.718cm-.2cm,y=11.974cm,rotate=0]{\fontsize{36}{36.14}\selectfont{ {0.07}}}
\pgfsetdash{}{0cm}
\pgfpathmoveto{\pgfqpoint{5.847cm}{15.137cm}}
\pgfpathlineto{\pgfqpoint{6.018cm}{15.137cm}}
\pgfusepath{stroke}
\pgfsetdash{}{0cm}
\pgfpathmoveto{\pgfqpoint{17.621cm}{15.137cm}}
\pgfpathlineto{\pgfqpoint{17.448cm}{15.137cm}}
\pgfusepath{stroke}
\pgftext[x=4.715cm-.2cm,y=15.105cm,rotate=0]{\fontsize{36}{36.14}\selectfont{ {0.08}}}
\pgfsetdash{}{0cm}
\pgfpathmoveto{\pgfqpoint{5.847cm}{18.271cm}}
\pgfpathlineto{\pgfqpoint{6.018cm}{18.271cm}}
\pgfusepath{stroke}
\pgfsetdash{}{0cm}
\pgfpathmoveto{\pgfqpoint{17.621cm}{18.271cm}}
\pgfpathlineto{\pgfqpoint{17.448cm}{18.271cm}}
\pgfusepath{stroke}
\pgftext[x=4.714cm-.2cm,y=18.239cm,rotate=0]{\fontsize{36}{36.14}\selectfont{ {0.09}}}
\pgfsetdash{}{0cm}
\pgfpathmoveto{\pgfqpoint{5.847cm}{21.405cm}}
\pgfpathlineto{\pgfqpoint{6.018cm}{21.405cm}}
\pgfusepath{stroke}
\pgfsetdash{}{0cm}
\pgfpathmoveto{\pgfqpoint{17.621cm}{21.405cm}}
\pgfpathlineto{\pgfqpoint{17.448cm}{21.405cm}}
\pgfusepath{stroke}
\pgftext[x=4.925cm-.2cm,y=21.373cm,rotate=0]{\fontsize{36}{36.14}\selectfont{ {0.1}}}
\pgfsetdash{}{0cm}
\pgfpathmoveto{\pgfqpoint{5.847cm}{5.739cm}}
\pgfpathlineto{\pgfqpoint{17.621cm}{5.739cm}}
\pgfusepath{stroke}
\pgfsetdash{}{0cm}
\pgfpathmoveto{\pgfqpoint{5.847cm}{22.969cm}}
\pgfpathlineto{\pgfqpoint{17.621cm}{22.969cm}}
\pgfusepath{stroke}
\pgfsetdash{}{0cm}
\pgfpathmoveto{\pgfqpoint{5.847cm}{5.739cm}}
\pgfpathlineto{\pgfqpoint{5.847cm}{22.969cm}}
\pgfusepath{stroke}
\pgfsetdash{}{0cm}
\pgfpathmoveto{\pgfqpoint{17.621cm}{5.739cm}}
\pgfpathlineto{\pgfqpoint{17.621cm}{22.969cm}}
\pgfusepath{stroke}
\begin{pgfscope}
\pgfpathmoveto{\pgfqpoint{5.847cm}{22.969cm}}
\pgfpathlineto{\pgfqpoint{17.624cm}{22.969cm}}
\pgfpathlineto{\pgfqpoint{17.624cm}{5.736cm}}
\pgfpathlineto{\pgfqpoint{5.847cm}{5.736cm}}
\pgfpathclose
\pgfusepath{clip}
\pgfsetdash{}{0cm}
\pgfsetlinewidth{1.058mm}
\pgfpathmoveto{\pgfqpoint{7.023cm}{21.434cm}}
\pgfpathlineto{\pgfqpoint{8.202cm}{21.434cm}}
\pgfpathlineto{\pgfqpoint{9.378cm}{21.434cm}}
\pgfpathlineto{\pgfqpoint{10.557cm}{21.434cm}}
\pgfpathlineto{\pgfqpoint{11.733cm}{21.434cm}}
\pgfpathlineto{\pgfqpoint{12.912cm}{21.434cm}}
\pgfpathlineto{\pgfqpoint{14.088cm}{21.434cm}}
\pgfpathlineto{\pgfqpoint{15.266cm}{21.434cm}}
\pgfpathlineto{\pgfqpoint{16.442cm}{21.434cm}}
\pgfpathlineto{\pgfqpoint{17.621cm}{21.434cm}}
\pgfusepath{stroke}
\end{pgfscope}
\pgfsetdash{}{0cm}
\pgfsetlinewidth{1.058mm}
\pgfpathmoveto{\pgfqpoint{7.164cm}{21.434cm}}
\pgfpathcurveto{\pgfqpoint{7.164cm}{21.356cm}}{\pgfqpoint{7.101cm}{21.293cm}}{\pgfqpoint{7.023cm}{21.293cm}}
\pgfpathcurveto{\pgfqpoint{6.945cm}{21.293cm}}{\pgfqpoint{6.882cm}{21.356cm}}{\pgfqpoint{6.882cm}{21.434cm}}
\pgfpathcurveto{\pgfqpoint{6.882cm}{21.512cm}}{\pgfqpoint{6.945cm}{21.575cm}}{\pgfqpoint{7.023cm}{21.575cm}}
\pgfpathcurveto{\pgfqpoint{7.101cm}{21.575cm}}{\pgfqpoint{7.164cm}{21.512cm}}{\pgfqpoint{7.164cm}{21.434cm}}
\pgfusepath{stroke}
\pgfsetdash{}{0cm}
\pgfpathmoveto{\pgfqpoint{8.343cm}{21.434cm}}
\pgfpathcurveto{\pgfqpoint{8.343cm}{21.356cm}}{\pgfqpoint{8.28cm}{21.293cm}}{\pgfqpoint{8.202cm}{21.293cm}}
\pgfpathcurveto{\pgfqpoint{8.124cm}{21.293cm}}{\pgfqpoint{8.061cm}{21.356cm}}{\pgfqpoint{8.061cm}{21.434cm}}
\pgfpathcurveto{\pgfqpoint{8.061cm}{21.512cm}}{\pgfqpoint{8.124cm}{21.575cm}}{\pgfqpoint{8.202cm}{21.575cm}}
\pgfpathcurveto{\pgfqpoint{8.28cm}{21.575cm}}{\pgfqpoint{8.343cm}{21.512cm}}{\pgfqpoint{8.343cm}{21.434cm}}
\pgfusepath{stroke}
\pgfsetdash{}{0cm}
\pgfpathmoveto{\pgfqpoint{9.519cm}{21.434cm}}
\pgfpathcurveto{\pgfqpoint{9.519cm}{21.356cm}}{\pgfqpoint{9.456cm}{21.293cm}}{\pgfqpoint{9.378cm}{21.293cm}}
\pgfpathcurveto{\pgfqpoint{9.3cm}{21.293cm}}{\pgfqpoint{9.237cm}{21.356cm}}{\pgfqpoint{9.237cm}{21.434cm}}
\pgfpathcurveto{\pgfqpoint{9.237cm}{21.512cm}}{\pgfqpoint{9.3cm}{21.575cm}}{\pgfqpoint{9.378cm}{21.575cm}}
\pgfpathcurveto{\pgfqpoint{9.456cm}{21.575cm}}{\pgfqpoint{9.519cm}{21.512cm}}{\pgfqpoint{9.519cm}{21.434cm}}
\pgfusepath{stroke}
\pgfsetdash{}{0cm}
\pgfpathmoveto{\pgfqpoint{10.698cm}{21.434cm}}
\pgfpathcurveto{\pgfqpoint{10.698cm}{21.356cm}}{\pgfqpoint{10.635cm}{21.293cm}}{\pgfqpoint{10.557cm}{21.293cm}}
\pgfpathcurveto{\pgfqpoint{10.479cm}{21.293cm}}{\pgfqpoint{10.416cm}{21.356cm}}{\pgfqpoint{10.416cm}{21.434cm}}
\pgfpathcurveto{\pgfqpoint{10.416cm}{21.512cm}}{\pgfqpoint{10.479cm}{21.575cm}}{\pgfqpoint{10.557cm}{21.575cm}}
\pgfpathcurveto{\pgfqpoint{10.635cm}{21.575cm}}{\pgfqpoint{10.698cm}{21.512cm}}{\pgfqpoint{10.698cm}{21.434cm}}
\pgfusepath{stroke}
\pgfsetdash{}{0cm}
\pgfpathmoveto{\pgfqpoint{11.874cm}{21.434cm}}
\pgfpathcurveto{\pgfqpoint{11.874cm}{21.356cm}}{\pgfqpoint{11.811cm}{21.293cm}}{\pgfqpoint{11.733cm}{21.293cm}}
\pgfpathcurveto{\pgfqpoint{11.655cm}{21.293cm}}{\pgfqpoint{11.592cm}{21.356cm}}{\pgfqpoint{11.592cm}{21.434cm}}
\pgfpathcurveto{\pgfqpoint{11.592cm}{21.512cm}}{\pgfqpoint{11.655cm}{21.575cm}}{\pgfqpoint{11.733cm}{21.575cm}}
\pgfpathcurveto{\pgfqpoint{11.811cm}{21.575cm}}{\pgfqpoint{11.874cm}{21.512cm}}{\pgfqpoint{11.874cm}{21.434cm}}
\pgfusepath{stroke}
\pgfsetdash{}{0cm}
\pgfpathmoveto{\pgfqpoint{13.053cm}{21.434cm}}
\pgfpathcurveto{\pgfqpoint{13.053cm}{21.356cm}}{\pgfqpoint{12.99cm}{21.293cm}}{\pgfqpoint{12.912cm}{21.293cm}}
\pgfpathcurveto{\pgfqpoint{12.834cm}{21.293cm}}{\pgfqpoint{12.771cm}{21.356cm}}{\pgfqpoint{12.771cm}{21.434cm}}
\pgfpathcurveto{\pgfqpoint{12.771cm}{21.512cm}}{\pgfqpoint{12.834cm}{21.575cm}}{\pgfqpoint{12.912cm}{21.575cm}}
\pgfpathcurveto{\pgfqpoint{12.99cm}{21.575cm}}{\pgfqpoint{13.053cm}{21.512cm}}{\pgfqpoint{13.053cm}{21.434cm}}
\pgfusepath{stroke}
\pgfsetdash{}{0cm}
\pgfpathmoveto{\pgfqpoint{14.229cm}{21.434cm}}
\pgfpathcurveto{\pgfqpoint{14.229cm}{21.356cm}}{\pgfqpoint{14.166cm}{21.293cm}}{\pgfqpoint{14.088cm}{21.293cm}}
\pgfpathcurveto{\pgfqpoint{14.01cm}{21.293cm}}{\pgfqpoint{13.946cm}{21.356cm}}{\pgfqpoint{13.946cm}{21.434cm}}
\pgfpathcurveto{\pgfqpoint{13.946cm}{21.512cm}}{\pgfqpoint{14.01cm}{21.575cm}}{\pgfqpoint{14.088cm}{21.575cm}}
\pgfpathcurveto{\pgfqpoint{14.166cm}{21.575cm}}{\pgfqpoint{14.229cm}{21.512cm}}{\pgfqpoint{14.229cm}{21.434cm}}
\pgfusepath{stroke}
\pgfsetdash{}{0cm}
\pgfpathmoveto{\pgfqpoint{15.408cm}{21.434cm}}
\pgfpathcurveto{\pgfqpoint{15.408cm}{21.356cm}}{\pgfqpoint{15.344cm}{21.293cm}}{\pgfqpoint{15.266cm}{21.293cm}}
\pgfpathcurveto{\pgfqpoint{15.189cm}{21.293cm}}{\pgfqpoint{15.125cm}{21.356cm}}{\pgfqpoint{15.125cm}{21.434cm}}
\pgfpathcurveto{\pgfqpoint{15.125cm}{21.512cm}}{\pgfqpoint{15.189cm}{21.575cm}}{\pgfqpoint{15.266cm}{21.575cm}}
\pgfpathcurveto{\pgfqpoint{15.344cm}{21.575cm}}{\pgfqpoint{15.408cm}{21.512cm}}{\pgfqpoint{15.408cm}{21.434cm}}
\pgfusepath{stroke}
\pgfsetdash{}{0cm}
\pgfpathmoveto{\pgfqpoint{16.583cm}{21.434cm}}
\pgfpathcurveto{\pgfqpoint{16.583cm}{21.356cm}}{\pgfqpoint{16.52cm}{21.293cm}}{\pgfqpoint{16.442cm}{21.293cm}}
\pgfpathcurveto{\pgfqpoint{16.364cm}{21.293cm}}{\pgfqpoint{16.301cm}{21.356cm}}{\pgfqpoint{16.301cm}{21.434cm}}
\pgfpathcurveto{\pgfqpoint{16.301cm}{21.512cm}}{\pgfqpoint{16.364cm}{21.575cm}}{\pgfqpoint{16.442cm}{21.575cm}}
\pgfpathcurveto{\pgfqpoint{16.52cm}{21.575cm}}{\pgfqpoint{16.583cm}{21.512cm}}{\pgfqpoint{16.583cm}{21.434cm}}
\pgfusepath{stroke}
\pgfsetdash{}{0cm}
\pgfpathmoveto{\pgfqpoint{17.762cm}{21.434cm}}
\pgfpathcurveto{\pgfqpoint{17.762cm}{21.356cm}}{\pgfqpoint{17.699cm}{21.293cm}}{\pgfqpoint{17.621cm}{21.293cm}}
\pgfpathcurveto{\pgfqpoint{17.543cm}{21.293cm}}{\pgfqpoint{17.48cm}{21.356cm}}{\pgfqpoint{17.48cm}{21.434cm}}
\pgfpathcurveto{\pgfqpoint{17.48cm}{21.512cm}}{\pgfqpoint{17.543cm}{21.575cm}}{\pgfqpoint{17.621cm}{21.575cm}}
\pgfpathcurveto{\pgfqpoint{17.699cm}{21.575cm}}{\pgfqpoint{17.762cm}{21.512cm}}{\pgfqpoint{17.762cm}{21.434cm}}
\pgfusepath{stroke}
\begin{pgfscope}
\pgfpathmoveto{\pgfqpoint{5.847cm}{22.969cm}}
\pgfpathlineto{\pgfqpoint{17.624cm}{22.969cm}}
\pgfpathlineto{\pgfqpoint{17.624cm}{5.736cm}}
\pgfpathlineto{\pgfqpoint{5.847cm}{5.736cm}}
\pgfpathclose
\pgfusepath{clip}
\pgfsetdash{}{0cm}
\definecolor{eps2pgf_color}{rgb}{0,0,1}\pgfsetstrokecolor{eps2pgf_color}\pgfsetfillcolor{eps2pgf_color}
\pgfpathmoveto{\pgfqpoint{7.023cm}{19.923cm}}
\pgfpathlineto{\pgfqpoint{8.202cm}{18.48cm}}
\pgfpathlineto{\pgfqpoint{9.378cm}{17.116cm}}
\pgfpathlineto{\pgfqpoint{10.557cm}{15.819cm}}
\pgfpathlineto{\pgfqpoint{11.733cm}{14.587cm}}
\pgfpathlineto{\pgfqpoint{12.912cm}{13.417cm}}
\pgfpathlineto{\pgfqpoint{14.088cm}{12.303cm}}
\pgfpathlineto{\pgfqpoint{15.266cm}{11.265cm}}
\pgfpathlineto{\pgfqpoint{16.442cm}{10.272cm}}
\pgfpathlineto{\pgfqpoint{17.621cm}{9.349cm}}
\pgfusepath{stroke}
\end{pgfscope}
\pgfsetdash{}{0cm}
\pgfsetmiterjoin
\definecolor{eps2pgf_color}{rgb}{0,0,1}\pgfsetstrokecolor{eps2pgf_color}\pgfsetfillcolor{eps2pgf_color}
\pgfpathmoveto{\pgfqpoint{6.912cm}{20.035cm}}
\pgfpathlineto{\pgfqpoint{7.135cm}{20.035cm}}
\pgfpathlineto{\pgfqpoint{7.135cm}{19.811cm}}
\pgfpathlineto{\pgfqpoint{6.912cm}{19.811cm}}
\pgfpathlineto{\pgfqpoint{6.912cm}{20.035cm}}
\pgfpathclose
\pgfusepath{stroke}
\pgfsetdash{}{0cm}
\pgfpathmoveto{\pgfqpoint{8.09cm}{18.591cm}}
\pgfpathlineto{\pgfqpoint{8.314cm}{18.591cm}}
\pgfpathlineto{\pgfqpoint{8.314cm}{18.368cm}}
\pgfpathlineto{\pgfqpoint{8.09cm}{18.368cm}}
\pgfpathlineto{\pgfqpoint{8.09cm}{18.591cm}}
\pgfpathclose
\pgfusepath{stroke}
\pgfsetdash{}{0cm}
\pgfpathmoveto{\pgfqpoint{9.266cm}{17.227cm}}
\pgfpathlineto{\pgfqpoint{9.49cm}{17.227cm}}
\pgfpathlineto{\pgfqpoint{9.49cm}{17.004cm}}
\pgfpathlineto{\pgfqpoint{9.266cm}{17.004cm}}
\pgfpathlineto{\pgfqpoint{9.266cm}{17.227cm}}
\pgfpathclose
\pgfusepath{stroke}
\pgfsetdash{}{0cm}
\pgfpathmoveto{\pgfqpoint{10.445cm}{15.931cm}}
\pgfpathlineto{\pgfqpoint{10.669cm}{15.931cm}}
\pgfpathlineto{\pgfqpoint{10.669cm}{15.707cm}}
\pgfpathlineto{\pgfqpoint{10.445cm}{15.707cm}}
\pgfpathlineto{\pgfqpoint{10.445cm}{15.931cm}}
\pgfpathclose
\pgfusepath{stroke}
\pgfsetdash{}{0cm}
\pgfpathmoveto{\pgfqpoint{11.621cm}{14.699cm}}
\pgfpathlineto{\pgfqpoint{11.845cm}{14.699cm}}
\pgfpathlineto{\pgfqpoint{11.845cm}{14.476cm}}
\pgfpathlineto{\pgfqpoint{11.621cm}{14.476cm}}
\pgfpathlineto{\pgfqpoint{11.621cm}{14.699cm}}
\pgfpathclose
\pgfusepath{stroke}
\pgfsetdash{}{0cm}
\pgfpathmoveto{\pgfqpoint{12.8cm}{13.529cm}}
\pgfpathlineto{\pgfqpoint{13.023cm}{13.529cm}}
\pgfpathlineto{\pgfqpoint{13.023cm}{13.306cm}}
\pgfpathlineto{\pgfqpoint{12.8cm}{13.306cm}}
\pgfpathlineto{\pgfqpoint{12.8cm}{13.529cm}}
\pgfpathclose
\pgfusepath{stroke}
\pgfsetdash{}{0cm}
\pgfpathmoveto{\pgfqpoint{13.976cm}{12.415cm}}
\pgfpathlineto{\pgfqpoint{14.199cm}{12.415cm}}
\pgfpathlineto{\pgfqpoint{14.199cm}{12.191cm}}
\pgfpathlineto{\pgfqpoint{13.976cm}{12.191cm}}
\pgfpathlineto{\pgfqpoint{13.976cm}{12.415cm}}
\pgfpathclose
\pgfusepath{stroke}
\pgfsetdash{}{0cm}
\pgfpathmoveto{\pgfqpoint{15.155cm}{11.377cm}}
\pgfpathlineto{\pgfqpoint{15.378cm}{11.377cm}}
\pgfpathlineto{\pgfqpoint{15.378cm}{11.154cm}}
\pgfpathlineto{\pgfqpoint{15.155cm}{11.154cm}}
\pgfpathlineto{\pgfqpoint{15.155cm}{11.377cm}}
\pgfpathclose
\pgfusepath{stroke}
\pgfsetdash{}{0cm}
\pgfpathmoveto{\pgfqpoint{16.331cm}{10.383cm}}
\pgfpathlineto{\pgfqpoint{16.554cm}{10.383cm}}
\pgfpathlineto{\pgfqpoint{16.554cm}{10.16cm}}
\pgfpathlineto{\pgfqpoint{16.331cm}{10.16cm}}
\pgfpathlineto{\pgfqpoint{16.331cm}{10.383cm}}
\pgfpathclose
\pgfusepath{stroke}
\pgfsetdash{}{0cm}
\pgfpathmoveto{\pgfqpoint{17.51cm}{9.46cm}}
\pgfpathlineto{\pgfqpoint{17.733cm}{9.46cm}}
\pgfpathlineto{\pgfqpoint{17.733cm}{9.237cm}}
\pgfpathlineto{\pgfqpoint{17.51cm}{9.237cm}}
\pgfpathlineto{\pgfqpoint{17.51cm}{9.46cm}}
\pgfpathclose
\pgfusepath{stroke}
\begin{pgfscope}
\pgfpathmoveto{\pgfqpoint{5.847cm}{22.969cm}}
\pgfpathlineto{\pgfqpoint{17.624cm}{22.969cm}}
\pgfpathlineto{\pgfqpoint{17.624cm}{5.736cm}}
\pgfpathlineto{\pgfqpoint{5.847cm}{5.736cm}}
\pgfpathclose
\pgfusepath{clip}
\pgfsetdash{}{0cm}
\definecolor{eps2pgf_color}{rgb}{0,1,0}\pgfsetstrokecolor{eps2pgf_color}\pgfsetfillcolor{eps2pgf_color}
\pgfpathmoveto{\pgfqpoint{7.023cm}{19.92cm}}
\pgfpathlineto{\pgfqpoint{8.202cm}{18.477cm}}
\pgfpathlineto{\pgfqpoint{9.378cm}{17.104cm}}
\pgfpathlineto{\pgfqpoint{10.557cm}{15.79cm}}
\pgfpathlineto{\pgfqpoint{11.733cm}{14.546cm}}
\pgfpathlineto{\pgfqpoint{12.912cm}{13.37cm}}
\pgfpathlineto{\pgfqpoint{14.088cm}{12.247cm}}
\pgfpathlineto{\pgfqpoint{15.266cm}{11.195cm}}
\pgfpathlineto{\pgfqpoint{16.442cm}{10.198cm}}
\pgfpathlineto{\pgfqpoint{17.621cm}{9.26cm}}
\pgfusepath{stroke}
\end{pgfscope}
\pgfsetdash{}{0cm}
\definecolor{eps2pgf_color}{rgb}{0,1,0}\pgfsetstrokecolor{eps2pgf_color}\pgfsetfillcolor{eps2pgf_color}
\pgfpathmoveto{\pgfqpoint{6.882cm}{19.92cm}}
\pgfpathlineto{\pgfqpoint{7.164cm}{19.92cm}}
\pgfusepath{stroke}
\pgfsetdash{}{0cm}
\pgfpathmoveto{\pgfqpoint{7.023cm}{20.061cm}}
\pgfpathlineto{\pgfqpoint{7.023cm}{19.779cm}}
\pgfusepath{stroke}
\pgfsetdash{}{0cm}
\pgfpathmoveto{\pgfqpoint{8.061cm}{18.477cm}}
\pgfpathlineto{\pgfqpoint{8.343cm}{18.477cm}}
\pgfusepath{stroke}
\pgfsetdash{}{0cm}
\pgfpathmoveto{\pgfqpoint{8.202cm}{18.618cm}}
\pgfpathlineto{\pgfqpoint{8.202cm}{18.336cm}}
\pgfusepath{stroke}
\pgfsetdash{}{0cm}
\pgfpathmoveto{\pgfqpoint{9.237cm}{17.104cm}}
\pgfpathlineto{\pgfqpoint{9.519cm}{17.104cm}}
\pgfusepath{stroke}
\pgfsetdash{}{0cm}
\pgfpathmoveto{\pgfqpoint{9.378cm}{17.245cm}}
\pgfpathlineto{\pgfqpoint{9.378cm}{16.963cm}}
\pgfusepath{stroke}
\pgfsetdash{}{0cm}
\pgfpathmoveto{\pgfqpoint{10.416cm}{15.79cm}}
\pgfpathlineto{\pgfqpoint{10.698cm}{15.79cm}}
\pgfusepath{stroke}
\pgfsetdash{}{0cm}
\pgfpathmoveto{\pgfqpoint{10.557cm}{15.931cm}}
\pgfpathlineto{\pgfqpoint{10.557cm}{15.649cm}}
\pgfusepath{stroke}
\pgfsetdash{}{0cm}
\pgfpathmoveto{\pgfqpoint{11.592cm}{14.546cm}}
\pgfpathlineto{\pgfqpoint{11.874cm}{14.546cm}}
\pgfusepath{stroke}
\pgfsetdash{}{0cm}
\pgfpathmoveto{\pgfqpoint{11.733cm}{14.687cm}}
\pgfpathlineto{\pgfqpoint{11.733cm}{14.405cm}}
\pgfusepath{stroke}
\pgfsetdash{}{0cm}
\pgfpathmoveto{\pgfqpoint{12.771cm}{13.37cm}}
\pgfpathlineto{\pgfqpoint{13.053cm}{13.37cm}}
\pgfusepath{stroke}
\pgfsetdash{}{0cm}
\pgfpathmoveto{\pgfqpoint{12.912cm}{13.511cm}}
\pgfpathlineto{\pgfqpoint{12.912cm}{13.229cm}}
\pgfusepath{stroke}
\pgfsetdash{}{0cm}
\pgfpathmoveto{\pgfqpoint{13.946cm}{12.247cm}}
\pgfpathlineto{\pgfqpoint{14.229cm}{12.247cm}}
\pgfusepath{stroke}
\pgfsetdash{}{0cm}
\pgfpathmoveto{\pgfqpoint{14.088cm}{12.388cm}}
\pgfpathlineto{\pgfqpoint{14.088cm}{12.106cm}}
\pgfusepath{stroke}
\pgfsetdash{}{0cm}
\pgfpathmoveto{\pgfqpoint{15.125cm}{11.195cm}}
\pgfpathlineto{\pgfqpoint{15.408cm}{11.195cm}}
\pgfusepath{stroke}
\pgfsetdash{}{0cm}
\pgfpathmoveto{\pgfqpoint{15.266cm}{11.336cm}}
\pgfpathlineto{\pgfqpoint{15.266cm}{11.054cm}}
\pgfusepath{stroke}
\pgfsetdash{}{0cm}
\pgfpathmoveto{\pgfqpoint{16.301cm}{10.198cm}}
\pgfpathlineto{\pgfqpoint{16.583cm}{10.198cm}}
\pgfusepath{stroke}
\pgfsetdash{}{0cm}
\pgfpathmoveto{\pgfqpoint{16.442cm}{10.339cm}}
\pgfpathlineto{\pgfqpoint{16.442cm}{10.057cm}}
\pgfusepath{stroke}
\pgfsetdash{}{0cm}
\pgfpathmoveto{\pgfqpoint{17.48cm}{9.26cm}}
\pgfpathlineto{\pgfqpoint{17.762cm}{9.26cm}}
\pgfusepath{stroke}
\pgfsetdash{}{0cm}
\pgfpathmoveto{\pgfqpoint{17.621cm}{9.402cm}}
\pgfpathlineto{\pgfqpoint{17.621cm}{9.119cm}}
\pgfusepath{stroke}
\pgfsetdash{}{0cm}
\pgfpathmoveto{\pgfqpoint{6.926cm}{20.017cm}}
\pgfpathlineto{\pgfqpoint{7.12cm}{19.823cm}}
\pgfusepath{stroke}
\pgfsetdash{}{0cm}
\pgfpathmoveto{\pgfqpoint{7.12cm}{20.017cm}}
\pgfpathlineto{\pgfqpoint{6.926cm}{19.823cm}}
\pgfusepath{stroke}
\pgfsetdash{}{0cm}
\pgfpathmoveto{\pgfqpoint{8.105cm}{18.574cm}}
\pgfpathlineto{\pgfqpoint{8.299cm}{18.38cm}}
\pgfusepath{stroke}
\pgfsetdash{}{0cm}
\pgfpathmoveto{\pgfqpoint{8.299cm}{18.574cm}}
\pgfpathlineto{\pgfqpoint{8.105cm}{18.38cm}}
\pgfusepath{stroke}
\pgfsetdash{}{0cm}
\pgfpathmoveto{\pgfqpoint{9.281cm}{17.201cm}}
\pgfpathlineto{\pgfqpoint{9.475cm}{17.007cm}}
\pgfusepath{stroke}
\pgfsetdash{}{0cm}
\pgfpathmoveto{\pgfqpoint{9.475cm}{17.201cm}}
\pgfpathlineto{\pgfqpoint{9.281cm}{17.007cm}}
\pgfusepath{stroke}
\pgfsetdash{}{0cm}
\pgfpathmoveto{\pgfqpoint{10.46cm}{15.887cm}}
\pgfpathlineto{\pgfqpoint{10.654cm}{15.693cm}}
\pgfusepath{stroke}
\pgfsetdash{}{0cm}
\pgfpathmoveto{\pgfqpoint{10.654cm}{15.887cm}}
\pgfpathlineto{\pgfqpoint{10.46cm}{15.693cm}}
\pgfusepath{stroke}
\pgfsetdash{}{0cm}
\pgfpathmoveto{\pgfqpoint{11.636cm}{14.643cm}}
\pgfpathlineto{\pgfqpoint{11.83cm}{14.449cm}}
\pgfusepath{stroke}
\pgfsetdash{}{0cm}
\pgfpathmoveto{\pgfqpoint{11.83cm}{14.643cm}}
\pgfpathlineto{\pgfqpoint{11.636cm}{14.449cm}}
\pgfusepath{stroke}
\pgfsetdash{}{0cm}
\pgfpathmoveto{\pgfqpoint{12.815cm}{13.467cm}}
\pgfpathlineto{\pgfqpoint{13.009cm}{13.273cm}}
\pgfusepath{stroke}
\pgfsetdash{}{0cm}
\pgfpathmoveto{\pgfqpoint{13.009cm}{13.467cm}}
\pgfpathlineto{\pgfqpoint{12.815cm}{13.273cm}}
\pgfusepath{stroke}
\pgfsetdash{}{0cm}
\pgfpathmoveto{\pgfqpoint{13.991cm}{12.344cm}}
\pgfpathlineto{\pgfqpoint{14.185cm}{12.15cm}}
\pgfusepath{stroke}
\pgfsetdash{}{0cm}
\pgfpathmoveto{\pgfqpoint{14.185cm}{12.344cm}}
\pgfpathlineto{\pgfqpoint{13.991cm}{12.15cm}}
\pgfusepath{stroke}
\pgfsetdash{}{0cm}
\pgfpathmoveto{\pgfqpoint{15.169cm}{11.292cm}}
\pgfpathlineto{\pgfqpoint{15.363cm}{11.098cm}}
\pgfusepath{stroke}
\pgfsetdash{}{0cm}
\pgfpathmoveto{\pgfqpoint{15.363cm}{11.292cm}}
\pgfpathlineto{\pgfqpoint{15.169cm}{11.098cm}}
\pgfusepath{stroke}
\pgfsetdash{}{0cm}
\pgfpathmoveto{\pgfqpoint{16.345cm}{10.295cm}}
\pgfpathlineto{\pgfqpoint{16.539cm}{10.101cm}}
\pgfusepath{stroke}
\pgfsetdash{}{0cm}
\pgfpathmoveto{\pgfqpoint{16.539cm}{10.295cm}}
\pgfpathlineto{\pgfqpoint{16.345cm}{10.101cm}}
\pgfusepath{stroke}
\pgfsetdash{}{0cm}
\pgfpathmoveto{\pgfqpoint{17.524cm}{9.357cm}}
\pgfpathlineto{\pgfqpoint{17.718cm}{9.163cm}}
\pgfusepath{stroke}
\pgfsetdash{}{0cm}
\pgfpathmoveto{\pgfqpoint{17.718cm}{9.357cm}}
\pgfpathlineto{\pgfqpoint{17.524cm}{9.163cm}}
\pgfusepath{stroke}
\begin{pgfscope}
\pgfpathmoveto{\pgfqpoint{5.847cm}{22.969cm}}
\pgfpathlineto{\pgfqpoint{17.624cm}{22.969cm}}
\pgfpathlineto{\pgfqpoint{17.624cm}{5.736cm}}
\pgfpathlineto{\pgfqpoint{5.847cm}{5.736cm}}
\pgfpathclose
\pgfusepath{clip}
\pgfsetdash{{0.018cm}{0.141cm}{0.212cm}{0.141cm}}{0cm}
\definecolor{eps2pgf_color}{rgb}{0,0,1}\pgfsetstrokecolor{eps2pgf_color}\pgfsetfillcolor{eps2pgf_color}
\pgfpathmoveto{\pgfqpoint{7.023cm}{19.838cm}}
\pgfpathlineto{\pgfqpoint{8.202cm}{18.306cm}}
\pgfpathlineto{\pgfqpoint{9.378cm}{16.854cm}}
\pgfpathlineto{\pgfqpoint{10.557cm}{15.472cm}}
\pgfpathlineto{\pgfqpoint{11.733cm}{14.152cm}}
\pgfpathlineto{\pgfqpoint{12.912cm}{12.897cm}}
\pgfpathlineto{\pgfqpoint{14.088cm}{11.695cm}}
\pgfpathlineto{\pgfqpoint{15.266cm}{10.569cm}}
\pgfpathlineto{\pgfqpoint{16.442cm}{9.49cm}}
\pgfpathlineto{\pgfqpoint{17.621cm}{8.478cm}}
\pgfusepath{stroke}
\end{pgfscope}
\pgfsetdash{}{0cm}
\definecolor{eps2pgf_color}{rgb}{0,0,1}\pgfsetstrokecolor{eps2pgf_color}\pgfsetfillcolor{eps2pgf_color}
\pgfpathmoveto{\pgfqpoint{6.912cm}{19.95cm}}
\pgfpathlineto{\pgfqpoint{7.135cm}{19.95cm}}
\pgfpathlineto{\pgfqpoint{7.135cm}{19.726cm}}
\pgfpathlineto{\pgfqpoint{6.912cm}{19.726cm}}
\pgfpathlineto{\pgfqpoint{6.912cm}{19.95cm}}
\pgfpathclose
\pgfusepath{stroke}
\pgfsetdash{}{0cm}
\pgfpathmoveto{\pgfqpoint{8.09cm}{18.418cm}}
\pgfpathlineto{\pgfqpoint{8.314cm}{18.418cm}}
\pgfpathlineto{\pgfqpoint{8.314cm}{18.195cm}}
\pgfpathlineto{\pgfqpoint{8.09cm}{18.195cm}}
\pgfpathlineto{\pgfqpoint{8.09cm}{18.418cm}}
\pgfpathclose
\pgfusepath{stroke}
\pgfsetdash{}{0cm}
\pgfpathmoveto{\pgfqpoint{9.266cm}{16.966cm}}
\pgfpathlineto{\pgfqpoint{9.49cm}{16.966cm}}
\pgfpathlineto{\pgfqpoint{9.49cm}{16.742cm}}
\pgfpathlineto{\pgfqpoint{9.266cm}{16.742cm}}
\pgfpathlineto{\pgfqpoint{9.266cm}{16.966cm}}
\pgfpathclose
\pgfusepath{stroke}
\pgfsetdash{}{0cm}
\pgfpathmoveto{\pgfqpoint{10.445cm}{15.584cm}}
\pgfpathlineto{\pgfqpoint{10.669cm}{15.584cm}}
\pgfpathlineto{\pgfqpoint{10.669cm}{15.361cm}}
\pgfpathlineto{\pgfqpoint{10.445cm}{15.361cm}}
\pgfpathlineto{\pgfqpoint{10.445cm}{15.584cm}}
\pgfpathclose
\pgfusepath{stroke}
\pgfsetdash{}{0cm}
\pgfpathmoveto{\pgfqpoint{11.621cm}{14.264cm}}
\pgfpathlineto{\pgfqpoint{11.845cm}{14.264cm}}
\pgfpathlineto{\pgfqpoint{11.845cm}{14.041cm}}
\pgfpathlineto{\pgfqpoint{11.621cm}{14.041cm}}
\pgfpathlineto{\pgfqpoint{11.621cm}{14.264cm}}
\pgfpathclose
\pgfusepath{stroke}
\pgfsetdash{}{0cm}
\pgfpathmoveto{\pgfqpoint{12.8cm}{13.009cm}}
\pgfpathlineto{\pgfqpoint{13.023cm}{13.009cm}}
\pgfpathlineto{\pgfqpoint{13.023cm}{12.785cm}}
\pgfpathlineto{\pgfqpoint{12.8cm}{12.785cm}}
\pgfpathlineto{\pgfqpoint{12.8cm}{13.009cm}}
\pgfpathclose
\pgfusepath{stroke}
\pgfsetdash{}{0cm}
\pgfpathmoveto{\pgfqpoint{13.976cm}{11.806cm}}
\pgfpathlineto{\pgfqpoint{14.199cm}{11.806cm}}
\pgfpathlineto{\pgfqpoint{14.199cm}{11.583cm}}
\pgfpathlineto{\pgfqpoint{13.976cm}{11.583cm}}
\pgfpathlineto{\pgfqpoint{13.976cm}{11.806cm}}
\pgfpathclose
\pgfusepath{stroke}
\pgfsetdash{}{0cm}
\pgfpathmoveto{\pgfqpoint{15.155cm}{10.68cm}}
\pgfpathlineto{\pgfqpoint{15.378cm}{10.68cm}}
\pgfpathlineto{\pgfqpoint{15.378cm}{10.457cm}}
\pgfpathlineto{\pgfqpoint{15.155cm}{10.457cm}}
\pgfpathlineto{\pgfqpoint{15.155cm}{10.68cm}}
\pgfpathclose
\pgfusepath{stroke}
\pgfsetdash{}{0cm}
\pgfpathmoveto{\pgfqpoint{16.331cm}{9.601cm}}
\pgfpathlineto{\pgfqpoint{16.554cm}{9.601cm}}
\pgfpathlineto{\pgfqpoint{16.554cm}{9.378cm}}
\pgfpathlineto{\pgfqpoint{16.331cm}{9.378cm}}
\pgfpathlineto{\pgfqpoint{16.331cm}{9.601cm}}
\pgfpathclose
\pgfusepath{stroke}
\pgfsetdash{}{0cm}
\pgfpathmoveto{\pgfqpoint{17.51cm}{8.59cm}}
\pgfpathlineto{\pgfqpoint{17.733cm}{8.59cm}}
\pgfpathlineto{\pgfqpoint{17.733cm}{8.367cm}}
\pgfpathlineto{\pgfqpoint{17.51cm}{8.367cm}}
\pgfpathlineto{\pgfqpoint{17.51cm}{8.59cm}}
\pgfpathclose
\pgfusepath{stroke}
\begin{pgfscope}
\pgfpathmoveto{\pgfqpoint{5.847cm}{22.969cm}}
\pgfpathlineto{\pgfqpoint{17.624cm}{22.969cm}}
\pgfpathlineto{\pgfqpoint{17.624cm}{5.736cm}}
\pgfpathlineto{\pgfqpoint{5.847cm}{5.736cm}}
\pgfpathclose
\pgfusepath{clip}
\end{pgfscope}
\definecolor{eps2pgf_color}{gray}{0}\pgfsetstrokecolor{eps2pgf_color}\pgfsetfillcolor{eps2pgf_color}
\pgftext[x=11.751cm,y=4.093cm,rotate=0]{\fontsize{36}{36.14}\selectfont{ {$\bmax$}}}
\pgftext[x=3.013cm,y=14.388cm,rotate=90]{\fontsize{36}{36.14}\selectfont{ { Average cost}}}
\pgftext[x=5.797cm,y=5.612cm,rotate=0]{\fontsize{10.04}{12.04}\selectfont{ { }}}
\pgftext[x=17.574cm,y=22.845cm,rotate=0]{\fontsize{10.04}{12.04}\selectfont{ { }}}
\pgftext[x=9.91cm+.1cm,y=10.24cm,rotate=0]{\fontsize{34}{36.14}\selectfont{ {No storage}}}
\begin{pgfscope}
\pgfpathmoveto{\pgfqpoint{6.024cm}{11.021cm}}
\pgfpathlineto{\pgfqpoint{13.429cm}{11.021cm}}
\pgfpathlineto{\pgfqpoint{13.429cm}{5.909cm}}
\pgfpathlineto{\pgfqpoint{6.024cm}{5.909cm}}
\pgfpathclose
\pgfusepath{clip}
\pgfsetdash{}{0cm}
\pgfpathmoveto{\pgfqpoint{6.229cm}{10.342cm}}
\pgfpathlineto{\pgfqpoint{7.258cm}{10.342cm}}
\pgfusepath{stroke}
\begin{pgfscope}
\pgfpathmoveto{\pgfqpoint{6.35cm}{10.736cm}}
\pgfpathlineto{\pgfqpoint{7.141cm}{10.736cm}}
\pgfpathlineto{\pgfqpoint{7.141cm}{9.945cm}}
\pgfpathlineto{\pgfqpoint{6.35cm}{9.945cm}}
\pgfpathclose
\pgfusepath{clip}
\pgfsetdash{}{0cm}
\pgfpathmoveto{\pgfqpoint{6.885cm}{10.342cm}}
\pgfpathcurveto{\pgfqpoint{6.885cm}{10.264cm}}{\pgfqpoint{6.822cm}{10.201cm}}{\pgfqpoint{6.744cm}{10.201cm}}
\pgfpathcurveto{\pgfqpoint{6.666cm}{10.201cm}}{\pgfqpoint{6.603cm}{10.264cm}}{\pgfqpoint{6.603cm}{10.342cm}}
\pgfpathcurveto{\pgfqpoint{6.603cm}{10.42cm}}{\pgfqpoint{6.666cm}{10.483cm}}{\pgfqpoint{6.744cm}{10.483cm}}
\pgfpathcurveto{\pgfqpoint{6.822cm}{10.483cm}}{\pgfqpoint{6.885cm}{10.42cm}}{\pgfqpoint{6.885cm}{10.342cm}}
\pgfusepath{stroke}
\end{pgfscope}
\end{pgfscope}
\pgftext[x=8.607cm+.2cm,y=9.096cm,rotate=0]{\fontsize{34}{36.14}\selectfont{ {OMG}}}
\begin{pgfscope}
\pgfpathmoveto{\pgfqpoint{6.024cm}{11.021cm}}
\pgfpathlineto{\pgfqpoint{13.429cm}{11.021cm}}
\pgfpathlineto{\pgfqpoint{13.429cm}{5.909cm}}
\pgfpathlineto{\pgfqpoint{6.024cm}{5.909cm}}
\pgfpathclose
\pgfusepath{clip}
\pgfsetdash{}{0cm}
\definecolor{eps2pgf_color}{rgb}{0,0,1}\pgfsetstrokecolor{eps2pgf_color}\pgfsetfillcolor{eps2pgf_color}
\pgfpathmoveto{\pgfqpoint{6.229cm}{9.093cm}}
\pgfpathlineto{\pgfqpoint{7.258cm}{9.093cm}}
\pgfusepath{stroke}
\begin{pgfscope}
\pgfpathmoveto{\pgfqpoint{6.35cm}{9.487cm}}
\pgfpathlineto{\pgfqpoint{7.141cm}{9.487cm}}
\pgfpathlineto{\pgfqpoint{7.141cm}{8.696cm}}
\pgfpathlineto{\pgfqpoint{6.35cm}{8.696cm}}
\pgfpathclose
\pgfusepath{clip}
\pgfsetdash{}{0cm}
\pgfpathmoveto{\pgfqpoint{6.632cm}{9.205cm}}
\pgfpathlineto{\pgfqpoint{6.856cm}{9.205cm}}
\pgfpathlineto{\pgfqpoint{6.856cm}{8.981cm}}
\pgfpathlineto{\pgfqpoint{6.632cm}{8.981cm}}
\pgfpathlineto{\pgfqpoint{6.632cm}{9.205cm}}
\pgfpathclose
\pgfusepath{stroke}
\end{pgfscope}
\end{pgfscope}
\pgftext[x=9.116cm,y=7.744cm,rotate=0]{\fontsize{34}{36.14}\selectfont{ {Greedy}}}
\begin{pgfscope}
\pgfpathmoveto{\pgfqpoint{6.024cm}{11.021cm}}
\pgfpathlineto{\pgfqpoint{13.429cm}{11.021cm}}
\pgfpathlineto{\pgfqpoint{13.429cm}{5.909cm}}
\pgfpathlineto{\pgfqpoint{6.024cm}{5.909cm}}
\pgfpathclose
\pgfusepath{clip}
\pgfsetdash{}{0cm}
\definecolor{eps2pgf_color}{rgb}{0,1,0}\pgfsetstrokecolor{eps2pgf_color}\pgfsetfillcolor{eps2pgf_color}
\pgfpathmoveto{\pgfqpoint{6.229cm}{7.843cm}}
\pgfpathlineto{\pgfqpoint{7.258cm}{7.843cm}}
\pgfusepath{stroke}
\begin{pgfscope}
\pgfpathmoveto{\pgfqpoint{6.35cm}{8.237cm}}
\pgfpathlineto{\pgfqpoint{7.141cm}{8.237cm}}
\pgfpathlineto{\pgfqpoint{7.141cm}{7.447cm}}
\pgfpathlineto{\pgfqpoint{6.35cm}{7.447cm}}
\pgfpathclose
\pgfusepath{clip}
\pgfsetdash{}{0cm}
\pgfpathmoveto{\pgfqpoint{6.603cm}{7.843cm}}
\pgfpathlineto{\pgfqpoint{6.885cm}{7.843cm}}
\pgfusepath{stroke}
\pgfsetdash{}{0cm}
\pgfpathmoveto{\pgfqpoint{6.744cm}{7.985cm}}
\pgfpathlineto{\pgfqpoint{6.744cm}{7.702cm}}
\pgfusepath{stroke}
\pgfsetdash{}{0cm}
\pgfpathmoveto{\pgfqpoint{6.647cm}{7.94cm}}
\pgfpathlineto{\pgfqpoint{6.841cm}{7.746cm}}
\pgfusepath{stroke}
\pgfsetdash{}{0cm}
\pgfpathmoveto{\pgfqpoint{6.841cm}{7.94cm}}
\pgfpathlineto{\pgfqpoint{6.647cm}{7.746cm}}
\pgfusepath{stroke}
\end{pgfscope}
\end{pgfscope}
\pgftext[x=10.342cm+.4cm,y=6.59cm,rotate=0]{\fontsize{34}{36.14}\selectfont{ {Lower bound }}}
\begin{pgfscope}
\pgfpathmoveto{\pgfqpoint{6.024cm}{11.021cm}}
\pgfpathlineto{\pgfqpoint{13.429cm}{11.021cm}}
\pgfpathlineto{\pgfqpoint{13.429cm}{5.909cm}}
\pgfpathlineto{\pgfqpoint{6.024cm}{5.909cm}}
\pgfpathclose
\pgfusepath{clip}
\pgfsetdash{{0.018cm}{0.141cm}{0.212cm}{0.141cm}}{0cm}
\definecolor{eps2pgf_color}{rgb}{0,0,1}\pgfsetstrokecolor{eps2pgf_color}\pgfsetfillcolor{eps2pgf_color}
\pgfpathmoveto{\pgfqpoint{6.229cm}{6.594cm}}
\pgfpathlineto{\pgfqpoint{7.258cm}{6.594cm}}
\pgfusepath{stroke}
\begin{pgfscope}
\pgfpathmoveto{\pgfqpoint{6.35cm}{6.988cm}}
\pgfpathlineto{\pgfqpoint{7.141cm}{6.988cm}}
\pgfpathlineto{\pgfqpoint{7.141cm}{6.197cm}}
\pgfpathlineto{\pgfqpoint{6.35cm}{6.197cm}}
\pgfpathclose
\pgfusepath{clip}
\pgfsetdash{}{0cm}
\pgfpathmoveto{\pgfqpoint{6.632cm}{6.706cm}}
\pgfpathlineto{\pgfqpoint{6.856cm}{6.706cm}}
\pgfpathlineto{\pgfqpoint{6.856cm}{6.482cm}}
\pgfpathlineto{\pgfqpoint{6.632cm}{6.482cm}}
\pgfpathlineto{\pgfqpoint{6.632cm}{6.706cm}}
\pgfpathclose
\pgfusepath{stroke}
\end{pgfscope}
\end{pgfscope}
\pgfsetdash{}{0cm}
\pgfsetlinewidth{0.176mm}
\definecolor{eps2pgf_color}{rgb}{0,0,1}\pgfsetstrokecolor{eps2pgf_color}\pgfsetfillcolor{eps2pgf_color}
\pgfusepath{stroke}
\end{pgfscope}
\end{pgfscope}
\end{pgfpicture}

%% file: fig/I2.tex
\scalebox{0.265}{\scalefont{2} \input{./fig/single_storage_I2.pgf}}

%% file: fig/single_storage_I2.pgf
% Created by Eps2pgf 0.7.0 (build on 2008-08-24) on Mon Apr 06 20:47:21 PDT 2015
\begin{pgfpicture}
\pgfpathmoveto{\pgfqpoint{2.293cm}{3.387cm}}
\pgfpathlineto{\pgfqpoint{19.226cm}{3.387cm}}
\pgfpathlineto{\pgfqpoint{19.226cm}{24.553cm}}
\pgfpathlineto{\pgfqpoint{2.293cm}{24.553cm}}
\pgfpathclose
\pgfusepath{clip}
\begin{pgfscope}
\begin{pgfscope}
\pgfpathmoveto{\pgfqpoint{2.293cm}{24.553cm}}
\pgfpathlineto{\pgfqpoint{19.229cm}{24.553cm}}
\pgfpathlineto{\pgfqpoint{19.229cm}{3.413cm}}
\pgfpathlineto{\pgfqpoint{2.293cm}{3.413cm}}
\pgfpathclose
\pgfusepath{clip}
\begin{pgfscope}
\definecolor{eps2pgf_color}{gray}{1}\pgfsetstrokecolor{eps2pgf_color}\pgfsetfillcolor{eps2pgf_color}
\pgfpathmoveto{\pgfqpoint{2.293cm}{24.553cm}}
\pgfpathlineto{\pgfqpoint{19.232cm}{24.553cm}}
\pgfpathlineto{\pgfqpoint{19.232cm}{3.41cm}}
\pgfpathlineto{\pgfqpoint{2.293cm}{3.41cm}}
\pgfpathclose
\pgfusepath{fill}
\end{pgfscope}
\definecolor{eps2pgf_color}{gray}{1}\pgfsetstrokecolor{eps2pgf_color}\pgfsetfillcolor{eps2pgf_color}
\pgfpathmoveto{\pgfqpoint{5.844cm}{5.739cm}}
\pgfpathlineto{\pgfqpoint{5.844cm}{22.969cm}}
\pgfpathlineto{\pgfqpoint{17.621cm}{22.969cm}}
\pgfpathlineto{\pgfqpoint{17.621cm}{5.739cm}}
\pgfpathclose
\pgfseteorule\pgfusepath{fill}\pgfsetnonzerorule
\pgfsetdash{}{0cm}
\pgfsetlinewidth{0.176mm}
\pgfsetroundjoin
\pgfpathmoveto{\pgfqpoint{5.844cm}{5.739cm}}
\pgfpathlineto{\pgfqpoint{5.844cm}{22.969cm}}
\pgfpathlineto{\pgfqpoint{17.621cm}{22.969cm}}
\pgfpathlineto{\pgfqpoint{17.621cm}{5.739cm}}
\pgfpathlineto{\pgfqpoint{5.844cm}{5.739cm}}
\pgfusepath{stroke}
\pgfsetdash{}{0cm}
\definecolor{eps2pgf_color}{gray}{0}\pgfsetstrokecolor{eps2pgf_color}\pgfsetfillcolor{eps2pgf_color}
\pgfpathmoveto{\pgfqpoint{5.844cm}{5.739cm}}
\pgfpathlineto{\pgfqpoint{17.621cm}{5.739cm}}
\pgfusepath{stroke}
\pgfsetdash{}{0cm}
\pgfpathmoveto{\pgfqpoint{5.844cm}{22.969cm}}
\pgfpathlineto{\pgfqpoint{17.621cm}{22.969cm}}
\pgfusepath{stroke}
\pgfsetdash{}{0cm}
\pgfpathmoveto{\pgfqpoint{5.844cm}{5.739cm}}
\pgfpathlineto{\pgfqpoint{5.844cm}{22.969cm}}
\pgfusepath{stroke}
\pgfsetdash{}{0cm}
\pgfpathmoveto{\pgfqpoint{17.621cm}{5.739cm}}
\pgfpathlineto{\pgfqpoint{17.621cm}{22.969cm}}
\pgfusepath{stroke}
\pgfsetdash{}{0cm}
\pgfpathmoveto{\pgfqpoint{5.844cm}{5.739cm}}
\pgfpathlineto{\pgfqpoint{17.621cm}{5.739cm}}
\pgfusepath{stroke}
\pgfsetdash{}{0cm}
\pgfpathmoveto{\pgfqpoint{5.844cm}{5.739cm}}
\pgfpathlineto{\pgfqpoint{5.844cm}{22.969cm}}
\pgfusepath{stroke}
\pgfsetdash{}{0cm}
\pgfpathmoveto{\pgfqpoint{5.844cm}{5.739cm}}
\pgfpathlineto{\pgfqpoint{5.844cm}{5.912cm}}
\pgfusepath{stroke}
\pgfsetdash{}{0cm}
\pgfpathmoveto{\pgfqpoint{5.844cm}{22.969cm}}
\pgfpathlineto{\pgfqpoint{5.844cm}{22.798cm}}
\pgfusepath{stroke}
\pgftext[x=5.845cm,y=5.016cm+.2cm,rotate=0]{\fontsize{36}{36.14}\selectfont{ {0}}}
\pgfsetdash{}{0cm}
\pgfpathmoveto{\pgfqpoint{11.733cm}{5.739cm}}
\pgfpathlineto{\pgfqpoint{11.733cm}{5.912cm}}
\pgfusepath{stroke}
\pgfsetdash{}{0cm}
\pgfpathmoveto{\pgfqpoint{11.733cm}{22.969cm}}
\pgfpathlineto{\pgfqpoint{11.733cm}{22.798cm}}
\pgfusepath{stroke}
\pgftext[x=11.731cm,y=5.016cm+.2cm,rotate=0]{\fontsize{36}{36.14}\selectfont{ {0.5}}}
\pgfsetdash{}{0cm}
\pgfpathmoveto{\pgfqpoint{17.621cm}{5.739cm}}
\pgfpathlineto{\pgfqpoint{17.621cm}{5.912cm}}
\pgfusepath{stroke}
\pgfsetdash{}{0cm}
\pgfpathmoveto{\pgfqpoint{17.621cm}{22.969cm}}
\pgfpathlineto{\pgfqpoint{17.621cm}{22.798cm}}
\pgfusepath{stroke}
\pgftext[x=17.571cm,y=5.026cm+.2cm,rotate=0]{\fontsize{36}{36.14}\selectfont{ {1}}}
\pgfsetdash{}{0cm}
\pgfpathmoveto{\pgfqpoint{5.844cm}{5.739cm}}
\pgfpathlineto{\pgfqpoint{6.018cm}{5.739cm}}
\pgfusepath{stroke}
\pgfsetdash{}{0cm}
\pgfpathmoveto{\pgfqpoint{17.621cm}{5.739cm}}
\pgfpathlineto{\pgfqpoint{17.448cm}{5.739cm}}
\pgfusepath{stroke}
\pgftext[x=4.714cm-.2cm,y=5.707cm,rotate=0]{\fontsize{36}{36.14}\selectfont{ {0.05}}}
\pgfsetdash{}{0cm}
\pgfpathmoveto{\pgfqpoint{5.844cm}{8.39cm}}
\pgfpathlineto{\pgfqpoint{6.018cm}{8.39cm}}
\pgfusepath{stroke}
\pgfsetdash{}{0cm}
\pgfpathmoveto{\pgfqpoint{17.621cm}{8.39cm}}
\pgfpathlineto{\pgfqpoint{17.448cm}{8.39cm}}
\pgfusepath{stroke}
\pgftext[x=4.716cm-.2cm,y=8.358cm,rotate=0]{\fontsize{36}{36.14}\selectfont{ {0.06}}}
\pgfsetdash{}{0cm}
\pgfpathmoveto{\pgfqpoint{5.844cm}{11.042cm}}
\pgfpathlineto{\pgfqpoint{6.018cm}{11.042cm}}
\pgfusepath{stroke}
\pgfsetdash{}{0cm}
\pgfpathmoveto{\pgfqpoint{17.621cm}{11.042cm}}
\pgfpathlineto{\pgfqpoint{17.448cm}{11.042cm}}
\pgfusepath{stroke}
\pgftext[x=4.718cm-.2cm,y=11.01cm,rotate=0]{\fontsize{36}{36.14}\selectfont{ {0.07}}}
\pgfsetdash{}{0cm}
\pgfpathmoveto{\pgfqpoint{5.844cm}{13.694cm}}
\pgfpathlineto{\pgfqpoint{6.018cm}{13.694cm}}
\pgfusepath{stroke}
\pgfsetdash{}{0cm}
\pgfpathmoveto{\pgfqpoint{17.621cm}{13.694cm}}
\pgfpathlineto{\pgfqpoint{17.448cm}{13.694cm}}
\pgfusepath{stroke}
\pgftext[x=4.715cm-.2cm,y=13.662cm,rotate=0]{\fontsize{36}{36.14}\selectfont{ {0.08}}}
\pgfsetdash{}{0cm}
\pgfpathmoveto{\pgfqpoint{5.844cm}{16.342cm}}
\pgfpathlineto{\pgfqpoint{6.018cm}{16.342cm}}
\pgfusepath{stroke}
\pgfsetdash{}{0cm}
\pgfpathmoveto{\pgfqpoint{17.621cm}{16.342cm}}
\pgfpathlineto{\pgfqpoint{17.448cm}{16.342cm}}
\pgfusepath{stroke}
\pgftext[x=4.714cm-.2cm,y=16.31cm,rotate=0]{\fontsize{36}{36.14}\selectfont{ {0.09}}}
\pgfsetdash{}{0cm}
\pgfpathmoveto{\pgfqpoint{5.844cm}{18.994cm}}
\pgfpathlineto{\pgfqpoint{6.018cm}{18.994cm}}
\pgfusepath{stroke}
\pgfsetdash{}{0cm}
\pgfpathmoveto{\pgfqpoint{17.621cm}{18.994cm}}
\pgfpathlineto{\pgfqpoint{17.448cm}{18.994cm}}
\pgfusepath{stroke}
\pgftext[x=4.925cm-.2cm,y=18.962cm,rotate=0]{\fontsize{36}{36.14}\selectfont{ {0.1}}}
\pgfsetdash{}{0cm}
\pgfpathmoveto{\pgfqpoint{5.844cm}{21.646cm}}
\pgfpathlineto{\pgfqpoint{6.018cm}{21.646cm}}
\pgfusepath{stroke}
\pgfsetdash{}{0cm}
\pgfpathmoveto{\pgfqpoint{17.621cm}{21.646cm}}
\pgfpathlineto{\pgfqpoint{17.448cm}{21.646cm}}
\pgfusepath{stroke}
\pgftext[x=4.632cm-.2cm,y=21.614cm,rotate=0]{\fontsize{36}{36.14}\selectfont{ {0.11}}}
\pgfsetdash{}{0cm}
\pgfpathmoveto{\pgfqpoint{5.844cm}{5.739cm}}
\pgfpathlineto{\pgfqpoint{17.621cm}{5.739cm}}
\pgfusepath{stroke}
\pgfsetdash{}{0cm}
\pgfpathmoveto{\pgfqpoint{5.844cm}{22.969cm}}
\pgfpathlineto{\pgfqpoint{17.621cm}{22.969cm}}
\pgfusepath{stroke}
\pgfsetdash{}{0cm}
\pgfpathmoveto{\pgfqpoint{5.844cm}{5.739cm}}
\pgfpathlineto{\pgfqpoint{5.844cm}{22.969cm}}
\pgfusepath{stroke}
\pgfsetdash{}{0cm}
\pgfpathmoveto{\pgfqpoint{17.621cm}{5.739cm}}
\pgfpathlineto{\pgfqpoint{17.621cm}{22.969cm}}
\pgfusepath{stroke}
\begin{pgfscope}
\pgfpathmoveto{\pgfqpoint{5.847cm}{22.969cm}}
\pgfpathlineto{\pgfqpoint{17.624cm}{22.969cm}}
\pgfpathlineto{\pgfqpoint{17.624cm}{5.736cm}}
\pgfpathlineto{\pgfqpoint{5.847cm}{5.736cm}}
\pgfpathclose
\pgfusepath{clip}
\pgfsetdash{}{0cm}
\pgfsetlinewidth{1.058mm}
\pgfpathmoveto{\pgfqpoint{7.023cm}{21.799cm}}
\pgfpathlineto{\pgfqpoint{8.202cm}{21.799cm}}
\pgfpathlineto{\pgfqpoint{9.378cm}{21.799cm}}
\pgfpathlineto{\pgfqpoint{10.557cm}{21.799cm}}
\pgfpathlineto{\pgfqpoint{11.733cm}{21.799cm}}
\pgfpathlineto{\pgfqpoint{12.912cm}{21.799cm}}
\pgfpathlineto{\pgfqpoint{14.088cm}{21.799cm}}
\pgfpathlineto{\pgfqpoint{15.266cm}{21.799cm}}
\pgfpathlineto{\pgfqpoint{16.442cm}{21.799cm}}
\pgfpathlineto{\pgfqpoint{17.621cm}{21.799cm}}
\pgfusepath{stroke}
\end{pgfscope}
\pgfsetdash{}{0cm}
\pgfsetlinewidth{1.058mm}
\pgfpathmoveto{\pgfqpoint{7.164cm}{21.799cm}}
\pgfpathcurveto{\pgfqpoint{7.164cm}{21.721cm}}{\pgfqpoint{7.101cm}{21.658cm}}{\pgfqpoint{7.023cm}{21.658cm}}
\pgfpathcurveto{\pgfqpoint{6.945cm}{21.658cm}}{\pgfqpoint{6.882cm}{21.721cm}}{\pgfqpoint{6.882cm}{21.799cm}}
\pgfpathcurveto{\pgfqpoint{6.882cm}{21.877cm}}{\pgfqpoint{6.945cm}{21.94cm}}{\pgfqpoint{7.023cm}{21.94cm}}
\pgfpathcurveto{\pgfqpoint{7.101cm}{21.94cm}}{\pgfqpoint{7.164cm}{21.877cm}}{\pgfqpoint{7.164cm}{21.799cm}}
\pgfusepath{stroke}
\pgfsetdash{}{0cm}
\pgfpathmoveto{\pgfqpoint{8.343cm}{21.799cm}}
\pgfpathcurveto{\pgfqpoint{8.343cm}{21.721cm}}{\pgfqpoint{8.28cm}{21.658cm}}{\pgfqpoint{8.202cm}{21.658cm}}
\pgfpathcurveto{\pgfqpoint{8.124cm}{21.658cm}}{\pgfqpoint{8.061cm}{21.721cm}}{\pgfqpoint{8.061cm}{21.799cm}}
\pgfpathcurveto{\pgfqpoint{8.061cm}{21.877cm}}{\pgfqpoint{8.124cm}{21.94cm}}{\pgfqpoint{8.202cm}{21.94cm}}
\pgfpathcurveto{\pgfqpoint{8.28cm}{21.94cm}}{\pgfqpoint{8.343cm}{21.877cm}}{\pgfqpoint{8.343cm}{21.799cm}}
\pgfusepath{stroke}
\pgfsetdash{}{0cm}
\pgfpathmoveto{\pgfqpoint{9.519cm}{21.799cm}}
\pgfpathcurveto{\pgfqpoint{9.519cm}{21.721cm}}{\pgfqpoint{9.456cm}{21.658cm}}{\pgfqpoint{9.378cm}{21.658cm}}
\pgfpathcurveto{\pgfqpoint{9.3cm}{21.658cm}}{\pgfqpoint{9.237cm}{21.721cm}}{\pgfqpoint{9.237cm}{21.799cm}}
\pgfpathcurveto{\pgfqpoint{9.237cm}{21.877cm}}{\pgfqpoint{9.3cm}{21.94cm}}{\pgfqpoint{9.378cm}{21.94cm}}
\pgfpathcurveto{\pgfqpoint{9.456cm}{21.94cm}}{\pgfqpoint{9.519cm}{21.877cm}}{\pgfqpoint{9.519cm}{21.799cm}}
\pgfusepath{stroke}
\pgfsetdash{}{0cm}
\pgfpathmoveto{\pgfqpoint{10.698cm}{21.799cm}}
\pgfpathcurveto{\pgfqpoint{10.698cm}{21.721cm}}{\pgfqpoint{10.635cm}{21.658cm}}{\pgfqpoint{10.557cm}{21.658cm}}
\pgfpathcurveto{\pgfqpoint{10.479cm}{21.658cm}}{\pgfqpoint{10.416cm}{21.721cm}}{\pgfqpoint{10.416cm}{21.799cm}}
\pgfpathcurveto{\pgfqpoint{10.416cm}{21.877cm}}{\pgfqpoint{10.479cm}{21.94cm}}{\pgfqpoint{10.557cm}{21.94cm}}
\pgfpathcurveto{\pgfqpoint{10.635cm}{21.94cm}}{\pgfqpoint{10.698cm}{21.877cm}}{\pgfqpoint{10.698cm}{21.799cm}}
\pgfusepath{stroke}
\pgfsetdash{}{0cm}
\pgfpathmoveto{\pgfqpoint{11.874cm}{21.799cm}}
\pgfpathcurveto{\pgfqpoint{11.874cm}{21.721cm}}{\pgfqpoint{11.811cm}{21.658cm}}{\pgfqpoint{11.733cm}{21.658cm}}
\pgfpathcurveto{\pgfqpoint{11.655cm}{21.658cm}}{\pgfqpoint{11.592cm}{21.721cm}}{\pgfqpoint{11.592cm}{21.799cm}}
\pgfpathcurveto{\pgfqpoint{11.592cm}{21.877cm}}{\pgfqpoint{11.655cm}{21.94cm}}{\pgfqpoint{11.733cm}{21.94cm}}
\pgfpathcurveto{\pgfqpoint{11.811cm}{21.94cm}}{\pgfqpoint{11.874cm}{21.877cm}}{\pgfqpoint{11.874cm}{21.799cm}}
\pgfusepath{stroke}
\pgfsetdash{}{0cm}
\pgfpathmoveto{\pgfqpoint{13.053cm}{21.799cm}}
\pgfpathcurveto{\pgfqpoint{13.053cm}{21.721cm}}{\pgfqpoint{12.99cm}{21.658cm}}{\pgfqpoint{12.912cm}{21.658cm}}
\pgfpathcurveto{\pgfqpoint{12.834cm}{21.658cm}}{\pgfqpoint{12.771cm}{21.721cm}}{\pgfqpoint{12.771cm}{21.799cm}}
\pgfpathcurveto{\pgfqpoint{12.771cm}{21.877cm}}{\pgfqpoint{12.834cm}{21.94cm}}{\pgfqpoint{12.912cm}{21.94cm}}
\pgfpathcurveto{\pgfqpoint{12.99cm}{21.94cm}}{\pgfqpoint{13.053cm}{21.877cm}}{\pgfqpoint{13.053cm}{21.799cm}}
\pgfusepath{stroke}
\pgfsetdash{}{0cm}
\pgfpathmoveto{\pgfqpoint{14.229cm}{21.799cm}}
\pgfpathcurveto{\pgfqpoint{14.229cm}{21.721cm}}{\pgfqpoint{14.166cm}{21.658cm}}{\pgfqpoint{14.088cm}{21.658cm}}
\pgfpathcurveto{\pgfqpoint{14.01cm}{21.658cm}}{\pgfqpoint{13.946cm}{21.721cm}}{\pgfqpoint{13.946cm}{21.799cm}}
\pgfpathcurveto{\pgfqpoint{13.946cm}{21.877cm}}{\pgfqpoint{14.01cm}{21.94cm}}{\pgfqpoint{14.088cm}{21.94cm}}
\pgfpathcurveto{\pgfqpoint{14.166cm}{21.94cm}}{\pgfqpoint{14.229cm}{21.877cm}}{\pgfqpoint{14.229cm}{21.799cm}}
\pgfusepath{stroke}
\pgfsetdash{}{0cm}
\pgfpathmoveto{\pgfqpoint{15.408cm}{21.799cm}}
\pgfpathcurveto{\pgfqpoint{15.408cm}{21.721cm}}{\pgfqpoint{15.344cm}{21.658cm}}{\pgfqpoint{15.266cm}{21.658cm}}
\pgfpathcurveto{\pgfqpoint{15.189cm}{21.658cm}}{\pgfqpoint{15.125cm}{21.721cm}}{\pgfqpoint{15.125cm}{21.799cm}}
\pgfpathcurveto{\pgfqpoint{15.125cm}{21.877cm}}{\pgfqpoint{15.189cm}{21.94cm}}{\pgfqpoint{15.266cm}{21.94cm}}
\pgfpathcurveto{\pgfqpoint{15.344cm}{21.94cm}}{\pgfqpoint{15.408cm}{21.877cm}}{\pgfqpoint{15.408cm}{21.799cm}}
\pgfusepath{stroke}
\pgfsetdash{}{0cm}
\pgfpathmoveto{\pgfqpoint{16.583cm}{21.799cm}}
\pgfpathcurveto{\pgfqpoint{16.583cm}{21.721cm}}{\pgfqpoint{16.52cm}{21.658cm}}{\pgfqpoint{16.442cm}{21.658cm}}
\pgfpathcurveto{\pgfqpoint{16.364cm}{21.658cm}}{\pgfqpoint{16.301cm}{21.721cm}}{\pgfqpoint{16.301cm}{21.799cm}}
\pgfpathcurveto{\pgfqpoint{16.301cm}{21.877cm}}{\pgfqpoint{16.364cm}{21.94cm}}{\pgfqpoint{16.442cm}{21.94cm}}
\pgfpathcurveto{\pgfqpoint{16.52cm}{21.94cm}}{\pgfqpoint{16.583cm}{21.877cm}}{\pgfqpoint{16.583cm}{21.799cm}}
\pgfusepath{stroke}
\pgfsetdash{}{0cm}
\pgfpathmoveto{\pgfqpoint{17.762cm}{21.799cm}}
\pgfpathcurveto{\pgfqpoint{17.762cm}{21.721cm}}{\pgfqpoint{17.699cm}{21.658cm}}{\pgfqpoint{17.621cm}{21.658cm}}
\pgfpathcurveto{\pgfqpoint{17.543cm}{21.658cm}}{\pgfqpoint{17.48cm}{21.721cm}}{\pgfqpoint{17.48cm}{21.799cm}}
\pgfpathcurveto{\pgfqpoint{17.48cm}{21.877cm}}{\pgfqpoint{17.543cm}{21.94cm}}{\pgfqpoint{17.621cm}{21.94cm}}
\pgfpathcurveto{\pgfqpoint{17.699cm}{21.94cm}}{\pgfqpoint{17.762cm}{21.877cm}}{\pgfqpoint{17.762cm}{21.799cm}}
\pgfusepath{stroke}
\begin{pgfscope}
\pgfpathmoveto{\pgfqpoint{5.847cm}{22.969cm}}
\pgfpathlineto{\pgfqpoint{17.624cm}{22.969cm}}
\pgfpathlineto{\pgfqpoint{17.624cm}{5.736cm}}
\pgfpathlineto{\pgfqpoint{5.847cm}{5.736cm}}
\pgfpathclose
\pgfusepath{clip}
\pgfsetdash{}{0cm}
\definecolor{eps2pgf_color}{rgb}{0,1,0}\pgfsetstrokecolor{eps2pgf_color}\pgfsetfillcolor{eps2pgf_color}
\pgfpathmoveto{\pgfqpoint{7.023cm}{21.202cm}}
\pgfpathlineto{\pgfqpoint{8.202cm}{20.617cm}}
\pgfpathlineto{\pgfqpoint{9.378cm}{20.294cm}}
\pgfpathlineto{\pgfqpoint{10.557cm}{20.393cm}}
\pgfpathlineto{\pgfqpoint{11.733cm}{19.897cm}}
\pgfpathlineto{\pgfqpoint{12.912cm}{19.556cm}}
\pgfpathlineto{\pgfqpoint{14.088cm}{19.362cm}}
\pgfpathlineto{\pgfqpoint{15.266cm}{19.238cm}}
\pgfpathlineto{\pgfqpoint{16.442cm}{19.112cm}}
\pgfpathlineto{\pgfqpoint{17.621cm}{18.965cm}}
\pgfusepath{stroke}
\end{pgfscope}
\pgfsetdash{}{0cm}
\definecolor{eps2pgf_color}{rgb}{0,1,0}\pgfsetstrokecolor{eps2pgf_color}\pgfsetfillcolor{eps2pgf_color}
\pgfpathmoveto{\pgfqpoint{6.882cm}{21.202cm}}
\pgfpathlineto{\pgfqpoint{7.164cm}{21.202cm}}
\pgfusepath{stroke}
\pgfsetdash{}{0cm}
\pgfpathmoveto{\pgfqpoint{7.023cm}{21.343cm}}
\pgfpathlineto{\pgfqpoint{7.023cm}{21.061cm}}
\pgfusepath{stroke}
\pgfsetdash{}{0cm}
\pgfpathmoveto{\pgfqpoint{8.061cm}{20.617cm}}
\pgfpathlineto{\pgfqpoint{8.343cm}{20.617cm}}
\pgfusepath{stroke}
\pgfsetdash{}{0cm}
\pgfpathmoveto{\pgfqpoint{8.202cm}{20.758cm}}
\pgfpathlineto{\pgfqpoint{8.202cm}{20.476cm}}
\pgfusepath{stroke}
\pgfsetdash{}{0cm}
\pgfpathmoveto{\pgfqpoint{9.237cm}{20.294cm}}
\pgfpathlineto{\pgfqpoint{9.519cm}{20.294cm}}
\pgfusepath{stroke}
\pgfsetdash{}{0cm}
\pgfpathmoveto{\pgfqpoint{9.378cm}{20.435cm}}
\pgfpathlineto{\pgfqpoint{9.378cm}{20.152cm}}
\pgfusepath{stroke}
\pgfsetdash{}{0cm}
\pgfpathmoveto{\pgfqpoint{10.416cm}{20.393cm}}
\pgfpathlineto{\pgfqpoint{10.698cm}{20.393cm}}
\pgfusepath{stroke}
\pgfsetdash{}{0cm}
\pgfpathmoveto{\pgfqpoint{10.557cm}{20.535cm}}
\pgfpathlineto{\pgfqpoint{10.557cm}{20.252cm}}
\pgfusepath{stroke}
\pgfsetdash{}{0cm}
\pgfpathmoveto{\pgfqpoint{11.592cm}{19.897cm}}
\pgfpathlineto{\pgfqpoint{11.874cm}{19.897cm}}
\pgfusepath{stroke}
\pgfsetdash{}{0cm}
\pgfpathmoveto{\pgfqpoint{11.733cm}{20.038cm}}
\pgfpathlineto{\pgfqpoint{11.733cm}{19.756cm}}
\pgfusepath{stroke}
\pgfsetdash{}{0cm}
\pgfpathmoveto{\pgfqpoint{12.771cm}{19.556cm}}
\pgfpathlineto{\pgfqpoint{13.053cm}{19.556cm}}
\pgfusepath{stroke}
\pgfsetdash{}{0cm}
\pgfpathmoveto{\pgfqpoint{12.912cm}{19.697cm}}
\pgfpathlineto{\pgfqpoint{12.912cm}{19.415cm}}
\pgfusepath{stroke}
\pgfsetdash{}{0cm}
\pgfpathmoveto{\pgfqpoint{13.946cm}{19.362cm}}
\pgfpathlineto{\pgfqpoint{14.229cm}{19.362cm}}
\pgfusepath{stroke}
\pgfsetdash{}{0cm}
\pgfpathmoveto{\pgfqpoint{14.088cm}{19.503cm}}
\pgfpathlineto{\pgfqpoint{14.088cm}{19.221cm}}
\pgfusepath{stroke}
\pgfsetdash{}{0cm}
\pgfpathmoveto{\pgfqpoint{15.125cm}{19.238cm}}
\pgfpathlineto{\pgfqpoint{15.408cm}{19.238cm}}
\pgfusepath{stroke}
\pgfsetdash{}{0cm}
\pgfpathmoveto{\pgfqpoint{15.266cm}{19.379cm}}
\pgfpathlineto{\pgfqpoint{15.266cm}{19.097cm}}
\pgfusepath{stroke}
\pgfsetdash{}{0cm}
\pgfpathmoveto{\pgfqpoint{16.301cm}{19.112cm}}
\pgfpathlineto{\pgfqpoint{16.583cm}{19.112cm}}
\pgfusepath{stroke}
\pgfsetdash{}{0cm}
\pgfpathmoveto{\pgfqpoint{16.442cm}{19.253cm}}
\pgfpathlineto{\pgfqpoint{16.442cm}{18.971cm}}
\pgfusepath{stroke}
\pgfsetdash{}{0cm}
\pgfpathmoveto{\pgfqpoint{17.48cm}{18.965cm}}
\pgfpathlineto{\pgfqpoint{17.762cm}{18.965cm}}
\pgfusepath{stroke}
\pgfsetdash{}{0cm}
\pgfpathmoveto{\pgfqpoint{17.621cm}{19.106cm}}
\pgfpathlineto{\pgfqpoint{17.621cm}{18.824cm}}
\pgfusepath{stroke}
\pgfsetdash{}{0cm}
\pgfpathmoveto{\pgfqpoint{6.926cm}{21.299cm}}
\pgfpathlineto{\pgfqpoint{7.12cm}{21.105cm}}
\pgfusepath{stroke}
\pgfsetdash{}{0cm}
\pgfpathmoveto{\pgfqpoint{7.12cm}{21.299cm}}
\pgfpathlineto{\pgfqpoint{6.926cm}{21.105cm}}
\pgfusepath{stroke}
\pgfsetdash{}{0cm}
\pgfpathmoveto{\pgfqpoint{8.105cm}{20.714cm}}
\pgfpathlineto{\pgfqpoint{8.299cm}{20.52cm}}
\pgfusepath{stroke}
\pgfsetdash{}{0cm}
\pgfpathmoveto{\pgfqpoint{8.299cm}{20.714cm}}
\pgfpathlineto{\pgfqpoint{8.105cm}{20.52cm}}
\pgfusepath{stroke}
\pgfsetdash{}{0cm}
\pgfpathmoveto{\pgfqpoint{9.281cm}{20.391cm}}
\pgfpathlineto{\pgfqpoint{9.475cm}{20.197cm}}
\pgfusepath{stroke}
\pgfsetdash{}{0cm}
\pgfpathmoveto{\pgfqpoint{9.475cm}{20.391cm}}
\pgfpathlineto{\pgfqpoint{9.281cm}{20.197cm}}
\pgfusepath{stroke}
\pgfsetdash{}{0cm}
\pgfpathmoveto{\pgfqpoint{10.46cm}{20.491cm}}
\pgfpathlineto{\pgfqpoint{10.654cm}{20.296cm}}
\pgfusepath{stroke}
\pgfsetdash{}{0cm}
\pgfpathmoveto{\pgfqpoint{10.654cm}{20.491cm}}
\pgfpathlineto{\pgfqpoint{10.46cm}{20.296cm}}
\pgfusepath{stroke}
\pgfsetdash{}{0cm}
\pgfpathmoveto{\pgfqpoint{11.636cm}{19.994cm}}
\pgfpathlineto{\pgfqpoint{11.83cm}{19.8cm}}
\pgfusepath{stroke}
\pgfsetdash{}{0cm}
\pgfpathmoveto{\pgfqpoint{11.83cm}{19.994cm}}
\pgfpathlineto{\pgfqpoint{11.636cm}{19.8cm}}
\pgfusepath{stroke}
\pgfsetdash{}{0cm}
\pgfpathmoveto{\pgfqpoint{12.815cm}{19.653cm}}
\pgfpathlineto{\pgfqpoint{13.009cm}{19.459cm}}
\pgfusepath{stroke}
\pgfsetdash{}{0cm}
\pgfpathmoveto{\pgfqpoint{13.009cm}{19.653cm}}
\pgfpathlineto{\pgfqpoint{12.815cm}{19.459cm}}
\pgfusepath{stroke}
\pgfsetdash{}{0cm}
\pgfpathmoveto{\pgfqpoint{13.991cm}{19.459cm}}
\pgfpathlineto{\pgfqpoint{14.185cm}{19.265cm}}
\pgfusepath{stroke}
\pgfsetdash{}{0cm}
\pgfpathmoveto{\pgfqpoint{14.185cm}{19.459cm}}
\pgfpathlineto{\pgfqpoint{13.991cm}{19.265cm}}
\pgfusepath{stroke}
\pgfsetdash{}{0cm}
\pgfpathmoveto{\pgfqpoint{15.169cm}{19.335cm}}
\pgfpathlineto{\pgfqpoint{15.363cm}{19.141cm}}
\pgfusepath{stroke}
\pgfsetdash{}{0cm}
\pgfpathmoveto{\pgfqpoint{15.363cm}{19.335cm}}
\pgfpathlineto{\pgfqpoint{15.169cm}{19.141cm}}
\pgfusepath{stroke}
\pgfsetdash{}{0cm}
\pgfpathmoveto{\pgfqpoint{16.345cm}{19.209cm}}
\pgfpathlineto{\pgfqpoint{16.539cm}{19.015cm}}
\pgfusepath{stroke}
\pgfsetdash{}{0cm}
\pgfpathmoveto{\pgfqpoint{16.539cm}{19.209cm}}
\pgfpathlineto{\pgfqpoint{16.345cm}{19.015cm}}
\pgfusepath{stroke}
\pgfsetdash{}{0cm}
\pgfpathmoveto{\pgfqpoint{17.524cm}{19.062cm}}
\pgfpathlineto{\pgfqpoint{17.718cm}{18.868cm}}
\pgfusepath{stroke}
\pgfsetdash{}{0cm}
\pgfpathmoveto{\pgfqpoint{17.718cm}{19.062cm}}
\pgfpathlineto{\pgfqpoint{17.524cm}{18.868cm}}
\pgfusepath{stroke}
\begin{pgfscope}
\pgfpathmoveto{\pgfqpoint{5.847cm}{22.969cm}}
\pgfpathlineto{\pgfqpoint{17.624cm}{22.969cm}}
\pgfpathlineto{\pgfqpoint{17.624cm}{5.736cm}}
\pgfpathlineto{\pgfqpoint{5.847cm}{5.736cm}}
\pgfpathclose
\pgfusepath{clip}
\pgfsetdash{}{0cm}
\definecolor{eps2pgf_color}{rgb}{0,0,1}\pgfsetstrokecolor{eps2pgf_color}\pgfsetfillcolor{eps2pgf_color}
\pgfpathmoveto{\pgfqpoint{7.023cm}{20.826cm}}
\pgfpathlineto{\pgfqpoint{8.202cm}{19.888cm}}
\pgfpathlineto{\pgfqpoint{9.378cm}{18.985cm}}
\pgfpathlineto{\pgfqpoint{10.557cm}{18.139cm}}
\pgfpathlineto{\pgfqpoint{11.733cm}{17.313cm}}
\pgfpathlineto{\pgfqpoint{12.912cm}{16.519cm}}
\pgfpathlineto{\pgfqpoint{14.088cm}{15.775cm}}
\pgfpathlineto{\pgfqpoint{15.266cm}{15.058cm}}
\pgfpathlineto{\pgfqpoint{16.442cm}{14.367cm}}
\pgfpathlineto{\pgfqpoint{17.621cm}{13.72cm}}
\pgfusepath{stroke}
\end{pgfscope}
\pgfsetdash{}{0cm}
\pgfsetmiterjoin
\definecolor{eps2pgf_color}{rgb}{0,0,1}\pgfsetstrokecolor{eps2pgf_color}\pgfsetfillcolor{eps2pgf_color}
\pgfpathmoveto{\pgfqpoint{6.912cm}{20.937cm}}
\pgfpathlineto{\pgfqpoint{7.135cm}{20.937cm}}
\pgfpathlineto{\pgfqpoint{7.135cm}{20.714cm}}
\pgfpathlineto{\pgfqpoint{6.912cm}{20.714cm}}
\pgfpathlineto{\pgfqpoint{6.912cm}{20.937cm}}
\pgfpathclose
\pgfusepath{stroke}
\pgfsetdash{}{0cm}
\pgfpathmoveto{\pgfqpoint{8.09cm}{20cm}}
\pgfpathlineto{\pgfqpoint{8.314cm}{20cm}}
\pgfpathlineto{\pgfqpoint{8.314cm}{19.776cm}}
\pgfpathlineto{\pgfqpoint{8.09cm}{19.776cm}}
\pgfpathlineto{\pgfqpoint{8.09cm}{20cm}}
\pgfpathclose
\pgfusepath{stroke}
\pgfsetdash{}{0cm}
\pgfpathmoveto{\pgfqpoint{9.266cm}{19.097cm}}
\pgfpathlineto{\pgfqpoint{9.49cm}{19.097cm}}
\pgfpathlineto{\pgfqpoint{9.49cm}{18.874cm}}
\pgfpathlineto{\pgfqpoint{9.266cm}{18.874cm}}
\pgfpathlineto{\pgfqpoint{9.266cm}{19.097cm}}
\pgfpathclose
\pgfusepath{stroke}
\pgfsetdash{}{0cm}
\pgfpathmoveto{\pgfqpoint{10.445cm}{18.25cm}}
\pgfpathlineto{\pgfqpoint{10.669cm}{18.25cm}}
\pgfpathlineto{\pgfqpoint{10.669cm}{18.027cm}}
\pgfpathlineto{\pgfqpoint{10.445cm}{18.027cm}}
\pgfpathlineto{\pgfqpoint{10.445cm}{18.25cm}}
\pgfpathclose
\pgfusepath{stroke}
\pgfsetdash{}{0cm}
\pgfpathmoveto{\pgfqpoint{11.621cm}{17.424cm}}
\pgfpathlineto{\pgfqpoint{11.845cm}{17.424cm}}
\pgfpathlineto{\pgfqpoint{11.845cm}{17.201cm}}
\pgfpathlineto{\pgfqpoint{11.621cm}{17.201cm}}
\pgfpathlineto{\pgfqpoint{11.621cm}{17.424cm}}
\pgfpathclose
\pgfusepath{stroke}
\pgfsetdash{}{0cm}
\pgfpathmoveto{\pgfqpoint{12.8cm}{16.631cm}}
\pgfpathlineto{\pgfqpoint{13.023cm}{16.631cm}}
\pgfpathlineto{\pgfqpoint{13.023cm}{16.407cm}}
\pgfpathlineto{\pgfqpoint{12.8cm}{16.407cm}}
\pgfpathlineto{\pgfqpoint{12.8cm}{16.631cm}}
\pgfpathclose
\pgfusepath{stroke}
\pgfsetdash{}{0cm}
\pgfpathmoveto{\pgfqpoint{13.976cm}{15.887cm}}
\pgfpathlineto{\pgfqpoint{14.199cm}{15.887cm}}
\pgfpathlineto{\pgfqpoint{14.199cm}{15.663cm}}
\pgfpathlineto{\pgfqpoint{13.976cm}{15.663cm}}
\pgfpathlineto{\pgfqpoint{13.976cm}{15.887cm}}
\pgfpathclose
\pgfusepath{stroke}
\pgfsetdash{}{0cm}
\pgfpathmoveto{\pgfqpoint{15.155cm}{15.169cm}}
\pgfpathlineto{\pgfqpoint{15.378cm}{15.169cm}}
\pgfpathlineto{\pgfqpoint{15.378cm}{14.946cm}}
\pgfpathlineto{\pgfqpoint{15.155cm}{14.946cm}}
\pgfpathlineto{\pgfqpoint{15.155cm}{15.169cm}}
\pgfpathclose
\pgfusepath{stroke}
\pgfsetdash{}{0cm}
\pgfpathmoveto{\pgfqpoint{16.331cm}{14.479cm}}
\pgfpathlineto{\pgfqpoint{16.554cm}{14.479cm}}
\pgfpathlineto{\pgfqpoint{16.554cm}{14.255cm}}
\pgfpathlineto{\pgfqpoint{16.331cm}{14.255cm}}
\pgfpathlineto{\pgfqpoint{16.331cm}{14.479cm}}
\pgfpathclose
\pgfusepath{stroke}
\pgfsetdash{}{0cm}
\pgfpathmoveto{\pgfqpoint{17.51cm}{13.832cm}}
\pgfpathlineto{\pgfqpoint{17.733cm}{13.832cm}}
\pgfpathlineto{\pgfqpoint{17.733cm}{13.608cm}}
\pgfpathlineto{\pgfqpoint{17.51cm}{13.608cm}}
\pgfpathlineto{\pgfqpoint{17.51cm}{13.832cm}}
\pgfpathclose
\pgfusepath{stroke}
\begin{pgfscope}
\pgfpathmoveto{\pgfqpoint{5.847cm}{22.969cm}}
\pgfpathlineto{\pgfqpoint{17.624cm}{22.969cm}}
\pgfpathlineto{\pgfqpoint{17.624cm}{5.736cm}}
\pgfpathlineto{\pgfqpoint{5.847cm}{5.736cm}}
\pgfpathclose
\pgfusepath{clip}
\pgfsetdash{}{0cm}
\definecolor{eps2pgf_color}{rgb}{0,1,1}\pgfsetstrokecolor{eps2pgf_color}\pgfsetfillcolor{eps2pgf_color}
\pgfpathmoveto{\pgfqpoint{7.023cm}{20.787cm}}
\pgfpathlineto{\pgfqpoint{8.202cm}{19.82cm}}
\pgfpathlineto{\pgfqpoint{9.378cm}{18.891cm}}
\pgfpathlineto{\pgfqpoint{10.557cm}{18.003cm}}
\pgfpathlineto{\pgfqpoint{11.733cm}{17.133cm}}
\pgfpathlineto{\pgfqpoint{12.912cm}{16.348cm}}
\pgfpathlineto{\pgfqpoint{14.088cm}{15.566cm}}
\pgfpathlineto{\pgfqpoint{15.266cm}{14.793cm}}
\pgfpathlineto{\pgfqpoint{16.442cm}{14.073cm}}
\pgfpathlineto{\pgfqpoint{17.621cm}{13.385cm}}
\pgfusepath{stroke}
\end{pgfscope}
\pgfsetdash{}{0cm}
\definecolor{eps2pgf_color}{rgb}{0,1,1}\pgfsetstrokecolor{eps2pgf_color}\pgfsetfillcolor{eps2pgf_color}
\pgfpathmoveto{\pgfqpoint{6.926cm}{20.884cm}}
\pgfpathlineto{\pgfqpoint{7.12cm}{20.69cm}}
\pgfusepath{stroke}
\pgfsetdash{}{0cm}
\pgfpathmoveto{\pgfqpoint{7.12cm}{20.884cm}}
\pgfpathlineto{\pgfqpoint{6.926cm}{20.69cm}}
\pgfusepath{stroke}
\pgfsetdash{}{0cm}
\pgfpathmoveto{\pgfqpoint{8.105cm}{19.917cm}}
\pgfpathlineto{\pgfqpoint{8.299cm}{19.723cm}}
\pgfusepath{stroke}
\pgfsetdash{}{0cm}
\pgfpathmoveto{\pgfqpoint{8.299cm}{19.917cm}}
\pgfpathlineto{\pgfqpoint{8.105cm}{19.723cm}}
\pgfusepath{stroke}
\pgfsetdash{}{0cm}
\pgfpathmoveto{\pgfqpoint{9.281cm}{18.988cm}}
\pgfpathlineto{\pgfqpoint{9.475cm}{18.794cm}}
\pgfusepath{stroke}
\pgfsetdash{}{0cm}
\pgfpathmoveto{\pgfqpoint{9.475cm}{18.988cm}}
\pgfpathlineto{\pgfqpoint{9.281cm}{18.794cm}}
\pgfusepath{stroke}
\pgfsetdash{}{0cm}
\pgfpathmoveto{\pgfqpoint{10.46cm}{18.1cm}}
\pgfpathlineto{\pgfqpoint{10.654cm}{17.906cm}}
\pgfusepath{stroke}
\pgfsetdash{}{0cm}
\pgfpathmoveto{\pgfqpoint{10.654cm}{18.1cm}}
\pgfpathlineto{\pgfqpoint{10.46cm}{17.906cm}}
\pgfusepath{stroke}
\pgfsetdash{}{0cm}
\pgfpathmoveto{\pgfqpoint{11.636cm}{17.23cm}}
\pgfpathlineto{\pgfqpoint{11.83cm}{17.036cm}}
\pgfusepath{stroke}
\pgfsetdash{}{0cm}
\pgfpathmoveto{\pgfqpoint{11.83cm}{17.23cm}}
\pgfpathlineto{\pgfqpoint{11.636cm}{17.036cm}}
\pgfusepath{stroke}
\pgfsetdash{}{0cm}
\pgfpathmoveto{\pgfqpoint{12.815cm}{16.445cm}}
\pgfpathlineto{\pgfqpoint{13.009cm}{16.251cm}}
\pgfusepath{stroke}
\pgfsetdash{}{0cm}
\pgfpathmoveto{\pgfqpoint{13.009cm}{16.445cm}}
\pgfpathlineto{\pgfqpoint{12.815cm}{16.251cm}}
\pgfusepath{stroke}
\pgfsetdash{}{0cm}
\pgfpathmoveto{\pgfqpoint{13.991cm}{15.663cm}}
\pgfpathlineto{\pgfqpoint{14.185cm}{15.469cm}}
\pgfusepath{stroke}
\pgfsetdash{}{0cm}
\pgfpathmoveto{\pgfqpoint{14.185cm}{15.663cm}}
\pgfpathlineto{\pgfqpoint{13.991cm}{15.469cm}}
\pgfusepath{stroke}
\pgfsetdash{}{0cm}
\pgfpathmoveto{\pgfqpoint{15.169cm}{14.89cm}}
\pgfpathlineto{\pgfqpoint{15.363cm}{14.696cm}}
\pgfusepath{stroke}
\pgfsetdash{}{0cm}
\pgfpathmoveto{\pgfqpoint{15.363cm}{14.89cm}}
\pgfpathlineto{\pgfqpoint{15.169cm}{14.696cm}}
\pgfusepath{stroke}
\pgfsetdash{}{0cm}
\pgfpathmoveto{\pgfqpoint{16.345cm}{14.17cm}}
\pgfpathlineto{\pgfqpoint{16.539cm}{13.976cm}}
\pgfusepath{stroke}
\pgfsetdash{}{0cm}
\pgfpathmoveto{\pgfqpoint{16.539cm}{14.17cm}}
\pgfpathlineto{\pgfqpoint{16.345cm}{13.976cm}}
\pgfusepath{stroke}
\pgfsetdash{}{0cm}
\pgfpathmoveto{\pgfqpoint{17.524cm}{13.482cm}}
\pgfpathlineto{\pgfqpoint{17.718cm}{13.288cm}}
\pgfusepath{stroke}
\pgfsetdash{}{0cm}
\pgfpathmoveto{\pgfqpoint{17.718cm}{13.482cm}}
\pgfpathlineto{\pgfqpoint{17.524cm}{13.288cm}}
\pgfusepath{stroke}
\begin{pgfscope}
\pgfpathmoveto{\pgfqpoint{5.847cm}{22.969cm}}
\pgfpathlineto{\pgfqpoint{17.624cm}{22.969cm}}
\pgfpathlineto{\pgfqpoint{17.624cm}{5.736cm}}
\pgfpathlineto{\pgfqpoint{5.847cm}{5.736cm}}
\pgfpathclose
\pgfusepath{clip}
\pgfsetdash{{0.212cm}}{0cm}
\definecolor{eps2pgf_color}{rgb}{0,0,1}\pgfsetstrokecolor{eps2pgf_color}\pgfsetfillcolor{eps2pgf_color}
\pgfpathmoveto{\pgfqpoint{7.023cm}{20.467cm}}
\pgfpathlineto{\pgfqpoint{8.202cm}{19.173cm}}
\pgfpathlineto{\pgfqpoint{9.378cm}{17.912cm}}
\pgfpathlineto{\pgfqpoint{10.557cm}{16.71cm}}
\pgfpathlineto{\pgfqpoint{11.733cm}{15.522cm}}
\pgfpathlineto{\pgfqpoint{12.912cm}{14.37cm}}
\pgfpathlineto{\pgfqpoint{14.088cm}{13.27cm}}
\pgfpathlineto{\pgfqpoint{15.266cm}{12.194cm}}
\pgfpathlineto{\pgfqpoint{16.442cm}{11.148cm}}
\pgfpathlineto{\pgfqpoint{17.621cm}{10.139cm}}
\pgfusepath{stroke}
\end{pgfscope}
\pgfsetdash{}{0cm}
\definecolor{eps2pgf_color}{rgb}{0,0,1}\pgfsetstrokecolor{eps2pgf_color}\pgfsetfillcolor{eps2pgf_color}
\pgfpathmoveto{\pgfqpoint{6.912cm}{20.579cm}}
\pgfpathlineto{\pgfqpoint{7.135cm}{20.579cm}}
\pgfpathlineto{\pgfqpoint{7.135cm}{20.355cm}}
\pgfpathlineto{\pgfqpoint{6.912cm}{20.355cm}}
\pgfpathlineto{\pgfqpoint{6.912cm}{20.579cm}}
\pgfpathclose
\pgfusepath{stroke}
\pgfsetdash{}{0cm}
\pgfpathmoveto{\pgfqpoint{8.09cm}{19.285cm}}
\pgfpathlineto{\pgfqpoint{8.314cm}{19.285cm}}
\pgfpathlineto{\pgfqpoint{8.314cm}{19.062cm}}
\pgfpathlineto{\pgfqpoint{8.09cm}{19.062cm}}
\pgfpathlineto{\pgfqpoint{8.09cm}{19.285cm}}
\pgfpathclose
\pgfusepath{stroke}
\pgfsetdash{}{0cm}
\pgfpathmoveto{\pgfqpoint{9.266cm}{18.024cm}}
\pgfpathlineto{\pgfqpoint{9.49cm}{18.024cm}}
\pgfpathlineto{\pgfqpoint{9.49cm}{17.801cm}}
\pgfpathlineto{\pgfqpoint{9.266cm}{17.801cm}}
\pgfpathlineto{\pgfqpoint{9.266cm}{18.024cm}}
\pgfpathclose
\pgfusepath{stroke}
\pgfsetdash{}{0cm}
\pgfpathmoveto{\pgfqpoint{10.445cm}{16.822cm}}
\pgfpathlineto{\pgfqpoint{10.669cm}{16.822cm}}
\pgfpathlineto{\pgfqpoint{10.669cm}{16.598cm}}
\pgfpathlineto{\pgfqpoint{10.445cm}{16.598cm}}
\pgfpathlineto{\pgfqpoint{10.445cm}{16.822cm}}
\pgfpathclose
\pgfusepath{stroke}
\pgfsetdash{}{0cm}
\pgfpathmoveto{\pgfqpoint{11.621cm}{15.634cm}}
\pgfpathlineto{\pgfqpoint{11.845cm}{15.634cm}}
\pgfpathlineto{\pgfqpoint{11.845cm}{15.411cm}}
\pgfpathlineto{\pgfqpoint{11.621cm}{15.411cm}}
\pgfpathlineto{\pgfqpoint{11.621cm}{15.634cm}}
\pgfpathclose
\pgfusepath{stroke}
\pgfsetdash{}{0cm}
\pgfpathmoveto{\pgfqpoint{12.8cm}{14.482cm}}
\pgfpathlineto{\pgfqpoint{13.023cm}{14.482cm}}
\pgfpathlineto{\pgfqpoint{13.023cm}{14.258cm}}
\pgfpathlineto{\pgfqpoint{12.8cm}{14.258cm}}
\pgfpathlineto{\pgfqpoint{12.8cm}{14.482cm}}
\pgfpathclose
\pgfusepath{stroke}
\pgfsetdash{}{0cm}
\pgfpathmoveto{\pgfqpoint{13.976cm}{13.382cm}}
\pgfpathlineto{\pgfqpoint{14.199cm}{13.382cm}}
\pgfpathlineto{\pgfqpoint{14.199cm}{13.159cm}}
\pgfpathlineto{\pgfqpoint{13.976cm}{13.159cm}}
\pgfpathlineto{\pgfqpoint{13.976cm}{13.382cm}}
\pgfpathclose
\pgfusepath{stroke}
\pgfsetdash{}{0cm}
\pgfpathmoveto{\pgfqpoint{15.155cm}{12.306cm}}
\pgfpathlineto{\pgfqpoint{15.378cm}{12.306cm}}
\pgfpathlineto{\pgfqpoint{15.378cm}{12.083cm}}
\pgfpathlineto{\pgfqpoint{15.155cm}{12.083cm}}
\pgfpathlineto{\pgfqpoint{15.155cm}{12.306cm}}
\pgfpathclose
\pgfusepath{stroke}
\pgfsetdash{}{0cm}
\pgfpathmoveto{\pgfqpoint{16.331cm}{11.259cm}}
\pgfpathlineto{\pgfqpoint{16.554cm}{11.259cm}}
\pgfpathlineto{\pgfqpoint{16.554cm}{11.036cm}}
\pgfpathlineto{\pgfqpoint{16.331cm}{11.036cm}}
\pgfpathlineto{\pgfqpoint{16.331cm}{11.259cm}}
\pgfpathclose
\pgfusepath{stroke}
\pgfsetdash{}{0cm}
\pgfpathmoveto{\pgfqpoint{17.51cm}{10.251cm}}
\pgfpathlineto{\pgfqpoint{17.733cm}{10.251cm}}
\pgfpathlineto{\pgfqpoint{17.733cm}{10.028cm}}
\pgfpathlineto{\pgfqpoint{17.51cm}{10.028cm}}
\pgfpathlineto{\pgfqpoint{17.51cm}{10.251cm}}
\pgfpathclose
\pgfusepath{stroke}
\begin{pgfscope}
\pgfpathmoveto{\pgfqpoint{5.847cm}{22.969cm}}
\pgfpathlineto{\pgfqpoint{17.624cm}{22.969cm}}
\pgfpathlineto{\pgfqpoint{17.624cm}{5.736cm}}
\pgfpathlineto{\pgfqpoint{5.847cm}{5.736cm}}
\pgfpathclose
\pgfusepath{clip}
\pgfsetdash{{0.212cm}}{0cm}
\definecolor{eps2pgf_color}{rgb}{0,1,1}\pgfsetstrokecolor{eps2pgf_color}\pgfsetfillcolor{eps2pgf_color}
\pgfpathmoveto{\pgfqpoint{7.023cm}{20.423cm}}
\pgfpathlineto{\pgfqpoint{8.202cm}{19.091cm}}
\pgfpathlineto{\pgfqpoint{9.378cm}{17.798cm}}
\pgfpathlineto{\pgfqpoint{10.557cm}{16.545cm}}
\pgfpathlineto{\pgfqpoint{11.733cm}{15.313cm}}
\pgfpathlineto{\pgfqpoint{12.912cm}{14.161cm}}
\pgfpathlineto{\pgfqpoint{14.088cm}{13.015cm}}
\pgfpathlineto{\pgfqpoint{15.266cm}{11.877cm}}
\pgfpathlineto{\pgfqpoint{16.442cm}{10.792cm}}
\pgfpathlineto{\pgfqpoint{17.621cm}{9.74cm}}
\pgfusepath{stroke}
\end{pgfscope}
\pgfsetdash{}{0cm}
\definecolor{eps2pgf_color}{rgb}{0,1,1}\pgfsetstrokecolor{eps2pgf_color}\pgfsetfillcolor{eps2pgf_color}
\pgfpathmoveto{\pgfqpoint{6.926cm}{20.52cm}}
\pgfpathlineto{\pgfqpoint{7.12cm}{20.326cm}}
\pgfusepath{stroke}
\pgfsetdash{}{0cm}
\pgfpathmoveto{\pgfqpoint{7.12cm}{20.52cm}}
\pgfpathlineto{\pgfqpoint{6.926cm}{20.326cm}}
\pgfusepath{stroke}
\pgfsetdash{}{0cm}
\pgfpathmoveto{\pgfqpoint{8.105cm}{19.188cm}}
\pgfpathlineto{\pgfqpoint{8.299cm}{18.994cm}}
\pgfusepath{stroke}
\pgfsetdash{}{0cm}
\pgfpathmoveto{\pgfqpoint{8.299cm}{19.188cm}}
\pgfpathlineto{\pgfqpoint{8.105cm}{18.994cm}}
\pgfusepath{stroke}
\pgfsetdash{}{0cm}
\pgfpathmoveto{\pgfqpoint{9.281cm}{17.895cm}}
\pgfpathlineto{\pgfqpoint{9.475cm}{17.701cm}}
\pgfusepath{stroke}
\pgfsetdash{}{0cm}
\pgfpathmoveto{\pgfqpoint{9.475cm}{17.895cm}}
\pgfpathlineto{\pgfqpoint{9.281cm}{17.701cm}}
\pgfusepath{stroke}
\pgfsetdash{}{0cm}
\pgfpathmoveto{\pgfqpoint{10.46cm}{16.642cm}}
\pgfpathlineto{\pgfqpoint{10.654cm}{16.448cm}}
\pgfusepath{stroke}
\pgfsetdash{}{0cm}
\pgfpathmoveto{\pgfqpoint{10.654cm}{16.642cm}}
\pgfpathlineto{\pgfqpoint{10.46cm}{16.448cm}}
\pgfusepath{stroke}
\pgfsetdash{}{0cm}
\pgfpathmoveto{\pgfqpoint{11.636cm}{15.411cm}}
\pgfpathlineto{\pgfqpoint{11.83cm}{15.216cm}}
\pgfusepath{stroke}
\pgfsetdash{}{0cm}
\pgfpathmoveto{\pgfqpoint{11.83cm}{15.411cm}}
\pgfpathlineto{\pgfqpoint{11.636cm}{15.216cm}}
\pgfusepath{stroke}
\pgfsetdash{}{0cm}
\pgfpathmoveto{\pgfqpoint{12.815cm}{14.258cm}}
\pgfpathlineto{\pgfqpoint{13.009cm}{14.064cm}}
\pgfusepath{stroke}
\pgfsetdash{}{0cm}
\pgfpathmoveto{\pgfqpoint{13.009cm}{14.258cm}}
\pgfpathlineto{\pgfqpoint{12.815cm}{14.064cm}}
\pgfusepath{stroke}
\pgfsetdash{}{0cm}
\pgfpathmoveto{\pgfqpoint{13.991cm}{13.112cm}}
\pgfpathlineto{\pgfqpoint{14.185cm}{12.918cm}}
\pgfusepath{stroke}
\pgfsetdash{}{0cm}
\pgfpathmoveto{\pgfqpoint{14.185cm}{13.112cm}}
\pgfpathlineto{\pgfqpoint{13.991cm}{12.918cm}}
\pgfusepath{stroke}
\pgfsetdash{}{0cm}
\pgfpathmoveto{\pgfqpoint{15.169cm}{11.974cm}}
\pgfpathlineto{\pgfqpoint{15.363cm}{11.78cm}}
\pgfusepath{stroke}
\pgfsetdash{}{0cm}
\pgfpathmoveto{\pgfqpoint{15.363cm}{11.974cm}}
\pgfpathlineto{\pgfqpoint{15.169cm}{11.78cm}}
\pgfusepath{stroke}
\pgfsetdash{}{0cm}
\pgfpathmoveto{\pgfqpoint{16.345cm}{10.889cm}}
\pgfpathlineto{\pgfqpoint{16.539cm}{10.695cm}}
\pgfusepath{stroke}
\pgfsetdash{}{0cm}
\pgfpathmoveto{\pgfqpoint{16.539cm}{10.889cm}}
\pgfpathlineto{\pgfqpoint{16.345cm}{10.695cm}}
\pgfusepath{stroke}
\pgfsetdash{}{0cm}
\pgfpathmoveto{\pgfqpoint{17.524cm}{9.837cm}}
\pgfpathlineto{\pgfqpoint{17.718cm}{9.643cm}}
\pgfusepath{stroke}
\pgfsetdash{}{0cm}
\pgfpathmoveto{\pgfqpoint{17.718cm}{9.837cm}}
\pgfpathlineto{\pgfqpoint{17.524cm}{9.643cm}}
\pgfusepath{stroke}
\begin{pgfscope}
\pgfpathmoveto{\pgfqpoint{5.847cm}{22.969cm}}
\pgfpathlineto{\pgfqpoint{17.624cm}{22.969cm}}
\pgfpathlineto{\pgfqpoint{17.624cm}{5.736cm}}
\pgfpathlineto{\pgfqpoint{5.847cm}{5.736cm}}
\pgfpathclose
\pgfusepath{clip}
\end{pgfscope}
\definecolor{eps2pgf_color}{gray}{0}\pgfsetstrokecolor{eps2pgf_color}\pgfsetfillcolor{eps2pgf_color}
\pgftext[x=11.751cm,y=4.093cm,rotate=0]{\fontsize{36}{36.14}\selectfont{ {$\bmax$}}}
\pgftext[x=3.013cm,y=14.388cm,rotate=90]{\fontsize{36}{36.14}\selectfont{ { Average cost}}}
\pgftext[x=5.797cm,y=5.612cm,rotate=0]{\fontsize{10.04}{12.04}\selectfont{ { }}}
\pgftext[x=17.574cm,y=22.845cm,rotate=0]{\fontsize{10.04}{12.04}\selectfont{ { }}}
\pgftext[x=9.924cm,y=12.733cm,rotate=0]{\fontsize{34}{36.14}\selectfont{ {No storage}}}
\begin{pgfscope}
\pgfpathmoveto{\pgfqpoint{6.024cm}{13.514cm}}
\pgfpathlineto{\pgfqpoint{17.13cm}{13.514cm}}
\pgfpathlineto{\pgfqpoint{17.13cm}{5.909cm}}
\pgfpathlineto{\pgfqpoint{6.024cm}{5.909cm}}
\pgfpathclose
\pgfusepath{clip}
\pgfsetdash{}{0cm}
\pgfpathmoveto{\pgfqpoint{6.229cm}{12.835cm}}
\pgfpathlineto{\pgfqpoint{7.27cm}{12.835cm}}
\pgfusepath{stroke}
\begin{pgfscope}
\pgfpathmoveto{\pgfqpoint{6.356cm}{13.229cm}}
\pgfpathlineto{\pgfqpoint{7.147cm}{13.229cm}}
\pgfpathlineto{\pgfqpoint{7.147cm}{12.438cm}}
\pgfpathlineto{\pgfqpoint{6.356cm}{12.438cm}}
\pgfpathclose
\pgfusepath{clip}
\pgfsetdash{}{0cm}
\pgfpathmoveto{\pgfqpoint{6.891cm}{12.835cm}}
\pgfpathcurveto{\pgfqpoint{6.891cm}{12.757cm}}{\pgfqpoint{6.828cm}{12.694cm}}{\pgfqpoint{6.75cm}{12.694cm}}
\pgfpathcurveto{\pgfqpoint{6.672cm}{12.694cm}}{\pgfqpoint{6.609cm}{12.757cm}}{\pgfqpoint{6.609cm}{12.835cm}}
\pgfpathcurveto{\pgfqpoint{6.609cm}{12.913cm}}{\pgfqpoint{6.672cm}{12.976cm}}{\pgfqpoint{6.75cm}{12.976cm}}
\pgfpathcurveto{\pgfqpoint{6.828cm}{12.976cm}}{\pgfqpoint{6.891cm}{12.913cm}}{\pgfqpoint{6.891cm}{12.835cm}}
\pgfusepath{stroke}
\end{pgfscope}
\end{pgfscope}
\pgftext[x=9.131cm,y=11.489cm,rotate=0]{\fontsize{34}{36.14}\selectfont{ {Greedy}}}
\begin{pgfscope}
\pgfpathmoveto{\pgfqpoint{6.024cm}{13.514cm}}
\pgfpathlineto{\pgfqpoint{17.13cm}{13.514cm}}
\pgfpathlineto{\pgfqpoint{17.13cm}{5.909cm}}
\pgfpathlineto{\pgfqpoint{6.024cm}{5.909cm}}
\pgfpathclose
\pgfusepath{clip}
\pgfsetdash{}{0cm}
\definecolor{eps2pgf_color}{rgb}{0,1,0}\pgfsetstrokecolor{eps2pgf_color}\pgfsetfillcolor{eps2pgf_color}
\pgfpathmoveto{\pgfqpoint{6.229cm}{11.586cm}}
\pgfpathlineto{\pgfqpoint{7.27cm}{11.586cm}}
\pgfusepath{stroke}
\begin{pgfscope}
\pgfpathmoveto{\pgfqpoint{6.356cm}{11.98cm}}
\pgfpathlineto{\pgfqpoint{7.147cm}{11.98cm}}
\pgfpathlineto{\pgfqpoint{7.147cm}{11.189cm}}
\pgfpathlineto{\pgfqpoint{6.356cm}{11.189cm}}
\pgfpathclose
\pgfusepath{clip}
\pgfsetdash{}{0cm}
\pgfpathmoveto{\pgfqpoint{6.609cm}{11.586cm}}
\pgfpathlineto{\pgfqpoint{6.891cm}{11.586cm}}
\pgfusepath{stroke}
\pgfsetdash{}{0cm}
\pgfpathmoveto{\pgfqpoint{6.75cm}{11.727cm}}
\pgfpathlineto{\pgfqpoint{6.75cm}{11.445cm}}
\pgfusepath{stroke}
\pgfsetdash{}{0cm}
\pgfpathmoveto{\pgfqpoint{6.653cm}{11.683cm}}
\pgfpathlineto{\pgfqpoint{6.847cm}{11.489cm}}
\pgfusepath{stroke}
\pgfsetdash{}{0cm}
\pgfpathmoveto{\pgfqpoint{6.847cm}{11.683cm}}
\pgfpathlineto{\pgfqpoint{6.653cm}{11.489cm}}
\pgfusepath{stroke}
\end{pgfscope}
\end{pgfscope}
\pgftext[x=10.183cm+.4cm,y=10.243cm,rotate=0]{\fontsize{34}{36.14}\selectfont{ {OMG(\texttt{minS})}}}
\begin{pgfscope}
\pgfpathmoveto{\pgfqpoint{6.024cm}{13.514cm}}
\pgfpathlineto{\pgfqpoint{17.13cm}{13.514cm}}
\pgfpathlineto{\pgfqpoint{17.13cm}{5.909cm}}
\pgfpathlineto{\pgfqpoint{6.024cm}{5.909cm}}
\pgfpathclose
\pgfusepath{clip}
\pgfsetdash{}{0cm}
\definecolor{eps2pgf_color}{rgb}{0,0,1}\pgfsetstrokecolor{eps2pgf_color}\pgfsetfillcolor{eps2pgf_color}
\pgfpathmoveto{\pgfqpoint{6.229cm}{10.339cm}}
\pgfpathlineto{\pgfqpoint{7.27cm}{10.339cm}}
\pgfusepath{stroke}
\begin{pgfscope}
\pgfpathmoveto{\pgfqpoint{6.356cm}{10.733cm}}
\pgfpathlineto{\pgfqpoint{7.147cm}{10.733cm}}
\pgfpathlineto{\pgfqpoint{7.147cm}{9.942cm}}
\pgfpathlineto{\pgfqpoint{6.356cm}{9.942cm}}
\pgfpathclose
\pgfusepath{clip}
\pgfsetdash{}{0cm}
\pgfpathmoveto{\pgfqpoint{6.638cm}{10.451cm}}
\pgfpathlineto{\pgfqpoint{6.862cm}{10.451cm}}
\pgfpathlineto{\pgfqpoint{6.862cm}{10.228cm}}
\pgfpathlineto{\pgfqpoint{6.638cm}{10.228cm}}
\pgfpathlineto{\pgfqpoint{6.638cm}{10.451cm}}
\pgfpathclose
\pgfusepath{stroke}
\end{pgfscope}
\end{pgfscope}
\pgftext[x=10.476cm+.15cm,y=8.994cm,rotate=0]{\fontsize{34}{36.14}\selectfont{ {OMG(\texttt{maxW})}}}
\begin{pgfscope}
\pgfpathmoveto{\pgfqpoint{6.024cm}{13.514cm}}
\pgfpathlineto{\pgfqpoint{17.13cm}{13.514cm}}
\pgfpathlineto{\pgfqpoint{17.13cm}{5.909cm}}
\pgfpathlineto{\pgfqpoint{6.024cm}{5.909cm}}
\pgfpathclose
\pgfusepath{clip}
\pgfsetdash{}{0cm}
\definecolor{eps2pgf_color}{rgb}{0,1,1}\pgfsetstrokecolor{eps2pgf_color}\pgfsetfillcolor{eps2pgf_color}
\pgfpathmoveto{\pgfqpoint{6.229cm}{9.09cm}}
\pgfpathlineto{\pgfqpoint{7.27cm}{9.09cm}}
\pgfusepath{stroke}
\begin{pgfscope}
\pgfpathmoveto{\pgfqpoint{6.356cm}{9.484cm}}
\pgfpathlineto{\pgfqpoint{7.147cm}{9.484cm}}
\pgfpathlineto{\pgfqpoint{7.147cm}{8.693cm}}
\pgfpathlineto{\pgfqpoint{6.356cm}{8.693cm}}
\pgfpathclose
\pgfusepath{clip}
\pgfsetdash{}{0cm}
\pgfpathmoveto{\pgfqpoint{6.653cm}{9.187cm}}
\pgfpathlineto{\pgfqpoint{6.847cm}{8.993cm}}
\pgfusepath{stroke}
\pgfsetdash{}{0cm}
\pgfpathmoveto{\pgfqpoint{6.847cm}{9.187cm}}
\pgfpathlineto{\pgfqpoint{6.653cm}{8.993cm}}
\pgfusepath{stroke}
\end{pgfscope}
\end{pgfscope}
\pgftext[x=11.909cm+.5cm,y=7.745cm,rotate=0]{\fontsize{34}{36.14}\selectfont{ {Lower bound(\texttt{minS})}}}
\begin{pgfscope}
\pgfpathmoveto{\pgfqpoint{6.024cm}{13.514cm}}
\pgfpathlineto{\pgfqpoint{17.13cm}{13.514cm}}
\pgfpathlineto{\pgfqpoint{17.13cm}{5.909cm}}
\pgfpathlineto{\pgfqpoint{6.024cm}{5.909cm}}
\pgfpathclose
\pgfusepath{clip}
\pgfsetdash{{0.212cm}}{0cm}
\definecolor{eps2pgf_color}{rgb}{0,0,1}\pgfsetstrokecolor{eps2pgf_color}\pgfsetfillcolor{eps2pgf_color}
\pgfpathmoveto{\pgfqpoint{6.229cm}{7.843cm}}
\pgfpathlineto{\pgfqpoint{7.27cm}{7.843cm}}
\pgfusepath{stroke}
\begin{pgfscope}
\pgfpathmoveto{\pgfqpoint{6.356cm}{8.237cm}}
\pgfpathlineto{\pgfqpoint{7.147cm}{8.237cm}}
\pgfpathlineto{\pgfqpoint{7.147cm}{7.447cm}}
\pgfpathlineto{\pgfqpoint{6.356cm}{7.447cm}}
\pgfpathclose
\pgfusepath{clip}
\pgfsetdash{}{0cm}
\pgfpathmoveto{\pgfqpoint{6.638cm}{7.955cm}}
\pgfpathlineto{\pgfqpoint{6.862cm}{7.955cm}}
\pgfpathlineto{\pgfqpoint{6.862cm}{7.732cm}}
\pgfpathlineto{\pgfqpoint{6.638cm}{7.732cm}}
\pgfpathlineto{\pgfqpoint{6.638cm}{7.955cm}}
\pgfpathclose
\pgfusepath{stroke}
\end{pgfscope}
\end{pgfscope}
\pgftext[x=12.203cm+.15cm,y=6.496cm,rotate=0]{\fontsize{34}{36.14}\selectfont{ {Lower bound(\texttt{maxW})}}}
\begin{pgfscope}
\pgfpathmoveto{\pgfqpoint{6.024cm}{13.514cm}}
\pgfpathlineto{\pgfqpoint{17.13cm}{13.514cm}}
\pgfpathlineto{\pgfqpoint{17.13cm}{5.909cm}}
\pgfpathlineto{\pgfqpoint{6.024cm}{5.909cm}}
\pgfpathclose
\pgfusepath{clip}
\pgfsetdash{{0.212cm}}{0cm}
\definecolor{eps2pgf_color}{rgb}{0,1,1}\pgfsetstrokecolor{eps2pgf_color}\pgfsetfillcolor{eps2pgf_color}
\pgfpathmoveto{\pgfqpoint{6.229cm}{6.594cm}}
\pgfpathlineto{\pgfqpoint{7.27cm}{6.594cm}}
\pgfusepath{stroke}
\begin{pgfscope}
\pgfpathmoveto{\pgfqpoint{6.356cm}{6.988cm}}
\pgfpathlineto{\pgfqpoint{7.147cm}{6.988cm}}
\pgfpathlineto{\pgfqpoint{7.147cm}{6.197cm}}
\pgfpathlineto{\pgfqpoint{6.356cm}{6.197cm}}
\pgfpathclose
\pgfusepath{clip}
\pgfsetdash{}{0cm}
\pgfpathmoveto{\pgfqpoint{6.653cm}{6.691cm}}
\pgfpathlineto{\pgfqpoint{6.847cm}{6.497cm}}
\pgfusepath{stroke}
\pgfsetdash{}{0cm}
\pgfpathmoveto{\pgfqpoint{6.847cm}{6.691cm}}
\pgfpathlineto{\pgfqpoint{6.653cm}{6.497cm}}
\pgfusepath{stroke}
\end{pgfscope}
\end{pgfscope}
\pgfsetdash{}{0cm}
\pgfsetlinewidth{0.176mm}
\definecolor{eps2pgf_color}{rgb}{0,1,1}\pgfsetstrokecolor{eps2pgf_color}\pgfsetfillcolor{eps2pgf_color}
\pgfusepath{stroke}
\end{pgfscope}
\end{pgfscope}
\end{pgfpicture}

%% file: fig/I3.tex
\scalebox{0.24}{\scalefont{2} \input{./fig/single_storage_costSaving_I1.pgf}}

%% file: fig/single_storage_costSaving_I1.pgf
% Created by Eps2pgf 0.7.0 (build on 2008-08-24) on Mon Apr 06 20:56:31 PDT 2015
\begin{pgfpicture}
\pgfpathmoveto{\pgfqpoint{2.293cm}{3.387cm}}
\pgfpathlineto{\pgfqpoint{19.226cm}{3.387cm}}
\pgfpathlineto{\pgfqpoint{19.226cm}{24.553cm}}
\pgfpathlineto{\pgfqpoint{2.293cm}{24.553cm}}
\pgfpathclose
\pgfusepath{clip}
\begin{pgfscope}
\begin{pgfscope}
\pgfpathmoveto{\pgfqpoint{2.293cm}{24.553cm}}
\pgfpathlineto{\pgfqpoint{19.229cm}{24.553cm}}
\pgfpathlineto{\pgfqpoint{19.229cm}{3.413cm}}
\pgfpathlineto{\pgfqpoint{2.293cm}{3.413cm}}
\pgfpathclose
\pgfusepath{clip}
\begin{pgfscope}
\definecolor{eps2pgf_color}{gray}{1}\pgfsetstrokecolor{eps2pgf_color}\pgfsetfillcolor{eps2pgf_color}
\pgfpathmoveto{\pgfqpoint{2.293cm}{24.553cm}}
\pgfpathlineto{\pgfqpoint{19.232cm}{24.553cm}}
\pgfpathlineto{\pgfqpoint{19.232cm}{3.41cm}}
\pgfpathlineto{\pgfqpoint{2.293cm}{3.41cm}}
\pgfpathclose
\pgfusepath{fill}
\end{pgfscope}
\definecolor{eps2pgf_color}{gray}{1}\pgfsetstrokecolor{eps2pgf_color}\pgfsetfillcolor{eps2pgf_color}
\pgfpathmoveto{\pgfqpoint{4.965cm}{5.739cm}}
\pgfpathlineto{\pgfqpoint{4.965cm}{22.969cm}}
\pgfpathlineto{\pgfqpoint{17.621cm}{22.969cm}}
\pgfpathlineto{\pgfqpoint{17.621cm}{5.739cm}}
\pgfpathclose
\pgfseteorule\pgfusepath{fill}\pgfsetnonzerorule
\pgfsetdash{}{0cm}
\pgfsetlinewidth{0.176mm}
\pgfsetroundjoin
\pgfpathmoveto{\pgfqpoint{4.965cm}{5.739cm}}
\pgfpathlineto{\pgfqpoint{4.965cm}{22.969cm}}
\pgfpathlineto{\pgfqpoint{17.621cm}{22.969cm}}
\pgfpathlineto{\pgfqpoint{17.621cm}{5.739cm}}
\pgfpathlineto{\pgfqpoint{4.965cm}{5.739cm}}
\pgfusepath{stroke}
\pgfsetdash{}{0cm}
\definecolor{eps2pgf_color}{gray}{0}\pgfsetstrokecolor{eps2pgf_color}\pgfsetfillcolor{eps2pgf_color}
\pgfpathmoveto{\pgfqpoint{4.965cm}{5.739cm}}
\pgfpathlineto{\pgfqpoint{17.621cm}{5.739cm}}
\pgfusepath{stroke}
\pgfsetdash{}{0cm}
\pgfpathmoveto{\pgfqpoint{4.965cm}{22.969cm}}
\pgfpathlineto{\pgfqpoint{17.621cm}{22.969cm}}
\pgfusepath{stroke}
\pgfsetdash{}{0cm}
\pgfpathmoveto{\pgfqpoint{4.965cm}{5.739cm}}
\pgfpathlineto{\pgfqpoint{4.965cm}{22.969cm}}
\pgfusepath{stroke}
\pgfsetdash{}{0cm}
\pgfpathmoveto{\pgfqpoint{17.621cm}{5.739cm}}
\pgfpathlineto{\pgfqpoint{17.621cm}{22.969cm}}
\pgfusepath{stroke}
\pgfsetdash{}{0cm}
\pgfpathmoveto{\pgfqpoint{4.965cm}{5.739cm}}
\pgfpathlineto{\pgfqpoint{17.621cm}{5.739cm}}
\pgfusepath{stroke}
\pgfsetdash{}{0cm}
\pgfpathmoveto{\pgfqpoint{4.965cm}{5.739cm}}
\pgfpathlineto{\pgfqpoint{4.965cm}{22.969cm}}
\pgfusepath{stroke}
\pgfsetdash{}{0cm}
\pgfpathmoveto{\pgfqpoint{4.965cm}{5.739cm}}
\pgfpathlineto{\pgfqpoint{4.965cm}{5.912cm}}
\pgfusepath{stroke}
\pgfsetdash{}{0cm}
\pgfpathmoveto{\pgfqpoint{4.965cm}{22.969cm}}
\pgfpathlineto{\pgfqpoint{4.965cm}{22.798cm}}
\pgfusepath{stroke}
\pgftext[x=4.966cm,y=5.016cm+.2cm,rotate=0]{\fontsize{36}{36.14}\selectfont{ {0}}}
\pgfsetdash{}{0cm}
\pgfpathmoveto{\pgfqpoint{11.292cm}{5.739cm}}
\pgfpathlineto{\pgfqpoint{11.292cm}{5.912cm}}
\pgfusepath{stroke}
\pgfsetdash{}{0cm}
\pgfpathmoveto{\pgfqpoint{11.292cm}{22.969cm}}
\pgfpathlineto{\pgfqpoint{11.292cm}{22.798cm}}
\pgfusepath{stroke}
\pgftext[x=11.29cm,y=5.016cm+.2cm,rotate=0]{\fontsize{36}{36.14}\selectfont{ {0.5}}}
\pgfsetdash{}{0cm}
\pgfpathmoveto{\pgfqpoint{17.621cm}{5.739cm}}
\pgfpathlineto{\pgfqpoint{17.621cm}{5.912cm}}
\pgfusepath{stroke}
\pgfsetdash{}{0cm}
\pgfpathmoveto{\pgfqpoint{17.621cm}{22.969cm}}
\pgfpathlineto{\pgfqpoint{17.621cm}{22.798cm}}
\pgfusepath{stroke}
\pgftext[x=17.571cm,y=5.026cm+.2cm,rotate=0]{\fontsize{36}{36.14}\selectfont{ {1}}}
\pgfsetdash{}{0cm}
\pgfpathmoveto{\pgfqpoint{4.965cm}{5.739cm}}
\pgfpathlineto{\pgfqpoint{5.136cm}{5.739cm}}
\pgfusepath{stroke}
\pgfsetdash{}{0cm}
\pgfpathmoveto{\pgfqpoint{17.621cm}{5.739cm}}
\pgfpathlineto{\pgfqpoint{17.448cm}{5.739cm}}
\pgfusepath{stroke}
\pgftext[x=4.569cm-.2cm,y=5.707cm,rotate=0]{\fontsize{36}{36.14}\selectfont{ {0}}}
\pgfsetdash{}{0cm}
\pgfpathmoveto{\pgfqpoint{4.965cm}{8.611cm}}
\pgfpathlineto{\pgfqpoint{5.136cm}{8.611cm}}
\pgfusepath{stroke}
\pgfsetdash{}{0cm}
\pgfpathmoveto{\pgfqpoint{17.621cm}{8.611cm}}
\pgfpathlineto{\pgfqpoint{17.448cm}{8.611cm}}
\pgfusepath{stroke}
\pgftext[x=4.309cm-.2cm,y=8.579cm,rotate=0]{\fontsize{36}{36.14}\selectfont{ {10}}}
\pgfsetdash{}{0cm}
\pgfpathmoveto{\pgfqpoint{4.965cm}{11.483cm}}
\pgfpathlineto{\pgfqpoint{5.136cm}{11.483cm}}
\pgfusepath{stroke}
\pgfsetdash{}{0cm}
\pgfpathmoveto{\pgfqpoint{17.621cm}{11.483cm}}
\pgfpathlineto{\pgfqpoint{17.448cm}{11.483cm}}
\pgfusepath{stroke}
\pgftext[x=4.269cm-.2cm,y=11.451cm,rotate=0]{\fontsize{36}{36.14}\selectfont{ {20}}}
\pgfsetdash{}{0cm}
\pgfpathmoveto{\pgfqpoint{4.965cm}{14.355cm}}
\pgfpathlineto{\pgfqpoint{5.136cm}{14.355cm}}
\pgfusepath{stroke}
\pgfsetdash{}{0cm}
\pgfpathmoveto{\pgfqpoint{17.621cm}{14.355cm}}
\pgfpathlineto{\pgfqpoint{17.448cm}{14.355cm}}
\pgfusepath{stroke}
\pgftext[x=4.273cm-.2cm,y=14.323cm,rotate=0]{\fontsize{36}{36.14}\selectfont{ {30}}}
\pgfsetdash{}{0cm}
\pgfpathmoveto{\pgfqpoint{4.965cm}{17.227cm}}
\pgfpathlineto{\pgfqpoint{5.136cm}{17.227cm}}
\pgfusepath{stroke}
\pgfsetdash{}{0cm}
\pgfpathmoveto{\pgfqpoint{17.621cm}{17.227cm}}
\pgfpathlineto{\pgfqpoint{17.448cm}{17.227cm}}
\pgfusepath{stroke}
\pgftext[x=4.269cm-.2cm,y=17.195cm,rotate=0]{\fontsize{36}{36.14}\selectfont{ {40}}}
\pgfsetdash{}{0cm}
\pgfpathmoveto{\pgfqpoint{4.965cm}{20.1cm}}
\pgfpathlineto{\pgfqpoint{5.136cm}{20.1cm}}
\pgfusepath{stroke}
\pgfsetdash{}{0cm}
\pgfpathmoveto{\pgfqpoint{17.621cm}{20.1cm}}
\pgfpathlineto{\pgfqpoint{17.448cm}{20.1cm}}
\pgfusepath{stroke}
\pgftext[x=4.272cm-.2cm,y=20.068cm,rotate=0]{\fontsize{36}{36.14}\selectfont{ {50}}}
\pgfsetdash{}{0cm}
\pgfpathmoveto{\pgfqpoint{4.965cm}{22.969cm}}
\pgfpathlineto{\pgfqpoint{5.136cm}{22.969cm}}
\pgfusepath{stroke}
\pgfsetdash{}{0cm}
\pgfpathmoveto{\pgfqpoint{17.621cm}{22.969cm}}
\pgfpathlineto{\pgfqpoint{17.448cm}{22.969cm}}
\pgfusepath{stroke}
\pgftext[x=4.275cm-.2cm,y=22.937cm,rotate=0]{\fontsize{36}{36.14}\selectfont{ {60}}}
\pgfsetdash{}{0cm}
\pgfpathmoveto{\pgfqpoint{4.965cm}{5.739cm}}
\pgfpathlineto{\pgfqpoint{17.621cm}{5.739cm}}
\pgfusepath{stroke}
\pgfsetdash{}{0cm}
\pgfpathmoveto{\pgfqpoint{4.965cm}{22.969cm}}
\pgfpathlineto{\pgfqpoint{17.621cm}{22.969cm}}
\pgfusepath{stroke}
\pgfsetdash{}{0cm}
\pgfpathmoveto{\pgfqpoint{4.965cm}{5.739cm}}
\pgfpathlineto{\pgfqpoint{4.965cm}{22.969cm}}
\pgfusepath{stroke}
\pgfsetdash{}{0cm}
\pgfpathmoveto{\pgfqpoint{17.621cm}{5.739cm}}
\pgfpathlineto{\pgfqpoint{17.621cm}{22.969cm}}
\pgfusepath{stroke}
\begin{pgfscope}
\pgfpathmoveto{\pgfqpoint{4.965cm}{22.969cm}}
\pgfpathlineto{\pgfqpoint{17.624cm}{22.969cm}}
\pgfpathlineto{\pgfqpoint{17.624cm}{5.736cm}}
\pgfpathlineto{\pgfqpoint{4.965cm}{5.736cm}}
\pgfpathclose
\pgfusepath{clip}
\pgfsetdash{}{0cm}
\pgfsetlinewidth{1.058mm}
\definecolor{eps2pgf_color}{rgb}{0,0,1}\pgfsetstrokecolor{eps2pgf_color}\pgfsetfillcolor{eps2pgf_color}
\pgfpathmoveto{\pgfqpoint{6.229cm}{7.123cm}}
\pgfpathlineto{\pgfqpoint{7.497cm}{8.446cm}}
\pgfpathlineto{\pgfqpoint{8.761cm}{9.696cm}}
\pgfpathlineto{\pgfqpoint{10.028cm}{10.883cm}}
\pgfpathlineto{\pgfqpoint{11.292cm}{12.009cm}}
\pgfpathlineto{\pgfqpoint{12.559cm}{13.079cm}}
\pgfpathlineto{\pgfqpoint{13.823cm}{14.102cm}}
\pgfpathlineto{\pgfqpoint{15.09cm}{15.052cm}}
\pgfpathlineto{\pgfqpoint{16.354cm}{15.96cm}}
\pgfpathlineto{\pgfqpoint{17.621cm}{16.807cm}}
\pgfusepath{stroke}
\end{pgfscope}
\pgfsetdash{}{0cm}
\pgfsetlinewidth{1.058mm}
\pgfsetmiterjoin
\definecolor{eps2pgf_color}{rgb}{0,0,1}\pgfsetstrokecolor{eps2pgf_color}\pgfsetfillcolor{eps2pgf_color}
\pgfpathmoveto{\pgfqpoint{6.118cm}{7.235cm}}
\pgfpathlineto{\pgfqpoint{6.341cm}{7.235cm}}
\pgfpathlineto{\pgfqpoint{6.341cm}{7.011cm}}
\pgfpathlineto{\pgfqpoint{6.118cm}{7.011cm}}
\pgfpathlineto{\pgfqpoint{6.118cm}{7.235cm}}
\pgfpathclose
\pgfusepath{stroke}
\pgfsetdash{}{0cm}
\pgfpathmoveto{\pgfqpoint{7.385cm}{8.558cm}}
\pgfpathlineto{\pgfqpoint{7.608cm}{8.558cm}}
\pgfpathlineto{\pgfqpoint{7.608cm}{8.334cm}}
\pgfpathlineto{\pgfqpoint{7.385cm}{8.334cm}}
\pgfpathlineto{\pgfqpoint{7.385cm}{8.558cm}}
\pgfpathclose
\pgfusepath{stroke}
\pgfsetdash{}{0cm}
\pgfpathmoveto{\pgfqpoint{8.649cm}{9.807cm}}
\pgfpathlineto{\pgfqpoint{8.872cm}{9.807cm}}
\pgfpathlineto{\pgfqpoint{8.872cm}{9.584cm}}
\pgfpathlineto{\pgfqpoint{8.649cm}{9.584cm}}
\pgfpathlineto{\pgfqpoint{8.649cm}{9.807cm}}
\pgfpathclose
\pgfusepath{stroke}
\pgfsetdash{}{0cm}
\pgfpathmoveto{\pgfqpoint{9.916cm}{10.995cm}}
\pgfpathlineto{\pgfqpoint{10.139cm}{10.995cm}}
\pgfpathlineto{\pgfqpoint{10.139cm}{10.771cm}}
\pgfpathlineto{\pgfqpoint{9.916cm}{10.771cm}}
\pgfpathlineto{\pgfqpoint{9.916cm}{10.995cm}}
\pgfpathclose
\pgfusepath{stroke}
\pgfsetdash{}{0cm}
\pgfpathmoveto{\pgfqpoint{11.18cm}{12.121cm}}
\pgfpathlineto{\pgfqpoint{11.404cm}{12.121cm}}
\pgfpathlineto{\pgfqpoint{11.404cm}{11.897cm}}
\pgfpathlineto{\pgfqpoint{11.18cm}{11.897cm}}
\pgfpathlineto{\pgfqpoint{11.18cm}{12.121cm}}
\pgfpathclose
\pgfusepath{stroke}
\pgfsetdash{}{0cm}
\pgfpathmoveto{\pgfqpoint{12.447cm}{13.191cm}}
\pgfpathlineto{\pgfqpoint{12.671cm}{13.191cm}}
\pgfpathlineto{\pgfqpoint{12.671cm}{12.968cm}}
\pgfpathlineto{\pgfqpoint{12.447cm}{12.968cm}}
\pgfpathlineto{\pgfqpoint{12.447cm}{13.191cm}}
\pgfpathclose
\pgfusepath{stroke}
\pgfsetdash{}{0cm}
\pgfpathmoveto{\pgfqpoint{13.711cm}{14.214cm}}
\pgfpathlineto{\pgfqpoint{13.935cm}{14.214cm}}
\pgfpathlineto{\pgfqpoint{13.935cm}{13.991cm}}
\pgfpathlineto{\pgfqpoint{13.711cm}{13.991cm}}
\pgfpathlineto{\pgfqpoint{13.711cm}{14.214cm}}
\pgfpathclose
\pgfusepath{stroke}
\pgfsetdash{}{0cm}
\pgfpathmoveto{\pgfqpoint{14.978cm}{15.164cm}}
\pgfpathlineto{\pgfqpoint{15.202cm}{15.164cm}}
\pgfpathlineto{\pgfqpoint{15.202cm}{14.94cm}}
\pgfpathlineto{\pgfqpoint{14.978cm}{14.94cm}}
\pgfpathlineto{\pgfqpoint{14.978cm}{15.164cm}}
\pgfpathclose
\pgfusepath{stroke}
\pgfsetdash{}{0cm}
\pgfpathmoveto{\pgfqpoint{16.242cm}{16.072cm}}
\pgfpathlineto{\pgfqpoint{16.466cm}{16.072cm}}
\pgfpathlineto{\pgfqpoint{16.466cm}{15.849cm}}
\pgfpathlineto{\pgfqpoint{16.242cm}{15.849cm}}
\pgfpathlineto{\pgfqpoint{16.242cm}{16.072cm}}
\pgfpathclose
\pgfusepath{stroke}
\pgfsetdash{}{0cm}
\pgfpathmoveto{\pgfqpoint{17.51cm}{16.919cm}}
\pgfpathlineto{\pgfqpoint{17.733cm}{16.919cm}}
\pgfpathlineto{\pgfqpoint{17.733cm}{16.695cm}}
\pgfpathlineto{\pgfqpoint{17.51cm}{16.695cm}}
\pgfpathlineto{\pgfqpoint{17.51cm}{16.919cm}}
\pgfpathclose
\pgfusepath{stroke}
\begin{pgfscope}
\pgfpathmoveto{\pgfqpoint{4.965cm}{22.969cm}}
\pgfpathlineto{\pgfqpoint{17.624cm}{22.969cm}}
\pgfpathlineto{\pgfqpoint{17.624cm}{5.736cm}}
\pgfpathlineto{\pgfqpoint{4.965cm}{5.736cm}}
\pgfpathclose
\pgfusepath{clip}
\pgfsetdash{}{0cm}
\definecolor{eps2pgf_color}{rgb}{0,1,0}\pgfsetstrokecolor{eps2pgf_color}\pgfsetfillcolor{eps2pgf_color}
\pgfpathmoveto{\pgfqpoint{6.229cm}{7.126cm}}
\pgfpathlineto{\pgfqpoint{7.497cm}{8.449cm}}
\pgfpathlineto{\pgfqpoint{8.761cm}{9.704cm}}
\pgfpathlineto{\pgfqpoint{10.028cm}{10.91cm}}
\pgfpathlineto{\pgfqpoint{11.292cm}{12.047cm}}
\pgfpathlineto{\pgfqpoint{12.559cm}{13.126cm}}
\pgfpathlineto{\pgfqpoint{13.823cm}{14.152cm}}
\pgfpathlineto{\pgfqpoint{15.09cm}{15.117cm}}
\pgfpathlineto{\pgfqpoint{16.354cm}{16.031cm}}
\pgfpathlineto{\pgfqpoint{17.621cm}{16.889cm}}
\pgfusepath{stroke}
\end{pgfscope}
\pgfsetdash{}{0cm}
\definecolor{eps2pgf_color}{rgb}{0,1,0}\pgfsetstrokecolor{eps2pgf_color}\pgfsetfillcolor{eps2pgf_color}
\pgfpathmoveto{\pgfqpoint{6.088cm}{7.126cm}}
\pgfpathlineto{\pgfqpoint{6.371cm}{7.126cm}}
\pgfusepath{stroke}
\pgfsetdash{}{0cm}
\pgfpathmoveto{\pgfqpoint{6.229cm}{7.267cm}}
\pgfpathlineto{\pgfqpoint{6.229cm}{6.985cm}}
\pgfusepath{stroke}
\pgfsetdash{}{0cm}
\pgfpathmoveto{\pgfqpoint{7.355cm}{8.449cm}}
\pgfpathlineto{\pgfqpoint{7.638cm}{8.449cm}}
\pgfusepath{stroke}
\pgfsetdash{}{0cm}
\pgfpathmoveto{\pgfqpoint{7.497cm}{8.59cm}}
\pgfpathlineto{\pgfqpoint{7.497cm}{8.308cm}}
\pgfusepath{stroke}
\pgfsetdash{}{0cm}
\pgfpathmoveto{\pgfqpoint{8.62cm}{9.704cm}}
\pgfpathlineto{\pgfqpoint{8.902cm}{9.704cm}}
\pgfusepath{stroke}
\pgfsetdash{}{0cm}
\pgfpathmoveto{\pgfqpoint{8.761cm}{9.845cm}}
\pgfpathlineto{\pgfqpoint{8.761cm}{9.563cm}}
\pgfusepath{stroke}
\pgfsetdash{}{0cm}
\pgfpathmoveto{\pgfqpoint{9.887cm}{10.91cm}}
\pgfpathlineto{\pgfqpoint{10.169cm}{10.91cm}}
\pgfusepath{stroke}
\pgfsetdash{}{0cm}
\pgfpathmoveto{\pgfqpoint{10.028cm}{11.051cm}}
\pgfpathlineto{\pgfqpoint{10.028cm}{10.769cm}}
\pgfusepath{stroke}
\pgfsetdash{}{0cm}
\pgfpathmoveto{\pgfqpoint{11.151cm}{12.047cm}}
\pgfpathlineto{\pgfqpoint{11.433cm}{12.047cm}}
\pgfusepath{stroke}
\pgfsetdash{}{0cm}
\pgfpathmoveto{\pgfqpoint{11.292cm}{12.188cm}}
\pgfpathlineto{\pgfqpoint{11.292cm}{11.906cm}}
\pgfusepath{stroke}
\pgfsetdash{}{0cm}
\pgfpathmoveto{\pgfqpoint{12.418cm}{13.126cm}}
\pgfpathlineto{\pgfqpoint{12.7cm}{13.126cm}}
\pgfusepath{stroke}
\pgfsetdash{}{0cm}
\pgfpathmoveto{\pgfqpoint{12.559cm}{13.267cm}}
\pgfpathlineto{\pgfqpoint{12.559cm}{12.985cm}}
\pgfusepath{stroke}
\pgfsetdash{}{0cm}
\pgfpathmoveto{\pgfqpoint{13.682cm}{14.152cm}}
\pgfpathlineto{\pgfqpoint{13.964cm}{14.152cm}}
\pgfusepath{stroke}
\pgfsetdash{}{0cm}
\pgfpathmoveto{\pgfqpoint{13.823cm}{14.293cm}}
\pgfpathlineto{\pgfqpoint{13.823cm}{14.011cm}}
\pgfusepath{stroke}
\pgfsetdash{}{0cm}
\pgfpathmoveto{\pgfqpoint{14.949cm}{15.117cm}}
\pgfpathlineto{\pgfqpoint{15.231cm}{15.117cm}}
\pgfusepath{stroke}
\pgfsetdash{}{0cm}
\pgfpathmoveto{\pgfqpoint{15.09cm}{15.258cm}}
\pgfpathlineto{\pgfqpoint{15.09cm}{14.975cm}}
\pgfusepath{stroke}
\pgfsetdash{}{0cm}
\pgfpathmoveto{\pgfqpoint{16.213cm}{16.031cm}}
\pgfpathlineto{\pgfqpoint{16.495cm}{16.031cm}}
\pgfusepath{stroke}
\pgfsetdash{}{0cm}
\pgfpathmoveto{\pgfqpoint{16.354cm}{16.172cm}}
\pgfpathlineto{\pgfqpoint{16.354cm}{15.89cm}}
\pgfusepath{stroke}
\pgfsetdash{}{0cm}
\pgfpathmoveto{\pgfqpoint{17.48cm}{16.889cm}}
\pgfpathlineto{\pgfqpoint{17.762cm}{16.889cm}}
\pgfusepath{stroke}
\pgfsetdash{}{0cm}
\pgfpathmoveto{\pgfqpoint{17.621cm}{17.03cm}}
\pgfpathlineto{\pgfqpoint{17.621cm}{16.748cm}}
\pgfusepath{stroke}
\pgfsetdash{}{0cm}
\pgfpathmoveto{\pgfqpoint{6.132cm}{7.223cm}}
\pgfpathlineto{\pgfqpoint{6.326cm}{7.029cm}}
\pgfusepath{stroke}
\pgfsetdash{}{0cm}
\pgfpathmoveto{\pgfqpoint{6.326cm}{7.223cm}}
\pgfpathlineto{\pgfqpoint{6.132cm}{7.029cm}}
\pgfusepath{stroke}
\pgfsetdash{}{0cm}
\pgfpathmoveto{\pgfqpoint{7.4cm}{8.546cm}}
\pgfpathlineto{\pgfqpoint{7.594cm}{8.352cm}}
\pgfusepath{stroke}
\pgfsetdash{}{0cm}
\pgfpathmoveto{\pgfqpoint{7.594cm}{8.546cm}}
\pgfpathlineto{\pgfqpoint{7.4cm}{8.352cm}}
\pgfusepath{stroke}
\pgfsetdash{}{0cm}
\pgfpathmoveto{\pgfqpoint{8.664cm}{9.801cm}}
\pgfpathlineto{\pgfqpoint{8.858cm}{9.607cm}}
\pgfusepath{stroke}
\pgfsetdash{}{0cm}
\pgfpathmoveto{\pgfqpoint{8.858cm}{9.801cm}}
\pgfpathlineto{\pgfqpoint{8.664cm}{9.607cm}}
\pgfusepath{stroke}
\pgfsetdash{}{0cm}
\pgfpathmoveto{\pgfqpoint{9.931cm}{11.007cm}}
\pgfpathlineto{\pgfqpoint{10.125cm}{10.813cm}}
\pgfusepath{stroke}
\pgfsetdash{}{0cm}
\pgfpathmoveto{\pgfqpoint{10.125cm}{11.007cm}}
\pgfpathlineto{\pgfqpoint{9.931cm}{10.813cm}}
\pgfusepath{stroke}
\pgfsetdash{}{0cm}
\pgfpathmoveto{\pgfqpoint{11.195cm}{12.144cm}}
\pgfpathlineto{\pgfqpoint{11.389cm}{11.95cm}}
\pgfusepath{stroke}
\pgfsetdash{}{0cm}
\pgfpathmoveto{\pgfqpoint{11.389cm}{12.144cm}}
\pgfpathlineto{\pgfqpoint{11.195cm}{11.95cm}}
\pgfusepath{stroke}
\pgfsetdash{}{0cm}
\pgfpathmoveto{\pgfqpoint{12.462cm}{13.223cm}}
\pgfpathlineto{\pgfqpoint{12.656cm}{13.029cm}}
\pgfusepath{stroke}
\pgfsetdash{}{0cm}
\pgfpathmoveto{\pgfqpoint{12.656cm}{13.223cm}}
\pgfpathlineto{\pgfqpoint{12.462cm}{13.029cm}}
\pgfusepath{stroke}
\pgfsetdash{}{0cm}
\pgfpathmoveto{\pgfqpoint{13.726cm}{14.249cm}}
\pgfpathlineto{\pgfqpoint{13.92cm}{14.055cm}}
\pgfusepath{stroke}
\pgfsetdash{}{0cm}
\pgfpathmoveto{\pgfqpoint{13.92cm}{14.249cm}}
\pgfpathlineto{\pgfqpoint{13.726cm}{14.055cm}}
\pgfusepath{stroke}
\pgfsetdash{}{0cm}
\pgfpathmoveto{\pgfqpoint{14.993cm}{15.214cm}}
\pgfpathlineto{\pgfqpoint{15.187cm}{15.02cm}}
\pgfusepath{stroke}
\pgfsetdash{}{0cm}
\pgfpathmoveto{\pgfqpoint{15.187cm}{15.214cm}}
\pgfpathlineto{\pgfqpoint{14.993cm}{15.02cm}}
\pgfusepath{stroke}
\pgfsetdash{}{0cm}
\pgfpathmoveto{\pgfqpoint{16.257cm}{16.128cm}}
\pgfpathlineto{\pgfqpoint{16.451cm}{15.934cm}}
\pgfusepath{stroke}
\pgfsetdash{}{0cm}
\pgfpathmoveto{\pgfqpoint{16.451cm}{16.128cm}}
\pgfpathlineto{\pgfqpoint{16.257cm}{15.934cm}}
\pgfusepath{stroke}
\pgfsetdash{}{0cm}
\pgfpathmoveto{\pgfqpoint{17.524cm}{16.986cm}}
\pgfpathlineto{\pgfqpoint{17.718cm}{16.792cm}}
\pgfusepath{stroke}
\pgfsetdash{}{0cm}
\pgfpathmoveto{\pgfqpoint{17.718cm}{16.986cm}}
\pgfpathlineto{\pgfqpoint{17.524cm}{16.792cm}}
\pgfusepath{stroke}
\begin{pgfscope}
\pgfpathmoveto{\pgfqpoint{4.965cm}{22.969cm}}
\pgfpathlineto{\pgfqpoint{17.624cm}{22.969cm}}
\pgfpathlineto{\pgfqpoint{17.624cm}{5.736cm}}
\pgfpathlineto{\pgfqpoint{4.965cm}{5.736cm}}
\pgfpathclose
\pgfusepath{clip}
\pgfsetdash{{0.018cm}{0.141cm}{0.212cm}{0.141cm}}{0cm}
\definecolor{eps2pgf_color}{rgb}{0,0,1}\pgfsetstrokecolor{eps2pgf_color}\pgfsetfillcolor{eps2pgf_color}
\pgfpathmoveto{\pgfqpoint{6.229cm}{7.203cm}}
\pgfpathlineto{\pgfqpoint{7.497cm}{8.605cm}}
\pgfpathlineto{\pgfqpoint{8.761cm}{9.934cm}}
\pgfpathlineto{\pgfqpoint{10.028cm}{11.201cm}}
\pgfpathlineto{\pgfqpoint{11.292cm}{12.409cm}}
\pgfpathlineto{\pgfqpoint{12.559cm}{13.558cm}}
\pgfpathlineto{\pgfqpoint{13.823cm}{14.661cm}}
\pgfpathlineto{\pgfqpoint{15.09cm}{15.69cm}}
\pgfpathlineto{\pgfqpoint{16.354cm}{16.678cm}}
\pgfpathlineto{\pgfqpoint{17.621cm}{17.604cm}}
\pgfusepath{stroke}
\end{pgfscope}
\pgfsetdash{}{0cm}
\definecolor{eps2pgf_color}{rgb}{0,0,1}\pgfsetstrokecolor{eps2pgf_color}\pgfsetfillcolor{eps2pgf_color}
\pgfpathmoveto{\pgfqpoint{6.118cm}{7.314cm}}
\pgfpathlineto{\pgfqpoint{6.341cm}{7.314cm}}
\pgfpathlineto{\pgfqpoint{6.341cm}{7.091cm}}
\pgfpathlineto{\pgfqpoint{6.118cm}{7.091cm}}
\pgfpathlineto{\pgfqpoint{6.118cm}{7.314cm}}
\pgfpathclose
\pgfusepath{stroke}
\pgfsetdash{}{0cm}
\pgfpathmoveto{\pgfqpoint{7.385cm}{8.717cm}}
\pgfpathlineto{\pgfqpoint{7.608cm}{8.717cm}}
\pgfpathlineto{\pgfqpoint{7.608cm}{8.493cm}}
\pgfpathlineto{\pgfqpoint{7.385cm}{8.493cm}}
\pgfpathlineto{\pgfqpoint{7.385cm}{8.717cm}}
\pgfpathclose
\pgfusepath{stroke}
\pgfsetdash{}{0cm}
\pgfpathmoveto{\pgfqpoint{8.649cm}{10.045cm}}
\pgfpathlineto{\pgfqpoint{8.872cm}{10.045cm}}
\pgfpathlineto{\pgfqpoint{8.872cm}{9.822cm}}
\pgfpathlineto{\pgfqpoint{8.649cm}{9.822cm}}
\pgfpathlineto{\pgfqpoint{8.649cm}{10.045cm}}
\pgfpathclose
\pgfusepath{stroke}
\pgfsetdash{}{0cm}
\pgfpathmoveto{\pgfqpoint{9.916cm}{11.312cm}}
\pgfpathlineto{\pgfqpoint{10.139cm}{11.312cm}}
\pgfpathlineto{\pgfqpoint{10.139cm}{11.089cm}}
\pgfpathlineto{\pgfqpoint{9.916cm}{11.089cm}}
\pgfpathlineto{\pgfqpoint{9.916cm}{11.312cm}}
\pgfpathclose
\pgfusepath{stroke}
\pgfsetdash{}{0cm}
\pgfpathmoveto{\pgfqpoint{11.18cm}{12.521cm}}
\pgfpathlineto{\pgfqpoint{11.404cm}{12.521cm}}
\pgfpathlineto{\pgfqpoint{11.404cm}{12.297cm}}
\pgfpathlineto{\pgfqpoint{11.18cm}{12.297cm}}
\pgfpathlineto{\pgfqpoint{11.18cm}{12.521cm}}
\pgfpathclose
\pgfusepath{stroke}
\pgfsetdash{}{0cm}
\pgfpathmoveto{\pgfqpoint{12.447cm}{13.67cm}}
\pgfpathlineto{\pgfqpoint{12.671cm}{13.67cm}}
\pgfpathlineto{\pgfqpoint{12.671cm}{13.447cm}}
\pgfpathlineto{\pgfqpoint{12.447cm}{13.447cm}}
\pgfpathlineto{\pgfqpoint{12.447cm}{13.67cm}}
\pgfpathclose
\pgfusepath{stroke}
\pgfsetdash{}{0cm}
\pgfpathmoveto{\pgfqpoint{13.711cm}{14.773cm}}
\pgfpathlineto{\pgfqpoint{13.935cm}{14.773cm}}
\pgfpathlineto{\pgfqpoint{13.935cm}{14.549cm}}
\pgfpathlineto{\pgfqpoint{13.711cm}{14.549cm}}
\pgfpathlineto{\pgfqpoint{13.711cm}{14.773cm}}
\pgfpathclose
\pgfusepath{stroke}
\pgfsetdash{}{0cm}
\pgfpathmoveto{\pgfqpoint{14.978cm}{15.802cm}}
\pgfpathlineto{\pgfqpoint{15.202cm}{15.802cm}}
\pgfpathlineto{\pgfqpoint{15.202cm}{15.578cm}}
\pgfpathlineto{\pgfqpoint{14.978cm}{15.578cm}}
\pgfpathlineto{\pgfqpoint{14.978cm}{15.802cm}}
\pgfpathclose
\pgfusepath{stroke}
\pgfsetdash{}{0cm}
\pgfpathmoveto{\pgfqpoint{16.242cm}{16.789cm}}
\pgfpathlineto{\pgfqpoint{16.466cm}{16.789cm}}
\pgfpathlineto{\pgfqpoint{16.466cm}{16.566cm}}
\pgfpathlineto{\pgfqpoint{16.242cm}{16.566cm}}
\pgfpathlineto{\pgfqpoint{16.242cm}{16.789cm}}
\pgfpathclose
\pgfusepath{stroke}
\pgfsetdash{}{0cm}
\pgfpathmoveto{\pgfqpoint{17.51cm}{17.715cm}}
\pgfpathlineto{\pgfqpoint{17.733cm}{17.715cm}}
\pgfpathlineto{\pgfqpoint{17.733cm}{17.492cm}}
\pgfpathlineto{\pgfqpoint{17.51cm}{17.492cm}}
\pgfpathlineto{\pgfqpoint{17.51cm}{17.715cm}}
\pgfpathclose
\pgfusepath{stroke}
\begin{pgfscope}
\pgfpathmoveto{\pgfqpoint{4.965cm}{22.969cm}}
\pgfpathlineto{\pgfqpoint{17.624cm}{22.969cm}}
\pgfpathlineto{\pgfqpoint{17.624cm}{5.736cm}}
\pgfpathlineto{\pgfqpoint{4.965cm}{5.736cm}}
\pgfpathclose
\pgfusepath{clip}
\end{pgfscope}
\definecolor{eps2pgf_color}{gray}{0}\pgfsetstrokecolor{eps2pgf_color}\pgfsetfillcolor{eps2pgf_color}
\pgftext[x=11.31cm,y=4.093cm,rotate=0]{\fontsize{36}{36.14}\selectfont{ {$\bmax$}}}
\pgftext[x=3.013cm,y=14.366cm,rotate=90]{\fontsize{36}{36.14}\selectfont{ {Percentage cost saving (\%)}}}
\pgftext[x=4.915cm,y=5.612cm,rotate=0]{\fontsize{10.04}{12.04}\selectfont{ { }}}
\pgftext[x=17.574cm,y=22.845cm,rotate=0]{\fontsize{10.04}{12.04}\selectfont{ { }}}
\pgftext[x=7.731cm+.2cm,y=22.117cm,rotate=0]{\fontsize{36}{36.14}\selectfont{ {OMG}}}
\begin{pgfscope}
\pgfpathmoveto{\pgfqpoint{5.142cm}{22.792cm}}
\pgfpathlineto{\pgfqpoint{12.571cm}{22.792cm}}
\pgfpathlineto{\pgfqpoint{12.571cm}{18.927cm}}
\pgfpathlineto{\pgfqpoint{5.142cm}{18.927cm}}
\pgfpathclose
\pgfusepath{clip}
\pgfsetdash{}{0cm}
\definecolor{eps2pgf_color}{rgb}{0,0,1}\pgfsetstrokecolor{eps2pgf_color}\pgfsetfillcolor{eps2pgf_color}
\pgfpathmoveto{\pgfqpoint{5.348cm}{22.113cm}}
\pgfpathlineto{\pgfqpoint{6.379cm}{22.113cm}}
\pgfusepath{stroke}
\begin{pgfscope}
\pgfpathmoveto{\pgfqpoint{5.468cm}{22.507cm}}
\pgfpathlineto{\pgfqpoint{6.259cm}{22.507cm}}
\pgfpathlineto{\pgfqpoint{6.259cm}{21.716cm}}
\pgfpathlineto{\pgfqpoint{5.468cm}{21.716cm}}
\pgfpathclose
\pgfusepath{clip}
\pgfsetdash{}{0cm}
\pgfpathmoveto{\pgfqpoint{5.75cm}{22.225cm}}
\pgfpathlineto{\pgfqpoint{5.974cm}{22.225cm}}
\pgfpathlineto{\pgfqpoint{5.974cm}{22.002cm}}
\pgfpathlineto{\pgfqpoint{5.75cm}{22.002cm}}
\pgfpathlineto{\pgfqpoint{5.75cm}{22.225cm}}
\pgfpathclose
\pgfusepath{stroke}
\end{pgfscope}
\end{pgfscope}
\pgftext[x=8.24cm,y=20.761cm,rotate=0]{\fontsize{36}{36.14}\selectfont{ {Greedy}}}
\begin{pgfscope}
\pgfpathmoveto{\pgfqpoint{5.142cm}{22.792cm}}
\pgfpathlineto{\pgfqpoint{12.571cm}{22.792cm}}
\pgfpathlineto{\pgfqpoint{12.571cm}{18.927cm}}
\pgfpathlineto{\pgfqpoint{5.142cm}{18.927cm}}
\pgfpathclose
\pgfusepath{clip}
\pgfsetdash{}{0cm}
\definecolor{eps2pgf_color}{rgb}{0,1,0}\pgfsetstrokecolor{eps2pgf_color}\pgfsetfillcolor{eps2pgf_color}
\pgfpathmoveto{\pgfqpoint{5.348cm}{20.861cm}}
\pgfpathlineto{\pgfqpoint{6.379cm}{20.861cm}}
\pgfusepath{stroke}
\begin{pgfscope}
\pgfpathmoveto{\pgfqpoint{5.468cm}{21.255cm}}
\pgfpathlineto{\pgfqpoint{6.259cm}{21.255cm}}
\pgfpathlineto{\pgfqpoint{6.259cm}{20.464cm}}
\pgfpathlineto{\pgfqpoint{5.468cm}{20.464cm}}
\pgfpathclose
\pgfusepath{clip}
\pgfsetdash{}{0cm}
\pgfpathmoveto{\pgfqpoint{5.721cm}{20.861cm}}
\pgfpathlineto{\pgfqpoint{6.003cm}{20.861cm}}
\pgfusepath{stroke}
\pgfsetdash{}{0cm}
\pgfpathmoveto{\pgfqpoint{5.862cm}{21.002cm}}
\pgfpathlineto{\pgfqpoint{5.862cm}{20.72cm}}
\pgfusepath{stroke}
\pgfsetdash{}{0cm}
\pgfpathmoveto{\pgfqpoint{5.765cm}{20.958cm}}
\pgfpathlineto{\pgfqpoint{5.959cm}{20.764cm}}
\pgfusepath{stroke}
\pgfsetdash{}{0cm}
\pgfpathmoveto{\pgfqpoint{5.959cm}{20.958cm}}
\pgfpathlineto{\pgfqpoint{5.765cm}{20.764cm}}
\pgfusepath{stroke}
\end{pgfscope}
\end{pgfscope}
\pgftext[x=9.468cm+.4cm,y=19.506cm,rotate=0]{\fontsize{36}{36.14}\selectfont{ {Upper bound }}}
\begin{pgfscope}
\pgfpathmoveto{\pgfqpoint{5.142cm}{22.792cm}}
\pgfpathlineto{\pgfqpoint{12.571cm}{22.792cm}}
\pgfpathlineto{\pgfqpoint{12.571cm}{18.927cm}}
\pgfpathlineto{\pgfqpoint{5.142cm}{18.927cm}}
\pgfpathclose
\pgfusepath{clip}
\pgfsetdash{{0.018cm}{0.141cm}{0.212cm}{0.141cm}}{0cm}
\definecolor{eps2pgf_color}{rgb}{0,0,1}\pgfsetstrokecolor{eps2pgf_color}\pgfsetfillcolor{eps2pgf_color}
\pgfpathmoveto{\pgfqpoint{5.348cm}{19.612cm}}
\pgfpathlineto{\pgfqpoint{6.379cm}{19.612cm}}
\pgfusepath{stroke}
\begin{pgfscope}
\pgfpathmoveto{\pgfqpoint{5.468cm}{20.005cm}}
\pgfpathlineto{\pgfqpoint{6.259cm}{20.005cm}}
\pgfpathlineto{\pgfqpoint{6.259cm}{19.215cm}}
\pgfpathlineto{\pgfqpoint{5.468cm}{19.215cm}}
\pgfpathclose
\pgfusepath{clip}
\pgfsetdash{}{0cm}
\pgfpathmoveto{\pgfqpoint{5.75cm}{19.723cm}}
\pgfpathlineto{\pgfqpoint{5.974cm}{19.723cm}}
\pgfpathlineto{\pgfqpoint{5.974cm}{19.5cm}}
\pgfpathlineto{\pgfqpoint{5.75cm}{19.5cm}}
\pgfpathlineto{\pgfqpoint{5.75cm}{19.723cm}}
\pgfpathclose
\pgfusepath{stroke}
\end{pgfscope}
\end{pgfscope}
\pgfsetdash{}{0cm}
\pgfsetlinewidth{0.176mm}
\definecolor{eps2pgf_color}{rgb}{0,0,1}\pgfsetstrokecolor{eps2pgf_color}\pgfsetfillcolor{eps2pgf_color}
\pgfusepath{stroke}
\end{pgfscope}
\end{pgfscope}
\end{pgfpicture}

%% file: fig/I4.tex
\scalebox{0.24}{\scalefont{2} \input{./fig/single_storage_costSaving_I2.pgf}}

%% file: fig/single_storage_costSaving_I2.pgf
% Created by Eps2pgf 0.7.0 (build on 2008-08-24) on Mon Apr 06 20:56:34 PDT 2015
\begin{pgfpicture}
\pgfpathmoveto{\pgfqpoint{2.293cm}{3.387cm}}
\pgfpathlineto{\pgfqpoint{19.226cm}{3.387cm}}
\pgfpathlineto{\pgfqpoint{19.226cm}{24.553cm}}
\pgfpathlineto{\pgfqpoint{2.293cm}{24.553cm}}
\pgfpathclose
\pgfusepath{clip}
\begin{pgfscope}
\begin{pgfscope}
\pgfpathmoveto{\pgfqpoint{2.293cm}{24.553cm}}
\pgfpathlineto{\pgfqpoint{19.229cm}{24.553cm}}
\pgfpathlineto{\pgfqpoint{19.229cm}{3.413cm}}
\pgfpathlineto{\pgfqpoint{2.293cm}{3.413cm}}
\pgfpathclose
\pgfusepath{clip}
\begin{pgfscope}
\definecolor{eps2pgf_color}{gray}{1}\pgfsetstrokecolor{eps2pgf_color}\pgfsetfillcolor{eps2pgf_color}
\pgfpathmoveto{\pgfqpoint{2.293cm}{24.553cm}}
\pgfpathlineto{\pgfqpoint{19.232cm}{24.553cm}}
\pgfpathlineto{\pgfqpoint{19.232cm}{3.41cm}}
\pgfpathlineto{\pgfqpoint{2.293cm}{3.41cm}}
\pgfpathclose
\pgfusepath{fill}
\end{pgfscope}
\definecolor{eps2pgf_color}{gray}{1}\pgfsetstrokecolor{eps2pgf_color}\pgfsetfillcolor{eps2pgf_color}
\pgfpathmoveto{\pgfqpoint{4.965cm}{5.739cm}}
\pgfpathlineto{\pgfqpoint{4.965cm}{22.969cm}}
\pgfpathlineto{\pgfqpoint{17.621cm}{22.969cm}}
\pgfpathlineto{\pgfqpoint{17.621cm}{5.739cm}}
\pgfpathclose
\pgfseteorule\pgfusepath{fill}\pgfsetnonzerorule
\pgfsetdash{}{0cm}
\pgfsetlinewidth{0.176mm}
\pgfsetroundjoin
\pgfpathmoveto{\pgfqpoint{4.965cm}{5.739cm}}
\pgfpathlineto{\pgfqpoint{4.965cm}{22.969cm}}
\pgfpathlineto{\pgfqpoint{17.621cm}{22.969cm}}
\pgfpathlineto{\pgfqpoint{17.621cm}{5.739cm}}
\pgfpathlineto{\pgfqpoint{4.965cm}{5.739cm}}
\pgfusepath{stroke}
\pgfsetdash{}{0cm}
\definecolor{eps2pgf_color}{gray}{0}\pgfsetstrokecolor{eps2pgf_color}\pgfsetfillcolor{eps2pgf_color}
\pgfpathmoveto{\pgfqpoint{4.965cm}{5.739cm}}
\pgfpathlineto{\pgfqpoint{17.621cm}{5.739cm}}
\pgfusepath{stroke}
\pgfsetdash{}{0cm}
\pgfpathmoveto{\pgfqpoint{4.965cm}{22.969cm}}
\pgfpathlineto{\pgfqpoint{17.621cm}{22.969cm}}
\pgfusepath{stroke}
\pgfsetdash{}{0cm}
\pgfpathmoveto{\pgfqpoint{4.965cm}{5.739cm}}
\pgfpathlineto{\pgfqpoint{4.965cm}{22.969cm}}
\pgfusepath{stroke}
\pgfsetdash{}{0cm}
\pgfpathmoveto{\pgfqpoint{17.621cm}{5.739cm}}
\pgfpathlineto{\pgfqpoint{17.621cm}{22.969cm}}
\pgfusepath{stroke}
\pgfsetdash{}{0cm}
\pgfpathmoveto{\pgfqpoint{4.965cm}{5.739cm}}
\pgfpathlineto{\pgfqpoint{17.621cm}{5.739cm}}
\pgfusepath{stroke}
\pgfsetdash{}{0cm}
\pgfpathmoveto{\pgfqpoint{4.965cm}{5.739cm}}
\pgfpathlineto{\pgfqpoint{4.965cm}{22.969cm}}
\pgfusepath{stroke}
\pgfsetdash{}{0cm}
\pgfpathmoveto{\pgfqpoint{4.965cm}{5.739cm}}
\pgfpathlineto{\pgfqpoint{4.965cm}{5.912cm}}
\pgfusepath{stroke}
\pgfsetdash{}{0cm}
\pgfpathmoveto{\pgfqpoint{4.965cm}{22.969cm}}
\pgfpathlineto{\pgfqpoint{4.965cm}{22.798cm}}
\pgfusepath{stroke}
\pgftext[x=4.966cm,y=5.016cm+.2cm,rotate=0]{\fontsize{36}{36.14}\selectfont{ {0}}}
\pgfsetdash{}{0cm}
\pgfpathmoveto{\pgfqpoint{11.292cm}{5.739cm}}
\pgfpathlineto{\pgfqpoint{11.292cm}{5.912cm}}
\pgfusepath{stroke}
\pgfsetdash{}{0cm}
\pgfpathmoveto{\pgfqpoint{11.292cm}{22.969cm}}
\pgfpathlineto{\pgfqpoint{11.292cm}{22.798cm}}
\pgfusepath{stroke}
\pgftext[x=11.29cm,y=5.016cm+.2cm,rotate=0]{\fontsize{36}{36.14}\selectfont{ {0.5}}}
\pgfsetdash{}{0cm}
\pgfpathmoveto{\pgfqpoint{17.621cm}{5.739cm}}
\pgfpathlineto{\pgfqpoint{17.621cm}{5.912cm}}
\pgfusepath{stroke}
\pgfsetdash{}{0cm}
\pgfpathmoveto{\pgfqpoint{17.621cm}{22.969cm}}
\pgfpathlineto{\pgfqpoint{17.621cm}{22.798cm}}
\pgfusepath{stroke}
\pgftext[x=17.571cm,y=5.026cm+.2cm,rotate=0]{\fontsize{36}{36.14}\selectfont{ {1}}}
\pgfsetdash{}{0cm}
\pgfpathmoveto{\pgfqpoint{4.965cm}{5.739cm}}
\pgfpathlineto{\pgfqpoint{5.136cm}{5.739cm}}
\pgfusepath{stroke}
\pgfsetdash{}{0cm}
\pgfpathmoveto{\pgfqpoint{17.621cm}{5.739cm}}
\pgfpathlineto{\pgfqpoint{17.448cm}{5.739cm}}
\pgfusepath{stroke}
\pgftext[x=4.569cm-.2cm,y=5.707cm,rotate=0]{\fontsize{36}{36.14}\selectfont{ {0}}}
\pgfsetdash{}{0cm}
\pgfpathmoveto{\pgfqpoint{4.965cm}{8.611cm}}
\pgfpathlineto{\pgfqpoint{5.136cm}{8.611cm}}
\pgfusepath{stroke}
\pgfsetdash{}{0cm}
\pgfpathmoveto{\pgfqpoint{17.621cm}{8.611cm}}
\pgfpathlineto{\pgfqpoint{17.448cm}{8.611cm}}
\pgfusepath{stroke}
\pgftext[x=4.309cm-.2cm,y=8.579cm,rotate=0]{\fontsize{36}{36.14}\selectfont{ {10}}}
\pgfsetdash{}{0cm}
\pgfpathmoveto{\pgfqpoint{4.965cm}{11.483cm}}
\pgfpathlineto{\pgfqpoint{5.136cm}{11.483cm}}
\pgfusepath{stroke}
\pgfsetdash{}{0cm}
\pgfpathmoveto{\pgfqpoint{17.621cm}{11.483cm}}
\pgfpathlineto{\pgfqpoint{17.448cm}{11.483cm}}
\pgfusepath{stroke}
\pgftext[x=4.269cm-.2cm,y=11.451cm,rotate=0]{\fontsize{36}{36.14}\selectfont{ {20}}}
\pgfsetdash{}{0cm}
\pgfpathmoveto{\pgfqpoint{4.965cm}{14.355cm}}
\pgfpathlineto{\pgfqpoint{5.136cm}{14.355cm}}
\pgfusepath{stroke}
\pgfsetdash{}{0cm}
\pgfpathmoveto{\pgfqpoint{17.621cm}{14.355cm}}
\pgfpathlineto{\pgfqpoint{17.448cm}{14.355cm}}
\pgfusepath{stroke}
\pgftext[x=4.273cm-.2cm,y=14.323cm,rotate=0]{\fontsize{36}{36.14}\selectfont{ {30}}}
\pgfsetdash{}{0cm}
\pgfpathmoveto{\pgfqpoint{4.965cm}{17.227cm}}
\pgfpathlineto{\pgfqpoint{5.136cm}{17.227cm}}
\pgfusepath{stroke}
\pgfsetdash{}{0cm}
\pgfpathmoveto{\pgfqpoint{17.621cm}{17.227cm}}
\pgfpathlineto{\pgfqpoint{17.448cm}{17.227cm}}
\pgfusepath{stroke}
\pgftext[x=4.269cm-.2cm,y=17.195cm,rotate=0]{\fontsize{36}{36.14}\selectfont{ {40}}}
\pgfsetdash{}{0cm}
\pgfpathmoveto{\pgfqpoint{4.965cm}{20.1cm}}
\pgfpathlineto{\pgfqpoint{5.136cm}{20.1cm}}
\pgfusepath{stroke}
\pgfsetdash{}{0cm}
\pgfpathmoveto{\pgfqpoint{17.621cm}{20.1cm}}
\pgfpathlineto{\pgfqpoint{17.448cm}{20.1cm}}
\pgfusepath{stroke}
\pgftext[x=4.272cm-.2cm,y=20.068cm,rotate=0]{\fontsize{36}{36.14}\selectfont{ {50}}}
\pgfsetdash{}{0cm}
\pgfpathmoveto{\pgfqpoint{4.965cm}{22.969cm}}
\pgfpathlineto{\pgfqpoint{5.136cm}{22.969cm}}
\pgfusepath{stroke}
\pgfsetdash{}{0cm}
\pgfpathmoveto{\pgfqpoint{17.621cm}{22.969cm}}
\pgfpathlineto{\pgfqpoint{17.448cm}{22.969cm}}
\pgfusepath{stroke}
\pgftext[x=4.275cm-.2cm,y=22.937cm,rotate=0]{\fontsize{36}{36.14}\selectfont{ {60}}}
\pgfsetdash{}{0cm}
\pgfpathmoveto{\pgfqpoint{4.965cm}{5.739cm}}
\pgfpathlineto{\pgfqpoint{17.621cm}{5.739cm}}
\pgfusepath{stroke}
\pgfsetdash{}{0cm}
\pgfpathmoveto{\pgfqpoint{4.965cm}{22.969cm}}
\pgfpathlineto{\pgfqpoint{17.621cm}{22.969cm}}
\pgfusepath{stroke}
\pgfsetdash{}{0cm}
\pgfpathmoveto{\pgfqpoint{4.965cm}{5.739cm}}
\pgfpathlineto{\pgfqpoint{4.965cm}{22.969cm}}
\pgfusepath{stroke}
\pgfsetdash{}{0cm}
\pgfpathmoveto{\pgfqpoint{17.621cm}{5.739cm}}
\pgfpathlineto{\pgfqpoint{17.621cm}{22.969cm}}
\pgfusepath{stroke}
\begin{pgfscope}
\pgfpathmoveto{\pgfqpoint{4.965cm}{22.969cm}}
\pgfpathlineto{\pgfqpoint{17.624cm}{22.969cm}}
\pgfpathlineto{\pgfqpoint{17.624cm}{5.736cm}}
\pgfpathlineto{\pgfqpoint{4.965cm}{5.736cm}}
\pgfpathclose
\pgfusepath{clip}
\pgfsetdash{}{0cm}
\pgfsetlinewidth{1.058mm}
\definecolor{eps2pgf_color}{rgb}{0,1,0}\pgfsetstrokecolor{eps2pgf_color}\pgfsetfillcolor{eps2pgf_color}
\pgfpathmoveto{\pgfqpoint{6.229cm}{6.324cm}}
\pgfpathlineto{\pgfqpoint{7.497cm}{6.897cm}}
\pgfpathlineto{\pgfqpoint{8.761cm}{7.214cm}}
\pgfpathlineto{\pgfqpoint{10.028cm}{7.117cm}}
\pgfpathlineto{\pgfqpoint{11.292cm}{7.602cm}}
\pgfpathlineto{\pgfqpoint{12.559cm}{7.938cm}}
\pgfpathlineto{\pgfqpoint{13.823cm}{8.126cm}}
\pgfpathlineto{\pgfqpoint{15.09cm}{8.249cm}}
\pgfpathlineto{\pgfqpoint{16.354cm}{8.373cm}}
\pgfpathlineto{\pgfqpoint{17.621cm}{8.517cm}}
\pgfusepath{stroke}
\end{pgfscope}
\pgfsetdash{}{0cm}
\pgfsetlinewidth{1.058mm}
\definecolor{eps2pgf_color}{rgb}{0,1,0}\pgfsetstrokecolor{eps2pgf_color}\pgfsetfillcolor{eps2pgf_color}
\pgfpathmoveto{\pgfqpoint{6.088cm}{6.324cm}}
\pgfpathlineto{\pgfqpoint{6.371cm}{6.324cm}}
\pgfusepath{stroke}
\pgfsetdash{}{0cm}
\pgfpathmoveto{\pgfqpoint{6.229cm}{6.465cm}}
\pgfpathlineto{\pgfqpoint{6.229cm}{6.182cm}}
\pgfusepath{stroke}
\pgfsetdash{}{0cm}
\pgfpathmoveto{\pgfqpoint{7.355cm}{6.897cm}}
\pgfpathlineto{\pgfqpoint{7.638cm}{6.897cm}}
\pgfusepath{stroke}
\pgfsetdash{}{0cm}
\pgfpathmoveto{\pgfqpoint{7.497cm}{7.038cm}}
\pgfpathlineto{\pgfqpoint{7.497cm}{6.756cm}}
\pgfusepath{stroke}
\pgfsetdash{}{0cm}
\pgfpathmoveto{\pgfqpoint{8.62cm}{7.214cm}}
\pgfpathlineto{\pgfqpoint{8.902cm}{7.214cm}}
\pgfusepath{stroke}
\pgfsetdash{}{0cm}
\pgfpathmoveto{\pgfqpoint{8.761cm}{7.355cm}}
\pgfpathlineto{\pgfqpoint{8.761cm}{7.073cm}}
\pgfusepath{stroke}
\pgfsetdash{}{0cm}
\pgfpathmoveto{\pgfqpoint{9.887cm}{7.117cm}}
\pgfpathlineto{\pgfqpoint{10.169cm}{7.117cm}}
\pgfusepath{stroke}
\pgfsetdash{}{0cm}
\pgfpathmoveto{\pgfqpoint{10.028cm}{7.258cm}}
\pgfpathlineto{\pgfqpoint{10.028cm}{6.976cm}}
\pgfusepath{stroke}
\pgfsetdash{}{0cm}
\pgfpathmoveto{\pgfqpoint{11.151cm}{7.602cm}}
\pgfpathlineto{\pgfqpoint{11.433cm}{7.602cm}}
\pgfusepath{stroke}
\pgfsetdash{}{0cm}
\pgfpathmoveto{\pgfqpoint{11.292cm}{7.743cm}}
\pgfpathlineto{\pgfqpoint{11.292cm}{7.461cm}}
\pgfusepath{stroke}
\pgfsetdash{}{0cm}
\pgfpathmoveto{\pgfqpoint{12.418cm}{7.938cm}}
\pgfpathlineto{\pgfqpoint{12.7cm}{7.938cm}}
\pgfusepath{stroke}
\pgfsetdash{}{0cm}
\pgfpathmoveto{\pgfqpoint{12.559cm}{8.079cm}}
\pgfpathlineto{\pgfqpoint{12.559cm}{7.796cm}}
\pgfusepath{stroke}
\pgfsetdash{}{0cm}
\pgfpathmoveto{\pgfqpoint{13.682cm}{8.126cm}}
\pgfpathlineto{\pgfqpoint{13.964cm}{8.126cm}}
\pgfusepath{stroke}
\pgfsetdash{}{0cm}
\pgfpathmoveto{\pgfqpoint{13.823cm}{8.267cm}}
\pgfpathlineto{\pgfqpoint{13.823cm}{7.985cm}}
\pgfusepath{stroke}
\pgfsetdash{}{0cm}
\pgfpathmoveto{\pgfqpoint{14.949cm}{8.249cm}}
\pgfpathlineto{\pgfqpoint{15.231cm}{8.249cm}}
\pgfusepath{stroke}
\pgfsetdash{}{0cm}
\pgfpathmoveto{\pgfqpoint{15.09cm}{8.39cm}}
\pgfpathlineto{\pgfqpoint{15.09cm}{8.108cm}}
\pgfusepath{stroke}
\pgfsetdash{}{0cm}
\pgfpathmoveto{\pgfqpoint{16.213cm}{8.373cm}}
\pgfpathlineto{\pgfqpoint{16.495cm}{8.373cm}}
\pgfusepath{stroke}
\pgfsetdash{}{0cm}
\pgfpathmoveto{\pgfqpoint{16.354cm}{8.514cm}}
\pgfpathlineto{\pgfqpoint{16.354cm}{8.231cm}}
\pgfusepath{stroke}
\pgfsetdash{}{0cm}
\pgfpathmoveto{\pgfqpoint{17.48cm}{8.517cm}}
\pgfpathlineto{\pgfqpoint{17.762cm}{8.517cm}}
\pgfusepath{stroke}
\pgfsetdash{}{0cm}
\pgfpathmoveto{\pgfqpoint{17.621cm}{8.658cm}}
\pgfpathlineto{\pgfqpoint{17.621cm}{8.376cm}}
\pgfusepath{stroke}
\pgfsetdash{}{0cm}
\pgfpathmoveto{\pgfqpoint{6.132cm}{6.421cm}}
\pgfpathlineto{\pgfqpoint{6.326cm}{6.227cm}}
\pgfusepath{stroke}
\pgfsetdash{}{0cm}
\pgfpathmoveto{\pgfqpoint{6.326cm}{6.421cm}}
\pgfpathlineto{\pgfqpoint{6.132cm}{6.227cm}}
\pgfusepath{stroke}
\pgfsetdash{}{0cm}
\pgfpathmoveto{\pgfqpoint{7.4cm}{6.994cm}}
\pgfpathlineto{\pgfqpoint{7.594cm}{6.8cm}}
\pgfusepath{stroke}
\pgfsetdash{}{0cm}
\pgfpathmoveto{\pgfqpoint{7.594cm}{6.994cm}}
\pgfpathlineto{\pgfqpoint{7.4cm}{6.8cm}}
\pgfusepath{stroke}
\pgfsetdash{}{0cm}
\pgfpathmoveto{\pgfqpoint{8.664cm}{7.311cm}}
\pgfpathlineto{\pgfqpoint{8.858cm}{7.117cm}}
\pgfusepath{stroke}
\pgfsetdash{}{0cm}
\pgfpathmoveto{\pgfqpoint{8.858cm}{7.311cm}}
\pgfpathlineto{\pgfqpoint{8.664cm}{7.117cm}}
\pgfusepath{stroke}
\pgfsetdash{}{0cm}
\pgfpathmoveto{\pgfqpoint{9.931cm}{7.214cm}}
\pgfpathlineto{\pgfqpoint{10.125cm}{7.02cm}}
\pgfusepath{stroke}
\pgfsetdash{}{0cm}
\pgfpathmoveto{\pgfqpoint{10.125cm}{7.214cm}}
\pgfpathlineto{\pgfqpoint{9.931cm}{7.02cm}}
\pgfusepath{stroke}
\pgfsetdash{}{0cm}
\pgfpathmoveto{\pgfqpoint{11.195cm}{7.699cm}}
\pgfpathlineto{\pgfqpoint{11.389cm}{7.505cm}}
\pgfusepath{stroke}
\pgfsetdash{}{0cm}
\pgfpathmoveto{\pgfqpoint{11.389cm}{7.699cm}}
\pgfpathlineto{\pgfqpoint{11.195cm}{7.505cm}}
\pgfusepath{stroke}
\pgfsetdash{}{0cm}
\pgfpathmoveto{\pgfqpoint{12.462cm}{8.035cm}}
\pgfpathlineto{\pgfqpoint{12.656cm}{7.84cm}}
\pgfusepath{stroke}
\pgfsetdash{}{0cm}
\pgfpathmoveto{\pgfqpoint{12.656cm}{8.035cm}}
\pgfpathlineto{\pgfqpoint{12.462cm}{7.84cm}}
\pgfusepath{stroke}
\pgfsetdash{}{0cm}
\pgfpathmoveto{\pgfqpoint{13.726cm}{8.223cm}}
\pgfpathlineto{\pgfqpoint{13.92cm}{8.029cm}}
\pgfusepath{stroke}
\pgfsetdash{}{0cm}
\pgfpathmoveto{\pgfqpoint{13.92cm}{8.223cm}}
\pgfpathlineto{\pgfqpoint{13.726cm}{8.029cm}}
\pgfusepath{stroke}
\pgfsetdash{}{0cm}
\pgfpathmoveto{\pgfqpoint{14.993cm}{8.346cm}}
\pgfpathlineto{\pgfqpoint{15.187cm}{8.152cm}}
\pgfusepath{stroke}
\pgfsetdash{}{0cm}
\pgfpathmoveto{\pgfqpoint{15.187cm}{8.346cm}}
\pgfpathlineto{\pgfqpoint{14.993cm}{8.152cm}}
\pgfusepath{stroke}
\pgfsetdash{}{0cm}
\pgfpathmoveto{\pgfqpoint{16.257cm}{8.47cm}}
\pgfpathlineto{\pgfqpoint{16.451cm}{8.276cm}}
\pgfusepath{stroke}
\pgfsetdash{}{0cm}
\pgfpathmoveto{\pgfqpoint{16.451cm}{8.47cm}}
\pgfpathlineto{\pgfqpoint{16.257cm}{8.276cm}}
\pgfusepath{stroke}
\pgfsetdash{}{0cm}
\pgfpathmoveto{\pgfqpoint{17.524cm}{8.614cm}}
\pgfpathlineto{\pgfqpoint{17.718cm}{8.42cm}}
\pgfusepath{stroke}
\pgfsetdash{}{0cm}
\pgfpathmoveto{\pgfqpoint{17.718cm}{8.614cm}}
\pgfpathlineto{\pgfqpoint{17.524cm}{8.42cm}}
\pgfusepath{stroke}
\begin{pgfscope}
\pgfpathmoveto{\pgfqpoint{4.965cm}{22.969cm}}
\pgfpathlineto{\pgfqpoint{17.624cm}{22.969cm}}
\pgfpathlineto{\pgfqpoint{17.624cm}{5.736cm}}
\pgfpathlineto{\pgfqpoint{4.965cm}{5.736cm}}
\pgfpathclose
\pgfusepath{clip}
\pgfsetdash{}{0cm}
\definecolor{eps2pgf_color}{rgb}{0,0,1}\pgfsetstrokecolor{eps2pgf_color}\pgfsetfillcolor{eps2pgf_color}
\pgfpathmoveto{\pgfqpoint{6.229cm}{6.694cm}}
\pgfpathlineto{\pgfqpoint{7.497cm}{7.611cm}}
\pgfpathlineto{\pgfqpoint{8.761cm}{8.496cm}}
\pgfpathlineto{\pgfqpoint{10.028cm}{9.322cm}}
\pgfpathlineto{\pgfqpoint{11.292cm}{10.134cm}}
\pgfpathlineto{\pgfqpoint{12.559cm}{10.913cm}}
\pgfpathlineto{\pgfqpoint{13.823cm}{11.642cm}}
\pgfpathlineto{\pgfqpoint{15.09cm}{12.344cm}}
\pgfpathlineto{\pgfqpoint{16.354cm}{13.017cm}}
\pgfpathlineto{\pgfqpoint{17.621cm}{13.655cm}}
\pgfusepath{stroke}
\end{pgfscope}
\pgfsetdash{}{0cm}
\pgfsetmiterjoin
\definecolor{eps2pgf_color}{rgb}{0,0,1}\pgfsetstrokecolor{eps2pgf_color}\pgfsetfillcolor{eps2pgf_color}
\pgfpathmoveto{\pgfqpoint{6.118cm}{6.806cm}}
\pgfpathlineto{\pgfqpoint{6.341cm}{6.806cm}}
\pgfpathlineto{\pgfqpoint{6.341cm}{6.582cm}}
\pgfpathlineto{\pgfqpoint{6.118cm}{6.582cm}}
\pgfpathlineto{\pgfqpoint{6.118cm}{6.806cm}}
\pgfpathclose
\pgfusepath{stroke}
\pgfsetdash{}{0cm}
\pgfpathmoveto{\pgfqpoint{7.385cm}{7.723cm}}
\pgfpathlineto{\pgfqpoint{7.608cm}{7.723cm}}
\pgfpathlineto{\pgfqpoint{7.608cm}{7.499cm}}
\pgfpathlineto{\pgfqpoint{7.385cm}{7.499cm}}
\pgfpathlineto{\pgfqpoint{7.385cm}{7.723cm}}
\pgfpathclose
\pgfusepath{stroke}
\pgfsetdash{}{0cm}
\pgfpathmoveto{\pgfqpoint{8.649cm}{8.608cm}}
\pgfpathlineto{\pgfqpoint{8.872cm}{8.608cm}}
\pgfpathlineto{\pgfqpoint{8.872cm}{8.384cm}}
\pgfpathlineto{\pgfqpoint{8.649cm}{8.384cm}}
\pgfpathlineto{\pgfqpoint{8.649cm}{8.608cm}}
\pgfpathclose
\pgfusepath{stroke}
\pgfsetdash{}{0cm}
\pgfpathmoveto{\pgfqpoint{9.916cm}{9.434cm}}
\pgfpathlineto{\pgfqpoint{10.139cm}{9.434cm}}
\pgfpathlineto{\pgfqpoint{10.139cm}{9.21cm}}
\pgfpathlineto{\pgfqpoint{9.916cm}{9.21cm}}
\pgfpathlineto{\pgfqpoint{9.916cm}{9.434cm}}
\pgfpathclose
\pgfusepath{stroke}
\pgfsetdash{}{0cm}
\pgfpathmoveto{\pgfqpoint{11.18cm}{10.245cm}}
\pgfpathlineto{\pgfqpoint{11.404cm}{10.245cm}}
\pgfpathlineto{\pgfqpoint{11.404cm}{10.022cm}}
\pgfpathlineto{\pgfqpoint{11.18cm}{10.022cm}}
\pgfpathlineto{\pgfqpoint{11.18cm}{10.245cm}}
\pgfpathclose
\pgfusepath{stroke}
\pgfsetdash{}{0cm}
\pgfpathmoveto{\pgfqpoint{12.447cm}{11.024cm}}
\pgfpathlineto{\pgfqpoint{12.671cm}{11.024cm}}
\pgfpathlineto{\pgfqpoint{12.671cm}{10.801cm}}
\pgfpathlineto{\pgfqpoint{12.447cm}{10.801cm}}
\pgfpathlineto{\pgfqpoint{12.447cm}{11.024cm}}
\pgfpathclose
\pgfusepath{stroke}
\pgfsetdash{}{0cm}
\pgfpathmoveto{\pgfqpoint{13.711cm}{11.753cm}}
\pgfpathlineto{\pgfqpoint{13.935cm}{11.753cm}}
\pgfpathlineto{\pgfqpoint{13.935cm}{11.53cm}}
\pgfpathlineto{\pgfqpoint{13.711cm}{11.53cm}}
\pgfpathlineto{\pgfqpoint{13.711cm}{11.753cm}}
\pgfpathclose
\pgfusepath{stroke}
\pgfsetdash{}{0cm}
\pgfpathmoveto{\pgfqpoint{14.978cm}{12.456cm}}
\pgfpathlineto{\pgfqpoint{15.202cm}{12.456cm}}
\pgfpathlineto{\pgfqpoint{15.202cm}{12.233cm}}
\pgfpathlineto{\pgfqpoint{14.978cm}{12.233cm}}
\pgfpathlineto{\pgfqpoint{14.978cm}{12.456cm}}
\pgfpathclose
\pgfusepath{stroke}
\pgfsetdash{}{0cm}
\pgfpathmoveto{\pgfqpoint{16.242cm}{13.129cm}}
\pgfpathlineto{\pgfqpoint{16.466cm}{13.129cm}}
\pgfpathlineto{\pgfqpoint{16.466cm}{12.906cm}}
\pgfpathlineto{\pgfqpoint{16.242cm}{12.906cm}}
\pgfpathlineto{\pgfqpoint{16.242cm}{13.129cm}}
\pgfpathclose
\pgfusepath{stroke}
\pgfsetdash{}{0cm}
\pgfpathmoveto{\pgfqpoint{17.51cm}{13.767cm}}
\pgfpathlineto{\pgfqpoint{17.733cm}{13.767cm}}
\pgfpathlineto{\pgfqpoint{17.733cm}{13.544cm}}
\pgfpathlineto{\pgfqpoint{17.51cm}{13.544cm}}
\pgfpathlineto{\pgfqpoint{17.51cm}{13.767cm}}
\pgfpathclose
\pgfusepath{stroke}
\begin{pgfscope}
\pgfpathmoveto{\pgfqpoint{4.965cm}{22.969cm}}
\pgfpathlineto{\pgfqpoint{17.624cm}{22.969cm}}
\pgfpathlineto{\pgfqpoint{17.624cm}{5.736cm}}
\pgfpathlineto{\pgfqpoint{4.965cm}{5.736cm}}
\pgfpathclose
\pgfusepath{clip}
\pgfsetdash{}{0cm}
\definecolor{eps2pgf_color}{rgb}{0,1,1}\pgfsetstrokecolor{eps2pgf_color}\pgfsetfillcolor{eps2pgf_color}
\pgfpathmoveto{\pgfqpoint{6.229cm}{6.729cm}}
\pgfpathlineto{\pgfqpoint{7.497cm}{7.676cm}}
\pgfpathlineto{\pgfqpoint{8.761cm}{8.587cm}}
\pgfpathlineto{\pgfqpoint{10.028cm}{9.457cm}}
\pgfpathlineto{\pgfqpoint{11.292cm}{10.31cm}}
\pgfpathlineto{\pgfqpoint{12.559cm}{11.077cm}}
\pgfpathlineto{\pgfqpoint{13.823cm}{11.847cm}}
\pgfpathlineto{\pgfqpoint{15.09cm}{12.603cm}}
\pgfpathlineto{\pgfqpoint{16.354cm}{13.309cm}}
\pgfpathlineto{\pgfqpoint{17.621cm}{13.982cm}}
\pgfusepath{stroke}
\end{pgfscope}
\pgfsetdash{}{0cm}
\definecolor{eps2pgf_color}{rgb}{0,1,1}\pgfsetstrokecolor{eps2pgf_color}\pgfsetfillcolor{eps2pgf_color}
\pgfpathmoveto{\pgfqpoint{6.132cm}{6.826cm}}
\pgfpathlineto{\pgfqpoint{6.326cm}{6.632cm}}
\pgfusepath{stroke}
\pgfsetdash{}{0cm}
\pgfpathmoveto{\pgfqpoint{6.326cm}{6.826cm}}
\pgfpathlineto{\pgfqpoint{6.132cm}{6.632cm}}
\pgfusepath{stroke}
\pgfsetdash{}{0cm}
\pgfpathmoveto{\pgfqpoint{7.4cm}{7.773cm}}
\pgfpathlineto{\pgfqpoint{7.594cm}{7.579cm}}
\pgfusepath{stroke}
\pgfsetdash{}{0cm}
\pgfpathmoveto{\pgfqpoint{7.594cm}{7.773cm}}
\pgfpathlineto{\pgfqpoint{7.4cm}{7.579cm}}
\pgfusepath{stroke}
\pgfsetdash{}{0cm}
\pgfpathmoveto{\pgfqpoint{8.664cm}{8.684cm}}
\pgfpathlineto{\pgfqpoint{8.858cm}{8.49cm}}
\pgfusepath{stroke}
\pgfsetdash{}{0cm}
\pgfpathmoveto{\pgfqpoint{8.858cm}{8.684cm}}
\pgfpathlineto{\pgfqpoint{8.664cm}{8.49cm}}
\pgfusepath{stroke}
\pgfsetdash{}{0cm}
\pgfpathmoveto{\pgfqpoint{9.931cm}{9.554cm}}
\pgfpathlineto{\pgfqpoint{10.125cm}{9.36cm}}
\pgfusepath{stroke}
\pgfsetdash{}{0cm}
\pgfpathmoveto{\pgfqpoint{10.125cm}{9.554cm}}
\pgfpathlineto{\pgfqpoint{9.931cm}{9.36cm}}
\pgfusepath{stroke}
\pgfsetdash{}{0cm}
\pgfpathmoveto{\pgfqpoint{11.195cm}{10.407cm}}
\pgfpathlineto{\pgfqpoint{11.389cm}{10.213cm}}
\pgfusepath{stroke}
\pgfsetdash{}{0cm}
\pgfpathmoveto{\pgfqpoint{11.389cm}{10.407cm}}
\pgfpathlineto{\pgfqpoint{11.195cm}{10.213cm}}
\pgfusepath{stroke}
\pgfsetdash{}{0cm}
\pgfpathmoveto{\pgfqpoint{12.462cm}{11.174cm}}
\pgfpathlineto{\pgfqpoint{12.656cm}{10.98cm}}
\pgfusepath{stroke}
\pgfsetdash{}{0cm}
\pgfpathmoveto{\pgfqpoint{12.656cm}{11.174cm}}
\pgfpathlineto{\pgfqpoint{12.462cm}{10.98cm}}
\pgfusepath{stroke}
\pgfsetdash{}{0cm}
\pgfpathmoveto{\pgfqpoint{13.726cm}{11.944cm}}
\pgfpathlineto{\pgfqpoint{13.92cm}{11.75cm}}
\pgfusepath{stroke}
\pgfsetdash{}{0cm}
\pgfpathmoveto{\pgfqpoint{13.92cm}{11.944cm}}
\pgfpathlineto{\pgfqpoint{13.726cm}{11.75cm}}
\pgfusepath{stroke}
\pgfsetdash{}{0cm}
\pgfpathmoveto{\pgfqpoint{14.993cm}{12.7cm}}
\pgfpathlineto{\pgfqpoint{15.187cm}{12.506cm}}
\pgfusepath{stroke}
\pgfsetdash{}{0cm}
\pgfpathmoveto{\pgfqpoint{15.187cm}{12.7cm}}
\pgfpathlineto{\pgfqpoint{14.993cm}{12.506cm}}
\pgfusepath{stroke}
\pgfsetdash{}{0cm}
\pgfpathmoveto{\pgfqpoint{16.257cm}{13.406cm}}
\pgfpathlineto{\pgfqpoint{16.451cm}{13.212cm}}
\pgfusepath{stroke}
\pgfsetdash{}{0cm}
\pgfpathmoveto{\pgfqpoint{16.451cm}{13.406cm}}
\pgfpathlineto{\pgfqpoint{16.257cm}{13.212cm}}
\pgfusepath{stroke}
\pgfsetdash{}{0cm}
\pgfpathmoveto{\pgfqpoint{17.524cm}{14.079cm}}
\pgfpathlineto{\pgfqpoint{17.718cm}{13.885cm}}
\pgfusepath{stroke}
\pgfsetdash{}{0cm}
\pgfpathmoveto{\pgfqpoint{17.718cm}{14.079cm}}
\pgfpathlineto{\pgfqpoint{17.524cm}{13.885cm}}
\pgfusepath{stroke}
\begin{pgfscope}
\pgfpathmoveto{\pgfqpoint{4.965cm}{22.969cm}}
\pgfpathlineto{\pgfqpoint{17.624cm}{22.969cm}}
\pgfpathlineto{\pgfqpoint{17.624cm}{5.736cm}}
\pgfpathlineto{\pgfqpoint{4.965cm}{5.736cm}}
\pgfpathclose
\pgfusepath{clip}
\pgfsetdash{{0.212cm}}{0cm}
\definecolor{eps2pgf_color}{rgb}{0,0,1}\pgfsetstrokecolor{eps2pgf_color}\pgfsetfillcolor{eps2pgf_color}
\pgfpathmoveto{\pgfqpoint{6.229cm}{7.044cm}}
\pgfpathlineto{\pgfqpoint{7.497cm}{8.314cm}}
\pgfpathlineto{\pgfqpoint{8.761cm}{9.549cm}}
\pgfpathlineto{\pgfqpoint{10.028cm}{10.727cm}}
\pgfpathlineto{\pgfqpoint{11.292cm}{11.889cm}}
\pgfpathlineto{\pgfqpoint{12.559cm}{13.017cm}}
\pgfpathlineto{\pgfqpoint{13.823cm}{14.096cm}}
\pgfpathlineto{\pgfqpoint{15.09cm}{15.149cm}}
\pgfpathlineto{\pgfqpoint{16.354cm}{16.175cm}}
\pgfpathlineto{\pgfqpoint{17.621cm}{17.16cm}}
\pgfusepath{stroke}
\end{pgfscope}
\pgfsetdash{}{0cm}
\definecolor{eps2pgf_color}{rgb}{0,0,1}\pgfsetstrokecolor{eps2pgf_color}\pgfsetfillcolor{eps2pgf_color}
\pgfpathmoveto{\pgfqpoint{6.118cm}{7.156cm}}
\pgfpathlineto{\pgfqpoint{6.341cm}{7.156cm}}
\pgfpathlineto{\pgfqpoint{6.341cm}{6.932cm}}
\pgfpathlineto{\pgfqpoint{6.118cm}{6.932cm}}
\pgfpathlineto{\pgfqpoint{6.118cm}{7.156cm}}
\pgfpathclose
\pgfusepath{stroke}
\pgfsetdash{}{0cm}
\pgfpathmoveto{\pgfqpoint{7.385cm}{8.426cm}}
\pgfpathlineto{\pgfqpoint{7.608cm}{8.426cm}}
\pgfpathlineto{\pgfqpoint{7.608cm}{8.202cm}}
\pgfpathlineto{\pgfqpoint{7.385cm}{8.202cm}}
\pgfpathlineto{\pgfqpoint{7.385cm}{8.426cm}}
\pgfpathclose
\pgfusepath{stroke}
\pgfsetdash{}{0cm}
\pgfpathmoveto{\pgfqpoint{8.649cm}{9.66cm}}
\pgfpathlineto{\pgfqpoint{8.872cm}{9.66cm}}
\pgfpathlineto{\pgfqpoint{8.872cm}{9.437cm}}
\pgfpathlineto{\pgfqpoint{8.649cm}{9.437cm}}
\pgfpathlineto{\pgfqpoint{8.649cm}{9.66cm}}
\pgfpathclose
\pgfusepath{stroke}
\pgfsetdash{}{0cm}
\pgfpathmoveto{\pgfqpoint{9.916cm}{10.839cm}}
\pgfpathlineto{\pgfqpoint{10.139cm}{10.839cm}}
\pgfpathlineto{\pgfqpoint{10.139cm}{10.616cm}}
\pgfpathlineto{\pgfqpoint{9.916cm}{10.616cm}}
\pgfpathlineto{\pgfqpoint{9.916cm}{10.839cm}}
\pgfpathclose
\pgfusepath{stroke}
\pgfsetdash{}{0cm}
\pgfpathmoveto{\pgfqpoint{11.18cm}{12cm}}
\pgfpathlineto{\pgfqpoint{11.404cm}{12cm}}
\pgfpathlineto{\pgfqpoint{11.404cm}{11.777cm}}
\pgfpathlineto{\pgfqpoint{11.18cm}{11.777cm}}
\pgfpathlineto{\pgfqpoint{11.18cm}{12cm}}
\pgfpathclose
\pgfusepath{stroke}
\pgfsetdash{}{0cm}
\pgfpathmoveto{\pgfqpoint{12.447cm}{13.129cm}}
\pgfpathlineto{\pgfqpoint{12.671cm}{13.129cm}}
\pgfpathlineto{\pgfqpoint{12.671cm}{12.906cm}}
\pgfpathlineto{\pgfqpoint{12.447cm}{12.906cm}}
\pgfpathlineto{\pgfqpoint{12.447cm}{13.129cm}}
\pgfpathclose
\pgfusepath{stroke}
\pgfsetdash{}{0cm}
\pgfpathmoveto{\pgfqpoint{13.711cm}{14.208cm}}
\pgfpathlineto{\pgfqpoint{13.935cm}{14.208cm}}
\pgfpathlineto{\pgfqpoint{13.935cm}{13.985cm}}
\pgfpathlineto{\pgfqpoint{13.711cm}{13.985cm}}
\pgfpathlineto{\pgfqpoint{13.711cm}{14.208cm}}
\pgfpathclose
\pgfusepath{stroke}
\pgfsetdash{}{0cm}
\pgfpathmoveto{\pgfqpoint{14.978cm}{15.261cm}}
\pgfpathlineto{\pgfqpoint{15.202cm}{15.261cm}}
\pgfpathlineto{\pgfqpoint{15.202cm}{15.037cm}}
\pgfpathlineto{\pgfqpoint{14.978cm}{15.037cm}}
\pgfpathlineto{\pgfqpoint{14.978cm}{15.261cm}}
\pgfpathclose
\pgfusepath{stroke}
\pgfsetdash{}{0cm}
\pgfpathmoveto{\pgfqpoint{16.242cm}{16.287cm}}
\pgfpathlineto{\pgfqpoint{16.466cm}{16.287cm}}
\pgfpathlineto{\pgfqpoint{16.466cm}{16.063cm}}
\pgfpathlineto{\pgfqpoint{16.242cm}{16.063cm}}
\pgfpathlineto{\pgfqpoint{16.242cm}{16.287cm}}
\pgfpathclose
\pgfusepath{stroke}
\pgfsetdash{}{0cm}
\pgfpathmoveto{\pgfqpoint{17.51cm}{17.271cm}}
\pgfpathlineto{\pgfqpoint{17.733cm}{17.271cm}}
\pgfpathlineto{\pgfqpoint{17.733cm}{17.048cm}}
\pgfpathlineto{\pgfqpoint{17.51cm}{17.048cm}}
\pgfpathlineto{\pgfqpoint{17.51cm}{17.271cm}}
\pgfpathclose
\pgfusepath{stroke}
\begin{pgfscope}
\pgfpathmoveto{\pgfqpoint{4.965cm}{22.969cm}}
\pgfpathlineto{\pgfqpoint{17.624cm}{22.969cm}}
\pgfpathlineto{\pgfqpoint{17.624cm}{5.736cm}}
\pgfpathlineto{\pgfqpoint{4.965cm}{5.736cm}}
\pgfpathclose
\pgfusepath{clip}
\pgfsetdash{{0.212cm}}{0cm}
\definecolor{eps2pgf_color}{rgb}{0,1,1}\pgfsetstrokecolor{eps2pgf_color}\pgfsetfillcolor{eps2pgf_color}
\pgfpathmoveto{\pgfqpoint{6.229cm}{7.088cm}}
\pgfpathlineto{\pgfqpoint{7.497cm}{8.39cm}}
\pgfpathlineto{\pgfqpoint{8.761cm}{9.66cm}}
\pgfpathlineto{\pgfqpoint{10.028cm}{10.886cm}}
\pgfpathlineto{\pgfqpoint{11.292cm}{12.094cm}}
\pgfpathlineto{\pgfqpoint{12.559cm}{13.22cm}}
\pgfpathlineto{\pgfqpoint{13.823cm}{14.346cm}}
\pgfpathlineto{\pgfqpoint{15.09cm}{15.46cm}}
\pgfpathlineto{\pgfqpoint{16.354cm}{16.522cm}}
\pgfpathlineto{\pgfqpoint{17.621cm}{17.554cm}}
\pgfusepath{stroke}
\end{pgfscope}
\pgfsetdash{}{0cm}
\definecolor{eps2pgf_color}{rgb}{0,1,1}\pgfsetstrokecolor{eps2pgf_color}\pgfsetfillcolor{eps2pgf_color}
\pgfpathmoveto{\pgfqpoint{6.132cm}{7.185cm}}
\pgfpathlineto{\pgfqpoint{6.326cm}{6.991cm}}
\pgfusepath{stroke}
\pgfsetdash{}{0cm}
\pgfpathmoveto{\pgfqpoint{6.326cm}{7.185cm}}
\pgfpathlineto{\pgfqpoint{6.132cm}{6.991cm}}
\pgfusepath{stroke}
\pgfsetdash{}{0cm}
\pgfpathmoveto{\pgfqpoint{7.4cm}{8.487cm}}
\pgfpathlineto{\pgfqpoint{7.594cm}{8.293cm}}
\pgfusepath{stroke}
\pgfsetdash{}{0cm}
\pgfpathmoveto{\pgfqpoint{7.594cm}{8.487cm}}
\pgfpathlineto{\pgfqpoint{7.4cm}{8.293cm}}
\pgfusepath{stroke}
\pgfsetdash{}{0cm}
\pgfpathmoveto{\pgfqpoint{8.664cm}{9.757cm}}
\pgfpathlineto{\pgfqpoint{8.858cm}{9.563cm}}
\pgfusepath{stroke}
\pgfsetdash{}{0cm}
\pgfpathmoveto{\pgfqpoint{8.858cm}{9.757cm}}
\pgfpathlineto{\pgfqpoint{8.664cm}{9.563cm}}
\pgfusepath{stroke}
\pgfsetdash{}{0cm}
\pgfpathmoveto{\pgfqpoint{9.931cm}{10.983cm}}
\pgfpathlineto{\pgfqpoint{10.125cm}{10.789cm}}
\pgfusepath{stroke}
\pgfsetdash{}{0cm}
\pgfpathmoveto{\pgfqpoint{10.125cm}{10.983cm}}
\pgfpathlineto{\pgfqpoint{9.931cm}{10.789cm}}
\pgfusepath{stroke}
\pgfsetdash{}{0cm}
\pgfpathmoveto{\pgfqpoint{11.195cm}{12.191cm}}
\pgfpathlineto{\pgfqpoint{11.389cm}{11.997cm}}
\pgfusepath{stroke}
\pgfsetdash{}{0cm}
\pgfpathmoveto{\pgfqpoint{11.389cm}{12.191cm}}
\pgfpathlineto{\pgfqpoint{11.195cm}{11.997cm}}
\pgfusepath{stroke}
\pgfsetdash{}{0cm}
\pgfpathmoveto{\pgfqpoint{12.462cm}{13.317cm}}
\pgfpathlineto{\pgfqpoint{12.656cm}{13.123cm}}
\pgfusepath{stroke}
\pgfsetdash{}{0cm}
\pgfpathmoveto{\pgfqpoint{12.656cm}{13.317cm}}
\pgfpathlineto{\pgfqpoint{12.462cm}{13.123cm}}
\pgfusepath{stroke}
\pgfsetdash{}{0cm}
\pgfpathmoveto{\pgfqpoint{13.726cm}{14.443cm}}
\pgfpathlineto{\pgfqpoint{13.92cm}{14.249cm}}
\pgfusepath{stroke}
\pgfsetdash{}{0cm}
\pgfpathmoveto{\pgfqpoint{13.92cm}{14.443cm}}
\pgfpathlineto{\pgfqpoint{13.726cm}{14.249cm}}
\pgfusepath{stroke}
\pgfsetdash{}{0cm}
\pgfpathmoveto{\pgfqpoint{14.993cm}{15.558cm}}
\pgfpathlineto{\pgfqpoint{15.187cm}{15.363cm}}
\pgfusepath{stroke}
\pgfsetdash{}{0cm}
\pgfpathmoveto{\pgfqpoint{15.187cm}{15.558cm}}
\pgfpathlineto{\pgfqpoint{14.993cm}{15.363cm}}
\pgfusepath{stroke}
\pgfsetdash{}{0cm}
\pgfpathmoveto{\pgfqpoint{16.257cm}{16.619cm}}
\pgfpathlineto{\pgfqpoint{16.451cm}{16.425cm}}
\pgfusepath{stroke}
\pgfsetdash{}{0cm}
\pgfpathmoveto{\pgfqpoint{16.451cm}{16.619cm}}
\pgfpathlineto{\pgfqpoint{16.257cm}{16.425cm}}
\pgfusepath{stroke}
\pgfsetdash{}{0cm}
\pgfpathmoveto{\pgfqpoint{17.524cm}{17.651cm}}
\pgfpathlineto{\pgfqpoint{17.718cm}{17.457cm}}
\pgfusepath{stroke}
\pgfsetdash{}{0cm}
\pgfpathmoveto{\pgfqpoint{17.718cm}{17.651cm}}
\pgfpathlineto{\pgfqpoint{17.524cm}{17.457cm}}
\pgfusepath{stroke}
\begin{pgfscope}
\pgfpathmoveto{\pgfqpoint{4.965cm}{22.969cm}}
\pgfpathlineto{\pgfqpoint{17.624cm}{22.969cm}}
\pgfpathlineto{\pgfqpoint{17.624cm}{5.736cm}}
\pgfpathlineto{\pgfqpoint{4.965cm}{5.736cm}}
\pgfpathclose
\pgfusepath{clip}
\end{pgfscope}
\definecolor{eps2pgf_color}{gray}{0}\pgfsetstrokecolor{eps2pgf_color}\pgfsetfillcolor{eps2pgf_color}
\pgftext[x=11.31cm,y=4.093cm,rotate=0]{\fontsize{36}{36.14}\selectfont{ {$\bmax$}}}
\pgftext[x=3.013cm,y=14.366cm,rotate=90]{\fontsize{36}{36.14}\selectfont{ {Percentage cost saving (\%)}}}
\pgftext[x=4.915cm,y=5.612cm,rotate=0]{\fontsize{10.04}{12.04}\selectfont{ { }}}
\pgftext[x=17.574cm,y=22.845cm,rotate=0]{\fontsize{10.04}{12.04}\selectfont{ { }}}
\pgftext[x=8.252cm-.2cm,y=22.014cm,rotate=0]{\fontsize{36}{36.14}\selectfont{ {Greedy}}}
\begin{pgfscope}
\pgfpathmoveto{\pgfqpoint{5.142cm}{22.792cm}}
\pgfpathlineto{\pgfqpoint{16.275cm}{22.792cm}}
\pgfpathlineto{\pgfqpoint{16.275cm}{16.434cm}}
\pgfpathlineto{\pgfqpoint{5.142cm}{16.434cm}}
\pgfpathclose
\pgfusepath{clip}
\pgfsetdash{}{0cm}
\definecolor{eps2pgf_color}{rgb}{0,1,0}\pgfsetstrokecolor{eps2pgf_color}\pgfsetfillcolor{eps2pgf_color}
\pgfpathmoveto{\pgfqpoint{5.348cm}{22.113cm}}
\pgfpathlineto{\pgfqpoint{6.391cm}{22.113cm}}
\pgfusepath{stroke}
\begin{pgfscope}
\pgfpathmoveto{\pgfqpoint{5.477cm}{22.507cm}}
\pgfpathlineto{\pgfqpoint{6.268cm}{22.507cm}}
\pgfpathlineto{\pgfqpoint{6.268cm}{21.716cm}}
\pgfpathlineto{\pgfqpoint{5.477cm}{21.716cm}}
\pgfpathclose
\pgfusepath{clip}
\pgfsetdash{}{0cm}
\pgfpathmoveto{\pgfqpoint{5.73cm}{22.113cm}}
\pgfpathlineto{\pgfqpoint{6.012cm}{22.113cm}}
\pgfusepath{stroke}
\pgfsetdash{}{0cm}
\pgfpathmoveto{\pgfqpoint{5.871cm}{22.254cm}}
\pgfpathlineto{\pgfqpoint{5.871cm}{21.972cm}}
\pgfusepath{stroke}
\pgfsetdash{}{0cm}
\pgfpathmoveto{\pgfqpoint{5.774cm}{22.21cm}}
\pgfpathlineto{\pgfqpoint{5.968cm}{22.016cm}}
\pgfusepath{stroke}
\pgfsetdash{}{0cm}
\pgfpathmoveto{\pgfqpoint{5.968cm}{22.21cm}}
\pgfpathlineto{\pgfqpoint{5.774cm}{22.016cm}}
\pgfusepath{stroke}
\end{pgfscope}
\end{pgfscope}
\pgftext[x=9.304cm+.4cm,y=20.768cm,rotate=0]{\fontsize{36}{36.14}\selectfont{ {OMG (\texttt{minS})}}}
\begin{pgfscope}
\pgfpathmoveto{\pgfqpoint{5.142cm}{22.792cm}}
\pgfpathlineto{\pgfqpoint{16.275cm}{22.792cm}}
\pgfpathlineto{\pgfqpoint{16.275cm}{16.434cm}}
\pgfpathlineto{\pgfqpoint{5.142cm}{16.434cm}}
\pgfpathclose
\pgfusepath{clip}
\pgfsetdash{}{0cm}
\definecolor{eps2pgf_color}{rgb}{0,0,1}\pgfsetstrokecolor{eps2pgf_color}\pgfsetfillcolor{eps2pgf_color}
\pgfpathmoveto{\pgfqpoint{5.348cm}{20.864cm}}
\pgfpathlineto{\pgfqpoint{6.391cm}{20.864cm}}
\pgfusepath{stroke}
\begin{pgfscope}
\pgfpathmoveto{\pgfqpoint{5.477cm}{21.258cm}}
\pgfpathlineto{\pgfqpoint{6.268cm}{21.258cm}}
\pgfpathlineto{\pgfqpoint{6.268cm}{20.467cm}}
\pgfpathlineto{\pgfqpoint{5.477cm}{20.467cm}}
\pgfpathclose
\pgfusepath{clip}
\pgfsetdash{}{0cm}
\pgfpathmoveto{\pgfqpoint{5.759cm}{20.976cm}}
\pgfpathlineto{\pgfqpoint{5.983cm}{20.976cm}}
\pgfpathlineto{\pgfqpoint{5.983cm}{20.752cm}}
\pgfpathlineto{\pgfqpoint{5.759cm}{20.752cm}}
\pgfpathlineto{\pgfqpoint{5.759cm}{20.976cm}}
\pgfpathclose
\pgfusepath{stroke}
\end{pgfscope}
\end{pgfscope}
\pgftext[x=9.597cm+.2cm,y=19.519cm,rotate=0]{\fontsize{36}{36.14}\selectfont{ {OMG (\texttt{maxW})}}}
\begin{pgfscope}
\pgfpathmoveto{\pgfqpoint{5.142cm}{22.792cm}}
\pgfpathlineto{\pgfqpoint{16.275cm}{22.792cm}}
\pgfpathlineto{\pgfqpoint{16.275cm}{16.434cm}}
\pgfpathlineto{\pgfqpoint{5.142cm}{16.434cm}}
\pgfpathclose
\pgfusepath{clip}
\pgfsetdash{}{0cm}
\definecolor{eps2pgf_color}{rgb}{0,1,1}\pgfsetstrokecolor{eps2pgf_color}\pgfsetfillcolor{eps2pgf_color}
\pgfpathmoveto{\pgfqpoint{5.348cm}{19.614cm}}
\pgfpathlineto{\pgfqpoint{6.391cm}{19.614cm}}
\pgfusepath{stroke}
\begin{pgfscope}
\pgfpathmoveto{\pgfqpoint{5.477cm}{20.008cm}}
\pgfpathlineto{\pgfqpoint{6.268cm}{20.008cm}}
\pgfpathlineto{\pgfqpoint{6.268cm}{19.218cm}}
\pgfpathlineto{\pgfqpoint{5.477cm}{19.218cm}}
\pgfpathclose
\pgfusepath{clip}
\pgfsetdash{}{0cm}
\pgfpathmoveto{\pgfqpoint{5.774cm}{19.711cm}}
\pgfpathlineto{\pgfqpoint{5.968cm}{19.517cm}}
\pgfusepath{stroke}
\pgfsetdash{}{0cm}
\pgfpathmoveto{\pgfqpoint{5.968cm}{19.711cm}}
\pgfpathlineto{\pgfqpoint{5.774cm}{19.517cm}}
\pgfusepath{stroke}
\end{pgfscope}
\end{pgfscope}
\pgftext[x=11.031cm+.6cm,y=18.272cm,rotate=0]{\fontsize{36}{36.14}\selectfont{ {Upper bound (\texttt{minS})}}}
\begin{pgfscope}
\pgfpathmoveto{\pgfqpoint{5.142cm}{22.792cm}}
\pgfpathlineto{\pgfqpoint{16.275cm}{22.792cm}}
\pgfpathlineto{\pgfqpoint{16.275cm}{16.434cm}}
\pgfpathlineto{\pgfqpoint{5.142cm}{16.434cm}}
\pgfpathclose
\pgfusepath{clip}
\pgfsetdash{{0.212cm}}{0cm}
\definecolor{eps2pgf_color}{rgb}{0,0,1}\pgfsetstrokecolor{eps2pgf_color}\pgfsetfillcolor{eps2pgf_color}
\pgfpathmoveto{\pgfqpoint{5.348cm}{18.368cm}}
\pgfpathlineto{\pgfqpoint{6.391cm}{18.368cm}}
\pgfusepath{stroke}
\begin{pgfscope}
\pgfpathmoveto{\pgfqpoint{5.477cm}{18.762cm}}
\pgfpathlineto{\pgfqpoint{6.268cm}{18.762cm}}
\pgfpathlineto{\pgfqpoint{6.268cm}{17.971cm}}
\pgfpathlineto{\pgfqpoint{5.477cm}{17.971cm}}
\pgfpathclose
\pgfusepath{clip}
\pgfsetdash{}{0cm}
\pgfpathmoveto{\pgfqpoint{5.759cm}{18.48cm}}
\pgfpathlineto{\pgfqpoint{5.983cm}{18.48cm}}
\pgfpathlineto{\pgfqpoint{5.983cm}{18.256cm}}
\pgfpathlineto{\pgfqpoint{5.759cm}{18.256cm}}
\pgfpathlineto{\pgfqpoint{5.759cm}{18.48cm}}
\pgfpathclose
\pgfusepath{stroke}
\end{pgfscope}
\end{pgfscope}
\pgftext[x=11.325cm+.3cm,y=17.021cm,rotate=0]{\fontsize{36}{36.14}\selectfont{ {Upper bound (\texttt{maxW})}}}
\begin{pgfscope}
\pgfpathmoveto{\pgfqpoint{5.142cm}{22.792cm}}
\pgfpathlineto{\pgfqpoint{16.275cm}{22.792cm}}
\pgfpathlineto{\pgfqpoint{16.275cm}{16.434cm}}
\pgfpathlineto{\pgfqpoint{5.142cm}{16.434cm}}
\pgfpathclose
\pgfusepath{clip}
\pgfsetdash{{0.212cm}}{0cm}
\definecolor{eps2pgf_color}{rgb}{0,1,1}\pgfsetstrokecolor{eps2pgf_color}\pgfsetfillcolor{eps2pgf_color}
\pgfpathmoveto{\pgfqpoint{5.348cm}{17.119cm}}
\pgfpathlineto{\pgfqpoint{6.391cm}{17.119cm}}
\pgfusepath{stroke}
\begin{pgfscope}
\pgfpathmoveto{\pgfqpoint{5.477cm}{17.512cm}}
\pgfpathlineto{\pgfqpoint{6.268cm}{17.512cm}}
\pgfpathlineto{\pgfqpoint{6.268cm}{16.722cm}}
\pgfpathlineto{\pgfqpoint{5.477cm}{16.722cm}}
\pgfpathclose
\pgfusepath{clip}
\pgfsetdash{}{0cm}
\pgfpathmoveto{\pgfqpoint{5.774cm}{17.216cm}}
\pgfpathlineto{\pgfqpoint{5.968cm}{17.022cm}}
\pgfusepath{stroke}
\pgfsetdash{}{0cm}
\pgfpathmoveto{\pgfqpoint{5.968cm}{17.216cm}}
\pgfpathlineto{\pgfqpoint{5.774cm}{17.022cm}}
\pgfusepath{stroke}
\end{pgfscope}
\end{pgfscope}
\pgfsetdash{}{0cm}
\pgfsetlinewidth{0.176mm}
\definecolor{eps2pgf_color}{rgb}{0,1,1}\pgfsetstrokecolor{eps2pgf_color}\pgfsetfillcolor{eps2pgf_color}
\pgfusepath{stroke}
\end{pgfscope}
\end{pgfscope}
\end{pgfpicture}

%% file: I4_p.tex
\scalebox{0.35}{\scalefont{1.28} \input{./fig/I4_p_short.pgf}}

%% file: fig/I4_p_short.pgf
% Created by Eps2pgf 0.7.0 (build on 2008-08-24) on Thu May 22 19:24:26 PDT 2014
\begin{pgfpicture}
\pgfpathmoveto{\pgfqpoint{-1.87cm}{9.737cm}}
\pgfpathlineto{\pgfqpoint{23.495cm}{9.737cm}}
\pgfpathlineto{\pgfqpoint{23.495cm}{18.203cm}}
\pgfpathlineto{\pgfqpoint{-1.87cm}{18.203cm}}
\pgfpathclose
\pgfusepath{clip}
\begin{pgfscope}
\begin{pgfscope}
\pgfpathmoveto{\pgfqpoint{-1.87cm}{18.203cm}}
\pgfpathlineto{\pgfqpoint{23.527cm}{18.203cm}}
\pgfpathlineto{\pgfqpoint{23.527cm}{9.745cm}}
\pgfpathlineto{\pgfqpoint{-1.87cm}{9.745cm}}
\pgfpathclose
\pgfusepath{clip}
\begin{pgfscope}
\definecolor{eps2pgf_color}{gray}{1}\pgfsetstrokecolor{eps2pgf_color}\pgfsetfillcolor{eps2pgf_color}
\pgfpathmoveto{\pgfqpoint{-1.87cm}{18.203cm}}
\pgfpathlineto{\pgfqpoint{23.53cm}{18.203cm}}
\pgfpathlineto{\pgfqpoint{23.53cm}{9.743cm}}
\pgfpathlineto{\pgfqpoint{-1.87cm}{9.743cm}}
\pgfpathclose
\pgfusepath{fill}
\end{pgfscope}
\definecolor{eps2pgf_color}{gray}{1}\pgfsetstrokecolor{eps2pgf_color}\pgfsetfillcolor{eps2pgf_color}
\pgfpathmoveto{\pgfqpoint{1.432cm}{11.168cm}}
\pgfpathlineto{\pgfqpoint{1.432cm}{17.992cm}}
\pgfpathlineto{\pgfqpoint{23.175cm}{17.992cm}}
\pgfpathlineto{\pgfqpoint{23.175cm}{11.168cm}}
\pgfpathclose
\pgfseteorule\pgfusepath{fill}\pgfsetnonzerorule
\pgfsetdash{}{0cm}
\pgfsetlinewidth{0.176mm}
\pgfsetroundjoin
\pgfpathmoveto{\pgfqpoint{1.432cm}{11.168cm}}
\pgfpathlineto{\pgfqpoint{1.432cm}{17.992cm}}
\pgfpathlineto{\pgfqpoint{23.175cm}{17.992cm}}
\pgfpathlineto{\pgfqpoint{23.175cm}{11.168cm}}
\pgfpathlineto{\pgfqpoint{1.432cm}{11.168cm}}
\pgfusepath{stroke}
\pgfsetdash{}{0cm}
\definecolor{eps2pgf_color}{gray}{0}\pgfsetstrokecolor{eps2pgf_color}\pgfsetfillcolor{eps2pgf_color}
\pgfpathmoveto{\pgfqpoint{1.432cm}{11.168cm}}
\pgfpathlineto{\pgfqpoint{23.175cm}{11.168cm}}
\pgfusepath{stroke}
\pgfsetdash{}{0cm}
\pgfpathmoveto{\pgfqpoint{1.432cm}{17.992cm}}
\pgfpathlineto{\pgfqpoint{23.175cm}{17.992cm}}
\pgfusepath{stroke}
\pgfsetdash{}{0cm}
\pgfpathmoveto{\pgfqpoint{1.432cm}{11.168cm}}
\pgfpathlineto{\pgfqpoint{1.432cm}{17.992cm}}
\pgfusepath{stroke}
\pgfsetdash{}{0cm}
\pgfpathmoveto{\pgfqpoint{23.175cm}{11.168cm}}
\pgfpathlineto{\pgfqpoint{23.175cm}{17.992cm}}
\pgfusepath{stroke}
\pgfsetdash{}{0cm}
\pgfpathmoveto{\pgfqpoint{1.432cm}{11.168cm}}
\pgfpathlineto{\pgfqpoint{23.175cm}{11.168cm}}
\pgfusepath{stroke}
\pgfsetdash{}{0cm}
\pgfpathmoveto{\pgfqpoint{1.432cm}{11.168cm}}
\pgfpathlineto{\pgfqpoint{1.432cm}{17.992cm}}
\pgfusepath{stroke}
\pgfsetdash{}{0cm}
\pgfpathmoveto{\pgfqpoint{1.432cm}{11.48cm}}
\pgfpathlineto{\pgfqpoint{1.646cm}{11.48cm}}
\pgfusepath{stroke}
\pgfsetdash{}{0cm}
\pgfpathmoveto{\pgfqpoint{23.175cm}{11.48cm}}
\pgfpathlineto{\pgfqpoint{22.957cm}{11.48cm}}
\pgfusepath{stroke}
\pgftext[x=1.094cm,y=11.455cm,rotate=0]{\fontsize{20}{7.06}\selectfont{{0}}}
\pgfsetdash{}{0cm}
\pgfpathmoveto{\pgfqpoint{1.432cm}{13.55cm}}
\pgfpathlineto{\pgfqpoint{1.646cm}{13.55cm}}
\pgfusepath{stroke}
\pgfsetdash{}{0cm}
\pgfpathmoveto{\pgfqpoint{23.175cm}{13.55cm}}
\pgfpathlineto{\pgfqpoint{22.957cm}{13.55cm}}
\pgfusepath{stroke}
\pgftext[x=0.857cm,y=13.525cm,rotate=0]{\fontsize{20}{7.06}\selectfont{{50}}}
\pgfsetdash{}{0cm}
\pgfpathmoveto{\pgfqpoint{1.432cm}{15.622cm}}
\pgfpathlineto{\pgfqpoint{1.646cm}{15.622cm}}
\pgfusepath{stroke}
\pgfsetdash{}{0cm}
\pgfpathmoveto{\pgfqpoint{23.175cm}{15.622cm}}
\pgfpathlineto{\pgfqpoint{22.957cm}{15.622cm}}
\pgfusepath{stroke}
\pgftext[x=0.651cm,y=15.597cm,rotate=0]{\fontsize{20}{7.06}\selectfont{{100}}}
\pgfsetdash{}{0cm}
\pgfpathmoveto{\pgfqpoint{1.432cm}{17.692cm}}
\pgfpathlineto{\pgfqpoint{1.646cm}{17.692cm}}
\pgfusepath{stroke}
\pgfsetdash{}{0cm}
\pgfpathmoveto{\pgfqpoint{23.175cm}{17.692cm}}
\pgfpathlineto{\pgfqpoint{22.957cm}{17.692cm}}
\pgfusepath{stroke}
\pgftext[x=0.651cm,y=17.667cm,rotate=0]{\fontsize{20}{7.06}\selectfont{{150}}}
\pgfsetdash{}{0cm}
\pgfpathmoveto{\pgfqpoint{1.432cm}{11.168cm}}
\pgfpathlineto{\pgfqpoint{23.175cm}{11.168cm}}
\pgfusepath{stroke}
\pgfsetdash{}{0cm}
\pgfpathmoveto{\pgfqpoint{1.432cm}{17.992cm}}
\pgfpathlineto{\pgfqpoint{23.175cm}{17.992cm}}
\pgfusepath{stroke}
\pgfsetdash{}{0cm}
\pgfpathmoveto{\pgfqpoint{1.432cm}{11.168cm}}
\pgfpathlineto{\pgfqpoint{1.432cm}{17.992cm}}
\pgfusepath{stroke}
\pgfsetdash{}{0cm}
\pgfpathmoveto{\pgfqpoint{23.175cm}{11.168cm}}
\pgfpathlineto{\pgfqpoint{23.175cm}{17.992cm}}
\pgfusepath{stroke}
\begin{pgfscope}
\pgfpathmoveto{\pgfqpoint{1.432cm}{17.992cm}}
\pgfpathlineto{\pgfqpoint{23.178cm}{17.992cm}}
\pgfpathlineto{\pgfqpoint{23.178cm}{11.165cm}}
\pgfpathlineto{\pgfqpoint{1.432cm}{11.165cm}}
\pgfpathclose
\pgfusepath{clip}
\pgfsetdash{{0.212cm}}{0cm}
\pgfpathmoveto{\pgfqpoint{1.884cm}{13.026cm}}
\pgfpathlineto{\pgfqpoint{1.884cm}{14.073cm}}
\pgfusepath{stroke}
\pgfsetdash{{0.212cm}}{0cm}
\pgfpathmoveto{\pgfqpoint{2.79cm}{13.026cm}}
\pgfpathlineto{\pgfqpoint{2.79cm}{14.073cm}}
\pgfusepath{stroke}
\pgfsetdash{{0.212cm}}{0cm}
\pgfpathmoveto{\pgfqpoint{3.695cm}{13.15cm}}
\pgfpathlineto{\pgfqpoint{3.695cm}{13.923cm}}
\pgfusepath{stroke}
\pgfsetdash{{0.212cm}}{0cm}
\pgfpathmoveto{\pgfqpoint{4.601cm}{13.15cm}}
\pgfpathlineto{\pgfqpoint{4.601cm}{13.923cm}}
\pgfusepath{stroke}
\pgfsetdash{{0.212cm}}{0cm}
\pgfpathmoveto{\pgfqpoint{5.506cm}{12.973cm}}
\pgfpathlineto{\pgfqpoint{5.506cm}{13.632cm}}
\pgfusepath{stroke}
\pgfsetdash{{0.212cm}}{0cm}
\pgfpathmoveto{\pgfqpoint{6.412cm}{12.973cm}}
\pgfpathlineto{\pgfqpoint{6.412cm}{13.632cm}}
\pgfusepath{stroke}
\pgfsetdash{{0.212cm}}{0cm}
\pgfpathmoveto{\pgfqpoint{7.32cm}{12.962cm}}
\pgfpathlineto{\pgfqpoint{7.32cm}{13.658cm}}
\pgfusepath{stroke}
\pgfsetdash{{0.212cm}}{0cm}
\pgfpathmoveto{\pgfqpoint{8.226cm}{12.962cm}}
\pgfpathlineto{\pgfqpoint{8.226cm}{13.658cm}}
\pgfusepath{stroke}
\pgfsetdash{{0.212cm}}{0cm}
\pgfpathmoveto{\pgfqpoint{9.131cm}{13.17cm}}
\pgfpathlineto{\pgfqpoint{9.131cm}{14.114cm}}
\pgfusepath{stroke}
\pgfsetdash{{0.212cm}}{0cm}
\pgfpathmoveto{\pgfqpoint{10.037cm}{13.17cm}}
\pgfpathlineto{\pgfqpoint{10.037cm}{14.114cm}}
\pgfusepath{stroke}
\pgfsetdash{{0.212cm}}{0cm}
\pgfpathmoveto{\pgfqpoint{10.942cm}{14.508cm}}
\pgfpathlineto{\pgfqpoint{10.942cm}{17.266cm}}
\pgfusepath{stroke}
\pgfsetdash{{0.212cm}}{0cm}
\pgfpathmoveto{\pgfqpoint{11.847cm}{14.508cm}}
\pgfpathlineto{\pgfqpoint{11.847cm}{17.266cm}}
\pgfusepath{stroke}
\pgfsetdash{{0.212cm}}{0cm}
\pgfpathmoveto{\pgfqpoint{12.756cm}{14.54cm}}
\pgfpathlineto{\pgfqpoint{12.756cm}{16.616cm}}
\pgfusepath{stroke}
\pgfsetdash{{0.212cm}}{0cm}
\pgfpathmoveto{\pgfqpoint{13.661cm}{14.54cm}}
\pgfpathlineto{\pgfqpoint{13.661cm}{16.616cm}}
\pgfusepath{stroke}
\pgfsetdash{{0.212cm}}{0cm}
\pgfpathmoveto{\pgfqpoint{14.567cm}{14.834cm}}
\pgfpathlineto{\pgfqpoint{14.567cm}{15.622cm}}
\pgfusepath{stroke}
\pgfsetdash{{0.212cm}}{0cm}
\pgfpathmoveto{\pgfqpoint{15.472cm}{14.834cm}}
\pgfpathlineto{\pgfqpoint{15.472cm}{15.622cm}}
\pgfusepath{stroke}
\pgfsetdash{{0.212cm}}{0cm}
\pgfpathmoveto{\pgfqpoint{16.378cm}{14.19cm}}
\pgfpathlineto{\pgfqpoint{16.378cm}{16.087cm}}
\pgfusepath{stroke}
\pgfsetdash{{0.212cm}}{0cm}
\pgfpathmoveto{\pgfqpoint{17.283cm}{14.19cm}}
\pgfpathlineto{\pgfqpoint{17.283cm}{16.087cm}}
\pgfusepath{stroke}
\pgfsetdash{{0.212cm}}{0cm}
\pgfpathmoveto{\pgfqpoint{18.192cm}{13.741cm}}
\pgfpathlineto{\pgfqpoint{18.192cm}{14.599cm}}
\pgfusepath{stroke}
\pgfsetdash{{0.212cm}}{0cm}
\pgfpathmoveto{\pgfqpoint{19.097cm}{13.741cm}}
\pgfpathlineto{\pgfqpoint{19.097cm}{14.599cm}}
\pgfusepath{stroke}
\pgfsetdash{{0.212cm}}{0cm}
\pgfpathmoveto{\pgfqpoint{20.002cm}{13.726cm}}
\pgfpathlineto{\pgfqpoint{20.002cm}{15.167cm}}
\pgfusepath{stroke}
\pgfsetdash{{0.212cm}}{0cm}
\pgfpathmoveto{\pgfqpoint{20.908cm}{13.726cm}}
\pgfpathlineto{\pgfqpoint{20.908cm}{15.167cm}}
\pgfusepath{stroke}
\pgfsetdash{{0.212cm}}{0cm}
\pgfpathmoveto{\pgfqpoint{21.813cm}{13.382cm}}
\pgfpathlineto{\pgfqpoint{21.813cm}{14.67cm}}
\pgfusepath{stroke}
\pgfsetdash{{0.212cm}}{0cm}
\pgfpathmoveto{\pgfqpoint{22.719cm}{13.382cm}}
\pgfpathlineto{\pgfqpoint{22.719cm}{14.67cm}}
\pgfusepath{stroke}
\pgfsetdash{{0.212cm}}{0cm}
\pgfpathmoveto{\pgfqpoint{1.884cm}{11.918cm}}
\pgfpathlineto{\pgfqpoint{1.884cm}{12.268cm}}
\pgfusepath{stroke}
\pgfsetdash{{0.212cm}}{0cm}
\pgfpathmoveto{\pgfqpoint{2.79cm}{11.918cm}}
\pgfpathlineto{\pgfqpoint{2.79cm}{12.268cm}}
\pgfusepath{stroke}
\pgfsetdash{{0.212cm}}{0cm}
\pgfpathmoveto{\pgfqpoint{3.695cm}{11.759cm}}
\pgfpathlineto{\pgfqpoint{3.695cm}{12.244cm}}
\pgfusepath{stroke}
\pgfsetdash{{0.212cm}}{0cm}
\pgfpathmoveto{\pgfqpoint{4.601cm}{11.759cm}}
\pgfpathlineto{\pgfqpoint{4.601cm}{12.244cm}}
\pgfusepath{stroke}
\pgfsetdash{{0.212cm}}{0cm}
\pgfpathmoveto{\pgfqpoint{5.506cm}{11.48cm}}
\pgfpathlineto{\pgfqpoint{5.506cm}{12.191cm}}
\pgfusepath{stroke}
\pgfsetdash{{0.212cm}}{0cm}
\pgfpathmoveto{\pgfqpoint{6.412cm}{11.48cm}}
\pgfpathlineto{\pgfqpoint{6.412cm}{12.191cm}}
\pgfusepath{stroke}
\pgfsetdash{{0.212cm}}{0cm}
\pgfpathmoveto{\pgfqpoint{7.32cm}{11.48cm}}
\pgfpathlineto{\pgfqpoint{7.32cm}{12.138cm}}
\pgfusepath{stroke}
\pgfsetdash{{0.212cm}}{0cm}
\pgfpathmoveto{\pgfqpoint{8.226cm}{11.48cm}}
\pgfpathlineto{\pgfqpoint{8.226cm}{12.138cm}}
\pgfusepath{stroke}
\pgfsetdash{{0.212cm}}{0cm}
\pgfpathmoveto{\pgfqpoint{9.131cm}{12.009cm}}
\pgfpathlineto{\pgfqpoint{9.131cm}{12.238cm}}
\pgfusepath{stroke}
\pgfsetdash{{0.212cm}}{0cm}
\pgfpathmoveto{\pgfqpoint{10.037cm}{12.009cm}}
\pgfpathlineto{\pgfqpoint{10.037cm}{12.238cm}}
\pgfusepath{stroke}
\pgfsetdash{{0.212cm}}{0cm}
\pgfpathmoveto{\pgfqpoint{10.942cm}{11.989cm}}
\pgfpathlineto{\pgfqpoint{10.942cm}{12.397cm}}
\pgfusepath{stroke}
\pgfsetdash{{0.212cm}}{0cm}
\pgfpathmoveto{\pgfqpoint{11.847cm}{11.989cm}}
\pgfpathlineto{\pgfqpoint{11.847cm}{12.397cm}}
\pgfusepath{stroke}
\pgfsetdash{{0.212cm}}{0cm}
\pgfpathmoveto{\pgfqpoint{12.756cm}{11.992cm}}
\pgfpathlineto{\pgfqpoint{12.756cm}{12.547cm}}
\pgfusepath{stroke}
\pgfsetdash{{0.212cm}}{0cm}
\pgfpathmoveto{\pgfqpoint{13.661cm}{11.992cm}}
\pgfpathlineto{\pgfqpoint{13.661cm}{12.547cm}}
\pgfusepath{stroke}
\pgfsetdash{{0.212cm}}{0cm}
\pgfpathmoveto{\pgfqpoint{14.567cm}{11.753cm}}
\pgfpathlineto{\pgfqpoint{14.567cm}{12.491cm}}
\pgfusepath{stroke}
\pgfsetdash{{0.212cm}}{0cm}
\pgfpathmoveto{\pgfqpoint{15.472cm}{11.753cm}}
\pgfpathlineto{\pgfqpoint{15.472cm}{12.491cm}}
\pgfusepath{stroke}
\pgfsetdash{{0.212cm}}{0cm}
\pgfpathmoveto{\pgfqpoint{16.378cm}{12.056cm}}
\pgfpathlineto{\pgfqpoint{16.378cm}{12.474cm}}
\pgfusepath{stroke}
\pgfsetdash{{0.212cm}}{0cm}
\pgfpathmoveto{\pgfqpoint{17.283cm}{12.056cm}}
\pgfpathlineto{\pgfqpoint{17.283cm}{12.474cm}}
\pgfusepath{stroke}
\pgfsetdash{{0.212cm}}{0cm}
\pgfpathmoveto{\pgfqpoint{18.192cm}{12.106cm}}
\pgfpathlineto{\pgfqpoint{18.192cm}{12.424cm}}
\pgfusepath{stroke}
\pgfsetdash{{0.212cm}}{0cm}
\pgfpathmoveto{\pgfqpoint{19.097cm}{12.106cm}}
\pgfpathlineto{\pgfqpoint{19.097cm}{12.424cm}}
\pgfusepath{stroke}
\pgfsetdash{{0.212cm}}{0cm}
\pgfpathmoveto{\pgfqpoint{20.002cm}{12.094cm}}
\pgfpathlineto{\pgfqpoint{20.002cm}{12.412cm}}
\pgfusepath{stroke}
\pgfsetdash{{0.212cm}}{0cm}
\pgfpathmoveto{\pgfqpoint{20.908cm}{12.094cm}}
\pgfpathlineto{\pgfqpoint{20.908cm}{12.412cm}}
\pgfusepath{stroke}
\pgfsetdash{{0.212cm}}{0cm}
\pgfpathmoveto{\pgfqpoint{21.813cm}{11.765cm}}
\pgfpathlineto{\pgfqpoint{21.813cm}{12.359cm}}
\pgfusepath{stroke}
\pgfsetdash{{0.212cm}}{0cm}
\pgfpathmoveto{\pgfqpoint{22.719cm}{11.765cm}}
\pgfpathlineto{\pgfqpoint{22.719cm}{12.359cm}}
\pgfusepath{stroke}
\pgfsetdash{}{0cm}
\pgfpathmoveto{\pgfqpoint{1.77cm}{14.073cm}}
\pgfpathlineto{\pgfqpoint{1.996cm}{14.073cm}}
\pgfusepath{stroke}
\pgfsetdash{}{0cm}
\pgfpathmoveto{\pgfqpoint{2.675cm}{14.073cm}}
\pgfpathlineto{\pgfqpoint{2.902cm}{14.073cm}}
\pgfusepath{stroke}
\pgfsetdash{}{0cm}
\pgfpathmoveto{\pgfqpoint{3.581cm}{13.923cm}}
\pgfpathlineto{\pgfqpoint{3.807cm}{13.923cm}}
\pgfusepath{stroke}
\pgfsetdash{}{0cm}
\pgfpathmoveto{\pgfqpoint{4.489cm}{13.923cm}}
\pgfpathlineto{\pgfqpoint{4.715cm}{13.923cm}}
\pgfusepath{stroke}
\pgfsetdash{}{0cm}
\pgfpathmoveto{\pgfqpoint{5.395cm}{13.632cm}}
\pgfpathlineto{\pgfqpoint{5.621cm}{13.632cm}}
\pgfusepath{stroke}
\pgfsetdash{}{0cm}
\pgfpathmoveto{\pgfqpoint{6.3cm}{13.632cm}}
\pgfpathlineto{\pgfqpoint{6.526cm}{13.632cm}}
\pgfusepath{stroke}
\pgfsetdash{}{0cm}
\pgfpathmoveto{\pgfqpoint{7.205cm}{13.658cm}}
\pgfpathlineto{\pgfqpoint{7.432cm}{13.658cm}}
\pgfusepath{stroke}
\pgfsetdash{}{0cm}
\pgfpathmoveto{\pgfqpoint{8.111cm}{13.658cm}}
\pgfpathlineto{\pgfqpoint{8.337cm}{13.658cm}}
\pgfusepath{stroke}
\pgfsetdash{}{0cm}
\pgfpathmoveto{\pgfqpoint{9.016cm}{14.114cm}}
\pgfpathlineto{\pgfqpoint{9.243cm}{14.114cm}}
\pgfusepath{stroke}
\pgfsetdash{}{0cm}
\pgfpathmoveto{\pgfqpoint{9.925cm}{14.114cm}}
\pgfpathlineto{\pgfqpoint{10.151cm}{14.114cm}}
\pgfusepath{stroke}
\pgfsetdash{}{0cm}
\pgfpathmoveto{\pgfqpoint{10.83cm}{17.266cm}}
\pgfpathlineto{\pgfqpoint{11.057cm}{17.266cm}}
\pgfusepath{stroke}
\pgfsetdash{}{0cm}
\pgfpathmoveto{\pgfqpoint{11.736cm}{17.266cm}}
\pgfpathlineto{\pgfqpoint{11.962cm}{17.266cm}}
\pgfusepath{stroke}
\pgfsetdash{}{0cm}
\pgfpathmoveto{\pgfqpoint{12.641cm}{16.616cm}}
\pgfpathlineto{\pgfqpoint{12.868cm}{16.616cm}}
\pgfusepath{stroke}
\pgfsetdash{}{0cm}
\pgfpathmoveto{\pgfqpoint{13.547cm}{16.616cm}}
\pgfpathlineto{\pgfqpoint{13.773cm}{16.616cm}}
\pgfusepath{stroke}
\pgfsetdash{}{0cm}
\pgfpathmoveto{\pgfqpoint{14.452cm}{15.622cm}}
\pgfpathlineto{\pgfqpoint{14.678cm}{15.622cm}}
\pgfusepath{stroke}
\pgfsetdash{}{0cm}
\pgfpathmoveto{\pgfqpoint{15.361cm}{15.622cm}}
\pgfpathlineto{\pgfqpoint{15.587cm}{15.622cm}}
\pgfusepath{stroke}
\pgfsetdash{}{0cm}
\pgfpathmoveto{\pgfqpoint{16.266cm}{16.087cm}}
\pgfpathlineto{\pgfqpoint{16.492cm}{16.087cm}}
\pgfusepath{stroke}
\pgfsetdash{}{0cm}
\pgfpathmoveto{\pgfqpoint{17.171cm}{16.087cm}}
\pgfpathlineto{\pgfqpoint{17.398cm}{16.087cm}}
\pgfusepath{stroke}
\pgfsetdash{}{0cm}
\pgfpathmoveto{\pgfqpoint{18.077cm}{14.599cm}}
\pgfpathlineto{\pgfqpoint{18.303cm}{14.599cm}}
\pgfusepath{stroke}
\pgfsetdash{}{0cm}
\pgfpathmoveto{\pgfqpoint{18.982cm}{14.599cm}}
\pgfpathlineto{\pgfqpoint{19.209cm}{14.599cm}}
\pgfusepath{stroke}
\pgfsetdash{}{0cm}
\pgfpathmoveto{\pgfqpoint{19.888cm}{15.167cm}}
\pgfpathlineto{\pgfqpoint{20.114cm}{15.167cm}}
\pgfusepath{stroke}
\pgfsetdash{}{0cm}
\pgfpathmoveto{\pgfqpoint{20.796cm}{15.167cm}}
\pgfpathlineto{\pgfqpoint{21.023cm}{15.167cm}}
\pgfusepath{stroke}
\pgfsetdash{}{0cm}
\pgfpathmoveto{\pgfqpoint{21.702cm}{14.67cm}}
\pgfpathlineto{\pgfqpoint{21.928cm}{14.67cm}}
\pgfusepath{stroke}
\pgfsetdash{}{0cm}
\pgfpathmoveto{\pgfqpoint{22.607cm}{14.67cm}}
\pgfpathlineto{\pgfqpoint{22.834cm}{14.67cm}}
\pgfusepath{stroke}
\pgfsetdash{}{0cm}
\pgfpathmoveto{\pgfqpoint{1.77cm}{11.918cm}}
\pgfpathlineto{\pgfqpoint{1.996cm}{11.918cm}}
\pgfusepath{stroke}
\pgfsetdash{}{0cm}
\pgfpathmoveto{\pgfqpoint{2.675cm}{11.918cm}}
\pgfpathlineto{\pgfqpoint{2.902cm}{11.918cm}}
\pgfusepath{stroke}
\pgfsetdash{}{0cm}
\pgfpathmoveto{\pgfqpoint{3.581cm}{11.759cm}}
\pgfpathlineto{\pgfqpoint{3.807cm}{11.759cm}}
\pgfusepath{stroke}
\pgfsetdash{}{0cm}
\pgfpathmoveto{\pgfqpoint{4.489cm}{11.759cm}}
\pgfpathlineto{\pgfqpoint{4.715cm}{11.759cm}}
\pgfusepath{stroke}
\pgfsetdash{}{0cm}
\pgfpathmoveto{\pgfqpoint{5.395cm}{11.48cm}}
\pgfpathlineto{\pgfqpoint{5.621cm}{11.48cm}}
\pgfusepath{stroke}
\pgfsetdash{}{0cm}
\pgfpathmoveto{\pgfqpoint{6.3cm}{11.48cm}}
\pgfpathlineto{\pgfqpoint{6.526cm}{11.48cm}}
\pgfusepath{stroke}
\pgfsetdash{}{0cm}
\pgfpathmoveto{\pgfqpoint{7.205cm}{11.48cm}}
\pgfpathlineto{\pgfqpoint{7.432cm}{11.48cm}}
\pgfusepath{stroke}
\pgfsetdash{}{0cm}
\pgfpathmoveto{\pgfqpoint{8.111cm}{11.48cm}}
\pgfpathlineto{\pgfqpoint{8.337cm}{11.48cm}}
\pgfusepath{stroke}
\pgfsetdash{}{0cm}
\pgfpathmoveto{\pgfqpoint{9.016cm}{12.009cm}}
\pgfpathlineto{\pgfqpoint{9.243cm}{12.009cm}}
\pgfusepath{stroke}
\pgfsetdash{}{0cm}
\pgfpathmoveto{\pgfqpoint{9.925cm}{12.009cm}}
\pgfpathlineto{\pgfqpoint{10.151cm}{12.009cm}}
\pgfusepath{stroke}
\pgfsetdash{}{0cm}
\pgfpathmoveto{\pgfqpoint{10.83cm}{11.989cm}}
\pgfpathlineto{\pgfqpoint{11.057cm}{11.989cm}}
\pgfusepath{stroke}
\pgfsetdash{}{0cm}
\pgfpathmoveto{\pgfqpoint{11.736cm}{11.989cm}}
\pgfpathlineto{\pgfqpoint{11.962cm}{11.989cm}}
\pgfusepath{stroke}
\pgfsetdash{}{0cm}
\pgfpathmoveto{\pgfqpoint{12.641cm}{11.992cm}}
\pgfpathlineto{\pgfqpoint{12.868cm}{11.992cm}}
\pgfusepath{stroke}
\pgfsetdash{}{0cm}
\pgfpathmoveto{\pgfqpoint{13.547cm}{11.992cm}}
\pgfpathlineto{\pgfqpoint{13.773cm}{11.992cm}}
\pgfusepath{stroke}
\pgfsetdash{}{0cm}
\pgfpathmoveto{\pgfqpoint{14.452cm}{11.753cm}}
\pgfpathlineto{\pgfqpoint{14.678cm}{11.753cm}}
\pgfusepath{stroke}
\pgfsetdash{}{0cm}
\pgfpathmoveto{\pgfqpoint{15.361cm}{11.753cm}}
\pgfpathlineto{\pgfqpoint{15.587cm}{11.753cm}}
\pgfusepath{stroke}
\pgfsetdash{}{0cm}
\pgfpathmoveto{\pgfqpoint{16.266cm}{12.056cm}}
\pgfpathlineto{\pgfqpoint{16.492cm}{12.056cm}}
\pgfusepath{stroke}
\pgfsetdash{}{0cm}
\pgfpathmoveto{\pgfqpoint{17.171cm}{12.056cm}}
\pgfpathlineto{\pgfqpoint{17.398cm}{12.056cm}}
\pgfusepath{stroke}
\pgfsetdash{}{0cm}
\pgfpathmoveto{\pgfqpoint{18.077cm}{12.106cm}}
\pgfpathlineto{\pgfqpoint{18.303cm}{12.106cm}}
\pgfusepath{stroke}
\pgfsetdash{}{0cm}
\pgfpathmoveto{\pgfqpoint{18.982cm}{12.106cm}}
\pgfpathlineto{\pgfqpoint{19.209cm}{12.106cm}}
\pgfusepath{stroke}
\pgfsetdash{}{0cm}
\pgfpathmoveto{\pgfqpoint{19.888cm}{12.094cm}}
\pgfpathlineto{\pgfqpoint{20.114cm}{12.094cm}}
\pgfusepath{stroke}
\pgfsetdash{}{0cm}
\pgfpathmoveto{\pgfqpoint{20.796cm}{12.094cm}}
\pgfpathlineto{\pgfqpoint{21.023cm}{12.094cm}}
\pgfusepath{stroke}
\pgfsetdash{}{0cm}
\pgfpathmoveto{\pgfqpoint{21.702cm}{11.765cm}}
\pgfpathlineto{\pgfqpoint{21.928cm}{11.765cm}}
\pgfusepath{stroke}
\pgfsetdash{}{0cm}
\pgfpathmoveto{\pgfqpoint{22.607cm}{11.765cm}}
\pgfpathlineto{\pgfqpoint{22.834cm}{11.765cm}}
\pgfusepath{stroke}
\pgfsetdash{}{0cm}
\definecolor{eps2pgf_color}{rgb}{0,0,1}\pgfsetstrokecolor{eps2pgf_color}\pgfsetfillcolor{eps2pgf_color}
\pgfpathmoveto{\pgfqpoint{1.658cm}{12.268cm}}
\pgfpathlineto{\pgfqpoint{1.658cm}{13.026cm}}
\pgfpathlineto{\pgfqpoint{2.111cm}{13.026cm}}
\pgfpathlineto{\pgfqpoint{2.111cm}{12.268cm}}
\pgfpathlineto{\pgfqpoint{1.658cm}{12.268cm}}
\pgfusepath{stroke}
\pgfsetdash{}{0cm}
\pgfpathmoveto{\pgfqpoint{2.564cm}{12.268cm}}
\pgfpathlineto{\pgfqpoint{2.564cm}{13.026cm}}
\pgfpathlineto{\pgfqpoint{3.016cm}{13.026cm}}
\pgfpathlineto{\pgfqpoint{3.016cm}{12.268cm}}
\pgfpathlineto{\pgfqpoint{2.564cm}{12.268cm}}
\pgfusepath{stroke}
\pgfsetdash{}{0cm}
\pgfpathmoveto{\pgfqpoint{3.469cm}{12.244cm}}
\pgfpathlineto{\pgfqpoint{3.469cm}{13.15cm}}
\pgfpathlineto{\pgfqpoint{3.922cm}{13.15cm}}
\pgfpathlineto{\pgfqpoint{3.922cm}{12.244cm}}
\pgfpathlineto{\pgfqpoint{3.469cm}{12.244cm}}
\pgfusepath{stroke}
\pgfsetdash{}{0cm}
\pgfpathmoveto{\pgfqpoint{4.374cm}{12.244cm}}
\pgfpathlineto{\pgfqpoint{4.374cm}{13.15cm}}
\pgfpathlineto{\pgfqpoint{4.827cm}{13.15cm}}
\pgfpathlineto{\pgfqpoint{4.827cm}{12.244cm}}
\pgfpathlineto{\pgfqpoint{4.374cm}{12.244cm}}
\pgfusepath{stroke}
\pgfsetdash{}{0cm}
\pgfpathmoveto{\pgfqpoint{5.28cm}{12.191cm}}
\pgfpathlineto{\pgfqpoint{5.28cm}{12.973cm}}
\pgfpathlineto{\pgfqpoint{5.733cm}{12.973cm}}
\pgfpathlineto{\pgfqpoint{5.733cm}{12.191cm}}
\pgfpathlineto{\pgfqpoint{5.28cm}{12.191cm}}
\pgfusepath{stroke}
\pgfsetdash{}{0cm}
\pgfpathmoveto{\pgfqpoint{6.185cm}{12.191cm}}
\pgfpathlineto{\pgfqpoint{6.185cm}{12.973cm}}
\pgfpathlineto{\pgfqpoint{6.638cm}{12.973cm}}
\pgfpathlineto{\pgfqpoint{6.638cm}{12.191cm}}
\pgfpathlineto{\pgfqpoint{6.185cm}{12.191cm}}
\pgfusepath{stroke}
\pgfsetdash{}{0cm}
\pgfpathmoveto{\pgfqpoint{7.094cm}{12.138cm}}
\pgfpathlineto{\pgfqpoint{7.094cm}{12.962cm}}
\pgfpathlineto{\pgfqpoint{7.547cm}{12.962cm}}
\pgfpathlineto{\pgfqpoint{7.547cm}{12.138cm}}
\pgfpathlineto{\pgfqpoint{7.094cm}{12.138cm}}
\pgfusepath{stroke}
\pgfsetdash{}{0cm}
\pgfpathmoveto{\pgfqpoint{7.999cm}{12.138cm}}
\pgfpathlineto{\pgfqpoint{7.999cm}{12.962cm}}
\pgfpathlineto{\pgfqpoint{8.452cm}{12.962cm}}
\pgfpathlineto{\pgfqpoint{8.452cm}{12.138cm}}
\pgfpathlineto{\pgfqpoint{7.999cm}{12.138cm}}
\pgfusepath{stroke}
\pgfsetdash{}{0cm}
\pgfpathmoveto{\pgfqpoint{8.905cm}{12.238cm}}
\pgfpathlineto{\pgfqpoint{8.905cm}{13.17cm}}
\pgfpathlineto{\pgfqpoint{9.357cm}{13.17cm}}
\pgfpathlineto{\pgfqpoint{9.357cm}{12.238cm}}
\pgfpathlineto{\pgfqpoint{8.905cm}{12.238cm}}
\pgfusepath{stroke}
\pgfsetdash{}{0cm}
\pgfpathmoveto{\pgfqpoint{9.81cm}{12.238cm}}
\pgfpathlineto{\pgfqpoint{9.81cm}{13.17cm}}
\pgfpathlineto{\pgfqpoint{10.263cm}{13.17cm}}
\pgfpathlineto{\pgfqpoint{10.263cm}{12.238cm}}
\pgfpathlineto{\pgfqpoint{9.81cm}{12.238cm}}
\pgfusepath{stroke}
\pgfsetdash{}{0cm}
\pgfpathmoveto{\pgfqpoint{10.716cm}{12.397cm}}
\pgfpathlineto{\pgfqpoint{10.716cm}{14.508cm}}
\pgfpathlineto{\pgfqpoint{11.168cm}{14.508cm}}
\pgfpathlineto{\pgfqpoint{11.168cm}{12.397cm}}
\pgfpathlineto{\pgfqpoint{10.716cm}{12.397cm}}
\pgfusepath{stroke}
\pgfsetdash{}{0cm}
\pgfpathmoveto{\pgfqpoint{11.621cm}{12.397cm}}
\pgfpathlineto{\pgfqpoint{11.621cm}{14.508cm}}
\pgfpathlineto{\pgfqpoint{12.074cm}{14.508cm}}
\pgfpathlineto{\pgfqpoint{12.074cm}{12.397cm}}
\pgfpathlineto{\pgfqpoint{11.621cm}{12.397cm}}
\pgfusepath{stroke}
\pgfsetdash{}{0cm}
\pgfpathmoveto{\pgfqpoint{12.529cm}{12.547cm}}
\pgfpathlineto{\pgfqpoint{12.529cm}{14.54cm}}
\pgfpathlineto{\pgfqpoint{12.982cm}{14.54cm}}
\pgfpathlineto{\pgfqpoint{12.982cm}{12.547cm}}
\pgfpathlineto{\pgfqpoint{12.529cm}{12.547cm}}
\pgfusepath{stroke}
\pgfsetdash{}{0cm}
\pgfpathmoveto{\pgfqpoint{13.435cm}{12.547cm}}
\pgfpathlineto{\pgfqpoint{13.435cm}{14.54cm}}
\pgfpathlineto{\pgfqpoint{13.888cm}{14.54cm}}
\pgfpathlineto{\pgfqpoint{13.888cm}{12.547cm}}
\pgfpathlineto{\pgfqpoint{13.435cm}{12.547cm}}
\pgfusepath{stroke}
\pgfsetdash{}{0cm}
\pgfpathmoveto{\pgfqpoint{14.34cm}{12.491cm}}
\pgfpathlineto{\pgfqpoint{14.34cm}{14.834cm}}
\pgfpathlineto{\pgfqpoint{14.793cm}{14.834cm}}
\pgfpathlineto{\pgfqpoint{14.793cm}{12.491cm}}
\pgfpathlineto{\pgfqpoint{14.34cm}{12.491cm}}
\pgfusepath{stroke}
\pgfsetdash{}{0cm}
\pgfpathmoveto{\pgfqpoint{15.246cm}{12.491cm}}
\pgfpathlineto{\pgfqpoint{15.246cm}{14.834cm}}
\pgfpathlineto{\pgfqpoint{15.699cm}{14.834cm}}
\pgfpathlineto{\pgfqpoint{15.699cm}{12.491cm}}
\pgfpathlineto{\pgfqpoint{15.246cm}{12.491cm}}
\pgfusepath{stroke}
\pgfsetdash{}{0cm}
\pgfpathmoveto{\pgfqpoint{16.151cm}{12.474cm}}
\pgfpathlineto{\pgfqpoint{16.151cm}{14.19cm}}
\pgfpathlineto{\pgfqpoint{16.604cm}{14.19cm}}
\pgfpathlineto{\pgfqpoint{16.604cm}{12.474cm}}
\pgfpathlineto{\pgfqpoint{16.151cm}{12.474cm}}
\pgfusepath{stroke}
\pgfsetdash{}{0cm}
\pgfpathmoveto{\pgfqpoint{17.057cm}{12.474cm}}
\pgfpathlineto{\pgfqpoint{17.057cm}{14.19cm}}
\pgfpathlineto{\pgfqpoint{17.51cm}{14.19cm}}
\pgfpathlineto{\pgfqpoint{17.51cm}{12.474cm}}
\pgfpathlineto{\pgfqpoint{17.057cm}{12.474cm}}
\pgfusepath{stroke}
\pgfsetdash{}{0cm}
\pgfpathmoveto{\pgfqpoint{17.965cm}{12.424cm}}
\pgfpathlineto{\pgfqpoint{17.965cm}{13.741cm}}
\pgfpathlineto{\pgfqpoint{18.418cm}{13.741cm}}
\pgfpathlineto{\pgfqpoint{18.418cm}{12.424cm}}
\pgfpathlineto{\pgfqpoint{17.965cm}{12.424cm}}
\pgfusepath{stroke}
\pgfsetdash{}{0cm}
\pgfpathmoveto{\pgfqpoint{18.871cm}{12.424cm}}
\pgfpathlineto{\pgfqpoint{18.871cm}{13.741cm}}
\pgfpathlineto{\pgfqpoint{19.323cm}{13.741cm}}
\pgfpathlineto{\pgfqpoint{19.323cm}{12.424cm}}
\pgfpathlineto{\pgfqpoint{18.871cm}{12.424cm}}
\pgfusepath{stroke}
\pgfsetdash{}{0cm}
\pgfpathmoveto{\pgfqpoint{19.776cm}{12.412cm}}
\pgfpathlineto{\pgfqpoint{19.776cm}{13.726cm}}
\pgfpathlineto{\pgfqpoint{20.229cm}{13.726cm}}
\pgfpathlineto{\pgfqpoint{20.229cm}{12.412cm}}
\pgfpathlineto{\pgfqpoint{19.776cm}{12.412cm}}
\pgfusepath{stroke}
\pgfsetdash{}{0cm}
\pgfpathmoveto{\pgfqpoint{20.682cm}{12.412cm}}
\pgfpathlineto{\pgfqpoint{20.682cm}{13.726cm}}
\pgfpathlineto{\pgfqpoint{21.134cm}{13.726cm}}
\pgfpathlineto{\pgfqpoint{21.134cm}{12.412cm}}
\pgfpathlineto{\pgfqpoint{20.682cm}{12.412cm}}
\pgfusepath{stroke}
\pgfsetdash{}{0cm}
\pgfpathmoveto{\pgfqpoint{21.587cm}{12.359cm}}
\pgfpathlineto{\pgfqpoint{21.587cm}{13.382cm}}
\pgfpathlineto{\pgfqpoint{22.04cm}{13.382cm}}
\pgfpathlineto{\pgfqpoint{22.04cm}{12.359cm}}
\pgfpathlineto{\pgfqpoint{21.587cm}{12.359cm}}
\pgfusepath{stroke}
\pgfsetdash{}{0cm}
\pgfpathmoveto{\pgfqpoint{22.493cm}{12.359cm}}
\pgfpathlineto{\pgfqpoint{22.493cm}{13.382cm}}
\pgfpathlineto{\pgfqpoint{22.945cm}{13.382cm}}
\pgfpathlineto{\pgfqpoint{22.945cm}{12.359cm}}
\pgfpathlineto{\pgfqpoint{22.493cm}{12.359cm}}
\pgfusepath{stroke}
\pgfsetdash{}{0cm}
\definecolor{eps2pgf_color}{rgb}{1,0,0}\pgfsetstrokecolor{eps2pgf_color}\pgfsetfillcolor{eps2pgf_color}
\pgfpathmoveto{\pgfqpoint{1.658cm}{12.441cm}}
\pgfpathlineto{\pgfqpoint{2.111cm}{12.441cm}}
\pgfusepath{stroke}
\pgfsetdash{}{0cm}
\pgfpathmoveto{\pgfqpoint{2.564cm}{12.441cm}}
\pgfpathlineto{\pgfqpoint{3.016cm}{12.441cm}}
\pgfusepath{stroke}
\pgfsetdash{}{0cm}
\pgfpathmoveto{\pgfqpoint{3.469cm}{12.527cm}}
\pgfpathlineto{\pgfqpoint{3.922cm}{12.527cm}}
\pgfusepath{stroke}
\pgfsetdash{}{0cm}
\pgfpathmoveto{\pgfqpoint{4.374cm}{12.527cm}}
\pgfpathlineto{\pgfqpoint{4.827cm}{12.527cm}}
\pgfusepath{stroke}
\pgfsetdash{}{0cm}
\pgfpathmoveto{\pgfqpoint{5.28cm}{12.465cm}}
\pgfpathlineto{\pgfqpoint{5.733cm}{12.465cm}}
\pgfusepath{stroke}
\pgfsetdash{}{0cm}
\pgfpathmoveto{\pgfqpoint{6.185cm}{12.465cm}}
\pgfpathlineto{\pgfqpoint{6.638cm}{12.465cm}}
\pgfusepath{stroke}
\pgfsetdash{}{0cm}
\pgfpathmoveto{\pgfqpoint{7.094cm}{12.43cm}}
\pgfpathlineto{\pgfqpoint{7.547cm}{12.43cm}}
\pgfusepath{stroke}
\pgfsetdash{}{0cm}
\pgfpathmoveto{\pgfqpoint{7.999cm}{12.43cm}}
\pgfpathlineto{\pgfqpoint{8.452cm}{12.43cm}}
\pgfusepath{stroke}
\pgfsetdash{}{0cm}
\pgfpathmoveto{\pgfqpoint{8.905cm}{12.529cm}}
\pgfpathlineto{\pgfqpoint{9.357cm}{12.529cm}}
\pgfusepath{stroke}
\pgfsetdash{}{0cm}
\pgfpathmoveto{\pgfqpoint{9.81cm}{12.529cm}}
\pgfpathlineto{\pgfqpoint{10.263cm}{12.529cm}}
\pgfusepath{stroke}
\pgfsetdash{}{0cm}
\pgfpathmoveto{\pgfqpoint{10.716cm}{13.047cm}}
\pgfpathlineto{\pgfqpoint{11.168cm}{13.047cm}}
\pgfusepath{stroke}
\pgfsetdash{}{0cm}
\pgfpathmoveto{\pgfqpoint{11.621cm}{13.047cm}}
\pgfpathlineto{\pgfqpoint{12.074cm}{13.047cm}}
\pgfusepath{stroke}
\pgfsetdash{}{0cm}
\pgfpathmoveto{\pgfqpoint{12.529cm}{13.855cm}}
\pgfpathlineto{\pgfqpoint{12.982cm}{13.855cm}}
\pgfusepath{stroke}
\pgfsetdash{}{0cm}
\pgfpathmoveto{\pgfqpoint{13.435cm}{13.855cm}}
\pgfpathlineto{\pgfqpoint{13.888cm}{13.855cm}}
\pgfusepath{stroke}
\pgfsetdash{}{0cm}
\pgfpathmoveto{\pgfqpoint{14.34cm}{13.888cm}}
\pgfpathlineto{\pgfqpoint{14.793cm}{13.888cm}}
\pgfusepath{stroke}
\pgfsetdash{}{0cm}
\pgfpathmoveto{\pgfqpoint{15.246cm}{13.888cm}}
\pgfpathlineto{\pgfqpoint{15.699cm}{13.888cm}}
\pgfusepath{stroke}
\pgfsetdash{}{0cm}
\pgfpathmoveto{\pgfqpoint{16.151cm}{13.388cm}}
\pgfpathlineto{\pgfqpoint{16.604cm}{13.388cm}}
\pgfusepath{stroke}
\pgfsetdash{}{0cm}
\pgfpathmoveto{\pgfqpoint{17.057cm}{13.388cm}}
\pgfpathlineto{\pgfqpoint{17.51cm}{13.388cm}}
\pgfusepath{stroke}
\pgfsetdash{}{0cm}
\pgfpathmoveto{\pgfqpoint{17.965cm}{13.153cm}}
\pgfpathlineto{\pgfqpoint{18.418cm}{13.153cm}}
\pgfusepath{stroke}
\pgfsetdash{}{0cm}
\pgfpathmoveto{\pgfqpoint{18.871cm}{13.153cm}}
\pgfpathlineto{\pgfqpoint{19.323cm}{13.153cm}}
\pgfusepath{stroke}
\pgfsetdash{}{0cm}
\pgfpathmoveto{\pgfqpoint{19.776cm}{12.744cm}}
\pgfpathlineto{\pgfqpoint{20.229cm}{12.744cm}}
\pgfusepath{stroke}
\pgfsetdash{}{0cm}
\pgfpathmoveto{\pgfqpoint{20.682cm}{12.744cm}}
\pgfpathlineto{\pgfqpoint{21.134cm}{12.744cm}}
\pgfusepath{stroke}
\pgfsetdash{}{0cm}
\pgfpathmoveto{\pgfqpoint{21.587cm}{12.941cm}}
\pgfpathlineto{\pgfqpoint{22.04cm}{12.941cm}}
\pgfusepath{stroke}
\pgfsetdash{}{0cm}
\pgfpathmoveto{\pgfqpoint{22.493cm}{12.941cm}}
\pgfpathlineto{\pgfqpoint{22.945cm}{12.941cm}}
\pgfusepath{stroke}
\begin{pgfscope}
\pgfpathmoveto{\pgfqpoint{1.596cm}{15.878cm}}
\pgfpathlineto{\pgfqpoint{2.175cm}{15.878cm}}
\pgfpathlineto{\pgfqpoint{2.175cm}{13.911cm}}
\pgfpathlineto{\pgfqpoint{1.596cm}{13.911cm}}
\pgfpathclose
\pgfusepath{clip}
\pgfsetdash{}{0cm}
\pgfsetlinewidth{0.706mm}
\pgfpathmoveto{\pgfqpoint{1.779cm}{14.202cm}}
\pgfpathlineto{\pgfqpoint{1.99cm}{14.202cm}}
\pgfusepath{stroke}
\pgfsetdash{}{0cm}
\pgfpathmoveto{\pgfqpoint{1.884cm}{14.308cm}}
\pgfpathlineto{\pgfqpoint{1.884cm}{14.096cm}}
\pgfusepath{stroke}
\pgfsetdash{}{0cm}
\pgfpathmoveto{\pgfqpoint{1.779cm}{15.59cm}}
\pgfpathlineto{\pgfqpoint{1.99cm}{15.59cm}}
\pgfusepath{stroke}
\pgfsetdash{}{0cm}
\pgfpathmoveto{\pgfqpoint{1.884cm}{15.696cm}}
\pgfpathlineto{\pgfqpoint{1.884cm}{15.484cm}}
\pgfusepath{stroke}
\end{pgfscope}
\begin{pgfscope}
\pgfpathmoveto{\pgfqpoint{2.502cm}{15.878cm}}
\pgfpathlineto{\pgfqpoint{3.081cm}{15.878cm}}
\pgfpathlineto{\pgfqpoint{3.081cm}{13.911cm}}
\pgfpathlineto{\pgfqpoint{2.502cm}{13.911cm}}
\pgfpathclose
\pgfusepath{clip}
\pgfsetdash{}{0cm}
\pgfsetlinewidth{0.706mm}
\pgfpathmoveto{\pgfqpoint{2.684cm}{14.202cm}}
\pgfpathlineto{\pgfqpoint{2.896cm}{14.202cm}}
\pgfusepath{stroke}
\pgfsetdash{}{0cm}
\pgfpathmoveto{\pgfqpoint{2.79cm}{14.308cm}}
\pgfpathlineto{\pgfqpoint{2.79cm}{14.096cm}}
\pgfusepath{stroke}
\pgfsetdash{}{0cm}
\pgfpathmoveto{\pgfqpoint{2.684cm}{15.59cm}}
\pgfpathlineto{\pgfqpoint{2.896cm}{15.59cm}}
\pgfusepath{stroke}
\pgfsetdash{}{0cm}
\pgfpathmoveto{\pgfqpoint{2.79cm}{15.696cm}}
\pgfpathlineto{\pgfqpoint{2.79cm}{15.484cm}}
\pgfusepath{stroke}
\end{pgfscope}
\begin{pgfscope}
\pgfpathmoveto{\pgfqpoint{5.218cm}{15.684cm}}
\pgfpathlineto{\pgfqpoint{5.797cm}{15.684cm}}
\pgfpathlineto{\pgfqpoint{5.797cm}{15.105cm}}
\pgfpathlineto{\pgfqpoint{5.218cm}{15.105cm}}
\pgfpathclose
\pgfusepath{clip}
\pgfsetdash{}{0cm}
\pgfsetlinewidth{0.706mm}
\pgfpathmoveto{\pgfqpoint{5.4cm}{15.396cm}}
\pgfpathlineto{\pgfqpoint{5.612cm}{15.396cm}}
\pgfusepath{stroke}
\pgfsetdash{}{0cm}
\pgfpathmoveto{\pgfqpoint{5.506cm}{15.502cm}}
\pgfpathlineto{\pgfqpoint{5.506cm}{15.29cm}}
\pgfusepath{stroke}
\end{pgfscope}
\begin{pgfscope}
\pgfpathmoveto{\pgfqpoint{6.124cm}{15.684cm}}
\pgfpathlineto{\pgfqpoint{6.703cm}{15.684cm}}
\pgfpathlineto{\pgfqpoint{6.703cm}{15.105cm}}
\pgfpathlineto{\pgfqpoint{6.124cm}{15.105cm}}
\pgfpathclose
\pgfusepath{clip}
\pgfsetdash{}{0cm}
\pgfsetlinewidth{0.706mm}
\pgfpathmoveto{\pgfqpoint{6.306cm}{15.396cm}}
\pgfpathlineto{\pgfqpoint{6.518cm}{15.396cm}}
\pgfusepath{stroke}
\pgfsetdash{}{0cm}
\pgfpathmoveto{\pgfqpoint{6.412cm}{15.502cm}}
\pgfpathlineto{\pgfqpoint{6.412cm}{15.29cm}}
\pgfusepath{stroke}
\end{pgfscope}
\begin{pgfscope}
\pgfpathmoveto{\pgfqpoint{7.032cm}{16.172cm}}
\pgfpathlineto{\pgfqpoint{7.611cm}{16.172cm}}
\pgfpathlineto{\pgfqpoint{7.611cm}{15.593cm}}
\pgfpathlineto{\pgfqpoint{7.032cm}{15.593cm}}
\pgfpathclose
\pgfusepath{clip}
\pgfsetdash{}{0cm}
\pgfsetlinewidth{0.706mm}
\pgfpathmoveto{\pgfqpoint{7.214cm}{15.884cm}}
\pgfpathlineto{\pgfqpoint{7.426cm}{15.884cm}}
\pgfusepath{stroke}
\pgfsetdash{}{0cm}
\pgfpathmoveto{\pgfqpoint{7.32cm}{15.99cm}}
\pgfpathlineto{\pgfqpoint{7.32cm}{15.778cm}}
\pgfusepath{stroke}
\end{pgfscope}
\begin{pgfscope}
\pgfpathmoveto{\pgfqpoint{7.938cm}{16.172cm}}
\pgfpathlineto{\pgfqpoint{8.517cm}{16.172cm}}
\pgfpathlineto{\pgfqpoint{8.517cm}{15.593cm}}
\pgfpathlineto{\pgfqpoint{7.938cm}{15.593cm}}
\pgfpathclose
\pgfusepath{clip}
\pgfsetdash{}{0cm}
\pgfsetlinewidth{0.706mm}
\pgfpathmoveto{\pgfqpoint{8.12cm}{15.884cm}}
\pgfpathlineto{\pgfqpoint{8.331cm}{15.884cm}}
\pgfusepath{stroke}
\pgfsetdash{}{0cm}
\pgfpathmoveto{\pgfqpoint{8.226cm}{15.99cm}}
\pgfpathlineto{\pgfqpoint{8.226cm}{15.778cm}}
\pgfusepath{stroke}
\end{pgfscope}
\begin{pgfscope}
\pgfpathmoveto{\pgfqpoint{8.843cm}{16.948cm}}
\pgfpathlineto{\pgfqpoint{9.422cm}{16.948cm}}
\pgfpathlineto{\pgfqpoint{9.422cm}{16.369cm}}
\pgfpathlineto{\pgfqpoint{8.843cm}{16.369cm}}
\pgfpathclose
\pgfusepath{clip}
\pgfsetdash{}{0cm}
\pgfsetlinewidth{0.706mm}
\pgfpathmoveto{\pgfqpoint{9.025cm}{16.66cm}}
\pgfpathlineto{\pgfqpoint{9.237cm}{16.66cm}}
\pgfusepath{stroke}
\pgfsetdash{}{0cm}
\pgfpathmoveto{\pgfqpoint{9.131cm}{16.766cm}}
\pgfpathlineto{\pgfqpoint{9.131cm}{16.554cm}}
\pgfusepath{stroke}
\end{pgfscope}
\begin{pgfscope}
\pgfpathmoveto{\pgfqpoint{9.748cm}{16.948cm}}
\pgfpathlineto{\pgfqpoint{10.328cm}{16.948cm}}
\pgfpathlineto{\pgfqpoint{10.328cm}{16.369cm}}
\pgfpathlineto{\pgfqpoint{9.748cm}{16.369cm}}
\pgfpathclose
\pgfusepath{clip}
\pgfsetdash{}{0cm}
\pgfsetlinewidth{0.706mm}
\pgfpathmoveto{\pgfqpoint{9.931cm}{16.66cm}}
\pgfpathlineto{\pgfqpoint{10.142cm}{16.66cm}}
\pgfusepath{stroke}
\pgfsetdash{}{0cm}
\pgfpathmoveto{\pgfqpoint{10.037cm}{16.766cm}}
\pgfpathlineto{\pgfqpoint{10.037cm}{16.554cm}}
\pgfusepath{stroke}
\end{pgfscope}
\begin{pgfscope}
\pgfpathmoveto{\pgfqpoint{17.903cm}{17.612cm}}
\pgfpathlineto{\pgfqpoint{18.483cm}{17.612cm}}
\pgfpathlineto{\pgfqpoint{18.483cm}{15.857cm}}
\pgfpathlineto{\pgfqpoint{17.903cm}{15.857cm}}
\pgfpathclose
\pgfusepath{clip}
\pgfsetdash{}{0cm}
\pgfsetlinewidth{0.706mm}
\pgfpathmoveto{\pgfqpoint{18.086cm}{16.148cm}}
\pgfpathlineto{\pgfqpoint{18.297cm}{16.148cm}}
\pgfusepath{stroke}
\pgfsetdash{}{0cm}
\pgfpathmoveto{\pgfqpoint{18.192cm}{16.254cm}}
\pgfpathlineto{\pgfqpoint{18.192cm}{16.043cm}}
\pgfusepath{stroke}
\pgfsetdash{}{0cm}
\pgfpathmoveto{\pgfqpoint{18.086cm}{17.324cm}}
\pgfpathlineto{\pgfqpoint{18.297cm}{17.324cm}}
\pgfusepath{stroke}
\pgfsetdash{}{0cm}
\pgfpathmoveto{\pgfqpoint{18.192cm}{17.43cm}}
\pgfpathlineto{\pgfqpoint{18.192cm}{17.218cm}}
\pgfusepath{stroke}
\end{pgfscope}
\begin{pgfscope}
\pgfpathmoveto{\pgfqpoint{18.809cm}{17.612cm}}
\pgfpathlineto{\pgfqpoint{19.388cm}{17.612cm}}
\pgfpathlineto{\pgfqpoint{19.388cm}{15.857cm}}
\pgfpathlineto{\pgfqpoint{18.809cm}{15.857cm}}
\pgfpathclose
\pgfusepath{clip}
\pgfsetdash{}{0cm}
\pgfsetlinewidth{0.706mm}
\pgfpathmoveto{\pgfqpoint{18.991cm}{16.148cm}}
\pgfpathlineto{\pgfqpoint{19.203cm}{16.148cm}}
\pgfusepath{stroke}
\pgfsetdash{}{0cm}
\pgfpathmoveto{\pgfqpoint{19.097cm}{16.254cm}}
\pgfpathlineto{\pgfqpoint{19.097cm}{16.043cm}}
\pgfusepath{stroke}
\pgfsetdash{}{0cm}
\pgfpathmoveto{\pgfqpoint{18.991cm}{17.324cm}}
\pgfpathlineto{\pgfqpoint{19.203cm}{17.324cm}}
\pgfusepath{stroke}
\pgfsetdash{}{0cm}
\pgfpathmoveto{\pgfqpoint{19.097cm}{17.43cm}}
\pgfpathlineto{\pgfqpoint{19.097cm}{17.218cm}}
\pgfusepath{stroke}
\end{pgfscope}
\begin{pgfscope}
\pgfpathmoveto{\pgfqpoint{19.714cm}{17.971cm}}
\pgfpathlineto{\pgfqpoint{20.294cm}{17.971cm}}
\pgfpathlineto{\pgfqpoint{20.294cm}{17.392cm}}
\pgfpathlineto{\pgfqpoint{19.714cm}{17.392cm}}
\pgfpathclose
\pgfusepath{clip}
\pgfsetdash{}{0cm}
\pgfsetlinewidth{0.706mm}
\pgfpathmoveto{\pgfqpoint{19.897cm}{17.683cm}}
\pgfpathlineto{\pgfqpoint{20.108cm}{17.683cm}}
\pgfusepath{stroke}
\pgfsetdash{}{0cm}
\pgfpathmoveto{\pgfqpoint{20.002cm}{17.789cm}}
\pgfpathlineto{\pgfqpoint{20.002cm}{17.577cm}}
\pgfusepath{stroke}
\end{pgfscope}
\begin{pgfscope}
\pgfpathmoveto{\pgfqpoint{20.62cm}{17.971cm}}
\pgfpathlineto{\pgfqpoint{21.199cm}{17.971cm}}
\pgfpathlineto{\pgfqpoint{21.199cm}{17.392cm}}
\pgfpathlineto{\pgfqpoint{20.62cm}{17.392cm}}
\pgfpathclose
\pgfusepath{clip}
\pgfsetdash{}{0cm}
\pgfsetlinewidth{0.706mm}
\pgfpathmoveto{\pgfqpoint{20.802cm}{17.683cm}}
\pgfpathlineto{\pgfqpoint{21.014cm}{17.683cm}}
\pgfusepath{stroke}
\pgfsetdash{}{0cm}
\pgfpathmoveto{\pgfqpoint{20.908cm}{17.789cm}}
\pgfpathlineto{\pgfqpoint{20.908cm}{17.577cm}}
\pgfusepath{stroke}
\end{pgfscope}
\begin{pgfscope}
\pgfpathmoveto{\pgfqpoint{21.525cm}{15.355cm}}
\pgfpathlineto{\pgfqpoint{22.104cm}{15.355cm}}
\pgfpathlineto{\pgfqpoint{22.104cm}{14.776cm}}
\pgfpathlineto{\pgfqpoint{21.525cm}{14.776cm}}
\pgfpathclose
\pgfusepath{clip}
\pgfsetdash{}{0cm}
\pgfsetlinewidth{0.706mm}
\pgfpathmoveto{\pgfqpoint{21.708cm}{15.067cm}}
\pgfpathlineto{\pgfqpoint{21.919cm}{15.067cm}}
\pgfusepath{stroke}
\pgfsetdash{}{0cm}
\pgfpathmoveto{\pgfqpoint{21.813cm}{15.172cm}}
\pgfpathlineto{\pgfqpoint{21.813cm}{14.961cm}}
\pgfusepath{stroke}
\end{pgfscope}
\begin{pgfscope}
\pgfpathmoveto{\pgfqpoint{22.431cm}{15.355cm}}
\pgfpathlineto{\pgfqpoint{23.01cm}{15.355cm}}
\pgfpathlineto{\pgfqpoint{23.01cm}{14.776cm}}
\pgfpathlineto{\pgfqpoint{22.431cm}{14.776cm}}
\pgfpathclose
\pgfusepath{clip}
\pgfsetdash{}{0cm}
\pgfsetlinewidth{0.706mm}
\pgfpathmoveto{\pgfqpoint{22.613cm}{15.067cm}}
\pgfpathlineto{\pgfqpoint{22.825cm}{15.067cm}}
\pgfusepath{stroke}
\pgfsetdash{}{0cm}
\pgfpathmoveto{\pgfqpoint{22.719cm}{15.172cm}}
\pgfpathlineto{\pgfqpoint{22.719cm}{14.961cm}}
\pgfusepath{stroke}
\end{pgfscope}
\end{pgfscope}
\pgftext[x=1.869cm,y=10.875cm,rotate=0]{\fontsize{20}{7.06}\selectfont{{1}}}
\pgftext[x=2.787cm,y=10.875cm,rotate=0]{\fontsize{20}{7.06}\selectfont{{2}}}
\pgftext[x=3.696cm,y=10.872cm,rotate=0]{\fontsize{20}{7.06}\selectfont{{3}}}
\pgftext[x=4.603cm,y=10.875cm,rotate=0]{\fontsize{20}{7.06}\selectfont{{4}}}
\pgftext[x=5.509cm,y=10.869cm,rotate=0]{\fontsize{20}{7.06}\selectfont{{5}}}
\pgftext[x=6.416cm,y=10.872cm,rotate=0]{\fontsize{20}{7.06}\selectfont{{6}}}
\pgftext[x=7.322cm,y=10.872cm,rotate=0]{\fontsize{20}{7.06}\selectfont{{7}}}
\pgftext[x=8.226cm,y=10.872cm,rotate=0]{\fontsize{20}{7.06}\selectfont{{8}}}
\pgftext[x=9.132cm,y=10.872cm,rotate=0]{\fontsize{20}{7.06}\selectfont{{9}}}
\pgftext[x=10.05cm,y=10.872cm,rotate=0]{\fontsize{20}{7.06}\selectfont{{10}}}
\pgftext[x=10.927cm,y=10.875cm,rotate=0]{\fontsize{20}{7.06}\selectfont{{11}}}
\pgftext[x=11.859cm,y=10.875cm,rotate=0]{\fontsize{20}{7.06}\selectfont{{12}}}
\pgftext[x=12.767cm,y=10.872cm,rotate=0]{\fontsize{20}{7.06}\selectfont{{13}}}
\pgftext[x=13.672cm,y=10.875cm,rotate=0]{\fontsize{20}{7.06}\selectfont{{14}}}
\pgftext[x=14.576cm,y=10.872cm,rotate=0]{\fontsize{20}{7.06}\selectfont{{15}}}
\pgftext[x=15.483cm,y=10.872cm,rotate=0]{\fontsize{20}{7.06}\selectfont{{16}}}
\pgftext[x=16.392cm,y=10.875cm,rotate=0]{\fontsize{20}{7.06}\selectfont{{17}}}
\pgftext[x=17.296cm,y=10.872cm,rotate=0]{\fontsize{20}{7.06}\selectfont{{18}}}
\pgftext[x=18.201cm,y=10.872cm,rotate=0]{\fontsize{20}{7.06}\selectfont{{19}}}
\pgftext[x=19.094cm,y=10.872cm,rotate=0]{\fontsize{20}{7.06}\selectfont{{20}}}
\pgftext[x=19.972cm,y=10.875cm,rotate=0]{\fontsize{20}{7.06}\selectfont{{21}}}
\pgftext[x=20.903cm,y=10.875cm,rotate=0]{\fontsize{20}{7.06}\selectfont{{22}}}
\pgftext[x=21.814cm,y=10.872cm,rotate=0]{\fontsize{20}{7.06}\selectfont{{23}}}
\pgftext[x=22.72cm,y=10.875cm,rotate=0]{\fontsize{20}{7.06}\selectfont{{24}}}
\pgftext[x=12.335cm,y=10.311cm-.2cm,rotate=0]{\fontsize{24}{7.06}\selectfont{{Hour}}}
\pgftext[x=-0.73cm,y=14.592cm,rotate=90]{\fontsize{24}{7.06}\selectfont{{Price (\$/MWh)}}}
\pgfsetdash{}{0cm}
\pgfusepath{stroke}
\end{pgfscope}
\end{pgfscope}
\end{pgfpicture}

%% file: I4_d.tex
\scalebox{0.35}{\scalefont{1.28} \input{./fig/I4_d_short.pgf}}

%% file: fig/I4_d_short.pgf
% Created by Eps2pgf 0.7.0 (build on 2008-08-24) on Thu May 22 19:24:07 PDT 2014
\begin{pgfpicture}
\pgfpathmoveto{\pgfqpoint{-1.87cm}{9.737cm}}
\pgfpathlineto{\pgfqpoint{23.495cm}{9.737cm}}
\pgfpathlineto{\pgfqpoint{23.495cm}{18.203cm}}
\pgfpathlineto{\pgfqpoint{-1.87cm}{18.203cm}}
\pgfpathclose
\pgfusepath{clip}
\begin{pgfscope}
\begin{pgfscope}
\pgfpathmoveto{\pgfqpoint{-1.87cm}{18.203cm}}
\pgfpathlineto{\pgfqpoint{23.527cm}{18.203cm}}
\pgfpathlineto{\pgfqpoint{23.527cm}{9.745cm}}
\pgfpathlineto{\pgfqpoint{-1.87cm}{9.745cm}}
\pgfpathclose
\pgfusepath{clip}
\begin{pgfscope}
\definecolor{eps2pgf_color}{gray}{1}\pgfsetstrokecolor{eps2pgf_color}\pgfsetfillcolor{eps2pgf_color}
\pgfpathmoveto{\pgfqpoint{-1.87cm}{18.203cm}}
\pgfpathlineto{\pgfqpoint{23.53cm}{18.203cm}}
\pgfpathlineto{\pgfqpoint{23.53cm}{9.743cm}}
\pgfpathlineto{\pgfqpoint{-1.87cm}{9.743cm}}
\pgfpathclose
\pgfusepath{fill}
\end{pgfscope}
\definecolor{eps2pgf_color}{gray}{1}\pgfsetstrokecolor{eps2pgf_color}\pgfsetfillcolor{eps2pgf_color}
\pgfpathmoveto{\pgfqpoint{1.432cm}{11.168cm}}
\pgfpathlineto{\pgfqpoint{1.432cm}{17.806cm}}
\pgfpathlineto{\pgfqpoint{23.175cm}{17.806cm}}
\pgfpathlineto{\pgfqpoint{23.175cm}{11.168cm}}
\pgfpathclose
\pgfseteorule\pgfusepath{fill}\pgfsetnonzerorule
\pgfsetdash{}{0cm}
\pgfsetlinewidth{0.176mm}
\pgfsetroundjoin
\pgfpathmoveto{\pgfqpoint{1.432cm}{11.168cm}}
\pgfpathlineto{\pgfqpoint{1.432cm}{17.806cm}}
\pgfpathlineto{\pgfqpoint{23.175cm}{17.806cm}}
\pgfpathlineto{\pgfqpoint{23.175cm}{11.168cm}}
\pgfpathlineto{\pgfqpoint{1.432cm}{11.168cm}}
\pgfusepath{stroke}
\pgfsetdash{}{0cm}
\definecolor{eps2pgf_color}{gray}{0}\pgfsetstrokecolor{eps2pgf_color}\pgfsetfillcolor{eps2pgf_color}
\pgfpathmoveto{\pgfqpoint{1.432cm}{11.168cm}}
\pgfpathlineto{\pgfqpoint{23.175cm}{11.168cm}}
\pgfusepath{stroke}
\pgfsetdash{}{0cm}
\pgfpathmoveto{\pgfqpoint{1.432cm}{17.806cm}}
\pgfpathlineto{\pgfqpoint{23.175cm}{17.806cm}}
\pgfusepath{stroke}
\pgfsetdash{}{0cm}
\pgfpathmoveto{\pgfqpoint{1.432cm}{11.168cm}}
\pgfpathlineto{\pgfqpoint{1.432cm}{17.806cm}}
\pgfusepath{stroke}
\pgfsetdash{}{0cm}
\pgfpathmoveto{\pgfqpoint{23.175cm}{11.168cm}}
\pgfpathlineto{\pgfqpoint{23.175cm}{17.806cm}}
\pgfusepath{stroke}
\pgfsetdash{}{0cm}
\pgfpathmoveto{\pgfqpoint{1.432cm}{11.168cm}}
\pgfpathlineto{\pgfqpoint{23.175cm}{11.168cm}}
\pgfusepath{stroke}
\pgfsetdash{}{0cm}
\pgfpathmoveto{\pgfqpoint{1.432cm}{11.168cm}}
\pgfpathlineto{\pgfqpoint{1.432cm}{17.806cm}}
\pgfusepath{stroke}
\pgfsetdash{}{0cm}
\pgfpathmoveto{\pgfqpoint{1.432cm}{11.515cm}}
\pgfpathlineto{\pgfqpoint{1.646cm}{11.515cm}}
\pgfusepath{stroke}
\pgfsetdash{}{0cm}
\pgfpathmoveto{\pgfqpoint{23.175cm}{11.515cm}}
\pgfpathlineto{\pgfqpoint{22.957cm}{11.515cm}}
\pgfusepath{stroke}
\pgftext[x=0.613cm,y=11.49cm,rotate=0]{ \fontsize{20}{7.06}\selectfont{{$-$50}}}
\pgfsetdash{}{0cm}
\pgfpathmoveto{\pgfqpoint{1.432cm}{14.623cm}}
\pgfpathlineto{\pgfqpoint{1.646cm}{14.623cm}}
\pgfusepath{stroke}
\pgfsetdash{}{0cm}
\pgfpathmoveto{\pgfqpoint{23.175cm}{14.623cm}}
\pgfpathlineto{\pgfqpoint{22.957cm}{14.623cm}}
\pgfusepath{stroke}
\pgftext[x=1.094cm,y=14.598cm,rotate=0]{ \fontsize{20}{7.06}\selectfont{{0}}}
\pgfsetdash{}{0cm}
\pgfpathmoveto{\pgfqpoint{1.432cm}{17.73cm}}
\pgfpathlineto{\pgfqpoint{1.646cm}{17.73cm}}
\pgfusepath{stroke}
\pgfsetdash{}{0cm}
\pgfpathmoveto{\pgfqpoint{23.175cm}{17.73cm}}
\pgfpathlineto{\pgfqpoint{22.957cm}{17.73cm}}
\pgfusepath{stroke}
\pgftext[x=0.857cm,y=17.705cm,rotate=0]{ \fontsize{20}{7.06}\selectfont{{50}}}
\pgfsetdash{}{0cm}
\pgfpathmoveto{\pgfqpoint{1.432cm}{11.168cm}}
\pgfpathlineto{\pgfqpoint{23.175cm}{11.168cm}}
\pgfusepath{stroke}
\pgfsetdash{}{0cm}
\pgfpathmoveto{\pgfqpoint{1.432cm}{17.806cm}}
\pgfpathlineto{\pgfqpoint{23.175cm}{17.806cm}}
\pgfusepath{stroke}
\pgfsetdash{}{0cm}
\pgfpathmoveto{\pgfqpoint{1.432cm}{11.168cm}}
\pgfpathlineto{\pgfqpoint{1.432cm}{17.806cm}}
\pgfusepath{stroke}
\pgfsetdash{}{0cm}
\pgfpathmoveto{\pgfqpoint{23.175cm}{11.168cm}}
\pgfpathlineto{\pgfqpoint{23.175cm}{17.806cm}}
\pgfusepath{stroke}
\begin{pgfscope}
\pgfpathmoveto{\pgfqpoint{1.432cm}{17.806cm}}
\pgfpathlineto{\pgfqpoint{23.178cm}{17.806cm}}
\pgfpathlineto{\pgfqpoint{23.178cm}{11.165cm}}
\pgfpathlineto{\pgfqpoint{1.432cm}{11.165cm}}
\pgfpathclose
\pgfusepath{clip}
\pgfsetdash{{0.212cm}}{0cm}
\pgfpathmoveto{\pgfqpoint{1.884cm}{14.928cm}}
\pgfpathlineto{\pgfqpoint{1.884cm}{16.686cm}}
\pgfusepath{stroke}
\pgfsetdash{{0.212cm}}{0cm}
\pgfpathmoveto{\pgfqpoint{2.79cm}{14.928cm}}
\pgfpathlineto{\pgfqpoint{2.79cm}{16.686cm}}
\pgfusepath{stroke}
\pgfsetdash{{0.212cm}}{0cm}
\pgfpathmoveto{\pgfqpoint{3.695cm}{15.09cm}}
\pgfpathlineto{\pgfqpoint{3.695cm}{17.183cm}}
\pgfusepath{stroke}
\pgfsetdash{{0.212cm}}{0cm}
\pgfpathmoveto{\pgfqpoint{4.601cm}{15.09cm}}
\pgfpathlineto{\pgfqpoint{4.601cm}{17.183cm}}
\pgfusepath{stroke}
\pgfsetdash{{0.212cm}}{0cm}
\pgfpathmoveto{\pgfqpoint{5.506cm}{15.305cm}}
\pgfpathlineto{\pgfqpoint{5.506cm}{17.363cm}}
\pgfusepath{stroke}
\pgfsetdash{{0.212cm}}{0cm}
\pgfpathmoveto{\pgfqpoint{6.412cm}{15.305cm}}
\pgfpathlineto{\pgfqpoint{6.412cm}{17.363cm}}
\pgfusepath{stroke}
\pgfsetdash{{0.212cm}}{0cm}
\pgfpathmoveto{\pgfqpoint{7.32cm}{15.405cm}}
\pgfpathlineto{\pgfqpoint{7.32cm}{17.507cm}}
\pgfusepath{stroke}
\pgfsetdash{{0.212cm}}{0cm}
\pgfpathmoveto{\pgfqpoint{8.226cm}{15.405cm}}
\pgfpathlineto{\pgfqpoint{8.226cm}{17.507cm}}
\pgfusepath{stroke}
\pgfsetdash{{0.212cm}}{0cm}
\pgfpathmoveto{\pgfqpoint{9.131cm}{15.455cm}}
\pgfpathlineto{\pgfqpoint{9.131cm}{17.474cm}}
\pgfusepath{stroke}
\pgfsetdash{{0.212cm}}{0cm}
\pgfpathmoveto{\pgfqpoint{10.037cm}{15.455cm}}
\pgfpathlineto{\pgfqpoint{10.037cm}{17.474cm}}
\pgfusepath{stroke}
\pgfsetdash{{0.212cm}}{0cm}
\pgfpathmoveto{\pgfqpoint{10.942cm}{15.164cm}}
\pgfpathlineto{\pgfqpoint{10.942cm}{16.631cm}}
\pgfusepath{stroke}
\pgfsetdash{{0.212cm}}{0cm}
\pgfpathmoveto{\pgfqpoint{11.847cm}{15.164cm}}
\pgfpathlineto{\pgfqpoint{11.847cm}{16.631cm}}
\pgfusepath{stroke}
\pgfsetdash{{0.212cm}}{0cm}
\pgfpathmoveto{\pgfqpoint{12.756cm}{15.169cm}}
\pgfpathlineto{\pgfqpoint{12.756cm}{16.939cm}}
\pgfusepath{stroke}
\pgfsetdash{{0.212cm}}{0cm}
\pgfpathmoveto{\pgfqpoint{13.661cm}{15.169cm}}
\pgfpathlineto{\pgfqpoint{13.661cm}{16.939cm}}
\pgfusepath{stroke}
\pgfsetdash{{0.212cm}}{0cm}
\pgfpathmoveto{\pgfqpoint{14.567cm}{15.202cm}}
\pgfpathlineto{\pgfqpoint{14.567cm}{15.996cm}}
\pgfusepath{stroke}
\pgfsetdash{{0.212cm}}{0cm}
\pgfpathmoveto{\pgfqpoint{15.472cm}{15.202cm}}
\pgfpathlineto{\pgfqpoint{15.472cm}{15.996cm}}
\pgfusepath{stroke}
\pgfsetdash{{0.212cm}}{0cm}
\pgfpathmoveto{\pgfqpoint{16.378cm}{15.02cm}}
\pgfpathlineto{\pgfqpoint{16.378cm}{15.884cm}}
\pgfusepath{stroke}
\pgfsetdash{{0.212cm}}{0cm}
\pgfpathmoveto{\pgfqpoint{17.283cm}{15.02cm}}
\pgfpathlineto{\pgfqpoint{17.283cm}{15.884cm}}
\pgfusepath{stroke}
\pgfsetdash{{0.212cm}}{0cm}
\pgfpathmoveto{\pgfqpoint{18.192cm}{15.031cm}}
\pgfpathlineto{\pgfqpoint{18.192cm}{15.71cm}}
\pgfusepath{stroke}
\pgfsetdash{{0.212cm}}{0cm}
\pgfpathmoveto{\pgfqpoint{19.097cm}{15.031cm}}
\pgfpathlineto{\pgfqpoint{19.097cm}{15.71cm}}
\pgfusepath{stroke}
\pgfsetdash{{0.212cm}}{0cm}
\pgfpathmoveto{\pgfqpoint{20.002cm}{15.125cm}}
\pgfpathlineto{\pgfqpoint{20.002cm}{15.978cm}}
\pgfusepath{stroke}
\pgfsetdash{{0.212cm}}{0cm}
\pgfpathmoveto{\pgfqpoint{20.908cm}{15.125cm}}
\pgfpathlineto{\pgfqpoint{20.908cm}{15.978cm}}
\pgfusepath{stroke}
\pgfsetdash{{0.212cm}}{0cm}
\pgfpathmoveto{\pgfqpoint{21.813cm}{14.97cm}}
\pgfpathlineto{\pgfqpoint{21.813cm}{16.954cm}}
\pgfusepath{stroke}
\pgfsetdash{{0.212cm}}{0cm}
\pgfpathmoveto{\pgfqpoint{22.719cm}{14.97cm}}
\pgfpathlineto{\pgfqpoint{22.719cm}{16.954cm}}
\pgfusepath{stroke}
\pgfsetdash{{0.212cm}}{0cm}
\pgfpathmoveto{\pgfqpoint{1.884cm}{11.471cm}}
\pgfpathlineto{\pgfqpoint{1.884cm}{13.076cm}}
\pgfusepath{stroke}
\pgfsetdash{{0.212cm}}{0cm}
\pgfpathmoveto{\pgfqpoint{2.79cm}{11.471cm}}
\pgfpathlineto{\pgfqpoint{2.79cm}{13.076cm}}
\pgfusepath{stroke}
\pgfsetdash{{0.212cm}}{0cm}
\pgfpathmoveto{\pgfqpoint{3.695cm}{12.191cm}}
\pgfpathlineto{\pgfqpoint{3.695cm}{13.603cm}}
\pgfusepath{stroke}
\pgfsetdash{{0.212cm}}{0cm}
\pgfpathmoveto{\pgfqpoint{4.601cm}{12.191cm}}
\pgfpathlineto{\pgfqpoint{4.601cm}{13.603cm}}
\pgfusepath{stroke}
\pgfsetdash{{0.212cm}}{0cm}
\pgfpathmoveto{\pgfqpoint{5.506cm}{13.032cm}}
\pgfpathlineto{\pgfqpoint{5.506cm}{13.591cm}}
\pgfusepath{stroke}
\pgfsetdash{{0.212cm}}{0cm}
\pgfpathmoveto{\pgfqpoint{6.412cm}{13.032cm}}
\pgfpathlineto{\pgfqpoint{6.412cm}{13.591cm}}
\pgfusepath{stroke}
\pgfsetdash{{0.212cm}}{0cm}
\pgfpathmoveto{\pgfqpoint{7.32cm}{12.403cm}}
\pgfpathlineto{\pgfqpoint{7.32cm}{13.15cm}}
\pgfusepath{stroke}
\pgfsetdash{{0.212cm}}{0cm}
\pgfpathmoveto{\pgfqpoint{8.226cm}{12.403cm}}
\pgfpathlineto{\pgfqpoint{8.226cm}{13.15cm}}
\pgfusepath{stroke}
\pgfsetdash{{0.212cm}}{0cm}
\pgfpathmoveto{\pgfqpoint{9.131cm}{11.771cm}}
\pgfpathlineto{\pgfqpoint{9.131cm}{13.367cm}}
\pgfusepath{stroke}
\pgfsetdash{{0.212cm}}{0cm}
\pgfpathmoveto{\pgfqpoint{10.037cm}{11.771cm}}
\pgfpathlineto{\pgfqpoint{10.037cm}{13.367cm}}
\pgfusepath{stroke}
\pgfsetdash{{0.212cm}}{0cm}
\pgfpathmoveto{\pgfqpoint{10.942cm}{11.906cm}}
\pgfpathlineto{\pgfqpoint{10.942cm}{13.188cm}}
\pgfusepath{stroke}
\pgfsetdash{{0.212cm}}{0cm}
\pgfpathmoveto{\pgfqpoint{11.847cm}{11.906cm}}
\pgfpathlineto{\pgfqpoint{11.847cm}{13.188cm}}
\pgfusepath{stroke}
\pgfsetdash{{0.212cm}}{0cm}
\pgfpathmoveto{\pgfqpoint{12.756cm}{12.424cm}}
\pgfpathlineto{\pgfqpoint{12.756cm}{13.785cm}}
\pgfusepath{stroke}
\pgfsetdash{{0.212cm}}{0cm}
\pgfpathmoveto{\pgfqpoint{13.661cm}{12.424cm}}
\pgfpathlineto{\pgfqpoint{13.661cm}{13.785cm}}
\pgfusepath{stroke}
\pgfsetdash{{0.212cm}}{0cm}
\pgfpathmoveto{\pgfqpoint{14.567cm}{11.75cm}}
\pgfpathlineto{\pgfqpoint{14.567cm}{13.529cm}}
\pgfusepath{stroke}
\pgfsetdash{{0.212cm}}{0cm}
\pgfpathmoveto{\pgfqpoint{15.472cm}{11.75cm}}
\pgfpathlineto{\pgfqpoint{15.472cm}{13.529cm}}
\pgfusepath{stroke}
\pgfsetdash{{0.212cm}}{0cm}
\pgfpathmoveto{\pgfqpoint{16.378cm}{12.324cm}}
\pgfpathlineto{\pgfqpoint{16.378cm}{13.87cm}}
\pgfusepath{stroke}
\pgfsetdash{{0.212cm}}{0cm}
\pgfpathmoveto{\pgfqpoint{17.283cm}{12.324cm}}
\pgfpathlineto{\pgfqpoint{17.283cm}{13.87cm}}
\pgfusepath{stroke}
\pgfsetdash{{0.212cm}}{0cm}
\pgfpathmoveto{\pgfqpoint{18.192cm}{12.465cm}}
\pgfpathlineto{\pgfqpoint{18.192cm}{13.547cm}}
\pgfusepath{stroke}
\pgfsetdash{{0.212cm}}{0cm}
\pgfpathmoveto{\pgfqpoint{19.097cm}{12.465cm}}
\pgfpathlineto{\pgfqpoint{19.097cm}{13.547cm}}
\pgfusepath{stroke}
\pgfsetdash{{0.212cm}}{0cm}
\pgfpathmoveto{\pgfqpoint{20.002cm}{12.318cm}}
\pgfpathlineto{\pgfqpoint{20.002cm}{13.573cm}}
\pgfusepath{stroke}
\pgfsetdash{{0.212cm}}{0cm}
\pgfpathmoveto{\pgfqpoint{20.908cm}{12.318cm}}
\pgfpathlineto{\pgfqpoint{20.908cm}{13.573cm}}
\pgfusepath{stroke}
\pgfsetdash{{0.212cm}}{0cm}
\pgfpathmoveto{\pgfqpoint{21.813cm}{11.789cm}}
\pgfpathlineto{\pgfqpoint{21.813cm}{13.479cm}}
\pgfusepath{stroke}
\pgfsetdash{{0.212cm}}{0cm}
\pgfpathmoveto{\pgfqpoint{22.719cm}{11.789cm}}
\pgfpathlineto{\pgfqpoint{22.719cm}{13.479cm}}
\pgfusepath{stroke}
\pgfsetdash{}{0cm}
\pgfpathmoveto{\pgfqpoint{1.77cm}{16.686cm}}
\pgfpathlineto{\pgfqpoint{1.996cm}{16.686cm}}
\pgfusepath{stroke}
\pgfsetdash{}{0cm}
\pgfpathmoveto{\pgfqpoint{2.675cm}{16.686cm}}
\pgfpathlineto{\pgfqpoint{2.902cm}{16.686cm}}
\pgfusepath{stroke}
\pgfsetdash{}{0cm}
\pgfpathmoveto{\pgfqpoint{3.581cm}{17.183cm}}
\pgfpathlineto{\pgfqpoint{3.807cm}{17.183cm}}
\pgfusepath{stroke}
\pgfsetdash{}{0cm}
\pgfpathmoveto{\pgfqpoint{4.489cm}{17.183cm}}
\pgfpathlineto{\pgfqpoint{4.715cm}{17.183cm}}
\pgfusepath{stroke}
\pgfsetdash{}{0cm}
\pgfpathmoveto{\pgfqpoint{5.395cm}{17.363cm}}
\pgfpathlineto{\pgfqpoint{5.621cm}{17.363cm}}
\pgfusepath{stroke}
\pgfsetdash{}{0cm}
\pgfpathmoveto{\pgfqpoint{6.3cm}{17.363cm}}
\pgfpathlineto{\pgfqpoint{6.526cm}{17.363cm}}
\pgfusepath{stroke}
\pgfsetdash{}{0cm}
\pgfpathmoveto{\pgfqpoint{7.205cm}{17.507cm}}
\pgfpathlineto{\pgfqpoint{7.432cm}{17.507cm}}
\pgfusepath{stroke}
\pgfsetdash{}{0cm}
\pgfpathmoveto{\pgfqpoint{8.111cm}{17.507cm}}
\pgfpathlineto{\pgfqpoint{8.337cm}{17.507cm}}
\pgfusepath{stroke}
\pgfsetdash{}{0cm}
\pgfpathmoveto{\pgfqpoint{9.016cm}{17.474cm}}
\pgfpathlineto{\pgfqpoint{9.243cm}{17.474cm}}
\pgfusepath{stroke}
\pgfsetdash{}{0cm}
\pgfpathmoveto{\pgfqpoint{9.925cm}{17.474cm}}
\pgfpathlineto{\pgfqpoint{10.151cm}{17.474cm}}
\pgfusepath{stroke}
\pgfsetdash{}{0cm}
\pgfpathmoveto{\pgfqpoint{10.83cm}{16.631cm}}
\pgfpathlineto{\pgfqpoint{11.057cm}{16.631cm}}
\pgfusepath{stroke}
\pgfsetdash{}{0cm}
\pgfpathmoveto{\pgfqpoint{11.736cm}{16.631cm}}
\pgfpathlineto{\pgfqpoint{11.962cm}{16.631cm}}
\pgfusepath{stroke}
\pgfsetdash{}{0cm}
\pgfpathmoveto{\pgfqpoint{12.641cm}{16.939cm}}
\pgfpathlineto{\pgfqpoint{12.868cm}{16.939cm}}
\pgfusepath{stroke}
\pgfsetdash{}{0cm}
\pgfpathmoveto{\pgfqpoint{13.547cm}{16.939cm}}
\pgfpathlineto{\pgfqpoint{13.773cm}{16.939cm}}
\pgfusepath{stroke}
\pgfsetdash{}{0cm}
\pgfpathmoveto{\pgfqpoint{14.452cm}{15.996cm}}
\pgfpathlineto{\pgfqpoint{14.678cm}{15.996cm}}
\pgfusepath{stroke}
\pgfsetdash{}{0cm}
\pgfpathmoveto{\pgfqpoint{15.361cm}{15.996cm}}
\pgfpathlineto{\pgfqpoint{15.587cm}{15.996cm}}
\pgfusepath{stroke}
\pgfsetdash{}{0cm}
\pgfpathmoveto{\pgfqpoint{16.266cm}{15.884cm}}
\pgfpathlineto{\pgfqpoint{16.492cm}{15.884cm}}
\pgfusepath{stroke}
\pgfsetdash{}{0cm}
\pgfpathmoveto{\pgfqpoint{17.171cm}{15.884cm}}
\pgfpathlineto{\pgfqpoint{17.398cm}{15.884cm}}
\pgfusepath{stroke}
\pgfsetdash{}{0cm}
\pgfpathmoveto{\pgfqpoint{18.077cm}{15.71cm}}
\pgfpathlineto{\pgfqpoint{18.303cm}{15.71cm}}
\pgfusepath{stroke}
\pgfsetdash{}{0cm}
\pgfpathmoveto{\pgfqpoint{18.982cm}{15.71cm}}
\pgfpathlineto{\pgfqpoint{19.209cm}{15.71cm}}
\pgfusepath{stroke}
\pgfsetdash{}{0cm}
\pgfpathmoveto{\pgfqpoint{19.888cm}{15.978cm}}
\pgfpathlineto{\pgfqpoint{20.114cm}{15.978cm}}
\pgfusepath{stroke}
\pgfsetdash{}{0cm}
\pgfpathmoveto{\pgfqpoint{20.796cm}{15.978cm}}
\pgfpathlineto{\pgfqpoint{21.023cm}{15.978cm}}
\pgfusepath{stroke}
\pgfsetdash{}{0cm}
\pgfpathmoveto{\pgfqpoint{21.702cm}{16.954cm}}
\pgfpathlineto{\pgfqpoint{21.928cm}{16.954cm}}
\pgfusepath{stroke}
\pgfsetdash{}{0cm}
\pgfpathmoveto{\pgfqpoint{22.607cm}{16.954cm}}
\pgfpathlineto{\pgfqpoint{22.834cm}{16.954cm}}
\pgfusepath{stroke}
\pgfsetdash{}{0cm}
\pgfpathmoveto{\pgfqpoint{1.77cm}{11.471cm}}
\pgfpathlineto{\pgfqpoint{1.996cm}{11.471cm}}
\pgfusepath{stroke}
\pgfsetdash{}{0cm}
\pgfpathmoveto{\pgfqpoint{2.675cm}{11.471cm}}
\pgfpathlineto{\pgfqpoint{2.902cm}{11.471cm}}
\pgfusepath{stroke}
\pgfsetdash{}{0cm}
\pgfpathmoveto{\pgfqpoint{3.581cm}{12.191cm}}
\pgfpathlineto{\pgfqpoint{3.807cm}{12.191cm}}
\pgfusepath{stroke}
\pgfsetdash{}{0cm}
\pgfpathmoveto{\pgfqpoint{4.489cm}{12.191cm}}
\pgfpathlineto{\pgfqpoint{4.715cm}{12.191cm}}
\pgfusepath{stroke}
\pgfsetdash{}{0cm}
\pgfpathmoveto{\pgfqpoint{5.395cm}{13.032cm}}
\pgfpathlineto{\pgfqpoint{5.621cm}{13.032cm}}
\pgfusepath{stroke}
\pgfsetdash{}{0cm}
\pgfpathmoveto{\pgfqpoint{6.3cm}{13.032cm}}
\pgfpathlineto{\pgfqpoint{6.526cm}{13.032cm}}
\pgfusepath{stroke}
\pgfsetdash{}{0cm}
\pgfpathmoveto{\pgfqpoint{7.205cm}{12.403cm}}
\pgfpathlineto{\pgfqpoint{7.432cm}{12.403cm}}
\pgfusepath{stroke}
\pgfsetdash{}{0cm}
\pgfpathmoveto{\pgfqpoint{8.111cm}{12.403cm}}
\pgfpathlineto{\pgfqpoint{8.337cm}{12.403cm}}
\pgfusepath{stroke}
\pgfsetdash{}{0cm}
\pgfpathmoveto{\pgfqpoint{9.016cm}{11.771cm}}
\pgfpathlineto{\pgfqpoint{9.243cm}{11.771cm}}
\pgfusepath{stroke}
\pgfsetdash{}{0cm}
\pgfpathmoveto{\pgfqpoint{9.925cm}{11.771cm}}
\pgfpathlineto{\pgfqpoint{10.151cm}{11.771cm}}
\pgfusepath{stroke}
\pgfsetdash{}{0cm}
\pgfpathmoveto{\pgfqpoint{10.83cm}{11.906cm}}
\pgfpathlineto{\pgfqpoint{11.057cm}{11.906cm}}
\pgfusepath{stroke}
\pgfsetdash{}{0cm}
\pgfpathmoveto{\pgfqpoint{11.736cm}{11.906cm}}
\pgfpathlineto{\pgfqpoint{11.962cm}{11.906cm}}
\pgfusepath{stroke}
\pgfsetdash{}{0cm}
\pgfpathmoveto{\pgfqpoint{12.641cm}{12.424cm}}
\pgfpathlineto{\pgfqpoint{12.868cm}{12.424cm}}
\pgfusepath{stroke}
\pgfsetdash{}{0cm}
\pgfpathmoveto{\pgfqpoint{13.547cm}{12.424cm}}
\pgfpathlineto{\pgfqpoint{13.773cm}{12.424cm}}
\pgfusepath{stroke}
\pgfsetdash{}{0cm}
\pgfpathmoveto{\pgfqpoint{14.452cm}{11.75cm}}
\pgfpathlineto{\pgfqpoint{14.678cm}{11.75cm}}
\pgfusepath{stroke}
\pgfsetdash{}{0cm}
\pgfpathmoveto{\pgfqpoint{15.361cm}{11.75cm}}
\pgfpathlineto{\pgfqpoint{15.587cm}{11.75cm}}
\pgfusepath{stroke}
\pgfsetdash{}{0cm}
\pgfpathmoveto{\pgfqpoint{16.266cm}{12.324cm}}
\pgfpathlineto{\pgfqpoint{16.492cm}{12.324cm}}
\pgfusepath{stroke}
\pgfsetdash{}{0cm}
\pgfpathmoveto{\pgfqpoint{17.171cm}{12.324cm}}
\pgfpathlineto{\pgfqpoint{17.398cm}{12.324cm}}
\pgfusepath{stroke}
\pgfsetdash{}{0cm}
\pgfpathmoveto{\pgfqpoint{18.077cm}{12.465cm}}
\pgfpathlineto{\pgfqpoint{18.303cm}{12.465cm}}
\pgfusepath{stroke}
\pgfsetdash{}{0cm}
\pgfpathmoveto{\pgfqpoint{18.982cm}{12.465cm}}
\pgfpathlineto{\pgfqpoint{19.209cm}{12.465cm}}
\pgfusepath{stroke}
\pgfsetdash{}{0cm}
\pgfpathmoveto{\pgfqpoint{19.888cm}{12.318cm}}
\pgfpathlineto{\pgfqpoint{20.114cm}{12.318cm}}
\pgfusepath{stroke}
\pgfsetdash{}{0cm}
\pgfpathmoveto{\pgfqpoint{20.796cm}{12.318cm}}
\pgfpathlineto{\pgfqpoint{21.023cm}{12.318cm}}
\pgfusepath{stroke}
\pgfsetdash{}{0cm}
\pgfpathmoveto{\pgfqpoint{21.702cm}{11.789cm}}
\pgfpathlineto{\pgfqpoint{21.928cm}{11.789cm}}
\pgfusepath{stroke}
\pgfsetdash{}{0cm}
\pgfpathmoveto{\pgfqpoint{22.607cm}{11.789cm}}
\pgfpathlineto{\pgfqpoint{22.834cm}{11.789cm}}
\pgfusepath{stroke}
\pgfsetdash{}{0cm}
\definecolor{eps2pgf_color}{rgb}{0,0,1}\pgfsetstrokecolor{eps2pgf_color}\pgfsetfillcolor{eps2pgf_color}
\pgfpathmoveto{\pgfqpoint{1.658cm}{13.076cm}}
\pgfpathlineto{\pgfqpoint{1.658cm}{14.928cm}}
\pgfpathlineto{\pgfqpoint{2.111cm}{14.928cm}}
\pgfpathlineto{\pgfqpoint{2.111cm}{13.076cm}}
\pgfpathlineto{\pgfqpoint{1.658cm}{13.076cm}}
\pgfusepath{stroke}
\pgfsetdash{}{0cm}
\pgfpathmoveto{\pgfqpoint{2.564cm}{13.076cm}}
\pgfpathlineto{\pgfqpoint{2.564cm}{14.928cm}}
\pgfpathlineto{\pgfqpoint{3.016cm}{14.928cm}}
\pgfpathlineto{\pgfqpoint{3.016cm}{13.076cm}}
\pgfpathlineto{\pgfqpoint{2.564cm}{13.076cm}}
\pgfusepath{stroke}
\pgfsetdash{}{0cm}
\pgfpathmoveto{\pgfqpoint{3.469cm}{13.603cm}}
\pgfpathlineto{\pgfqpoint{3.469cm}{15.09cm}}
\pgfpathlineto{\pgfqpoint{3.922cm}{15.09cm}}
\pgfpathlineto{\pgfqpoint{3.922cm}{13.603cm}}
\pgfpathlineto{\pgfqpoint{3.469cm}{13.603cm}}
\pgfusepath{stroke}
\pgfsetdash{}{0cm}
\pgfpathmoveto{\pgfqpoint{4.374cm}{13.603cm}}
\pgfpathlineto{\pgfqpoint{4.374cm}{15.09cm}}
\pgfpathlineto{\pgfqpoint{4.827cm}{15.09cm}}
\pgfpathlineto{\pgfqpoint{4.827cm}{13.603cm}}
\pgfpathlineto{\pgfqpoint{4.374cm}{13.603cm}}
\pgfusepath{stroke}
\pgfsetdash{}{0cm}
\pgfpathmoveto{\pgfqpoint{5.28cm}{13.591cm}}
\pgfpathlineto{\pgfqpoint{5.28cm}{15.305cm}}
\pgfpathlineto{\pgfqpoint{5.733cm}{15.305cm}}
\pgfpathlineto{\pgfqpoint{5.733cm}{13.591cm}}
\pgfpathlineto{\pgfqpoint{5.28cm}{13.591cm}}
\pgfusepath{stroke}
\pgfsetdash{}{0cm}
\pgfpathmoveto{\pgfqpoint{6.185cm}{13.591cm}}
\pgfpathlineto{\pgfqpoint{6.185cm}{15.305cm}}
\pgfpathlineto{\pgfqpoint{6.638cm}{15.305cm}}
\pgfpathlineto{\pgfqpoint{6.638cm}{13.591cm}}
\pgfpathlineto{\pgfqpoint{6.185cm}{13.591cm}}
\pgfusepath{stroke}
\pgfsetdash{}{0cm}
\pgfpathmoveto{\pgfqpoint{7.094cm}{13.15cm}}
\pgfpathlineto{\pgfqpoint{7.094cm}{15.405cm}}
\pgfpathlineto{\pgfqpoint{7.547cm}{15.405cm}}
\pgfpathlineto{\pgfqpoint{7.547cm}{13.15cm}}
\pgfpathlineto{\pgfqpoint{7.094cm}{13.15cm}}
\pgfusepath{stroke}
\pgfsetdash{}{0cm}
\pgfpathmoveto{\pgfqpoint{7.999cm}{13.15cm}}
\pgfpathlineto{\pgfqpoint{7.999cm}{15.405cm}}
\pgfpathlineto{\pgfqpoint{8.452cm}{15.405cm}}
\pgfpathlineto{\pgfqpoint{8.452cm}{13.15cm}}
\pgfpathlineto{\pgfqpoint{7.999cm}{13.15cm}}
\pgfusepath{stroke}
\pgfsetdash{}{0cm}
\pgfpathmoveto{\pgfqpoint{8.905cm}{13.367cm}}
\pgfpathlineto{\pgfqpoint{8.905cm}{15.455cm}}
\pgfpathlineto{\pgfqpoint{9.357cm}{15.455cm}}
\pgfpathlineto{\pgfqpoint{9.357cm}{13.367cm}}
\pgfpathlineto{\pgfqpoint{8.905cm}{13.367cm}}
\pgfusepath{stroke}
\pgfsetdash{}{0cm}
\pgfpathmoveto{\pgfqpoint{9.81cm}{13.367cm}}
\pgfpathlineto{\pgfqpoint{9.81cm}{15.455cm}}
\pgfpathlineto{\pgfqpoint{10.263cm}{15.455cm}}
\pgfpathlineto{\pgfqpoint{10.263cm}{13.367cm}}
\pgfpathlineto{\pgfqpoint{9.81cm}{13.367cm}}
\pgfusepath{stroke}
\pgfsetdash{}{0cm}
\pgfpathmoveto{\pgfqpoint{10.716cm}{13.188cm}}
\pgfpathlineto{\pgfqpoint{10.716cm}{15.164cm}}
\pgfpathlineto{\pgfqpoint{11.168cm}{15.164cm}}
\pgfpathlineto{\pgfqpoint{11.168cm}{13.188cm}}
\pgfpathlineto{\pgfqpoint{10.716cm}{13.188cm}}
\pgfusepath{stroke}
\pgfsetdash{}{0cm}
\pgfpathmoveto{\pgfqpoint{11.621cm}{13.188cm}}
\pgfpathlineto{\pgfqpoint{11.621cm}{15.164cm}}
\pgfpathlineto{\pgfqpoint{12.074cm}{15.164cm}}
\pgfpathlineto{\pgfqpoint{12.074cm}{13.188cm}}
\pgfpathlineto{\pgfqpoint{11.621cm}{13.188cm}}
\pgfusepath{stroke}
\pgfsetdash{}{0cm}
\pgfpathmoveto{\pgfqpoint{12.529cm}{13.785cm}}
\pgfpathlineto{\pgfqpoint{12.529cm}{15.169cm}}
\pgfpathlineto{\pgfqpoint{12.982cm}{15.169cm}}
\pgfpathlineto{\pgfqpoint{12.982cm}{13.785cm}}
\pgfpathlineto{\pgfqpoint{12.529cm}{13.785cm}}
\pgfusepath{stroke}
\pgfsetdash{}{0cm}
\pgfpathmoveto{\pgfqpoint{13.435cm}{13.785cm}}
\pgfpathlineto{\pgfqpoint{13.435cm}{15.169cm}}
\pgfpathlineto{\pgfqpoint{13.888cm}{15.169cm}}
\pgfpathlineto{\pgfqpoint{13.888cm}{13.785cm}}
\pgfpathlineto{\pgfqpoint{13.435cm}{13.785cm}}
\pgfusepath{stroke}
\pgfsetdash{}{0cm}
\pgfpathmoveto{\pgfqpoint{14.34cm}{13.529cm}}
\pgfpathlineto{\pgfqpoint{14.34cm}{15.202cm}}
\pgfpathlineto{\pgfqpoint{14.793cm}{15.202cm}}
\pgfpathlineto{\pgfqpoint{14.793cm}{13.529cm}}
\pgfpathlineto{\pgfqpoint{14.34cm}{13.529cm}}
\pgfusepath{stroke}
\pgfsetdash{}{0cm}
\pgfpathmoveto{\pgfqpoint{15.246cm}{13.529cm}}
\pgfpathlineto{\pgfqpoint{15.246cm}{15.202cm}}
\pgfpathlineto{\pgfqpoint{15.699cm}{15.202cm}}
\pgfpathlineto{\pgfqpoint{15.699cm}{13.529cm}}
\pgfpathlineto{\pgfqpoint{15.246cm}{13.529cm}}
\pgfusepath{stroke}
\pgfsetdash{}{0cm}
\pgfpathmoveto{\pgfqpoint{16.151cm}{13.87cm}}
\pgfpathlineto{\pgfqpoint{16.151cm}{15.02cm}}
\pgfpathlineto{\pgfqpoint{16.604cm}{15.02cm}}
\pgfpathlineto{\pgfqpoint{16.604cm}{13.87cm}}
\pgfpathlineto{\pgfqpoint{16.151cm}{13.87cm}}
\pgfusepath{stroke}
\pgfsetdash{}{0cm}
\pgfpathmoveto{\pgfqpoint{17.057cm}{13.87cm}}
\pgfpathlineto{\pgfqpoint{17.057cm}{15.02cm}}
\pgfpathlineto{\pgfqpoint{17.51cm}{15.02cm}}
\pgfpathlineto{\pgfqpoint{17.51cm}{13.87cm}}
\pgfpathlineto{\pgfqpoint{17.057cm}{13.87cm}}
\pgfusepath{stroke}
\pgfsetdash{}{0cm}
\pgfpathmoveto{\pgfqpoint{17.965cm}{13.547cm}}
\pgfpathlineto{\pgfqpoint{17.965cm}{15.031cm}}
\pgfpathlineto{\pgfqpoint{18.418cm}{15.031cm}}
\pgfpathlineto{\pgfqpoint{18.418cm}{13.547cm}}
\pgfpathlineto{\pgfqpoint{17.965cm}{13.547cm}}
\pgfusepath{stroke}
\pgfsetdash{}{0cm}
\pgfpathmoveto{\pgfqpoint{18.871cm}{13.547cm}}
\pgfpathlineto{\pgfqpoint{18.871cm}{15.031cm}}
\pgfpathlineto{\pgfqpoint{19.323cm}{15.031cm}}
\pgfpathlineto{\pgfqpoint{19.323cm}{13.547cm}}
\pgfpathlineto{\pgfqpoint{18.871cm}{13.547cm}}
\pgfusepath{stroke}
\pgfsetdash{}{0cm}
\pgfpathmoveto{\pgfqpoint{19.776cm}{13.573cm}}
\pgfpathlineto{\pgfqpoint{19.776cm}{15.125cm}}
\pgfpathlineto{\pgfqpoint{20.229cm}{15.125cm}}
\pgfpathlineto{\pgfqpoint{20.229cm}{13.573cm}}
\pgfpathlineto{\pgfqpoint{19.776cm}{13.573cm}}
\pgfusepath{stroke}
\pgfsetdash{}{0cm}
\pgfpathmoveto{\pgfqpoint{20.682cm}{13.573cm}}
\pgfpathlineto{\pgfqpoint{20.682cm}{15.125cm}}
\pgfpathlineto{\pgfqpoint{21.134cm}{15.125cm}}
\pgfpathlineto{\pgfqpoint{21.134cm}{13.573cm}}
\pgfpathlineto{\pgfqpoint{20.682cm}{13.573cm}}
\pgfusepath{stroke}
\pgfsetdash{}{0cm}
\pgfpathmoveto{\pgfqpoint{21.587cm}{13.479cm}}
\pgfpathlineto{\pgfqpoint{21.587cm}{14.97cm}}
\pgfpathlineto{\pgfqpoint{22.04cm}{14.97cm}}
\pgfpathlineto{\pgfqpoint{22.04cm}{13.479cm}}
\pgfpathlineto{\pgfqpoint{21.587cm}{13.479cm}}
\pgfusepath{stroke}
\pgfsetdash{}{0cm}
\pgfpathmoveto{\pgfqpoint{22.493cm}{13.479cm}}
\pgfpathlineto{\pgfqpoint{22.493cm}{14.97cm}}
\pgfpathlineto{\pgfqpoint{22.945cm}{14.97cm}}
\pgfpathlineto{\pgfqpoint{22.945cm}{13.479cm}}
\pgfpathlineto{\pgfqpoint{22.493cm}{13.479cm}}
\pgfusepath{stroke}
\pgfsetdash{}{0cm}
\definecolor{eps2pgf_color}{rgb}{1,0,0}\pgfsetstrokecolor{eps2pgf_color}\pgfsetfillcolor{eps2pgf_color}
\pgfpathmoveto{\pgfqpoint{1.658cm}{14.437cm}}
\pgfpathlineto{\pgfqpoint{2.111cm}{14.437cm}}
\pgfusepath{stroke}
\pgfsetdash{}{0cm}
\pgfpathmoveto{\pgfqpoint{2.564cm}{14.437cm}}
\pgfpathlineto{\pgfqpoint{3.016cm}{14.437cm}}
\pgfusepath{stroke}
\pgfsetdash{}{0cm}
\pgfpathmoveto{\pgfqpoint{3.469cm}{14.626cm}}
\pgfpathlineto{\pgfqpoint{3.922cm}{14.626cm}}
\pgfusepath{stroke}
\pgfsetdash{}{0cm}
\pgfpathmoveto{\pgfqpoint{4.374cm}{14.626cm}}
\pgfpathlineto{\pgfqpoint{4.827cm}{14.626cm}}
\pgfusepath{stroke}
\pgfsetdash{}{0cm}
\pgfpathmoveto{\pgfqpoint{5.28cm}{14.373cm}}
\pgfpathlineto{\pgfqpoint{5.733cm}{14.373cm}}
\pgfusepath{stroke}
\pgfsetdash{}{0cm}
\pgfpathmoveto{\pgfqpoint{6.185cm}{14.373cm}}
\pgfpathlineto{\pgfqpoint{6.638cm}{14.373cm}}
\pgfusepath{stroke}
\pgfsetdash{}{0cm}
\pgfpathmoveto{\pgfqpoint{7.094cm}{14.214cm}}
\pgfpathlineto{\pgfqpoint{7.547cm}{14.214cm}}
\pgfusepath{stroke}
\pgfsetdash{}{0cm}
\pgfpathmoveto{\pgfqpoint{7.999cm}{14.214cm}}
\pgfpathlineto{\pgfqpoint{8.452cm}{14.214cm}}
\pgfusepath{stroke}
\pgfsetdash{}{0cm}
\pgfpathmoveto{\pgfqpoint{8.905cm}{14.243cm}}
\pgfpathlineto{\pgfqpoint{9.357cm}{14.243cm}}
\pgfusepath{stroke}
\pgfsetdash{}{0cm}
\pgfpathmoveto{\pgfqpoint{9.81cm}{14.243cm}}
\pgfpathlineto{\pgfqpoint{10.263cm}{14.243cm}}
\pgfusepath{stroke}
\pgfsetdash{}{0cm}
\pgfpathmoveto{\pgfqpoint{10.716cm}{14.205cm}}
\pgfpathlineto{\pgfqpoint{11.168cm}{14.205cm}}
\pgfusepath{stroke}
\pgfsetdash{}{0cm}
\pgfpathmoveto{\pgfqpoint{11.621cm}{14.205cm}}
\pgfpathlineto{\pgfqpoint{12.074cm}{14.205cm}}
\pgfusepath{stroke}
\pgfsetdash{}{0cm}
\pgfpathmoveto{\pgfqpoint{12.529cm}{14.499cm}}
\pgfpathlineto{\pgfqpoint{12.982cm}{14.499cm}}
\pgfusepath{stroke}
\pgfsetdash{}{0cm}
\pgfpathmoveto{\pgfqpoint{13.435cm}{14.499cm}}
\pgfpathlineto{\pgfqpoint{13.888cm}{14.499cm}}
\pgfusepath{stroke}
\pgfsetdash{}{0cm}
\pgfpathmoveto{\pgfqpoint{14.34cm}{14.437cm}}
\pgfpathlineto{\pgfqpoint{14.793cm}{14.437cm}}
\pgfusepath{stroke}
\pgfsetdash{}{0cm}
\pgfpathmoveto{\pgfqpoint{15.246cm}{14.437cm}}
\pgfpathlineto{\pgfqpoint{15.699cm}{14.437cm}}
\pgfusepath{stroke}
\pgfsetdash{}{0cm}
\pgfpathmoveto{\pgfqpoint{16.151cm}{14.47cm}}
\pgfpathlineto{\pgfqpoint{16.604cm}{14.47cm}}
\pgfusepath{stroke}
\pgfsetdash{}{0cm}
\pgfpathmoveto{\pgfqpoint{17.057cm}{14.47cm}}
\pgfpathlineto{\pgfqpoint{17.51cm}{14.47cm}}
\pgfusepath{stroke}
\pgfsetdash{}{0cm}
\pgfpathmoveto{\pgfqpoint{17.965cm}{14.408cm}}
\pgfpathlineto{\pgfqpoint{18.418cm}{14.408cm}}
\pgfusepath{stroke}
\pgfsetdash{}{0cm}
\pgfpathmoveto{\pgfqpoint{18.871cm}{14.408cm}}
\pgfpathlineto{\pgfqpoint{19.323cm}{14.408cm}}
\pgfusepath{stroke}
\pgfsetdash{}{0cm}
\pgfpathmoveto{\pgfqpoint{19.776cm}{14.57cm}}
\pgfpathlineto{\pgfqpoint{20.229cm}{14.57cm}}
\pgfusepath{stroke}
\pgfsetdash{}{0cm}
\pgfpathmoveto{\pgfqpoint{20.682cm}{14.57cm}}
\pgfpathlineto{\pgfqpoint{21.134cm}{14.57cm}}
\pgfusepath{stroke}
\pgfsetdash{}{0cm}
\pgfpathmoveto{\pgfqpoint{21.587cm}{14.393cm}}
\pgfpathlineto{\pgfqpoint{22.04cm}{14.393cm}}
\pgfusepath{stroke}
\pgfsetdash{}{0cm}
\pgfpathmoveto{\pgfqpoint{22.493cm}{14.393cm}}
\pgfpathlineto{\pgfqpoint{22.945cm}{14.393cm}}
\pgfusepath{stroke}
\begin{pgfscope}
\pgfpathmoveto{\pgfqpoint{3.407cm}{17.745cm}}
\pgfpathlineto{\pgfqpoint{3.986cm}{17.745cm}}
\pgfpathlineto{\pgfqpoint{3.986cm}{17.166cm}}
\pgfpathlineto{\pgfqpoint{3.407cm}{17.166cm}}
\pgfpathclose
\pgfusepath{clip}
\pgfsetdash{}{0cm}
\pgfsetlinewidth{0.706mm}
\pgfpathmoveto{\pgfqpoint{3.59cm}{17.457cm}}
\pgfpathlineto{\pgfqpoint{3.801cm}{17.457cm}}
\pgfusepath{stroke}
\pgfsetdash{}{0cm}
\pgfpathmoveto{\pgfqpoint{3.695cm}{17.562cm}}
\pgfpathlineto{\pgfqpoint{3.695cm}{17.351cm}}
\pgfusepath{stroke}
\end{pgfscope}
\begin{pgfscope}
\pgfpathmoveto{\pgfqpoint{4.313cm}{17.745cm}}
\pgfpathlineto{\pgfqpoint{4.892cm}{17.745cm}}
\pgfpathlineto{\pgfqpoint{4.892cm}{17.166cm}}
\pgfpathlineto{\pgfqpoint{4.313cm}{17.166cm}}
\pgfpathclose
\pgfusepath{clip}
\pgfsetdash{}{0cm}
\pgfsetlinewidth{0.706mm}
\pgfpathmoveto{\pgfqpoint{4.495cm}{17.457cm}}
\pgfpathlineto{\pgfqpoint{4.707cm}{17.457cm}}
\pgfusepath{stroke}
\pgfsetdash{}{0cm}
\pgfpathmoveto{\pgfqpoint{4.601cm}{17.562cm}}
\pgfpathlineto{\pgfqpoint{4.601cm}{17.351cm}}
\pgfusepath{stroke}
\end{pgfscope}
\begin{pgfscope}
\pgfpathmoveto{\pgfqpoint{12.468cm}{11.939cm}}
\pgfpathlineto{\pgfqpoint{13.047cm}{11.939cm}}
\pgfpathlineto{\pgfqpoint{13.047cm}{11.218cm}}
\pgfpathlineto{\pgfqpoint{12.468cm}{11.218cm}}
\pgfpathclose
\pgfusepath{clip}
\pgfsetdash{}{0cm}
\pgfsetlinewidth{0.706mm}
\pgfpathmoveto{\pgfqpoint{12.65cm}{11.509cm}}
\pgfpathlineto{\pgfqpoint{12.862cm}{11.509cm}}
\pgfusepath{stroke}
\pgfsetdash{}{0cm}
\pgfpathmoveto{\pgfqpoint{12.756cm}{11.615cm}}
\pgfpathlineto{\pgfqpoint{12.756cm}{11.404cm}}
\pgfusepath{stroke}
\pgfsetdash{}{0cm}
\pgfpathmoveto{\pgfqpoint{12.65cm}{11.65cm}}
\pgfpathlineto{\pgfqpoint{12.862cm}{11.65cm}}
\pgfusepath{stroke}
\pgfsetdash{}{0cm}
\pgfpathmoveto{\pgfqpoint{12.756cm}{11.756cm}}
\pgfpathlineto{\pgfqpoint{12.756cm}{11.545cm}}
\pgfusepath{stroke}
\end{pgfscope}
\begin{pgfscope}
\pgfpathmoveto{\pgfqpoint{13.373cm}{11.939cm}}
\pgfpathlineto{\pgfqpoint{13.952cm}{11.939cm}}
\pgfpathlineto{\pgfqpoint{13.952cm}{11.218cm}}
\pgfpathlineto{\pgfqpoint{13.373cm}{11.218cm}}
\pgfpathclose
\pgfusepath{clip}
\pgfsetdash{}{0cm}
\pgfsetlinewidth{0.706mm}
\pgfpathmoveto{\pgfqpoint{13.555cm}{11.509cm}}
\pgfpathlineto{\pgfqpoint{13.767cm}{11.509cm}}
\pgfusepath{stroke}
\pgfsetdash{}{0cm}
\pgfpathmoveto{\pgfqpoint{13.661cm}{11.615cm}}
\pgfpathlineto{\pgfqpoint{13.661cm}{11.404cm}}
\pgfusepath{stroke}
\pgfsetdash{}{0cm}
\pgfpathmoveto{\pgfqpoint{13.555cm}{11.65cm}}
\pgfpathlineto{\pgfqpoint{13.767cm}{11.65cm}}
\pgfusepath{stroke}
\pgfsetdash{}{0cm}
\pgfpathmoveto{\pgfqpoint{13.661cm}{11.756cm}}
\pgfpathlineto{\pgfqpoint{13.661cm}{11.545cm}}
\pgfusepath{stroke}
\end{pgfscope}
\end{pgfscope}
\pgftext[x=1.869cm,y=10.875cm,rotate=0]{\fontsize{20}{7.06}\selectfont{{1}}}
\pgftext[x=2.787cm,y=10.875cm,rotate=0]{ \fontsize{20}{7.06}\selectfont{{2}}}
\pgftext[x=3.696cm,y=10.872cm,rotate=0]{ \fontsize{20}{7.06}\selectfont{{3}}}
\pgftext[x=4.603cm,y=10.875cm,rotate=0]{ \fontsize{20}{7.06}\selectfont{{4}}}
\pgftext[x=5.509cm,y=10.869cm,rotate=0]{ \fontsize{20}{7.06}\selectfont{{5}}}
\pgftext[x=6.416cm,y=10.872cm,rotate=0]{ \fontsize{20}{7.06}\selectfont{{6}}}
\pgftext[x=7.322cm,y=10.872cm,rotate=0]{ \fontsize{20}{7.06}\selectfont{{7}}}
\pgftext[x=8.226cm,y=10.872cm,rotate=0]{ \fontsize{20}{7.06}\selectfont{{8}}}
\pgftext[x=9.132cm,y=10.872cm,rotate=0]{ \fontsize{20}{7.06}\selectfont{{9}}}
\pgftext[x=10.05cm,y=10.872cm,rotate=0]{ \fontsize{20}{7.06}\selectfont{{10}}}
\pgftext[x=10.927cm,y=10.875cm,rotate=0]{ \fontsize{20}{7.06}\selectfont{{11}}}
\pgftext[x=11.859cm,y=10.875cm,rotate=0]{ \fontsize{20}{7.06}\selectfont{{12}}}
\pgftext[x=12.767cm,y=10.872cm,rotate=0]{ \fontsize{20}{7.06}\selectfont{{13}}}
\pgftext[x=13.672cm,y=10.875cm,rotate=0]{ \fontsize{20}{7.06}\selectfont{{14}}}
\pgftext[x=14.576cm,y=10.872cm,rotate=0]{ \fontsize{20}{7.06}\selectfont{{15}}}
\pgftext[x=15.483cm,y=10.872cm,rotate=0]{ \fontsize{20}{7.06}\selectfont{{16}}}
\pgftext[x=16.392cm,y=10.875cm,rotate=0]{ \fontsize{20}{7.06}\selectfont{{17}}}
\pgftext[x=17.296cm,y=10.872cm,rotate=0]{ \fontsize{20}{7.06}\selectfont{{18}}}
\pgftext[x=18.201cm,y=10.872cm,rotate=0]{ \fontsize{20}{7.06}\selectfont{{19}}}
\pgftext[x=19.094cm,y=10.872cm,rotate=0]{ \fontsize{20}{7.06}\selectfont{{20}}}
\pgftext[x=19.972cm,y=10.875cm,rotate=0]{ \fontsize{20}{7.06}\selectfont{{21}}}
\pgftext[x=20.903cm,y=10.875cm,rotate=0]{ \fontsize{20}{7.06}\selectfont{{22}}}
\pgftext[x=21.814cm,y=10.872cm,rotate=0]{ \fontsize{20}{7.06}\selectfont{{23}}}
\pgftext[x=22.72cm,y=10.875cm,rotate=0]{ \fontsize{20}{7.06}\selectfont{{24}}}
\pgftext[x=12.335cm,y=10.311cm-.2cm,rotate=0]{ \fontsize{24}{7.06}\selectfont{{Hour}}}
\pgftext[x=-0.753cm,y=14.533cm-.5cm,rotate=90]{ \fontsize{24}{7.06}\selectfont{{Forecast Error (MWh)}}}
\pgfsetdash{}{0cm}
\pgfusepath{stroke}
\end{pgfscope}
\end{pgfscope}
\end{pgfpicture}

%% file: appendix_1bus.tex
\section{Proofs of IID Case}
\label{sec:app:1bus}
We will prove the results in Section \ref{sec:onlineALG} by constructing a sequence of auxiliary optimization problems \textbf{P1} to \textbf{P3}. First, define
\[
\bar \u \defeq \lim_{T \to \infty} \frac{1}{T} \expec\left[ \sum_{t=1}^T \u_t\right],\,\,\bar \b \defeq \lim_{T \to \infty} \frac{1}{T} \expec\left[ \sum_{t=1}^T \b_t\right].
\]
Note that for $\b_1\in[\bmin,\bmax]$,
\[
\bar \u=\lim_{T \to \infty}\frac{1}{T}\expec\left[\sum_{t=1}^T\b_{t+1}-\lambda\b_t\right]=(1-\lambda)\bar\b.
\]
As $\b_t\in[\bmin,\bmax]$ for all $t\geq 0$, the above expression implies
\[
(1-\lambda)\bmin\leq \bar \u\leq (1-\lambda)\bmax.
\]
Then, problem~\eqref{prob:singleBusGeneral} can be equivalently written as follows
\begin{subequations}\label{prob:P1}
\begin{align}
\opttag{P1:} \minimize  &  \lim_{T \to \infty} \frac{1}{T} \expec \left[\sum_{t=1}^T \g_t \right]\\
\st  & \b_{t+1} = \la \b_t + \u_t, \label{P1:dynamics}\\
        %& \bmin \le \b_t \le \bmax, \\
        & \bmin - \la \b_t \le \u_t \le \bmax - \la \b_t, \label{P1:bbounds-u}\\
        & \umin \le \u_t \le \umax, \label{P1:ubounds}\\
        & f_t \in \F,\\
        &(1-\lambda)\bmin\leq \bar \u\leq (1-\lambda)\bmax \label{P1:ubar},
\end{align}
\end{subequations}
where bounds on $\b_t$ are replaced by \eqref{P1:bbounds-u}, and \eqref{P1:ubar} is added without loss of optimality.

\begin{figure}
\centering
\scalebox{0.6}{\scalefont{1.2} \input{./fig/proof2.pgf}}
\caption{An illustration of the proof procedure as relations between three problems considered. Here $\mathcal{S}$ denotes the sub-optimality bound.}
\label{fig:proofdiag}
\end{figure}
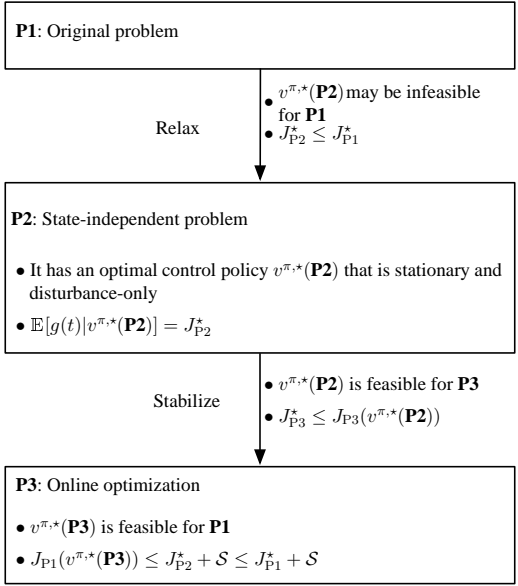

The proof procedure is depicted in the diagram shown in Figure~\ref{fig:proofdiag}. Here we use  \ylmod{$\J_\mathrm{P1}(\vpi)$ to denote the objective value of {\bf P1} with control policy sequence $\vpi = \{\upi, \fpi\}$, where $\upi$ and $\fpi$ are abbreviations of $\{\upi_t: t\ge 1\}$ and $\{\fpi_t: t \ge 1 \}$ respectively; $\vpistar(\mathrm{\bf P1})$ denotes an optimal control policy sequence for {\bf P1}, $\J^\star_\mathrm{P1} \defeq  \J_\mathrm{P1}(\vpistar(\mathrm{\bf P1}))$, and we define similar quantities for {\bf P2} and {\bf P3}. It is obvious that $\J_\mathrm{P1}(\vpi)=J(\vpi)$ and $\J^\star_\mathrm{P1}=J^\star$.} Here {\bf P2} is an auxilliary problem we construct to bridge the infinite horizon storage control problem {\bf P1} to online Lyapunov optimization problems {\bf P3} in \eqref{prob:P3}. It has the following form
\begin{subequations}\label{prob:P2}
\begin{align}
\opttag{P2:} \minimize  &  \lim_{T \to \infty} \frac{1}{T} \expec \left[\sum_{t=1}^T \g_t \right]\\
\st  %& \b_{t+1} = \la \b_t + \u_t, \label{P1:dynamics}\\
        %& \bmin \le \b_t \le \bmax, \\
        %& \bmin - \la \b_t \le \u_t \le \bmax - \la \b_t, \label{P1:bbounds-u}\\
        & \umin \le \u_t \le \umax, \label{P2:ubounds}\\
        & f_t \in \F,\\
        &(1-\lambda)\bmin\leq \bar \u\leq (1-\lambda)\bmax \label{P2:ubar}.
\end{align}
\end{subequations}
Notice that it has the same objective as \Pk{1}, and evidently it is a relaxation of \Pk{1}. This implies that \ylmod{$\vpistar(\Pk{2})$ (in particular $\upistar(\Pk{2})$)} may not be feasible for \Pk{1}, and
\begin{equation}
\Jk{2}^\star = \ylmod{\Jk{1}(\vpistar(\Pk{2}))} \le \Jk{1}^\star.
\end{equation}
The reason for the removal of state-dependent constraints \eqref{P1:bbounds-u} (and hence \eqref{P1:dynamics} as the sequence $\{\b_t: t\ge 1\}$ becomes irrelevant to the optimization of \ylmod{$\{\u_t: t\ge 1\}$)} in \Pk{2} is that the state-independent problem \Pk{2} has easy-to-characterize optimal stationary control policies. In particular, from the theory of stochastic network optimization \cite{NeelyBook}, the following result holds.
\begin{lemma}[Optimal Stationary Disturbance-Only Policies]\label{lemma:neely}\label{lem:stat}Under Assumption~\ref{ass:1}
there exists a stationary disturbance-only\footnote{The policy is a pure function (possibly randomized) of the current disturbances $\d_t$ and $\p_t$.} policy $\vpistat = (\upistat, \fpistat)$, satisfying \eqref{P2:ubounds} and \eqref{P2:ubar},
and providing the following guarantees for all $t$:
\begin{align}
  & (1-\la)\bmin\leq \expec[\ustat_t ] \leq (1-\la)\bmax,  \\
  & \expec[\g_t | \v_t = \vstat_t]  = \Jk{2}^\star, \label{eq:statCostLink}
\end{align}
where $\vstat_t = (\ustat_t, \fstat_t)$ is the control action induced by control policy $\vpistat$ at time $t$ and the expectation is taken over the randomization of $\d_t$, $\p_t$, and (possibly) $\vpistat$.
%\begin{equation*}\label{eq:singleBusGeneralCost_stat}
%\begin{split}
%\gstat_t=&\sum_{\ell = 1}^L \p(t, \ell) \Bigg(\alI(\ell) \d_t - \alC(\ell)\hC\left( \pos{\ustat_t}\right) \\
%&\quad + \alD(\ell) \hD\left(\neg{\ustat_t}\right) + \alConst(\ell)\Bigg)^+.
%\end{split}
%\end{equation*}
\end{lemma}
\begin{remark}\label{rk:neely}
Lemma~\ref{lemma:neely} holds for many non-i.i.d. disturbance processes as well.  One can generalize the results in Lemma~\ref{lemma:neely} to other stationary processes by invoking Theorem 4.5 of~\cite{NeelyBook}. Generalizing to the case without stationary assumptions is also possible; see~\cite{2010arXiv1003.3396N} and references therein for more details. 
\end{remark}

Equation~\eqref{eq:statCostLink} not only assures the storage operation induced by the stationary disturbance-only policy achieves the optimal cost, but also
guarantees that the expected stage-wise cost is a constant across time periods and equal to the optimal time average cost. This fact will later be exploited in order to establish the performance guarantee of our online algorithm. %By the merits of this lemma, in the sequel, we \ylmod{denote with $\v^\star(\Pk{2})$ }the control sequence obtained from an optimal stationary disturbance-only policy.

%\begin{figure}
%\centerline{
%\includegraphics[scale = .75]{../fig/proof-diag}}
%\caption{An illustration of the proof procedure as relations between problems considered in Section~\ref{sec:proofSingleBus}. Here $\mathcal{S}$ denotes the sub-optimality bound.}
%\label{fig:proofdiag}
%\end{figure}
%
%\begin{figure}
%\centering
%\scalebox{0.7}{\scalefont{1.2} \input{../fig/proof.pgf}}
%\caption{An illustration of the proof procedure as relations between problems considered in Section~\ref{sec:proofSingleBus}. Here $\mathcal{S}$ denotes the sub-optimality bound.}
%\label{fig:proofdiag}
%\end{figure}

An issue \ylmod{which arises in the application of control policy $\vpistar(\Pk{2})$ to} the original problem is that $\vpistar(\Pk{2})$ may not be feasible for \Pk{1}. To have the $\{\b_t: t\ge 1\}$ sequence induced by the storage operation sequence lie in the interval $[\bmin, \bmax]$, we construct a virtual queue related to $\b_t$ and use techniques from Lyapunov optimization to ``stabilize'' such a queue. Let the queueing state be a shifted version of the storage level:
\begin{equation}\label{eq:shift_queue}
\bs_t=\b_t+\ks,
\end{equation}
where the shift constant $\ks$ satisfies conditions~\eqref{eq:ksbounds}. We wish to minimize the stage-wise cost $\g_t$ and at the same time to maintain the queueing state close to zero.
This motivates us to consider solving the following optimization online (\ie, at the beginning of each time period $t$ after the realizations of stochastic parameters $\p_t$ and $\d_t$ have been observed)
\begin{subequations}\label{prob:P3}
\begin{align}
\opttag{P3:} \minimize  &  \la \bs_t \u_t + \W \tilde \g_t\label{eq:p3_obj}\\
\st  %& \b_{t+1} = \la \b_t + \u_t, \label{P1:dynamics}\\
        %& \bmin \le \b_t \le \bmax, \\
        %& \bmin - \la \b_t \le \u_t \le \bmax - \la \b_t, \label{P1:bbounds-u}\\
        & \umin \le \u_t \le \umax, \label{P3:ubounds}\\
        & \f_t \in \F, 
\end{align}
\end{subequations}
where the optimization variables are $\u_t$ and $\f_t$, 
%the stochastic parameters in $\g_t$ are replaced with their observed realizations, 
and $\W > 0$ is the weight parameter satisfying conditions~\eqref{eq:Wbounds}. 
%Note that the objective here is a weighted combination of the stage-wise cost and a linear term of $\u_t$, whose coefficient is positive when $\b_t$ is large, and negative when $\b_t$ is small.
%Meanwhile, different from \Pk{1} and \Pk{2}, which are specified for the entire horizon, \Pk{3} is specified and solved for each time period $t$.
We use the notations $\vhat_t$ for the solution to \Pk{3} at time period $t$, $\v^\star(\Pk{3})$ for the sequence $\{\vhat_t: t\ge 1\}$, $\Jkt{3}{t}(\v_t)$ for the objective function of \Pk{3} at time period $t$, and $\Jkt{3}{t}^\star$ for the corresponding optimal cost. Note that \Pk{3} is implemented in the online phase of Algorithm \ref{alg} (see the optimization problem in \eqref{prob:P3_alg}) and $\v^\star(\Pk{3})=\{\v_t,\,t\geq 1\}$ where $\v_t$ is the solution of problem \eqref{prob:P3_alg} at time $t$. \ylmod{Furthermore, denote with $\vpistar(\Pk{3})$ the corresponding control policy defined by the online optimization (which generates $\v^\star(\Pk{3})$).} We also define the corresponding quantities for $\u$ and $\f$. 

We break the proof of Theorem~\ref{thm:perf_lyap} into two parts -- feasibility and performance. In order to prove the feasibility of control policy $\upistar(\Pk{3})$ (and hence $\vpistar(\Pk{3})$), the following technical lemma is needed.

\begin{lemma}[Structural Properties of Online Optimization]\label{lem:ana_sol}
Let $\uhat_t$ be the optimal storage operation obtained via solving \eqref{prob:P3_alg} at time $t$. The following statements hold: 
\begin{enumerate}
\item if $\la (\b_t+\ks) + \W \Dl g  \geq 0$, then $\uhat_t = \umin$;
\item if $\la (\b_t+\ks) + \W \Du g \leq 0$, then $\uhat_t = \umax$.
\end{enumerate}
\end{lemma}
\begin{IEEEproof}
Let $J(\u,\f)=\la (\b_t+\ks) \u + \W \g_t(\u, \f, \drl_t, \prl_t)$ be the objective function of \eqref{prob:P3_alg} after the stochastic parameters $\drl_t$ and $\prl_t$ are realized. Recall $\phi_t(u, y) \defeq \g_t(\u, f, \drl_t, \prl_t)$ where $y=(f, \drl_t, \prl_t)$ and let $J_t(u)=\sup_{y\in\mathcal Y}\phi_t(u, y)$.
To show the set of sufficient conditions for $\uhat_t$ takes $\umax$ (or $\umin$), notice that the condition
\[
\lambda (\b_t+\ks) \leq-\W\Du \g
\]
 implies $ \partial_u J_t(u)\vert_{u=u_t}\subseteq(-\infty,0]$, for any given $y\in\mathcal Y$. Thus, for every given $\u\in[\umin,\umax]$, if $\beta$ is a constant such that
 \[
 J_t(v)-J_t(u)\geq \beta\cdot(v-\u),\,\, \forall v\in[\umin,\umax],
  \]
then the sub-differential condition implies that $\beta\leq 0$. Now, by substituting $\u=\umax$ in the above expression, one obtains $\beta\cdot(v-\u)\geq 0$ and $J_t(v)\geq J_t(\umax)$,  for all $v\in[\umin,\umax]$. Therefore, one concludes that $\u_t = \umax$ attains an optimal solution in \eqref{prob:P3_alg}.
 Similarly, the condition
 \[
 \lambda\tilde \b_t \geq-\W\Dl g
 \]
 implies $\partial_u J_t(u)\vert_{u=u_t}\subseteq[0,\infty) $. Based on analogous arguments, one concludes that $\u_t = \umin$ attains an optimal solution in \eqref{prob:P3_alg}.
\end{IEEEproof}

Now, we are in position to prove that \ylmod{the control policy $\vpistar(\Pk{3})$}  is a feasible solution to $\textbf{P1}$ (and the stochastic control problem in \eqref{prob:singleBusGeneral}).
\begin{IEEEproof}[Proof of Theorem~\ref{thm:perf_lyap}, Feasibility]
%\noindent{\bf Proof of Theorem~\ref{thm:feas_single}}
We first validate that the intervals of $\ks$ and $\W$ are non-empty. Note that from Assumption \ref{ass:1}, $\Wmax> 0$, thus it remains to show $\ksmax \ge \ksmin$.
%, this directly implies that the bound for $W$ is non-empty. Now, we turn to prove that the bound of $\ks$ is non-empty by showing $\ksmax\geq \ksmin$ whenever $0\leq \W  \le \Wmax$.
Based on \eqref{eq:W_max}, $W>0$, and $\Du\g\geq \Dl\g$, one obtains
\[
\W(\Du \g - \Dl \g)\leq [(\bmax - \bmin) - (\umax-\umin)].
\]
Re-arranging terms results in
\[
-\W \Dl g + \umax - \bmax \le - \W \Du g - \bmin + \umin,
\]
which further implies $\ksmax\geq \ksmin$.

We proceed to show that
\begin{equation}\label{exp:MI}
\bmin\leq  \b_t\leq \bmax,
\end{equation}
for $t= 1,2,\dots$, when \ylmod{control action $\u^\star(\Pk{3})$} is implemented.
The base case holds by assumption.
Let the inductive hypothesis be that (\ref{exp:MI}) holds at time $t$.
The storage level at $t+1$ is  then
$\b_{t+1}=\lambda \b_t+\uhat_t. $
%Recall from expression (\ref{eq:shift_K}) that $\ks\in[\ksmin,\ksmax]$.
We show  (\ref{exp:MI}) holds at $t+1$ by considering the following three cases.

\noindent{\bf Case 1.}  $- \W \Dl\g\leq\la\bs_t\leq\la ( \bmax+\ks)$. \\
First, it is easy to verify that the above interval for $\la \bs_t$ is non-empty using \eqref{eq:ineq_1} and $\ks \ge \ksmin$. Next, based on Lemma \ref{lem:ana_sol}, one obtains  $\uhat_t=\umin\le 0$ in this case. Therefore
\[
\b_{t+1} = \la \b_t + \umin \le  \la \bmax + \umin \le \bmax,
\]
where the last inequality follows from the feasibility assumption in Definition~\ref{assume:feas}. On the other hand,
\begin{align*}
\b_{t+1}& = \la \b_t + \umin \ge -\W\Dl \g  - \la \ks + \umin \\
\ge & -\W\Dl \g  - \la \ksmax + \umin \\
\ge  &  W [\Du g - \Dl g] + \bmin
\ge    \bmin,
\end{align*}
where the third inequality follows from the definition of $\ksmax$, and the fourth inequality used $\Du \g \ge \Dl \g$.
%This implies that when $\uhat_{k}$ is implemented,
%\[
%\begin{split}
%\tilde \b(k+1)&<\lambda\tilde \b(k)+(1-\lambda)\ks+\umin\\
%&\leq \lambda\bmax+\ks+\umin\leq \bmax+\ks.
%\end{split}
%\]
%The last inequality is due to the second statement in Assumption \ref{assume:feas}.  On the other hand,
%\[
%\begin{split}
%\tilde \b(k+1)\geq&\lambda\tilde \b(k)+(1-\lambda)\ks+\umin\\
%\geq & - \W \Dl\g+(1-\lambda)\ks+\umin\\
%\geq&\ks+\bmin ,
%\end{split}
%\]
%the last inequality follows from expression (\ref{eq:ineq_2}) and the fact that
%\[\begin{split}
%- \W \Dl\g\geq& - \W \Du\g\\
%=&{\lambda( \bmin+\ksmax)}+ {((1-\lambda)\bmin-\umin)^+}\\
%\geq &{\lambda( \bmin+\ks)}+ {((1-\lambda)\bmin-\umin)^+}\\
%\geq&\bmin+\lambda \ks-\umin.
%\end{split}\]

\noindent{\bf Case 2.}
 $\la( \bmin+\ks) \le  \la\bs_t  \le - \W \Du\g$.  \\
The above interval for $\la \bs_t$ is non-empty by \eqref{eq:ineq_2} and $\ks \le \ksmax$. Lemma \ref{lem:ana_sol} implies  $\uhat_t=\umax \ge 0$ in this case.
Therefore, by the feasibility assumption,
\[
\b_{t+1} = \la \b_t + \umax \ge \la \bmin + \umax \ge \bmin.
\]
On the other hand,
\begin{align*}
\!\!\b_{t+1} &= \la \b_t + \umax \le -\W\Du g - \la \ks + \umax \\
\le &  -\W\Du g - \la \ksmin + \umax \\
\le &  -W [\Du g - \Dl g] + \bmax
\le \bmax,
\end{align*}
where the third inequality used the definition of $\ksmin$, and the fourth inequality again is by $\Du \g \ge \Dl \g$.

%First, based on expression (\ref{eq:ineq_1}), we know that the above interval for $\lambda\tilde{\b}(k)$ is non-empty. Next, based on Lemma \ref{coro:tech_res}, one obtains that the solution of \textbf{P3} equals $\uhat(k)=\umax>0$. This implies that when $\uhat(k)$ is implemented,
%\[
%\begin{split}
%\tilde \b(k+1)&\geq\lambda\tilde \b(k)+\umax+(1-\lambda)\ks\\
%&\geq \lambda\bmin+\umax+\ks\geq \bmin+\ks.
%\end{split}
%\]
%The last inequality is due to the first statement in Assumption \ref{assume:feas}. On the other hand,
%\[
%\begin{split}
%\tilde \b(k+1)\leq&\lambda\tilde \b(k)+(1-\lambda)\ks+\umax\\
%\leq& -\W \Du\g+(1-\lambda)\ks+\umax\\
%\leq&\ks+\bmax ,
%\end{split}
%\]
%the last inequality follows from expression (\ref{eq:ineq_1}) and the fact that
%\[\begin{split}
%-\W \Du\g\leq& -\W \Dl\g\\
%=&{-(\umax-(1-\lambda)\bmax)^+}+{\lambda (\bmax+\ksmin)}\\
%\leq&{-(\umax-(1-\lambda)\bmax)^+}+{\lambda (\bmax+\ks)}\\
%\leq&\bmax+\lambda \ks-\umax.
%\end{split}\]
\noindent{\bf Case 3.} $-\W \Du\g<\la\bs_t<- \W \Dl\g$.  \\
%Based on expression (\ref{eq:ineq_1}) and expression (\ref{eq:ineq_2}), one first obtains
%\[
%- \W \Dl\g\geq {\lambda( \bmin+\ks)},\,\,
%- \W \Du\g\leq {\lambda (\bmax+\ks)}.
%\]
By  $\umin \le \uhat_t\le\umax$, one obtains
\begin{align*}
\b_{t+1} &=  \la \b_t + \uhat_t \le \la \b_t + \umax\\
< & - \W \Dl\g - \la \ks + \umax\\
\le &  - \W \Dl\g - \la \ksmin + \umax \le \bmax,
\end{align*}
where the last inequality is by the definition of $\ksmin$.
On the other hand,
\begin{align*}
\b_{t+1} &=  \la \b_t + \uhat_t \ge \la \b_t + \umin\\
> & - \W \Du\g - \la \ks + \umin\\
\ge &  - \W \Du\g - \la \ksmax + \umin\ge    \bmin,
\end{align*}
where the last inequality follows from the definition of $\ksmax$. 

Combining these three cases, and by mathematical induction, we conclude \eqref{exp:MI} holds for all $t = 1,2,\dots$.
%\[
%\begin{split}
%\tilde \b(k+1)&\leq \lambda\tilde \b(k)+\umax+(1-\lambda)\ks\\
%& < -\W \Dl\g+(1-\lambda)\ks+\umax\\
%&\leq \bmax+\ks,
%\end{split}
%\]
%the last inequality follows from expression (\ref{eq:ineq_3}) and the fact that
%\[
%\begin{split}
%\W \Dl\g\leq&{-(\umax-(1-\lambda)\bmax)^+}+{\lambda (\bmax+\ks)}\\
%\leq&\bmax+\lambda \ks-\umax.
%\end{split}
%\]
%On the other hand, based on similar arguments, one obtains
%\[
%\begin{split}
%\tilde \b(k+1)&\geq \lambda\tilde \b(k)+\umin+(1-\lambda)\ks\\
%& >- \W \Du\g +\umin+(1-\lambda)\ks\\
%&\geq \bmin+\ks,
%\end{split}
%\]
%the last inequality follows from expression (\ref{eq:ineq_1}) and the fact that
%\[
%\begin{split}
%-\W \Du\g\geq &{\lambda( \bmin+\ks)}+ {((1-\lambda)\bmin-\umin)^+}\\
%\geq&\bmin+\lambda \ks-\umin.
%\end{split}\]
%Thus, based on the above cases, we have shown that equation (\ref{exp:MI}) holds at $t=k+1$. By induction, this in turns implies that equation (\ref{exp:MI}) holds for any $t\geq 1$.
%
%Finally from equation (\ref{eq:shift_queue}) that $ \b_t=\tilde \b_t-\ks$. When the sequence of control actions is given by the the solution of \textbf{P3}, expression (\ref{exp:MI}) holds for any $t\in\{1,2,\ldots,\}$. This implies $\bmin\leq \b_t$ and $\b_t\leq \bmax$, for any $t\in\{1,2,\ldots,\}$, which completes the second part of the IEEEproof.
\end{IEEEproof}
%Next, by using the feasibility result and exploiting the relationships between the objective functions of \textbf{P1} to \textbf{P3}, one can prove the sub-optimality of $u^\star(\textbf{P3})$  with respect to \textbf{P1} using the ``drift--plus--penalty" Lyapunov optimization technique. This is also equivalent to the sub-optimality of storage operation sequence $\{\un_t,\,\t\geq 1\}$ generated by Algorithm \ref{alg} with respect to the stochastic control problem in \eqref{prob:singleBusGeneral}.

We proceed to prove the sub-optimality of \ylmod{control policy $\vpistar(\textbf{P3})$}.
\begin{IEEEproof}[Proof of Theorem~\ref{thm:perf_lyap}, Performance]
%\noindent{\bf Proof of Theorem~\ref{thm:perf_lyap}}
Consider a quadratic Lyapunov function $\L(\b) = \b^2/2$.
Let the corresponding Lyapunov drift be
\[
\Ls( \bs_t) = \expec\left[ \L(\bs_{t+1}) -  \L(\bs_t) \vert  \bs_t \right].
\]
Recall that
$
\bs_{t+1} = \b_{t+1} + \ks = \la \bs_t   + \u_t + (1-\la)\ks,
$
and so
\begin{align}
\Ls( \bs_t)%& = \expec\left[  \half\u_t^2 + \half(1-\la)^2\ks^2 -  \half (1-\la^2)\bs_t^2 + \la \bs_t \u_t + \la (1-\la) \bs_t \ks + (1-\la) \u_t\ks \vert  \bs_t \right]\\
& = \expec\big[ (1/2)(\u_t + (1-\la)\ks)^2  -  (1/2) (1-\la^2)\bs_t^2 \nn\\
                &\quad\quad\quad + \la \bs_t \u_t + \la (1-\la) \bs_t \ks \vert  \bs_t \big]\nn\\
& \le \Mone(\ks)    -  (1/2) (1-\la^2)\bs_t^2\nn \\
& \quad+\expec\big[ \la \bs_t \u_t + \la (1-\la) \bs_t \ks \vert  \bs_t \big]\nn\\
& \le\Mone(\ks) + \expec\left[ \la \bs_t (\u_t+(1-\la)\ks)\vert \bs_t \right].\label{eq:1busdrift}
\end{align}
It follows that, with arbitrary control action $\v_t$,
\begin{align}
& \Ls(\bs_t) + \W \expec [\g_t| \bs_t] \label{eq:ME}\\
\le & \Mone(\ks) + \la (1-\la)\bs_t\ks+  \expec\big[ \Jkt{3}{t}(\v_t) | \bs_t],\nn
\end{align}
%
%\begin{equation}
%\begin{split}
%&\Delta(\tilde \b_t) + \W \expec[\g_t \vert \tilde \b_t ] \\
%\leq& \Mone(\ks) +  \expec\big[ \la\tilde \b_t (\u_t+(1-\la)\ks)+ \W\g_t \vert \tilde \b_t \big],\\
%= &  \Mone(\ks) + \expec\big[
% \label{eq:ME}
%\end{split}
%\end{equation}
where it is clear that minimizing the right hand side of the above inequality over $\v_t$ is equivalent to minimizing the objective of \Pk{3}. Given that $\vstat_t$, \ylmod{the control action induced by disturbance-only stationary policy $\vpistat$ of \Pk{2}}  described in Lemma~\ref{lem:stat}, is feasible for \Pk{3}, the above inequality implies\footnote{The notation $\expec [\g_t|\vpistat]$ is an abbreviation for $\expec [\g_t|v_t = \vstat_t]$. Similar abbreviation appears in Appendix B.}
\begin{align}
& \Ls(\bs_t) + \W \expec [\g_t| \bs_t, \v_t = \vhat_t] \label{eq:ME}\\
\le & \Mone(\ks) + \la (1-\la)\bs_t\ks+  \expec\big[ \Jkt{3}{t}^\star | \bs_t] \nn \\
\le & \Mone(\ks) + \la (1-\la)\bs_t\ks+  \expec\big[ \Jkt{3}{t}(\vstat_t) | \bs_t] \nn\\
\stackrel{(a)}{=} & \Mone(\ks) +\la\bs_t\expec\left[ \ustat_t +(1-\la)\ks \right]+\W\expec [\g_t|\vpistat_t]\nn \\
\stackrel{(b)}{\le} & \M(\ks) + \W\expec[\g_t|\vstat_t]\stackrel{(c)}{\le}  \M(\ks) + \W \Jk{1}^\star.\nn
\end{align}
Here $(a)$ uses the fact that $\ustat_t$ is induced by a disturbance-only stationary policy;  $(b)$ follows from inequalities
$ |\bs_t|\leq \left(\max\left( (\bmax+\ks)^2,(\bmin+\ks)^2\right)\right)^{1/2} $ and
$\left|\expec\left[ \ustat_t\right]+(1-\la)\ks\right|
\le (1-\la) (\max( (\bmax+\ks)^2,(\bmin+\ks)^2))^{1/2};
$
and $(c)$ used $\expec[\g_t|\vstat_t]  = \Jk{2}^\star$ in Lemma~\ref{lem:stat} and $\Jk{2}^\star \le \Jk{1}^\star$. Taking expectation over $\bs_t$ on both sides gives\begin{align}
&  \expec\left[ \L(\bs_{t+1}) -  \L( \bs_t) \right] + \W\expec \left[\g_t|\v_t = \vhat_t \right]    \nn\\
  \leq & \M(\ks)+\W \Jk{1}^\star. \label{eq:diff_ly}
\end{align}
Summing expression \eqref{eq:diff_ly} over $t$ from $1$ to $T$, dividing both sides by $\W T$, taking the limit $T\rightarrow \infty$ and noting that $\J^\star_\mathrm{P1}=J^\star$, we obtain the performance bound in expression \eqref{eq:1busiidperf_bdd}.

\end{IEEEproof}
\begin{remark}[Finite Termination]\label{rk:ft}
In the above proof, one notes that with a finite $T$, we get the bound
\[
\frac{1}{T} \sum_{t=1}^T J_{\mathrm{P3},t}^\star \le  \Jk{1}^\star + \frac{\M(\ks)}{\W} + \frac{1}{WT} \expec [L(\bs_1) - L(\bs_{T+1})],
\]
where the last term could serve as a proxy for estimating the error in the performance bounds in Theorem~\ref{thm:perf_lyap} if a finite $T$ is used.
\end{remark}

Finally, Lemma~\ref{SDP_P3_PO} can be easily proved using the Schur complement as follows.
\begin{IEEEproof}[Proof of Lemma~\ref{SDP_P3_PO}]
%\noindent{\bf Proof of Lemma~\ref{SDP_P3_PO}}
Based on the following re-parametrizations
\[
\None=\Mone(\ks)/\W,\,\, \Ntwo=\Mtwo(\ks)/\W,
\]
(since $\W> 0$) one can easily show that problem \textbf{PO} has the same solution as the following optimization problem:
\begin{subequations}
\begin{align*}
\minimize \quad & \None+\la(1-\la)\Ntwo\\
\st \quad & \ksmin\le  \ks \le \ksmax,0 <  \W  \le \Wmax, \\
&2\None \W\geq  \left(\umin+(1-\la)\ks\right)^2,\\
&2\None \W\geq \left(\umax+(1-\la)\ks\right)^2,\\
&\Ntwo \W\geq  \left(\bmin+\ks\right)^2,\Ntwo \W \geq \left(\bmax+\ks\right)^2.
\end{align*}
\end{subequations}
The proof is completed by applying Schur complement on the last four constraints of the above optimization.
\end{IEEEproof}

%% file: fig/proof2.pgf
% Created by Eps2pgf 0.7.0 (build on 2008-08-24) on Tue May 20 23:10:13 PDT 2014
\begin{pgfpicture}
\pgfpathmoveto{\pgfqpoint{0cm}{0cm}}
\pgfpathlineto{\pgfqpoint{11.43cm}{0cm}}
\pgfpathlineto{\pgfqpoint{11.43cm}{13.053cm}}
\pgfpathlineto{\pgfqpoint{0cm}{13.053cm}}
\pgfpathclose
\pgfusepath{clip}
\begin{pgfscope}
\end{pgfscope}
\begin{pgfscope}
\begin{pgfscope}
\pgfpathmoveto{\pgfqpoint{0cm}{0cm}}
\pgfpathlineto{\pgfqpoint{11.43cm}{0cm}}
\pgfpathlineto{\pgfqpoint{11.43cm}{13.053cm}}
\pgfpathlineto{\pgfqpoint{0cm}{13.053cm}}
\pgfpathclose
\pgfusepath{clip}
\definecolor{eps2pgf_color}{rgb}{1,1,1}\pgfsetstrokecolor{eps2pgf_color}\pgfsetfillcolor{eps2pgf_color}
\pgfpathmoveto{\pgfqpoint{-2.487cm}{16.034cm}}
\pgfpathlineto{\pgfqpoint{17.833cm}{16.034cm}}
\pgfpathlineto{\pgfqpoint{17.833cm}{-9.825cm}}
\pgfpathlineto{\pgfqpoint{-2.487cm}{-9.825cm}}
\pgfpathclose
\pgfusepath{fill}
\pgfsetdash{}{0cm}
\pgfsetlinewidth{0.353mm}
\pgfsetroundcap
\pgfsetroundjoin
\definecolor{eps2pgf_color}{rgb}{0,0,0}\pgfsetstrokecolor{eps2pgf_color}\pgfsetfillcolor{eps2pgf_color}
\pgfpathmoveto{\pgfqpoint{0.071cm}{12.982cm}}
\pgfpathlineto{\pgfqpoint{11.359cm}{12.982cm}}
\pgfpathlineto{\pgfqpoint{11.359cm}{11.501cm}}
\pgfpathlineto{\pgfqpoint{0.071cm}{11.501cm}}
\pgfpathclose
\pgfusepath{stroke}
\definecolor{eps2pgf_color}{rgb}{0,0,0}\pgfsetstrokecolor{eps2pgf_color}\pgfsetfillcolor{eps2pgf_color}
\pgftext[x=0.2cm,y=12.209cm,rotate=0,left,base]{  \selectfont{{\bf P1}: Original problem}}
\pgfsetdash{}{0cm}
\definecolor{eps2pgf_color}{rgb}{0,0,0}\pgfsetstrokecolor{eps2pgf_color}\pgfsetfillcolor{eps2pgf_color}
\pgfpathmoveto{\pgfqpoint{0.071cm}{8.978cm}}
\pgfpathlineto{\pgfqpoint{11.359cm}{8.978cm}}
\pgfpathlineto{\pgfqpoint{11.359cm}{5.203cm}}
\pgfpathlineto{\pgfqpoint{0.071cm}{5.203cm}}
\pgfpathclose
\pgfusepath{stroke}
\definecolor{eps2pgf_color}{rgb}{0,0,0}\pgfsetstrokecolor{eps2pgf_color}\pgfsetfillcolor{eps2pgf_color}
\pgftext[x=0.2cm,y=8.041cm,rotate=0,left,base]{{\bf P2}: State-independent problem}
\pgfsetdash{}{0cm}
\pgftext[x=0.2cm,y=6.90cm,rotate=0,left,base]{ {$\bullet$ It has an optimal control policy $\vpistar(\Pk{2})$ that is stationary and }}
\pgftext[x=0.5cm,y=6.40cm,rotate=0,left,base]{ {disturbance-only}}
\pgftext[x=0.2cm,y=5.7cm,rotate=0,left,base]{ {$\bullet$ $\expec[\g(t)|\vpistar(\Pk{2})] = \Jk{2}^\star$}}
\pgfsetdash{}{0cm}
\definecolor{eps2pgf_color}{rgb}{0,0,0}\pgfsetstrokecolor{eps2pgf_color}\pgfsetfillcolor{eps2pgf_color}
\pgfpathmoveto{\pgfqpoint{0.071cm}{2.681cm}}
\pgfpathlineto{\pgfqpoint{11.359cm}{2.681cm}}
\pgfpathlineto{\pgfqpoint{11.359cm}{0.088cm}}
\pgfpathlineto{\pgfqpoint{0.071cm}{0.088cm}}
\pgfpathclose
\pgfusepath{stroke}
\definecolor{eps2pgf_color}{rgb}{0,0,0}\pgfsetstrokecolor{eps2pgf_color}\pgfsetfillcolor{eps2pgf_color}
\pgftext[x=0.2cm,y=2.14cm,rotate=0,left,base]{  \selectfont{{\bf P3}: Online optimization}}
\pgfsetdash{}{0cm}
\pgftext[x=0.2cm,y=1.2cm,rotate=0,left,base]{ $\bullet$ $\vpistar(\Pk{3})$ is feasible for \Pk{1}}
\pgftext[x=0.2cm,y=0.5cm,rotate=0,left,base]{ $\bullet$ $\Jk{1}(\vpistar(\Pk{3})) \le \Jk{2}^\star +\mathcal S \le \Jk{1}^\star + \mathcal S$}
\definecolor{eps2pgf_color}{rgb}{0,0,0}\pgfsetstrokecolor{eps2pgf_color}\pgfsetfillcolor{eps2pgf_color}
\pgfpathmoveto{\pgfqpoint{5.724cm}{11.501cm}}
\pgfpathlineto{\pgfqpoint{5.72cm}{9.345cm}}
\pgfusepath{stroke}
\pgfpathmoveto{\pgfqpoint{5.719cm}{9.063cm}}
\pgfpathlineto{\pgfqpoint{5.825cm}{9.345cm}}
\pgfpathlineto{\pgfqpoint{5.614cm}{9.345cm}}
\pgfpathclose
\pgfusepath{fill}
\pgfsetdash{}{0cm}
\pgfsetbuttcap
\pgfsetmiterjoin
\pgfpathmoveto{\pgfqpoint{5.719cm}{9.063cm}}
\pgfpathlineto{\pgfqpoint{5.825cm}{9.345cm}}
\pgfpathlineto{\pgfqpoint{5.614cm}{9.345cm}}
\pgfpathclose
\pgfusepath{stroke}
\pgfsetdash{}{0cm}
\pgfsetroundcap
\pgfsetroundjoin
\pgfpathmoveto{\pgfqpoint{5.715cm}{5.186cm}}
\pgfpathlineto{\pgfqpoint{5.715cm}{3.048cm}}
\pgfusepath{stroke}
\pgfpathmoveto{\pgfqpoint{5.715cm}{2.766cm}}
\pgfpathlineto{\pgfqpoint{5.821cm}{3.048cm}}
\pgfpathlineto{\pgfqpoint{5.609cm}{3.048cm}}
\pgfpathclose
\pgfusepath{fill}
\pgfsetdash{}{0cm}
\pgfsetbuttcap
\pgfsetmiterjoin
\pgfpathmoveto{\pgfqpoint{5.715cm}{2.766cm}}
\pgfpathlineto{\pgfqpoint{5.821cm}{3.048cm}}
\pgfpathlineto{\pgfqpoint{5.609cm}{3.048cm}}
\pgfpathclose
\pgfusepath{stroke}
\definecolor{eps2pgf_color}{rgb}{0,0,0}\pgfsetstrokecolor{eps2pgf_color}\pgfsetfillcolor{eps2pgf_color}
\pgftext[x=5.7cm,y=10.657cm,rotate=0,left,base]{ $\bullet$ {\parbox{15em}{$\vpistar(\Pk{2})\!$ may be infeasible \\for \Pk{1}}}}
\pgftext[x=5.7cm,y=9.957cm,rotate=0,left,base]{ $\bullet$ {$\Jk{2}^\star \le  \Jk{1}^\star$}}
\pgftext[x=5.7cm,y=4.38cm,rotate=0,left,base]{ $\bullet$ {$\vpistar(\Pk{2})$ is feasible for \Pk{3}}}
\pgftext[x=5.7cm,y=3.68cm,rotate=0,left,base]{ $\bullet$ {$\Jk{3}^\star \le  \Jk{3}(\vpistar(\Pk{2}))$}}
\pgftext[x=3.843cm,y=10.202cm,rotate=0]{ {Relax}}
\pgftext[x=4.04cm,y=4.138cm,rotate=0]{ {Stabilize}}
\end{pgfscope}
\end{pgfscope}
\end{pgfpicture}

%% file: appendix_markov_process.tex
\def\R{\mathcal{R}}
\section{Generalization to Non-IID Cases}\label{appendix_markov_process}
Markov models are widely used in the power system applications for the modeling of stochastic demand, renewable generation, and price processes (cf. \cite{1388488, 6545387, 5453026}).
We demonstrate how our results can be generalized to non-i.i.d. cases by  establishing similar performance bounds for \emph{ergodic} Markov chains. The proof technique is based on the well-known method of analyzing regenerative cycles of the underlying disturbance process. 

We consider the following particular disturbance model. 
Suppose that the uncertain parameter vector  $(\d_t,\p_t)$ is some deterministic function of the system stochastic state $\omega_t$, where $\omega_t$ follows a finite state ergodic Markov Chain, supported on $\Omega$. Here by ergodic, we mean $\{\omega_t: t\ge 1\}$ is stationary, positive recurrent and irreducible. 
Let $\omega^\mathrm{R}\in\Omega$ be the initial state of $\omega_t$. 
Since $\omega_t$ is an ergodic Markov chain, there exists a sequence of finite random return time $1 = T_1 < T_2< \dots < T_r < T_{r+1}< \dots$, for $r = 1, 2,\dots$, such that $\omega_{t}$ visits $\omega^\mathrm{R}$ for the $r$-th time at time $t = T_r$. %Define $\R_{t}$ as the number of visits of $\omega_0$ at time $t$. Specifically, $\R_{t}\defeq\max\{r:T_r\leq t\}$. 
From this sequence of return times, we define the $r-$th epoch as $[T_r,T_{r+1}-1]$ and the length of this epoch is defined as $\DT_r=T_{r+1}-T_r$.
Apparently, the sequence of $\{\DT_r: r \ge 1\}$ is i.i.d.. Let  $\DT$ be a random variable distributed as $\DT_1$ and independent with all $\DT_r$, $r \ge 1$. The positive recurrence assumption implies that $\expec[\DT]<\infty$. We also assume that the second moment of $\DT$ is bounded, \ie, $\expec\left[\DT^2\right]<\infty$.
% and define the mean return rate of state $w_0$ as $\renewrate = {1}/{\expec[\DT]}$. 

As the proof of the feasibility of the OMG algorithm does not depend on the assumptions on the disturbance process, we focus on the performance analysis in the remaining of this appendix.

\begin{theorem} [Performance]\label{thm:perf_lyap_renew}
The sub-optimality of storage operation \ylmod{control policy $\vpistar(\Pk{3})$} is bounded by $\Mfour(\ks)/\W$ with probability one, that is
\begin{equation}\label{eq:perf_bdd_renewal}
\Jk{1}^\star \le  \ylmod{\Jk{1}(\vpistar(\Pk{3}))} \le \Jk{1}^\star + \Mfour(\ks)/\W
\end{equation}
with probability one, where
\begin{equation}\label{eq:M4}
\Mfour(\ks)= \frac{\expec [\DT^2]}{\expec [\DT]} \Mone(\ks)+\frac{\la (1-\expec \left[\la^{\DT} \right])}{ \expec [\DT]} \Mtwo(\ks),
\end{equation}
and $\Mtwo(\ks)$ and $\Mone(\ks)$ are  defined in Theorem \ref{thm:perf_lyap}.
\end{theorem}
\begin{IEEEproof}
Consider a quadratic Lyapunov function $\L(\b) = \b^2/2$ and the corresponding Lyapunov drift $\Ls( \bs_t) = \expec\left[ \L(\bs_{t+1}) -  \L(\bs_t) \vert  \bs_t \right]$.
Based on the analysis in expression (\ref{eq:1busdrift}), we have that
\[
\Ls( \bs_t) \le\Mone(\ks) + \expec\left[ \la \bs_t (u_t+(1-\la)\ks)\vert \bs_t \right]
\]
holds for any $t$. Consider the $r$-th epoch $[T_r, T_{r+1}-1]$. For this analysis, we will first treat $T_r$ and $T_{r+1}$ as fixed deterministic quantities, and then consider that they are in fact random and take expectation over them.  Applying above inequality gives
\begin{align}
&\expec \left[ \sum_{t = T_r}^{T_{r+1}-1} \Ls(\bs_t) + W g_t \middle| \bs_{T_r}\right] \label{eq:mcp1}\\
& \le \DT_r \Mone(\ks) + \expec\left[ \sum_{t = T_r}^{T_{r+1}-1}\la \bs_t (u_t+(1-\la)\ks) + Wg_t\middle| \bs_{T_r} \right]. \nonumber
\end{align}
Using the tower property of iterative conditional expectation, one recognizes that the last term of the right hand side of~\eqref{eq:mcp1} is the same as the sum of the objectives of $\Pk{3}$ for $t = T_r, \dots, T_{r+1}-1$, apart from a constant term. As~\eqref{eq:mcp1} holds for arbitrary control policy, and the stationary disturbance only policy in Lemma~\ref{lemma:neely}, \ie, the solution of \Pk{2}, is feasible for $\Pk{3}$, we have
\begin{align*}
&\expec \left[ \sum_{t = T_r}^{T_{r+1}-1} \Ls(\bs_t) + W g_t \middle| \bs_{T_r}, \vpistar(\Pk{3})\right] \label{eq:mcp1}\\
& \le\!\! \DT_r \Mone(\ks)\!\! +\!\! \expec\!\!\left[\! \sum_{t = T_r}^{T_{r+1}-1}\!\!\!\!\la \bs_t (\u_t\!+\!(1\!-\!\la)\ks)\!\! + \!\!W\g_t\middle| \bs_{T_r},\!\vpistar(\Pk{3})\! \right]\\
& \le\!\! \DT_r \Mone(\ks)\!\! +\!\! \expec\!\!\left[\! \sum_{t = T_r}^{T_{r+1}-1}\!\!\!\!\la \bs_t (\u_t\!+\!(1\!-\!\la)\ks)\!\! + \!\!W\g_t\middle| \bs_{T_r},\!\vpistar(\Pk{2})\! \right]\\
& =\!\! \DT_r (\Mone(\ks)+ W\Jk{2}^\star)\!\! +\!\! \expec\!\!\left[\! \sum_{t = T_r}^{T_{r+1}-1}\!\!\!\!\la \bs_t (\ustat_t\!+(1-\la)\ks) \middle| \bs_{T_r}\! \right]\!,
\end{align*}
where the last identity is by Lemma~\ref{lemma:neely} (see Remark~\ref{rk:neely} for the applicability in this case). The fact that the disturbance process is Markov makes the one step bound for $\la \bs_t (\ustat_t\!+(1-\la)\ks) $ no longer directly applicable here. Instead, we bound the last term of the right hand side of the last inequality as follows:
\begin{align*}
 &\expec\!\!\left[\! \sum_{t = T_r}^{T_{r+1}-1}\!\!\!\!\la \bs_t (\ustat_t\!+(1-\la)\ks) \middle| \bs_{T_r}\! \right] \\
&\le \underbrace{\expec \left[\sum_{t = T_r}^{T_{r+1}-1}\la (\bs_t - \la^{t - T_r}\bs_{T_r}) (\ustat_t\!+(1-\la)\ks) \middle|\bs_{T_r}\right]}_{\mathcal{B}_1} \\
& \quad +  \underbrace{\bs_{T_r} \expec\left[\sum_{t = T_r}^{T_{r+1}-1}\la^{t - T_r+1}(\ustat_t\!+(1-\la)\ks)\middle|\bs_{T_r}\right]}_{\mathcal{B}_2} ,
\end{align*}
where by the same arguments proving \eqref{eq:ME}, 
\[
\mathcal B_2 \le \la (1-\la^{\DT_r}) \Mtwo(\ks).
\]
On the other hand,
\begin{align*}
&\mathcal B_1 \le  \expec \left[\sum_{t = T_r+1}^{T_{r+1}-1}\la |\bs_t - \la^{t - T_r}\bs_{T_r}| |\ustat_t\!+(1-\la)\ks| \middle|\bs_{T_r}\right]\\
& \le  \expec \left[\sum_{t = T_r+1}^{T_{r+1}-1} \sum_{\ell=1}^{t - T_r} \la^\ell |u_{t-\ell} + (1-\la) \ks  |  |\ustat_t\!+(1-\la)\ks| \middle|\bs_{T_r}\right]\\
%& \le  \expec \left[(1-\la)\sum_{t = T_r}^{T_{r+1}-1} \sum_{\ell=1}^{t - T_r} \la^\ell \Mtwo(\ks)^{1/2}  |\ustat_t\!+(1-\la)\ks| \right]\\
& \le 2  \Mone(\ks) \sum_{t = T_r+1}^{T_{r+1}-1} \sum_{\ell=1}^{t - T_r} \la^\ell  \le \Mone(\ks) \DT_r (\DT_r - 1), 
\end{align*}
where the first term (for $t = T_r$) in the summation that appeared in the definition of $\mathcal B_1$  is removed as it is zero, and the second inequality is due to the fact that
\[
\bs_{t_2} = \la^{t_2 - t_1} \bs_{t_1} + \sum_{\ell=1}^{t_2 - t_1} \la^{\ell - 1} (u_{t_2 - \ell} + (1- \la) \ks)
\]
for any $t_2 > t_1>0$.%\footnote{Alternative bounds that are tighter for certain combinations of the storage parameters $\Sm$ and the distribution of the return time $\DT$ may be developed for $\mathcal B_1$ based on bounding $\bs_t$ directly. However, we have chosen to present a bounding procedure which is mathematically simple and consistent with that in the i.i.d. case.} 
Thus for the $r$-th epoch, we have that
\begin{align*}
&\expec \left[ L(\bs_{T_{r+1}}) - L(\bs_{T_r}) + \sum_{t = T_r}^{T_{r+1}-1}  W g_t\middle|\bs_{T_r}, \vpistar(\Pk{3})\ \right]\\
&=\expec \left[ \sum_{t = T_r}^{T_{r+1}-1} \Ls(\bs_t) + W g_t \middle| \bs_{T_r}, \vpistar(\Pk{3})\right]\\
& \le   \Mone(\ks)\DT_r^2 + \la (1-\la^{\DT_r}) \Mtwo(\ks) +\DT_r W\Jk{2}^\star.
\end{align*}
Taking expectation over  the return times and $s_{T_r}$, and summing over epochs $1,\dots, R$ gives
\begin{align*}
&\expec \left[ L(\bs_{T_{R}}) - L(\bs_{1}) + \sum_{t = 1}^{T_{R}}  W g_t\middle|\bs_{1}, \vpistar(\Pk{3})\ \right]\\
&\le \!\!  R\expec [\DT^2]\Mone(\ks) \!\! + \!\!  R \la (1 \!\! - \!\! \expec \left[\la^{\DT} \right]) \Mtwo(\ks)  \!\! + \!\!  R\expec[\DT] W\Jk{2}^\star.
\end{align*}
Dividing both sides by $W R\expec [\DT]$ and sending $R\to \infty$ yields
\[
\Jk{1}(\vpistar (\Pk{3})) \le \Jk{2}^\star + \Mfour(\ks)/W \le \Jk{1}^\star + \Mfour(\ks)/W,
\]
where we have used the fact that, by elementary renewal theorem, $T_R/R \to \expec [\DT]$ with probability one, and that $\Jk{2}^\star \le \Jk{1}^\star$. 
 \end{IEEEproof}
\begin{remark}[Beyond Stationary Models]
The technique above can be easily generalized to other stationary processes of regenerative natures. Under suitable technical conditions, bootstrapping this analysis to processes that are not initially stationary, but converge to a limiting/stationary distribution,  such as many Markov models and martingales, is a standard excise in probability theory. Extending to processes that are fundamentally non-stationary requires a new analysis. Most importantly, the ``equilibrium''  notions of optimality may no longer apply. Interested readers are referred to \cite[Section 4.9.2]{NeelyBook} for the use of the so-called ``$T$-slot lookahead metric'' for establishing performance guarantees in non-stationary contexts.
\end{remark}

%% file: Main_singleStor_TPS_OMG.bbl
\begin{thebibliography}{10}
\providecommand{\url}[1]{#1}
\csname url@samestyle\endcsname
\providecommand{\newblock}{\relax}
\providecommand{\bibinfo}[2]{#2}
\providecommand{\BIBentrySTDinterwordspacing}{\spaceskip=0pt\relax}
\providecommand{\BIBentryALTinterwordstretchfactor}{4}
\providecommand{\BIBentryALTinterwordspacing}{\spaceskip=\fontdimen2\font plus
\BIBentryALTinterwordstretchfactor\fontdimen3\font minus
  \fontdimen4\font\relax}
\providecommand{\BIBforeignlanguage}[2]{{%
\expandafter\ifx\csname l@#1\endcsname\relax
\typeout{** WARNING: IEEEtran.bst: No hyphenation pattern has been}%
\typeout{** loaded for the language `#1'. Using the pattern for}%
\typeout{** the default language instead.}%
\else
\language=\csname l@#1\endcsname
\fi
#2}}
\providecommand{\BIBdecl}{\relax}
\BIBdecl

\bibitem{Denholm2010}
\BIBentryALTinterwordspacing
{National Renewable Energy Laboratory}. (2010) {The Role of Energy Storage with
  Renewable Electricity Generation}. [Online]. Available:
  \url{http://www.nrel.gov/wind/pdfs/47187.pdf}
\BIBentrySTDinterwordspacing

\bibitem{thermalStor1993}
B.~Daryanian and R.~E. Bohn, ``{Sizing of Electric Thermal Storage under Real
  Time Pricing},'' \emph{{IEEE} Trans. on Power Systems}, vol.~8, no.~1, pp.
  35--43, 1993.

\bibitem{ThatteXieStorValue2012s}
A.~A. Thatte and L.~Xie, ``{Towards a Unified Operational Value Index of Energy
  Storage in Smart Grid Environment},'' \emph{IEEE Trans. on Smart Grid},
  vol.~3, no.~3, pp. 1418--1426, 2012.

\bibitem{ObRACC2013}
G.~O'Brien and R.~Rajagopal, ``{A Method for Automatically Scheduling Notified
  Deferrable Loads},'' in \emph{Proc. of American Control Conference (ACC)},
  2013, pp. 5080--5085.

\bibitem{Callaway2009}
D.~S. Callaway, ``{Tapping the Energy Storage Potential in Electric Loads to
  Deliver Load Following and Regulation, with Application to Wind Energy},''
  \emph{Energy Conversion and Management}, vol.~50, no.~5, pp. 1389 -- 1400,
  2009.

\bibitem{HaoSanandajiPoollaVincent2013}
H.~Hao, B.~Sanandaji, K.~Poolla, and T.~Vincent, ``{Aggregate Flexibility of
  Thermostatically Controlled Loads},'' \emph{IEEE Trans. on Power Systems},
  vol.~30, no.~1, pp. 189--198, Jan 2015.

\bibitem{QinPESGM12}
J.~Qin, R.~Sevlian, D.~Varodayan, and R.~Rajagopal, ``{Optimal Electric Energy
  Storage Operation},'' in \emph{Proc. of IEEE Power and Energy Society General
  Meeting}, July 2012, pp. 1--6.

\bibitem{MITrampStor}
A.~{Faghih}, M.~{Roozbehani}, and M.~A. {Dahleh}, ``{On the Economic Value and
  Price-Responsiveness of Ramp-Constrained Storage},'' \emph{{ArXiv} e-prints},
  2012.

\bibitem{SuEGTPS}
H.~I. Su and A.~El~Gamal, ``{Modeling and Analysis of the Role of Energy
  Storage for Renewable Integration: Power Balancing},'' \emph{IEEE Trans. on
  Power Systems}, vol.~28, no.~4, pp. 4109--4117, 2013.

\bibitem{RLDSACC}
J.~Qin, H.~I. Su, and R.~Rajagopal, ``{Storage in Risk Limiting Dispatch:
  Control and Approximation},'' in \emph{Proc. of American Control Conference
  (ACC)}, 2013, pp. 4202--4208.

\bibitem{BitarRACC_colocated}
E.~Bitar, R.~Rajagopal, P.~Khargonekar, and K.~Poolla, ``{The Role of
  Co-Located Storage for Wind Power Producers in Conventional Electricity
  Markets},'' in \emph{Proc. of American Control Conference (ACC)}, 2011, pp.
  3886--3891.

\bibitem{Powell}
J.~H. Kim and W.~B. Powell, ``{Optimal Energy Commitments with Storage and
  Intermittent Supply},'' \emph{Operations Research}, vol.~59, no.~6, pp.
  1347--1360, 2011.

\bibitem{IBMload}
P.~M. {van de Ven}, N.~{Hegde}, L.~{Massoulie}, and T.~{Salonidis}, ``{Optimal
  Control of End-User Energy Storage},'' \emph{{ArXiv} e-prints}, 2012.

\bibitem{DataCenter}
R.~Urgaonkar, B.~Urgaonkar, M.~J. Neely, and A.~Sivasubramaniam, ``{Optimal
  Power Cost Management Using Stored Energy in Data Centers},'' in \emph{Proc.
  of the {ACM} SIGMETRICS Joint International Conference on Measurement and
  Modeling of Computer Systems}, 2011, pp. 221--232.

\bibitem{StorDRLongbo}
L.~Huang, J.~Walrand, and K.~Ramchandran, ``{Optimal Demand Response with
  Energy Storage Management},'' in \emph{Proc. of IEEE Third International
  Conference on Smart Grid Communications (SmartGridComm)}, 2012, pp. 61--66.

\bibitem{chow1991optimal}
C.~S. Chow and J.~N. Tsitsiklis, ``{An Optimal One-Way Multigrid Algorithm for
  Discrete-Time Stochastic Control},'' \emph{IEEE Trans. on Automatic Control},
  vol.~36, no.~8, pp. 898--914, 1991.

\bibitem{Jain:2010:SOM:1836310.1836495}
R.~Jain and P.~P. Varaiya, ``{Simulation-Based Optimization of Markov Decision
  Processes: An Empirical Process Theory Approach},'' \emph{Automatica},
  vol.~46, no.~8, pp. 1297--1304, Aug. 2010.

\bibitem{XieEtAlWindStorMPC}
L.~Xie, Y.~Gu, A.~Eskandari, and M.~Ehsani, ``{Fast MPC-Based Coordination of
  Wind Power and Battery Energy Storage Systems},'' \emph{Journal of Energy
  Engineering}, vol. 138, no.~2, pp. 43--53, 2012.

\bibitem{NRELStorValue2013}
\BIBentryALTinterwordspacing
{National Renewable Energy Laboratory}. (2013) {The Value of Energy Storage for
  Grid Applications}. [Online]. Available:
  \url{http://www.nrel.gov/docs/fy13osti/58465.pdf}
\BIBentrySTDinterwordspacing

\bibitem{NeelyBook}
M.~J. Neely, ``{Stochastic Network Optimization with Application to
  Communication and Queueing Systems},'' \emph{Synthesis Lectures on
  Communication Networks}, vol.~3, no.~1, pp. 1--211, 2010.

\bibitem{lyap1}
S.~Chen, P.~Sinha, and N.~Shroff, ``{Scheduling Heterogeneous Delay Tolerant
  Tasks in Smart Grid with Renewable Energy},'' in \emph{Proc. of IEEE 51st
  Annual Conference on Decision and Control (CDC)}, Dec 2012, pp. 1130--1135.

\bibitem{lyap2}
Q.~{Li}, T.~{Cui}, R.~{Negi}, F.~{Franchetti}, and M.~D. {Ilic}, ``{On-line
  Decentralized Charging of Plug-In Electric Vehicles in Power Systems},''
  \emph{ArXiv e-prints}, Jun. 2011.

\bibitem{QCYR:acm}
J.~Qin, Y.~Chow, J.~Yang, and R.~Rajagopal, ``{Modeling and Online Control of
  Generalized Energy Storage Networks},'' in \emph{Proc. of the 5th
  International Conference on Future Energy Systems (ACM e-Energy '14)}.\hskip
  1em plus 0.5em minus 0.4em\relax ACM, June 2014.

\bibitem{QCYR:report}
\BIBentryALTinterwordspacing
------, ``{Control of Generalized Energy Storage Networks},'' \emph{{Stanford
  S3L Report}}, 2014. [Online]. Available:
  \url{http://arxiv.org/abs/1504.05661}
\BIBentrySTDinterwordspacing

\bibitem{DecisionHorizon1988}
C.~Bes and S.~P. Sethi, ``{Concepts of Forecast and Decision Horizons:
  Applications to Dynamic Stochastic Optimization Problems},''
  \emph{Mathematics of Operations Research}, vol.~13, no.~2, pp. 295--310,
  1988.

\bibitem{bertsekas2007dynamic}
D.~P. Bertsekas, \emph{{Dynamic Programming and Optimal Control, Two-Volume
  Set}}.\hskip 1em plus 0.5em minus 0.4em\relax Athena Scientific, 2007.

\bibitem{bertsimas2009constructing}
D.~Bertsimas and D.~B. Brown, ``{Constructing Uncertainty Sets for Robust
  Linear Optimization},'' \emph{Operations research}, vol.~57, no.~6, pp.
  1483--1495, 2009.

\bibitem{NREL2010}
{National Renewable Energy Laboratory}, ``{Eastern Wind Integration and
  Transmission Study},'' Tech. Rep.

\bibitem{7101867}
J.~Qin, Y.~Chow, J.~Yang, and R.~Rajagopal, ``{Distributed Online Modified
  Greedy Algorithm for Networked Storage Operation Under Uncertainty},''
  \emph{IEEE Trans. on Smart Grid}, vol.~PP, no.~99, pp. 1--1, 2015.

\bibitem{2010arXiv1003.3396N}
M.~J. {Neely}, ``{Stability and Capacity Regions for Discrete Time Queueing
  Networks},'' \emph{ArXiv e-prints}, Mar. 2010.

\bibitem{1388488}
A.~Gonzalez, A.~Roque, and J.~Garcia-Gonzalez, ``{Modeling and Forecasting
  Electricity Prices with Input/Output Hidden Markov Models},'' \emph{IEEE
  Trans. on Power Systems}, vol.~20, no.~1, pp. 13--24, Feb 2005.

\bibitem{6545387}
A.~Albert and R.~Rajagopal, ``{Smart Meter Driven Segmentation: What Your
  Consumption Says About You},'' \emph{IEEE Trans. on Power Systems}, vol.~28,
  no.~4, pp. 4019--4030, Nov 2013.

\bibitem{5453026}
Y.-Y. Hong and K.-L. Pen, ``{Optimal VAR Planning Considering Intermittent Wind
  Power Using Markov Model and Quantum Evolutionary Algorithm},'' \emph{IEEE
  Trans. on Power Delivery}, vol.~25, no.~4, pp. 2987--2996, Oct 2010.

\end{thebibliography}
